# Discrete Graphical Models — An Optimization Perspective

Bogdan Savchynskyy

*Heidelberg University; bogdan.savchynskyy@iwr.uni-heidelberg.de*

---

ABSTRACT

This monograph is about combinatorial optimization. More precisely, about a special class of combinatorial problems known as *energy minimization* or *maximum a posteriori (MAP) inference in graphical models*, closely related to *weighted* and *valued constraint satisfaction problems* and having tight connections to *Markov random fields* and *quadratic pseudo-boolean optimization*. What distinguishes this monograph from a number of other monographs on graphical models is its focus: It considers graphical models, or, more precisely, MAP-inference for graphical models, purely as a combinatorial optimization problem. Modeling, applications, probabilistic interpretations and many other aspects are either ignored here or find their place in examples and remarks only.

Combinatorial optimization as a field is largely based on five fundamental topics: (i) integer linear programming and polyhedral optimization; (ii) totally unimodular matrices and the class of min-cost-flow problems; (iii) Lagrange decompositions and relaxations; (iv) dynamic programming and (v) submodularity, matroids and greedy algorithms. Each of these topics found its place in this monograph, although to a variable extent. The covering of each respective topic

---

reflects its importance for the considered MAP-inference problem.

Since optimization is the primary topic of this monograph, we mostly stick to the terminology widely used in optimization and where it was possible we tried to avoid the graphical models community-specific technical terms. The latter differ from sub-community to sub-community and, in our view, significantly complicate the information exchange between them.

The same holds also for the presentation of material in this monograph. If there is a choice when introducing mathematical constructs or proving statements, we prefer more general mathematical tools applicable in the whole field of operations research rather than to stick to graphical model-specific constructions. We additionally provide the graphical model-specific constructions if it turns out to be easier than the more general one. This way of presentation has two advantages. A reader familiar with a more general technique can grasp the new material faster. On the other hand, the monograph may serve as an introduction to combinatorial optimization for readers unfamiliar with this subject. To make the monograph even more suitable for both categories of readers we explicitly separate the mathematical optimization background chapters from those specific to graphical models.

We believe, therefore, that the monograph can be useful for undergraduate and graduate students studying optimization or graphical models, as well as for experts in optimization who want to have a look into graphical models. Moreover, we believe that even experts in graphical models can find new views on the known facts as well as a novel presentation of less known results in the monograph. These are for instance (i) a simple and general proof of equivalence of different acyclic Lagrange decompositions of a graphical model; (ii) a general scheme for analyzing convergence of dual block-coordinate descent methods; (iii) a short and self-contained analysis of a linear programming relaxation for binary graphical models, its persistency properties and its relation to quadratic pseudo-boolean optimization.

The present monograph is based on lectures given to undergraduate students of Technical University of Dresden and Heidelberg University. The selection of material is done in a way that it may serve as a basis for a semester course.

# Notation

To simplify reading of the monograph, some frequently used notations are collected here. Some of them, which we assume to be quite standard, are used without additional notice in the text. Others, typically more specialized, are introduced in the monograph. For those we point out the section and the page they are defined in.

## Standard notation

| | |
|---|---|
| $\mathbb{N}$ | the set of natural numbers |
| $\mathbb{Z}$ | the set of integer numbers |
| $\mathbb{R}$ | the set of real numbers |
| $\mathbb{R}^n$ | an $n$-dimensional vector space over the field of real numbers |
| $\mathbb{R}^n_+$ | the set of vectors with non-negative coordinates in $\mathbb{R}^n$, i.e. $\{x \in \mathbb{R}^n \colon x_i \geq 0, \ i = 1, 2, \ldots, n\}$; for $n = 1$ the notation simplifies to $\mathbb{R}_+$ |
| $x \in \mathcal{A}^{\mathcal{B}}$ | For any set $A$ and any finite set $\mathcal{B}$, this notation stands for a vector $x$ with $|\mathcal{B}|$ coordinates indexed by elements of $\mathcal{B}$, for each $b \in \mathcal{B}$ it holds that $x_b \in \mathcal{A}$. The only exception from this rule is the notation $\Delta^{\mathcal{B}}$, see below. |
| $x \geq y$ | comparison operations are applied coordinate-wise to vectors and point-wise to functions |
| $[\![\cdot]\!]$ | denotes the Iverson brackets, that is, for any predicate $A$ it holds that $[\![A]\!] = 1$ if $A$ is true, otherwise $[\![A]\!] = 0$ |
| $\langle c, x \rangle$ | the inner product, i.e. $\langle c, x \rangle = \sum_{i=1}^n c_i x_i$ |
| $\nabla f$ | gradient of the function $f$ |
| $O(\cdot)$ | for two functions $f \colon \mathbb{N} \to \mathbb{N}$ and $g \colon \mathbb{N} \to \mathbb{N}$ one writes $f = O(g)$, if there is a constant $c > 0$ and a number $n_0 \in \mathbb{N}$ such that $f(n) \leq c \cdot g(n)$ for all $n \geq n_0$ |

## Standard abbreviations

| | |
|---|---|
| w.r.t. | with respect to |
| w.l.o.g. | without loss of generality |
| s.t. | subject to |

## Notation defined in the monograph

| | |
|---|---|
| $\mathcal{N}_b(u)$ | set of graph vertexes incident to vertex $u$, see §1.1, page 9 |
| $\delta_{\mathcal{G}}(y)$, $\delta(y)$ | binary representation of the labeling $y$, i.e. a binary vector with non-zero coordinates corresponding to the labels $y_u$, $u \in \mathcal{Y}_u$, and label pairs $(y_u, y_v)$, $uv \in \mathcal{E}$, see page 49; |
| $\mathcal{I}$ | the set of indexes of the cost vector of a graphical model; $|\mathcal{I}|$ is equal to the number of coordinates of the cost vector, see §1.1 on page 11 |
| $\mathrm{vrtx}(P)$ | set of vertexes of the polyhedron $P$, see Definition 3.19 on page 32 |
| $\Delta^n$ | $n$-dimensional simplex, see Definition 3.21 on page 33 |
| $\Delta^{\mathcal{X}}$ | $|\mathcal{X}|$-dimensional simplex, with coordinates indexed by elements of $\mathcal{X}$, see Definition 3.21 on page 33 |
| $\mathrm{conv}(X)$ | convex hull of $X$, see Definition 3.28 on page 34 |
| $\mathrm{mi}[\theta_w]$ | binary vector with non-zero coordinates corresponding to locally minimal values of the cost vector $\theta_w$, see page 95 |
| $\mathrm{nz}[\mu]$ | binary vector with non-zero coordinates corresponding to the non-zero coordinates of $\mu$, see page 95 |
| $\mathrm{cl}(\xi)$ | arc-consistency closure of a binary vector $\xi$, see Definition 6.11 on page 99 |
| $\mathcal{J}$ | the set of indexes of the Lagrange dual vector for the MAP-inference problem; $|\mathcal{J}|$ is equal to the number of coordinates in the dual vector, see §6.1 on page 85 |
| $\left\langle \frac{1}{2} \right\rangle$, $\langle 0.7 \rangle$ | angular brackets are used in figures for coordinates of primal relaxed solutions, see e.g. Figure 4.1, 4.2, 6.4, 12.3 |

# 1

---

# Introduction to Inference for Graphical Models

---

There are many problems in computer science, which can be formulated in the form of so-called *Discrete Graphical Models*. Examples can be found in bio-informatics, communication theory, statistical physics, computer vision, signal processing, information retrieval and machine learning.

Discrete graphical models as a modeling tool naturally appear when

- the target object (the object we model) consists of *many small parts*,

- each part must be labeled by a label from a *finite set*, and

- parts (and, therefore, their labels) are *mutually dependent*.

**Example 1.1** (Image segmentation)**.** Image segmentation is a classical image analysis task: Each pixel of an input image must be assigned a label of an object visible in the image. For instance, if we consider images of street scenes, these labels could belong to the set {pedestrian, car, tree, building}.

The target object is an image, i.e. a two-dimensional array of pixels. Each pixel constitutes an elementary part of the image and must be labeled with a label from a finite set. The simplest assumption about image segments, i.e. groups of pixels having the same label, is the so called "compactness assumption". It states that it is more probable that neighboring pixels are labeled with the same label than with different ones.

**Example 1.2** (Depth reconstruction). Depth reconstruction is another important image analysis problem. In the classical setting there are (at least) two images taken from different viewpoints. The task is to match pixels from these two images to each other. Assuming the positions of cameras and their focal lengths are known, this allows us to estimate depth of the scene, which was photographed with the cameras.

As in the previous example, the target object is a two-dimensional pixel array, where each pixel constitutes an elementary part of the object and must be labeled with a label from a finite set. Here, the meaning of the labels is different: Each label represents depth information of the associated pixel in an image, i.e. how far the depicted observation is placed from the camera. Usually the set of labels is chosen as natural numbers in a given interval, for instance, $\{0, 1, \ldots, 255\}$.

Assuming that the observed surface is smooth, one would expect the difference $|s - s'|$ between labels $s$ and $s'$ in neighboring pixels to be small. The opposite would mean a jump in depth, or, in other words, non-smoothness of the surface.

**Example 1.3** (Cell tracking problem in bio-imaging). Given is a sequence of images that show the development of a living organism from an early embryo consisting of only a few cells to a fully grown animal. During this sequence, the images show moving and splitting cells.

Under the assumption that the image is already *pre-segmented*, i.e. the cells were already found in each image, the task at hand is to track each individual cell and its descendants from the first to the last frame.

The cells are the elementary parts of the considered object. Each cell in a given image frame is labeled with pointers to one or two cells in the next frame. One pointer means that the cell only moved, and two pointers correspond to a cell division. The simplest tracking model forbids two different cells to have the same descendants. This rule defines dependencies between object parts.

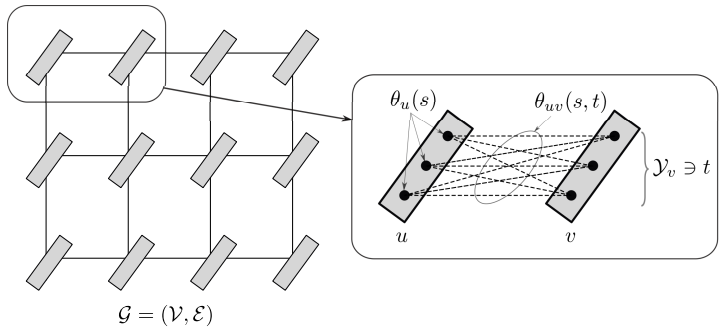

**Figure 1.1:** Example of a graphical model with grid structure. On the left, graph nodes are denoted with inclined rectangles, lines connecting nodes correspond to graph edges. On the right, two neighboring nodes are shown. Black circles inside rectangles correspond to the labels $s$ in the node $u$ (left rectangle) and $t$ in the node $v$ (right rectangle). Dashed lines correspond to each label pair $s, t$ with an assigned pairwise cost $\theta_{uv}(s, t)$.

## 1.1 Basic definitions

**Graph**  Let $\mathcal{G} = (\mathcal{V}, \mathcal{E})$ be an undirected graph consisting of a finite set of *nodes* $\mathcal{V}$ and a set of edges $\mathcal{E} \subseteq \binom{\mathcal{V}}{2}$. The set $\mathcal{E}$ will also be called *a neighborhood structure* of $\mathcal{V}$. For convenience, we will typically use lower case letters $u$ and $v$ for nodes of the graph, and write $uv$ to denote an edge $\{u, v\} \in \mathcal{E}$ connecting $u$ and $v$. Since the graph is undirected, $uv$ and $vu$ denote the same edge. The notation $\mathcal{N}_b(u)$ will be used for the set of nodes $\{v \mid uv \in \mathcal{E}\}$ connected to the node $u$.

The graph $\mathcal{G}$ is considered as a model of the considered target object, where the nodes represent the elementary object parts and edges stand for mutually dependencies between them.

In Examples 1.1 and 1.2 the graph $\mathcal{G}$ may have the grid structure of the underlying two-dimensional pixel array. In Example 1.3 cells of one image frame are neighbors, since their labels depend on each other.

**Labels and unary costs**  A finite *set of labels* $\mathcal{Y}_u$ is associated with each node $u \in \mathcal{V}$. Our preference for each label is expressed by the *unary cost function* $\theta_u \colon \mathcal{Y}_u \to \mathbb{R}$, which is defined for each node $u \in \mathcal{V}$. The value $\theta_u(s)$ determines the *cost*, which we pay for assigning label $s \in \mathcal{Y}_u$ to the node $u$. Sometimes we will use very high costs to implicitly forbid certain labels oder label pairs. The notation $\infty$ will be used to denote such high costs.

Unary costs are usually defined by what is known from observation. In Example 1.1, typically, the color distribution in the vicinity of a given pixel defines the cost of each possible label. The difference between color distributions from two or more images of the same scene taken from different viewpoints determines the unary costs in Example 1.2. Unary costs are often called the "data term" to emphasize that they depend on the input data or observation.

**Dependence and pairwise costs**      Dependencies between labels assigned to different graph nodes are modeled with *pairwise cost functions* $\theta_{uv} \colon \mathcal{Y}_u \times \mathcal{Y}_v \to \mathbb{R}$, which are defined for each edge $uv \in \mathcal{E}$ of the graph.

A simple (although not always the best) way to model the compactness assumption in Example 1.1 is to assign

$$\theta_{uv}(s,t) = \begin{cases} 0, & s = t \\ \alpha, & s \neq t \end{cases} \qquad (1.1)$$

for any pair of labels $(s,t) \in \mathcal{Y}_u \times \mathcal{Y}_v$ with some $\alpha > 0$. A simple way to model a smooth surface in depth reconstruction in Example 1.2 is to assign

$$\theta_{uv}(s,t) = |s - t| \,, \qquad (1.2)$$

to penalize large differences between depth in the neighboring nodes.

In the cell tracking example the pairwise costs should forbid the same labels to be assigned to neighboring nodes when no cell division happens:

$$\theta_{uv}(s,t) = \begin{cases} 0, & s \neq t \\ \infty, & s = t \,. \end{cases} \qquad (1.3)$$

This disallows that cells $u$ and $v$ "glue" to the same "parent" cell $s = t$. In case of cell division, this pairwise cost function can be extended in a natural way to disallow intersection of cell descendants.

These examples show that pairwise costs often incorporate the prior information about a considered object, therefore, they are often collectively referred to as the *prior*. However, this is not always the case. For instance, much better segmentation results can be obtained if the parameter $\alpha$ in (1.1) depends on the color distribution of the input image, i.e. on $uv$.

Costs and cost functions are also called *potentials* and *potential functions*. We prefer the term *cost* since it is more widely used in general optimization literature.

Since unary and pairwise costs are functions of discrete variables, they can be seen as vectors. Therefore we can treat the unary cost function $\theta_u$ as a *unary cost* vector $(\theta_u(s),\ s \in \mathcal{Y}_u)$. Similar reasoning holds also for each pairwise cost function, which can be considered as a *pairwise cost* vector $\theta_{uv} = (\theta_{uv}(s,t),\ (s,t) \in \mathcal{Y}_u \times \mathcal{Y}_v)$. Unless we use the word *vector* or *function*, the context will determine whether we refer to a vector or a function $\theta_u$ (or $\theta_{uv}$). All unary vectors stacked together form the vector of all unary costs $\theta_{\mathcal{V}} = (\theta_u,\ u \in \mathcal{V})$. The vector $\theta_{\mathcal{E}}$ of all pairwise costs is defined similarly as $(\theta_{uv},\ uv \in \mathcal{E})$. Stacking together the latter two results in a long *cost vector* $\theta = (\theta_{\mathcal{V}}, \theta_{\mathcal{E}})$ with dimension $\mathcal{I} := \sum_{u \in \mathcal{V}} |\mathcal{Y}_u| + \sum_{uv \in \mathcal{E}} |\mathcal{Y}_{uv}|$.

**Labeling** In the following, we will often use the notation $\mathcal{Y}_{\mathcal{A}}$ for all possible label assignments to a subset of nodes $\mathcal{A} \subseteq \mathcal{V}$. Formally, $\mathcal{Y}_{\mathcal{A}}$ stands for the Cartesian product $\prod_{u \in \mathcal{A}} \mathcal{Y}_u$. In particular, $\mathcal{Y}_{uv}$ denotes $\mathcal{Y}_u \times \mathcal{Y}_v$ and is the set of all possible pairs of labels in nodes $u$ and $v$. A vector $y \in \mathcal{Y}_{\mathcal{V}}$ of labels assigned to *all* nodes of the graph is called *labeling*. We will refer to coordinates of this vector with the node index, i.e. $y_u$ stands for the label assigned to the node $u$. One may also speak about *partial labelings*, if only a subset $\mathcal{A}$ of the nodes is labeled.

**Definition 1.4** (Graphical model). The triple $(\mathcal{G}, \mathcal{Y}_{\mathcal{V}}, \theta)$ consisting of a graph $\mathcal{G}$, discrete space of all labelings $\mathcal{Y}_{\mathcal{V}}$ and a corresponding cost vector $\theta$, is called *a graphical model*.

**Definition 1.5** (Energy minimization problem). The problem

$$y^* = \underset{y \in \mathcal{Y}_{\mathcal{V}}}{\arg\min} \left[ E(y; \theta) := \sum_{u \in \mathcal{V}} \theta_u(y_u) + \sum_{uv \in \mathcal{E}} \theta_{uv}(y_u, y_v) \right] \qquad (1.4)$$

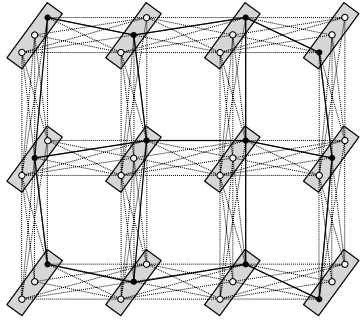

**Figure 1.2:** Labeling of the graphical model from Figure 1.1. Selected labels are marked as black circles and connected with solid lines. Each black circle corresponds to a unary cost and each solid line to a pairwise cost in the sum in the energy minimization problem (1.4).

of finding a labeling $y^*$ with minimal total cost will be called *energy minimization* or *maximum a posteriori (MAP) inference* problem for the graphical model $(\mathcal{G}, \mathcal{Y}_\mathcal{V}, \theta)$.

For the sake of notation we will sometimes use the short form of (1.4)

$$y^* = \underset{y \in \mathcal{Y}_\mathcal{V}}{\arg\min} \left[ E(y; \theta) := \sum_{w \in \mathcal{V} \cup \mathcal{E}} \theta_w(y_w) \right] \qquad (1.5)$$

with $y_w$ being equal to $y_u$, if $w$ corresponds to a node, i.e. $w = u \in \mathcal{V}$, and $y_{uv}$, if $w$ corresponds to an edge, i.e. $w = uv \in \mathcal{E}$.

Problems equivalent or very closely related to (1.4) have also other names depending on the corresponding community they are studied in: *maximum likelihood explanation (MLE) inference* (machine learning, natural language processing community), *weighted/valued/partial constraint satisfaction problem* (constraint satisfaction community).

## 1.2   Probabilistic interpretation

The name *MAP-inference* stems from the probabilistic interpretation of the problem (1.4). With the energy $E(y; \theta)$ one typically associates the exponential probability distribution

$$p(y) = \frac{1}{Z(\theta)} \exp\left(-E(y; \theta)\right), \qquad (1.6)$$

where the normalizer $Z(\theta)$ is known as *partition function*.

According to the distribution (1.6), problem (1.4) is equivalent to finding *the most probable* labeling $y$, i.e. the one maximizing $p(y)$. Since $E$ has the separable form (1.5) the expression (1.6) takes the form of the *Gibbs distribution*

$$p(y) = \frac{1}{Z(\theta)} \prod_{w \in \mathcal{V} \cup \mathcal{E}} \Theta_w(y_w) \qquad (1.7)$$

with $\Theta_w = \exp(-\theta_w)$. This explains the also frequently used name "factors" for the cost functions and their exponents $\Theta_w$.

The probabilistic interpretation (1.6) gives rise to several other *probabilistic inference* problems motivated by Bayesian statistical decision making theory. One computational problem, often referred to as *marginalization inference*, consists of computing *marginal* distributions

$$\hat{p}_u(s) := \sum_{y \in \mathcal{Y}_\mathcal{V} \,:\, y_u = s} p(y) \qquad (1.8)$$

for each node $u$ and label $s$ of a graphical model. These kinds of problems, although closely related to MAP-inference, are beyond the scope of this monograph.

## 1.3 Combinatorial complexity of MAP-inference

The number of possible labelings $y$ in (1.4) grows exponentially with the cardinality of $\mathcal{V}$, as it is equals $\prod_{v \in \mathcal{V}} |\mathcal{Y}_v|$. It results in $L^{|\mathcal{V}|}$ in case all nodes have the same number of labels $|\mathcal{Y}_u| = L$, $\forall u \in \mathcal{V}$.

However, an exponentially large set of solutions is not sufficient for polynomial $\mathcal{NP}$-hardness of a problem. For example, the shortest path between two nodes in a directed graph with positive edge weights has an exponentially large set of solutions, but is polynomially solvable by Dijkstra's algorithm.

Below we show that the MAP-inference (1.4) is, indeed, $\mathcal{NP}$-hard. To do so, it is sufficient to show that some $\mathcal{NP}$-complete decision problem is polynomially reducible to MAP-inference.

In the following construction we will show that the Hamiltonian cycle problem reduces to MAP-inference in polynomial time.

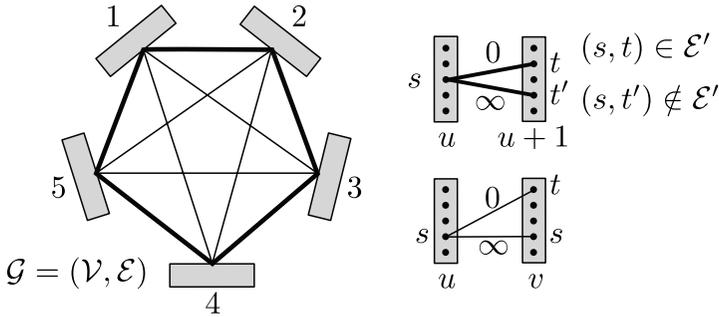

**Figure 1.3:** Illustration of the reduction of the Hamiltonian cycle problem to MAP-inference for a graph with 5 nodes. Edges of the MAP-inference graph $\mathcal{G}$ are divided into two groups: between nodes $u$ and $u+1$ (bold edges) and all others (thin edges). The corresponding pairwise costs are illustrated on the right (see also main text).

**Definition 1.6** (Hamiltonian cycle). A *Hamiltonian cycle* in a graph $\mathcal{G}$ is a cycle which visits each node exactly once.

The problem of deciding whether a given directed graph has a Hamiltonian cycle is known to be $\mathcal{NP}$-complete. To show $\mathcal{NP}$-hardness of MAP-inference, it is sufficient to reduce the Hamiltonian cycle problem to it.

Let $\mathcal{G}' = (\mathcal{V}', \mathcal{E}')$ be the graph for which one should solve the Hamiltonian cycle problem. Let us construct the following graphical model (see Figure 1.3): For the graph $\mathcal{G} = (\mathcal{V}, \mathcal{E})$ it holds that $\mathcal{V} = \mathcal{V}'$ and $\mathcal{E} = \binom{\mathcal{V}}{2}$. In other words, graph $\mathcal{G}$ contains the same nodes as graph $\mathcal{G}'$ and is fully connected. Moreover, we will order all nodes of the graph $\mathcal{G}$, i.e. $\mathcal{V} = \{1, 2, \ldots, |\mathcal{V}|\}$. This order is the order of nodes in the Hamiltonian cycle we are searching for. We will assume the operation $u+1$ to be defined modulo $|\mathcal{V}|$, i.e. $u+1$ defines the next element of the Hamiltonian cycle. In other words, if $u < |\mathcal{V}|$ then $u+1$ is the next natural number after $u$ and for $u = |\mathcal{V}|$ the element $u+1$ is equal to 1.

The set of labels $\mathcal{Y}_u := \mathcal{V}'$ is the same for each node $u \in \mathcal{V}$. Its elements index nodes of the graph $\mathcal{G}'$. A label $s$ assigned to a node $u \in \mathcal{V}$ encodes that the $u$-th node in the Hamiltonian cycle corresponds to the node $s$ of the graph $\mathcal{G}'$.

Unary costs are equal to 0. Pairwise costs are split into two groups. For a pair of nodes $\{u, u+1\} \in \mathcal{E}$ the cost reads

$$\theta_{u,u+1}(s,t) = \begin{cases} 0, & (s,t) \in \mathcal{E}' \\ \infty, & (s,t) \notin \mathcal{E}' . \end{cases} \tag{1.9}$$

It guarantees that two neighboring nodes of the Hamiltonian cycle are connected by an edge in the graph $\mathcal{G}'$.

To guarantee that no node is included twice in the Hamiltonian cycle, we set up other pairwise costs for $v \neq u+1$ and $u \neq v+1$ as follows:

$$\theta_{uv}(s,t) = \begin{cases} 0, & s \neq t \\ \infty, & s = t \end{cases} . \tag{1.10}$$

Such type of pairwise costs is sometimes called the *uniqueness constraints*, since these costs enforce that each label is selected at most ones.

Let $y$ be some labeling of the graphical model $\mathcal{G}$ such that $E(y, \theta) < \infty$. Then the sequence $(y_1, y_2, \ldots, y_{|\mathcal{V}|})$ is the Hamiltonian cycle by construction: there is an edge between $y_u$ and $y_{u+1}$ in $\mathcal{G}'$, and the set $\{y_1, y_2, \ldots, y_{|\mathcal{V}|}\}$ is exactly the set $\mathcal{V}'$.

All labelings have either value 0 or $\infty$. Therefore, the solution of the MAP-inference problem answers the question whether there is a labeling $y$ such that $E(y, \theta) < \infty$, and, therefore, whether there is a Hamiltonian cycle in the graph $\mathcal{G}'$.

Note that the same reduction of the Hamiltonian cycle problem to the MAP-inference could have also be done without using the infinite costs. Instead, any positive finite cost (e.g. 1) could be used in place of infinities. In this case the solution of the MAP-inference problem answers the question whether there is a labeling $y$ such that $E(y, \theta) = 0$, which is equivalent to the existence of a Hamiltonian cycle in the graph $\mathcal{G}'$.

## 1.4 Bibliography and further reading

For further examples of applications of graphical models in computer vision and image processing we refer to the collection [10]. Books [77, 47] can be recommended to learn more about the probabilistic view on graphical models. The monograph [135] concentrates on the exponential family (1.6) of distributions and its relation to graphical models.

A classical source to learn about the computational complexity of combinatorial problems is [27], a modern exposition is given in [5]. The most recent and comprehensive analysis of complexity of the MAP-inference problem is provided in [70].

The text books [25] and [103] can be recommended to learn about Bayesian decision theory.

The reduction of the Hamiltonian cycle problem to MAP-inference is reproduced from the lectures on structural pattern recognition given by Prof. Michail Schlesinger at National Technical University of Ukraine "Igor Sikorsky Kyiv Polytechnic Institute" where the author studied mathematics and computer science.

# 2

---

## Acyclic Graphical Models

---

In this chapter we will concentrate on an important subclass of graphical models, namely those having an acyclic structure, i.e. whose graph $\mathcal{G}$ is acyclic. This type of problem is important by itself, because it naturally allows us to model (i) time-series, as chain-structured graphs, and (ii) hierarchies, by tree-structured graphical models. As we show below, the MAP-inference problem on acyclic graphs can be solved very efficiently. Therefore, this type of model is also widely used as a building block for approximate algorithms. We will consider one class of such approximations in Chapters 9 and 10.

The algorithm that allows for an efficient solution of the MAP-inference problem for acyclic graphical models is the dynamic programming algorithm. Most importantly, for acyclic models the computational complexity of this algorithm grows only linearly with the problem size which is the best possible complexity one can expect.

The chapter is split into four parts. The first two are aggregated into a single section, and introduce the dynamic programming algorithm for chain- and general tree-structured graphs, respectively. The third part covers a related topic, the computation of the so-called *min-marginals* for chain-structured graphical models. This part shows its real importance when approximate solvers as in Chapter 10 are constructed. There the algorithm introduced below plays a crucial role. The fourth part comments on the use of dynamic programming for models whose graphs contain cycles.

## 2.1   MAP-Inference with dynamic programming

### 2.1.1   Chain-structured graphical models

To simplify our presentation we start with the special case where the graph $\mathcal{G}$ is chain-structured. In this case the set of nodes can be totally ordered, i.e. $\mathcal{V} = \{1, \ldots, n\}$. The set of edges contains pairs of neighboring nodes in the chain: $\mathcal{E} = \{(i, i+1) \colon i = 1, \ldots, n-1\}$.

Consider Figure 2.1, where a chain-structured graphical model is depicted. An optimal labeling starts in the first (left-most) node and ends in the last (right-most) one. To select its last label optimally assume we have computed the function $F_n \colon \mathcal{Y}_n \to \mathbb{R}$ such that $F_n(s) + \theta_n(s)$ denotes the cost of the best (minimal-cost) labeling with the very last label being $s$. In other words, the value $F_n(s)$ is the cost of the labeling without the unary cost $\theta_n(s)$. Then the minimization $y_n = \arg\min_{s \in \mathcal{Y}_n}(F_n(s) + \theta_n(s))$, which can be performed by enumeration, allows us to determine the last label of an optimal labeling $y$.

To find other labels the same considerations can be applied given that the functions $F_i \colon \mathcal{Y}_i \to \mathbb{R}$ are computed such that $F_i(s)$ equals the cost of the best partial labeling in the nodes $1, \ldots, i$ with the label $s$ assigned to the $i$-th node .

Functions $F_i$ can be recursively computed as

$$F_i(s) := \min_{t \in \mathcal{Y}_{i-1}} \left( F_{i-1}(t) + \theta_{i-1}(t) + \theta_{i-1,i}(t, s) \right) , \qquad (2.1)$$

with $F_1(s) = 0$ for all $s \in \mathcal{Y}_1$.

These considerations can be derived from the following recursive representation of the energy of a chain-structured graphical model, where functions $F_i \colon \mathcal{Y}_i \to \mathbb{R}$ for $i \in \mathcal{V}$ are introduced "on the fly":

$$\min_{y \in \mathcal{Y}_{\mathcal{V}}} E(y; \theta) = \min_{y_1, \ldots, y_n} \left( \sum_{i=1}^{n-1} \left( \theta_i(y_i) + \theta_{i,i+1}(y_i, y_{i+1}) \right) + \theta_n(y_n) \right)$$

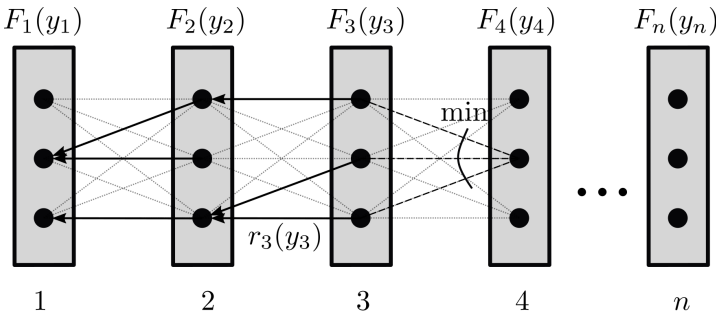

**Figure 2.1:** Illustration of dynamic programming on a chain: recursive expression (2.1) and Algorithm 1. Dash lines between nodes 4 and 3 connect the label $y_4$, which is currently computed, with all possible labels $y_3 \in \mathcal{Y}_3$. Minimization is performed over $y_3$ according to $F_4(y_4) = \min_{y_3 \in \mathcal{Y}_3} \left( F_3(y_3) + \theta_3(y_3) + \theta_{3,4}(y_3, y_4) \right)$. The arrows point to the result of this optimization during the previous steps of the algorithm and correspond to the values of $r_i(y_i)$. To reconstruct an optimal labeling one has to follow the arrows back from the last node to the first one. The starting label in the last node is chosen as $y_n = \min_{s \in \mathcal{Y}_n}(F_n(s) + \theta_n(s))$.

$$
\begin{aligned}
&= \min_{y_2,\ldots,y_n} \left( \underbrace{\min_{y_1 \in \mathcal{Y}_1} \left( \theta_1(y_1) + \theta_{1,2}(y_1, y_2) \right)}_{F_2(y_2)} \right. \\
&\quad \left. + \sum_{i=2}^{n-1} \left( \theta_i(y_i) + \theta_{i,i+1}(y_i, y_{i+1}) \right) + \theta_n(y_n) \right) \\
&= \min_{y_3,\ldots,y_n} \left( \underbrace{\min_{y_2 \in \mathcal{Y}_2} \left( F_2(y_2) + \theta_2(y_2) + \theta_{2,3}(y_2, y_3) \right)}_{F_3(y_3)} \right. \\
&\quad \left. + \sum_{i=3}^{n-1} \left( \theta_i(y_i) + \theta_{i,i+1}(y_i, y_{i+1}) \right) + \theta_n(y_n) \right) \\
&\quad \ldots \\
&= \min_{y_n \in \mathcal{Y}_n} \left( F_n(y_n) + \theta_n(y_n) \right).
\end{aligned}
\tag{2.2}
$$

The transformations (2.2) confirm that the energy of the optimal labeling can be computed as $\min_{y \in \mathcal{Y}_\mathcal{V}} E(y) = \min_{s \in \mathcal{Y}_n}(F_n(s) + \theta_n(s))$ and values $F_i(s)$ can be computed recursively as given by (2.1). We summarize these observations in Algorithm 1.

---

**Algorithm 1** Dynamic programming algorithm for chain-structured graphical models

---

1: Given: $(\mathcal{G}, \mathcal{Y}_{\mathcal{V}}, \theta)$ - a chain-structured graphical model
2: $F_1(s) := 0$, for all $s \in \mathcal{Y}_1$
3: **for** $i = 2$ **to** $n$ **do**
4:     $F_i(s) := \min_{t \in \mathcal{Y}_{i-1}} (F_{i-1}(t) + \theta_{i-1}(t) + \theta_{i-1,i}(t,s))$ for all $s \in \mathcal{Y}_i$

5:     $r_i(s) := \arg\min_{t \in \mathcal{Y}_{i-1}} (F_{i-1}(t) + \theta_{i-1}(t) + \theta_{i-1,i}(t,s))$
6: **end for**
7: $E^* = \min_{s \in \mathcal{Y}_n}(F_n(s) + \theta_n(s))$
8: **return** $E^*$, $\{r_i(s) \colon i = 2, \ldots, n, \ s \in \mathcal{Y}_i\}$

---

Values $F_i(s)$ are often called *forward min-marginals* or *forward messages*, since they represent an optimal labeling from the first node to the label $s$ in the node $i$, and all unary/pairwise costs with indexes smaller than $i$ are *marginalized out*.

Values $r_i(s)$ are pointers needed to recover the optimal labeling, when its energy is computed, as given by Algorithm 2. Note that $F_i(s)$ and $r_i(s)$ are obtained by the same minimization, therefore it must be performed only once.

---

**Algorithm 2** Reconstructing an optimal labeling

---

1: Given: $r_i(s) \colon i = 2, \ldots, n, \ s \in \mathcal{Y}_i$
2: $y_n = \arg\min_{s \in \mathcal{Y}_n}(F_n(s) + \theta_n(s))$
3: **for** $i = n$ **to** $2$ **do**
4:     $y_{i-1} = r_i(y_i)$
5: **end for**
6: **return** $y$

---

**Computation complexity** of each minimization in line 4 of Algorithm 1 is $O(|\mathcal{Y}_{i-1}|)$. One has to perform $|\mathcal{Y}_i|$ such minimizations on each iteration. Therefore, complexity of the algorithm reads

$$O\left(\sum_{i=2}^{n} |\mathcal{Y}_{i-1}||\mathcal{Y}_i|\right),$$

which results in $O(L^2|\mathcal{V}|) = O(L^2|\mathcal{E}|)$ if we assume that all nodes have an equal number $L$ of labels. This is exactly the memory, which is required to define a chain graphical model: $|\mathcal{V}|$ unary cost vectors of the size $L$ each plus $|\mathcal{E}| - 1$ pairwise costs of the size $L^2$ each. In other words, dynamic programming for chain-structured graphical models has linear complexity with respect to the model size.

The functions $F_i$ are known as *Bellman* functions, in honor of the author of dynamic programming. Besides the original name, Algorithm 1 is often refereed to as the *Viterbi algorithm*.

Algorithm 1 can also be run the opposite direction. In this case one speaks about *backward min-marginals* recursively defined as

$$\begin{aligned} B_n(s) &= 0, \;\; s \in \mathcal{Y}_n; \\ B_i(s) &:= \min_{t \in \mathcal{Y}_{i+1}} \left( B_{i+1}(t) + \theta_{i+1}(t) + \theta_{i,i+1}(s,t) \right), \;\; s \in \mathcal{Y}_i \end{aligned} \tag{2.3}$$

for $i = n-1, \ldots, 1$. Backward min-marginals can be computed similarly to forward min-marginals, by processing the nodes of a chain graph in an inverse order. In §2.2 we show why computing both forward and backward min-marginals can be useful.

### 2.1.2 Tree-structured graphical models

Dynamic programming is also applicable to MAP-inference in general acyclic graphs. Without loss of generality we will assume the graph to be connected, i.e. it must represent a single tree. For graphs with multiple connected components the MAP-inference problem decomposes into independent subproblems for each component.

Two important properties of a tree-structured graph are:

- If a tree has more than one node, there always exists a node $u$ which has a single incident edge, i.e. $|\mathcal{N}_b(u)| = 1$. Such nodes are called *leaves* of the tree.

- Removing any leaf and its incident edge turns a tree-structured graph again into a tree-structured graph.

Note that a chain is a tree, and in Algorithm 1 we remove iteratively the leave with the smallest index. This observation leads to the generalization of Algorithm 1 to arbitrary tree-structured graphs. For convenience we slightly change the definition of functions $F_v$ compared to the case of chain-structured models.

---

**Algorithm 3** Dynamic programming algorithm for tree-structured graphical models

---

1: Given: $(\mathcal{G}, \mathcal{Y}_\mathcal{V}, \theta)$ - a tree-structured graphical model
2: $F_v(s) := \theta_v(s)$, for all $v \in \mathcal{V}, \ s \in \mathcal{Y}_v$
3: **while** $|\mathcal{V}| > 1$ **do**
4:      Find a leaf $v \in \mathcal{V} \colon \mathcal{N}_b(v) = \{u\}$ for some $u \in \mathcal{V}$
5:      $F_u(s) := F_u(s) + \min_{t \in \mathcal{Y}_v} (F_v(t) + \theta_{uv}(s,t))$
6:      $\mathcal{V} := \mathcal{V} \backslash \{v\}$
7: **end while**
8: Let $v$ be the remaining node, i.e. $\mathcal{V} = \{v\}$
9: **return** $E^* = \min_{s \in \mathcal{Y}_v} F_v(s).$

---

Algorithm 3 iteratively searches for a leave node $v$ and eliminates it by recomputing $F_u$ for its single neighbor $u$. To prove correctness of the algorithm its derivation in an algebraic form analogous to (2.2) can be constructed.

On each iteration of Algorithm 3 one edge $uv$ (incident to a leave node $v$) is eliminated in line 5. This elimination has a complexity of $O(|\mathcal{Y}_u||\mathcal{Y}_v|)$. It has to be performed for all edges, which results in the overall complexity of Algorithm 3 being $O(\sum_{uv \in \mathcal{E}} |\mathcal{Y}_u||\mathcal{Y}_v|)$. The latter simplifies to $O(|\mathcal{E}|L^2)$ in case all nodes are associated with the same number of labels $L$. Similarly to the case of a chain-structured model, the complexity is linear in the model size.

In the same way as Algorithm 1, Algorithm 3 can be augmented with additional information to reconstruct an optimal labeling when the optimal energy is found.

## 2.2 Computation of min-marginals

Besides computing an optimal labeling for the problem as a whole, dynamic programming can be very efficiently applied to compute the energy of an optimal labeling containing a specified label in a particular graph node. Moreover, computing energies of *all such labelings*, i.e. for all labels in all nodes at once, has basically the same complexity as the MAP-inference problem. Here we will again assume that the graphical model is chain-structured, although the computations generalize straightforwardly to any acyclic graph.

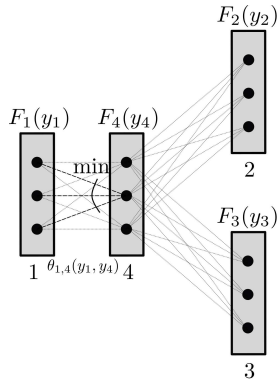

**Figure 2.2:** Illustration to Algorithm 3 of dynamic programming on a tree. Nodes are enumerated in a possible processing order. Dash lines between nodes 4 and 1 connect the label $y_4$, which is currently computed, with a possible label $y_1 \in \mathcal{Y}_1$. Minimization is performed over $y_1$ according to $F_4(y_4) = F_4(y_4) + \min_{y_1 \in \mathcal{Y}_1} \left( F_1(y_1) + \theta_{1,4}(y_1, y_4) \right)$.

**Definition 2.1.** Let $E_w : \mathcal{Y}_w \to \mathbb{R}$ be defined for all $w \in \mathcal{V} \cup \mathcal{E}$ as

$$E_w(s) = \min_{y \in \mathcal{Y}_\mathcal{V} : y_w = s} E(y; \theta), \qquad (2.4)$$

that is, $E_w(s)$ is the energy of the optimal labeling containing label (label pair) $s$ in node (edge) $w$. Values $E_w(s)$ are called *min-marginals*. We will distinguish *node min-marginals* $E_v$, $v \in \mathcal{V}$, and *edge min-marginals* $E_{uv}$, $uv \in \mathcal{E}$.

Using the definition of $F_i$ as in Algorithm 1, and $B_i$ as in (2.3), it is easy to see that

$$E_i(s) = F_i(s) + B_i(s) + \theta_i(s), \ i \in \mathcal{V}, \ s \in \mathcal{Y}_i, \qquad (2.5)$$

$$E_{i,i+1}(s,t) = F_i(s) + B_{i+1}(t) + \theta_i(s) + \theta_{i+1}(t) + \theta_{i,i+1}(s,t),$$
$$i = 1, \dots, n-1, \ (s,t) \in \mathcal{Y}_{i,i+1}. \qquad (2.6)$$

The algorithm, which computes min-marginals by first computing functions $F_i$ and $B_i$ and then applying (2.5) and/or (2.6), is often referred to as *forward-backward dynamic programming* or *forward-backward belief propagation* on a chain. The name comes from the fact that one has to perform computations in the order of increasing node indexes for $F_i$ ("forward move of dynamic programming") and in order of decreasing node indexes for $B_i$ ("backward move of dynamic programming").

## 2.3 On dynamic programming for cyclic graphs

**Graphs with low tree-width** The dynamic programming algorithm can be extended to cyclic graphical models. However, the complexity of the algorithm grows exponentially with the so called *tree width* of the graph. Loosely speaking, nodes and edges of the initial graph must be grouped together into "hyper-nodes" and "hyper-edges" in such a way that the "hyper-graph" defined by this grouping is acyclic. Tree-width is defined as the maximal number of nodes in a "hyper-node" of this graph. Unfortunately, an $n \times n$ grid has tree-width $n$, and so has the fully-connected graph with $n$ nodes. This implies an exponential complexity of the dynamic programming algorithm for the corresponding graphical models.

Therefore, typically dynamic programming can be used to simplify the problem by reducing it to a subgraph with high tree-width. To deal with this remaining problem, alternative approaches must be considered, which we will explore in later chapters of this monograph.

**Belief propagation** Another, heuristic, way to use the dynamic programming algorithm for cyclic graphs is to define an order on the set of its nodes and iteratively perform dynamic programming updates as in Algorithm 1 according to this order. This type of algorithm is usually referred to as *(loopy min-sum) belief propagation*. However, such algorithms have very weak convergence properties, they do not converge in general, and lack optimality guarantees even when they converge.

In Chapters 9 and 10 we will follow the more theory based Lagrangean decomposition approach which splits the graphical model into acyclic subgraphs, and enforces consistency of their solutions. Interestingly, some of the resulting algorithms turn out to be very similar to belief propagation, but contrary to it they (i) enjoy much stronger convergence guarantees, (ii) provide optimality bounds for the optimal energy, and (iii) perform better in numerous applications.

## 2.4    Bibliography and further reading

A historical reference for dynamic programming is the book [7]. Andrew Viterbi proposed a method for finding the most probable sequence of a hidden Markovian chain (which can be seen as a chain graphical model) in his seminal work [132]. Acyclic graphical models and generalizations of dynamic programming are presented in great detail in the excellent text-book [103] which we recommend for further reading on the topic. The definition and a detailed overview of tree-width computation is given in [11, 12], and the corresponding generalization of dynamic programming is presented in [67, p.228] as *message passing in junction trees*. Belief propagation was introduced in [80]. The work [63] is also widely cited for this algorithm.

# 3

---

# Background: (Integer) Linear Programs and Their Geometry

---

In the following, we will often treat MAP-inference for graphical models as an *integer linear program*. The theory of this mathematical object is well-developed, and plays an important role for constructing exact and approximative optimization algorithms. Therefore, we will give basic concepts and results of this theory which are necessary to describe the most efficient existing optimization methods for the considered MAP-inference problem.

We start this chapter with a general notion of a (constrained) optimization problem and its relaxations. Further, we will focus on a specific, polyhedral constraints, which will lead us naturally to the linear programs. Finally, we will consider linear programs with integer constraints, known as integer linear programs.

## 3.1 Optimization problems

In the optimization literature the notation

$$\min_{x \in X} f(x) \tag{3.1}$$

is adopted for a problem of *minimizing* a function $f \colon X \to \mathbb{R}$ on a set $X \subseteq \mathbb{R}^n$. More precisely, it means that *the optimal value* of $f$ defined as $f^* = \inf_{x \in X} f(x)$ must be found. If $f^* = -\infty$, that is, there exists a sequence $x^t \in X$ such that $f(x^t) \stackrel{t \to \infty}{\Longrightarrow} -\infty$, the problem is called *unbounded*. Another special case appears when $X$ is an empty set. Then the problem is called *infeasible*, the notation $f^* = \infty$ is adopted to this case. The set $X$ is called the *feasible set* of the problem (3.1), and $x \in X$ is a *feasible point*. The fact that the problem (3.6) is neither infeasible nor unbounded, in other words, *feasible* and *bounded*, is denoted as $-\infty < f^* < \infty$. A vector $x^* \in X$ such that $f(x^*) = f^*$ is called an *(optimal) solution* or *optimal point* of the problem. In general, an optimal solution is not unique, or may not exist, even if $-\infty < f^* < \infty$.

**Example 3.1.** Although the optimal values of the problems

$$\min_{x \geq 1} e^{\frac{1}{x}}; \quad \min_{x > 1} x \tag{3.2}$$

are finite and equal to 1, their optimal points do not exist.

Indeed, assume that a feasible point $x \geq 1$ is an optimal solution for the first problem. However, since $e^{\frac{1}{x+1}} < e^{\frac{1}{x}}$ and $x + 1 \geq 1$ is feasible, the point $x$ is not an optimal solution.

Similarly for the second problem, its objective value in $(1 + \frac{x-1}{2})$ is strictly smaller than in $x$ for any feasible $x$, and $1 + \frac{x-1}{2}$ is feasible as soon as $x$ is feasible.

In the following, we will mostly deal with problems where the existence of an optimal value implies also the existence of an optimal point.

A feasible point $x$ such that $f(x)$ is approximately equal to $f^*$ (denoted as $f(x) \approx f^*$) is often referred to as an *approximate solution* of the minimization problem (3.1).

The problem of the form (3.1) is referred to as a *constraint minimization* problem. The *constrained maximization* problem $\max_{x \in X} f(x)$ is defined by substituting min with max, inf with sup, and $-\infty$ with $\infty$ and the other way around. Minimization and maximization problems together are also called *optimization* problems.

**Example 3.2.** The following optimization problems are infeasible, as their respective feasible sets are empty:

$$\min_{\substack{x \geq 1 \\ x \leq 0}} f(x); \quad \max_{\substack{x \in \{0,1\} \\ x \geq 2}} f(x); \quad \min_{\substack{x \in \{0,1\} \\ x \in [0.2, 0.9]}} f(x). \tag{3.3}$$

**Example 3.3.** The following optimization problems are unbounded:

$$\min_{x \le 0} x + 1; \qquad \min_{x \ge 0} 2 - x; \qquad \max_{x \le 0} x^2 . \tag{3.4}$$

**Example 3.4.** The following optimization problems have multiple optimal solutions

$$\min_{x \in [-1,1]} 1 - x^2; \qquad \min_{x \in \{0,1\}} |x - 0.5|; \qquad \max_{x \in [0,1]} 1 . \tag{3.5}$$

In the rest of the monograph we will deal mostly with *linear objectives*, i.e. $f(x) = \langle c, x \rangle$ for some vector $c \in \mathbb{R}^n$. The feasible set $X$ will usually be either finite or *polyhedral*. The definition of the latter will be given in §3.2.

### 3.1.1 Relaxations of optimization problems

As we already learned in the example for the MAP-inference problem in Chapter 1, important optimization problems are often computationally intractable. Therefore, it is natural to consider tractable approximations of these problems to obtain an approximate solution of the problem or at least a bound on its optimal value. In the following, we define one of the most important types of such approximations:

**Definition 3.5.** The optimization problem

$$\min_{x \in X'} g(x) \tag{3.6}$$

is called a *relaxation* of the problem (3.1) if $X' \supseteq X$ and $g(x) \le f(x)$ for any $x \in X$. The solution of the relaxed problem is often referred to as a *relaxed solution*.

**Example 3.6.** The problems $\min_{x \ge 0} f(x)$ and $\min_{x \in [0,1]^n} f(x)$ are relaxations of the problem $\min_{x \in \{0,1\}^n} f(x)$.

**Example 3.7.** The problem $\min_{x \in [0,1]} x^2$ is a relaxation of the problem $\min_{x \in [0,1]} |x|$. The advantage of this relaxation is the differentiability of $x^2$, which may often lead to faster optimization algorithms.

**Example 3.8.** The problems $\min_{x \in X \subset \mathbb{R}^n} \langle c, x \rangle$ and $\min_{x \,:\, Ax = b} \langle c, x \rangle$ are relaxations of $\min_{\substack{x \in X \subset \mathbb{R}^n \\ Ax = b}} \langle c, x \rangle$.

The most important property of a relaxation is provided by the following proposition:

**Proposition 3.9.** Any relaxation constitutes a lower bound, i.e.

$$\min_{x \in X'} g(x) \leq \min_{x \in X} f(x),$$

if $X \subseteq X'$ and $g(x) \leq f(x)$ for $x \in X$. Moreover, if $x' \in X$ and $f(x') = g(x')$ holds for a relaxed solution $x' \in \arg\min_{x \in X'} g(x)$, then $x'$ is an optimal solution of the non-relaxed problem $\min_{x \in X} f(x)$ as well. In this case, one says that the lower bound provided by the relaxation is *tight*.

*Proof.* The first statement of the proposition is implied by the following inequalities $\min_{x \in X'} g(x) \leq \min_{x \in X} g(x) \leq \min_{x \in X} f(x)$, which follow from Definition 3.5.

The second statement follows from $f(x') = g(x') = \min_{x \in X'} g(x) \leq \min_{x \in X} f(x) \leq f(x')$. □

Relaxations, which provide tight lower bounds for *all* instances of optimization problems from a certain class, are called *tight relaxations*.

**Example 3.10.** Problem $\min\limits_{x \in [0,1]} cx$ is a tight relaxation for $\min\limits_{x \in \{0,1\}} cx$ for any constant $c$.

**Definition 3.11.** The relaxation $\min_{x \in X'} g(x)$ is called *tighter* than the relaxation $\min_{x \in \hat{X}'} \hat{g}(x)$, if it provides better lower bound, i.e. $\min_{x \in X'} g(x) \geq \min_{x \in \hat{X}'} \hat{g}(x)$.

The first part of Proposition 3.9, in particular, implies that the $\min_{x \in X'} g(x)$ is tighter than $\min_{x \in \hat{X}'} \hat{g}(x)$ if the latter is a relaxation of the former one.

**Definition 3.12.** Two relaxations are called *equivalent*, if their bounds coincide.

The following simple proposition shows that uniqueness of the solution of a tight relaxation implies uniqueness of the solution of the non-relaxed problem:

**Proposition 3.13.** Let $\min_{x \in X} f(x)$ be an optimization problem and $\min_{x \in X'} g(x)$ be its relaxation. If $x' \in X$ is the unique solution of the relaxed problem and $g(x') = f(x')$, then it is also the unique solution of the non-relaxed one.

*Proof.* Since $x'$ is the unique relaxed solution, it holds that $f(x') = g(x') < g(x)$ for $x \in X' \backslash \{x'\}$. Therefore, $f(x') < g(x) \le f(x)$ for all $x \in (X \backslash \{x'\}) \subseteq (X' \backslash \{x'\})$, which implies that $x'$ is the unique solution of the non-relaxed problem. $\square$

## 3.2    Linear constraints and polyhedra

In this section we consider a special type of constraint sets known as polyhedra. Polyhedra are the main components of linear and integer linear programs.

Let $a \ne \bar{0}$ be a vector in $\mathbb{R}^n$. The set $\{x \in \mathbb{R}^n \colon \langle a, x \rangle = b\}$ is called a *hyperplane* in $\mathbb{R}^n$, and $\{x \in \mathbb{R}^n \mid \langle a, x \rangle \le b\}$ is a *half-space* defined by this *hyperplane*.

**Definition 3.14.** A *polyhedron* $P$ in $\mathbb{R}^n$ is a set which can be represented by a finite set of linear inequalities, i.e. $P = \{x \in \mathbb{R}^n \colon Ax \le b\}$ for a matrix $A \in \mathbb{R}^{m \times n}$ and a vector $b \in \mathbb{R}^m$. A bounded polyhedron is called a *polytope*.

**Example 3.15.** The set $[0, 1]^n$, called an *n-dimensional cube*, is a polytope, since it can be represented as

$$\{x \in \mathbb{R}^n \mid -x_i \le 0, \ x_i \le 1, \ i = 1 \dots, n\}$$

and it is bounded.

Constraints $x \in [0, 1]^n$ constitute a special case of *box constraints*. The latter in general may specify a different feasible interval for each coordinate of $x$.

Polyhedra can be represented also by a mixture of equalities and inequalities:

**Proposition 3.16.** The set of the type $\{x \in \mathbb{R}^n \colon Ax = b, \ Bx \le d\}$, where $A$ and $B$ are arbitrary matrices of a suitable dimension, is a polyhedron in $\mathbb{R}^n$.

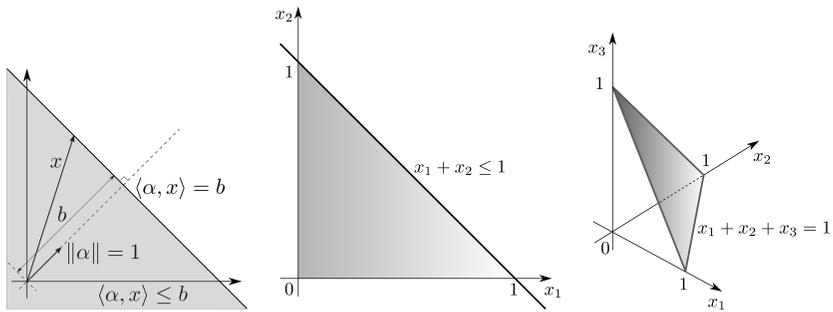

**Figure 3.1:** Examples of different polyhedra: (left) Hyperspace $\langle \alpha, x \rangle \leq b$; (middle) Polytope $\{x \geq 0, \ x_2 \geq 0, \ x_1 + x_2 \leq 1\}$; (right) Polytope (three-dimensional simplex) $\Delta^3 = \{x \in \mathbb{R}^3 : x_1 + x_2 + x_3 = 1, \ x_1, x_2, x_3 \geq 0\}$.

*Proof.* The proof follows from the fact that the equality constraints $Ax = b$ can be equivalently represented as inequalities:

$$\{Ax \leq b, \ -Ax \leq -b\}. \qquad \square$$

**Example 3.17.** A hyperplane $\langle a, x \rangle = 0$ is a polyhedron, and so is an intersection of multiple hyperplanes $Ax = b$. Indeed, as an intersection of hyperplanes $Ax = b$ defines a linear subspace, which is unbounded unless it is empty or consists of a single point.

The set of polyhedra, as well as the set of polytopes are closed with respect to intersection:

**Proposition 3.18.** The intersection of two polyhedra is a polyhedron itself. Moreover, if one of the polyhedra is a polytope, the intersection is a polytope as well.

*Proof.* The proof of the first statement follows from the fact that an intersection of two polyhedra given by $Ax \leq b$ and $Bx \leq d$, respectively, is the set $\{x \mid Ax \leq b, \ Bx \leq d\}$ constrained by both these sets of inequalities. The second statement is trivial as the intersection is a subset of the polytope at hand. $\qquad \square$

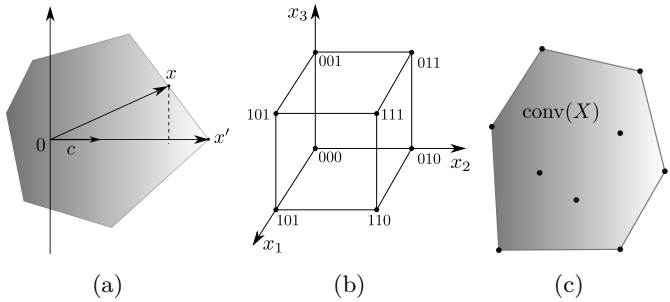

**Figure 3.2: (a)** Illustration of Definition 3.19 of a polyhedron vertex. The value of the scalar product $\langle c, x \rangle$ is the length of the orthogonal projection of a point $x$ of the polyhedron onto the direction of the vector $c$. The length of the projection of the vertex $x'$ is the longest one. **(b)** Illustration of the statement $\{0, 1\}^n = \text{vrtx}([0, 1]^n)$. **(c)** Illustration of the statement $\text{vrtx}(\text{conv}(X)) \subseteq X$. Each point in the finite set $X$ is denoted by a dot.

### 3.2.1    Convex hulls and vertices of polyhedra

There is an important one-to-one relation between a polytope and the set of its vertices, which plays an important role in the analysis of (integer) linear programs. Below, we define this relation.

**Definition 3.19.** A vector $x' \in \mathbb{R}^n$ is called a *vertex* of a polyhedron $P$ if there exists a vector $c \in \mathbb{R}^n$ such that the linear function $\langle c, x \rangle$ attains a unique maximum at $x'$ on $P$.

In other words, $x'$ is a vertex, if it is the unique solution of

$$\max_{x \in P} \langle c, x \rangle .$$

Definition 3.19 is illustrated in Figure 3.2. We will use the notation $\text{vrtx}(P)$ for the set of vertices of a polyhedron $P$.

Note that since $\max_{x \in P} \langle c, x \rangle = -\min_{x \in P}(-\langle c, x \rangle)$, the maximization problem $\max_{x \in P} \langle c, x \rangle$ can be equivalently exchanged with

$$\min_{x \in P} \langle c, x \rangle$$

in Definition 3.19.

**Example 3.20.** Let us show that $\{0,1\}^n \subseteq \text{vrtx}([0,1]^n)$. Indeed, consider $z \in \{0,1\}^n$ and a linear function $\langle c, x \rangle$ with $c_i = 2z_i - 1$. Then

$$\langle c, x \rangle = 2\sum_{i=1}^{n} z_i x_i - \sum_{i=1}^{n} x_i = 2\sum_{\substack{i=1,\ldots,n: \\ z_i=1}} x_i - \sum_{i=1}^{n} x_i$$

$$= \sum_{\substack{i=1,\ldots,n: \\ z_i=1}} x_i - \sum_{\substack{i=1,\ldots,n: \\ z_i=0}}^{n} x_i. \qquad (3.7)$$

Therefore, $\max_{x \in [0,1]^n} \langle c, x \rangle = \max_{x \in [0,1]^n} \sum_{\substack{i=1,\ldots,n: \\ z_i=1}} x_i - \sum_{\substack{i=1,\ldots,n: \\ z_i=0}} x_i = \sum_{i=1}^{n} z_i$, and the maximum is attained in the point $z$ only. Therefore, according to Definition 3.19, $z$ is a vertex of $[0,1]^n$.

Let $\mathbb{R}_+ := \{x \in \mathbb{R} \mid x \geq 0\}$ denote the set of non-negative real numbers and $\mathbb{R}_+^n = \prod_{i=1}^{n} \mathbb{R}_+$ stand for the non-negative cone of the $n$-dimensional vector space.

**Definition 3.21.** The polytope $\Delta^n := \{x \in \mathbb{R}_+^n : \sum_{i=1}^{n} x_i = 1\}$ is called *$n$-dimensional (probability) simplex.*[1] For finite $\mathcal{X}$ we also use the notation $\Delta^{\mathcal{X}}$ for vectors from $\Delta^{|\mathcal{X}|}$ whose coordinates are indexed by elements of the set $\mathcal{X}$.

A three-dimensional simplex is shown in Figure 3.1.

**Example 3.22.** Let us show that vectors $x^i \in \mathbb{R}^n$, $i = 1, \ldots, n$ such that $x_j^i = [\![ i = j ]\!]$ are vertices of the simplex $\Delta^n$. Indeed, for $c = x^i$ the maximum $\max_{x \in \Delta^n} \langle c, x \rangle$ is attained in a single point $x^i$. Therefore, $x^i \in \text{vrtx}(\Delta^n)$.

Note also that $x^i \in \mathbb{R}^n$, $i = 1, \ldots, n$, are the only binary (with coordinates 0 and 1) vectors in $\Delta^n$. In other words, $\{x^i \in \mathbb{R}^n \mid i = 1, \ldots, n\} = \Delta^n \cap \{0,1\}^n$. This implies that $\Delta^n \cap \{0,1\}^n \subseteq \text{vrtx}(\Delta^n)$. Moreover, later we will show that these sets are equal, i.e. $\Delta^n \cap \{0,1\}^n = \text{vrtx}(\Delta^n)$.

---

[1]Although dimensionality of the $n$-dimensional simplex as a manifold is $n - 1$, we stick to the above name. For the purpose of this book the dimensionality of the space plays more important role than the dimensionality of the simplex itself.

**Definition 3.23.** For any finite number of points $x^i \in \mathbb{R}^N$, $i = 1, \ldots, n$ and any $p \in \Delta^n$ the point $\sum_{i=1}^n p_i x^i$ is called *a convex combination* of $x^i$, $i = 1, \ldots, n$.

**Definition 3.24.** A set $X \subseteq \mathbb{R}^N$ is called *convex*, if for any $\alpha \in [0, 1]$ and any two points $x, z \in X$ their convex combination $\alpha x + (1 - \alpha)z$ also belongs to $X$.

It is sometimes easier to use an equivalent definition:

**Definition 3.25.** A set $X \subseteq \mathbb{R}^N$ is called *convex*, if for any finite subset $x^i \in X$, $i = 1, \ldots, n$ of its points and any $p \in \Delta^n$ it holds that $(\sum_{i=1}^n p_i x^i) \in X$.

Convexity is closed w.r.t. the intersection:

**Lemma 3.26.** Let $X$ be representable as the intersection of an arbitrary number of convex sets. Then $X$ is convex as well.

*Proof.* Let $x$ and $z$ belong to the intersection. Then $x$ and $z$ belong to each of the sets. Consider now one of these sets. Then, for any $\alpha \in [0, 1]$ the point $\hat{x} = \alpha x + (1 - \alpha)z$ belongs to this set as well due to its convexity. Since this consideration holds for any set involved in the intersection, $\hat{x}$ belongs to the intersection as well. $\qquad\square$

**Proposition 3.27.** Polyhedra are convex sets.

*Proof.* It is necessary to prove that the set $X = \{x \colon Ax \le b\}$ is convex, where $A$ is a matrix and $b$ a vector of suitable dimensions. Let $p \in \Delta^n$. The proof follows from

$$A \sum_{i=1}^n p_i x^i = \sum_{i=1}^n p_i A x^i \le \sum_{i=1}^n p_i b = b \sum_{i=1}^n p_i = b, \qquad (3.8)$$

where the inequality holds since $p_i \ge 0$, and the last equality holds due to $\sum_{i=1}^n p_i = 1$. $\qquad\square$

**Definition 3.28.** For $X \subset \mathbb{R}^n$ the set $\{s \in \mathbb{R}^n \mid \exists N \colon s = \sum_{i=1}^N p_i x^i, \ x^i \in X, \ p \in \Delta^N\}$ of points representable as a convex combination of a finite number of points of $X$ is called the *convex hull* of $X$ and will be denoted as $\mathrm{conv}(X)$.

**Exercise 3.29.** Prove that $X$ is a convex set if and only if $\operatorname{conv}(X) = X$.

**Exercise 3.30.** Prove that $\operatorname{conv}(X)$ is a minimum convex set containing $X$, that is, for any convex set $C \supseteq X$ it holds that $\operatorname{conv}(X) \subseteq C$.

**Proposition 3.31.** Let $Z \subseteq \mathbb{R}^n$ and $X \subseteq Z$. Then $\operatorname{conv}(X) \subseteq \operatorname{conv}(Z)$.

*Proof.* Any $x \in \operatorname{conv}(X)$ is representable as $\sum_{i=1}^{N} p_i x^i$ with some natural $N$, $p \in \Delta^N$ and $x^i \in X$. Since $x^i \in Z$, the point $x = \sum_{i=1}^{N} p_i x^i$ also belongs to $\operatorname{conv}(Z)$. $\qquad\square$

**Example 3.32.** Let us show that simplex $\Delta^n$ is a convex hull of vectors $\{x^i \in \mathbb{R}^n \mid i = 1, \ldots, n, \ x^i_j = [\![i = j]\!]\} = \Delta^n \cap \{0,1\}^n$.

Indeed, any vector $p \in \Delta^n$ is representable as $\sum_{i=1}^{n} p_i x^i$. The other way around, all vectors of the form $\sum_{i=1}^{n} p_i x^i$ belong to $\Delta^n$.

A simple but important technical lemma simplifies the computation of convex hulls in multidimensional spaces:

**Lemma 3.33.** Let $X = \prod_{i=1}^{n} X^i$ be the Cartesian product of $X^i \subseteq \mathbb{R}^{d_i}$, where $d_i$ is the dimension of the subspace. Then $\operatorname{conv}(X) = \prod_{i=1}^{n} \operatorname{conv}(X^i)$.

*Proof.* It is sufficient to give the proof for $n = 2$, since $\prod_{i=1}^{n} X^i = X' \times X^n$, where $X' = \prod_{i=1}^{n-1} X^i \subseteq \mathbb{R}^{\sum_{i=1}^{n-1} d_i}$, and the statement of the lemma can be applied recursively. Therefore, let $X = Y \times Z$, $Y \in \mathbb{R}^{d_1}$, $Z \in \mathbb{R}^{d_2}$.

Proof of $\operatorname{conv}(X) \subseteq (\operatorname{conv}(Y) \times \operatorname{conv}(Z))$: Elements of $\operatorname{conv}(X)$ are representable as $x = (\sum_{i=1}^{N} p_i y^i, \sum_{i=1}^{N} p_i z^i)$ with $p \in \Delta^N$. By definition, $\sum_{i=1}^{N} p_i y^i \in \operatorname{conv}(Y)$ and $\sum_{i=1}^{N} p_i z^i \in \operatorname{conv}(Z)$. Therefore, $x \in \operatorname{conv}(Y) \times \operatorname{conv}(Z)$.

Proof of $(\operatorname{conv}(Y) \times \operatorname{conv}(Z)) \subseteq \operatorname{conv}(X)$: Let $p \in \Delta^k$ and $q \in \Delta^m$. Elements of $\operatorname{conv}(Y) \times \operatorname{conv}(Z)$ are representable as

$$x = (\sum_{i=1}^{k} p_i y^i, \sum_{j=1}^{m} q_j z^j) = \sum_{i=1}^{k} \sum_{j=1}^{m} p_i q_j (y^i, z^j). \qquad (3.9)$$

Note that $p_i q_j \geq 0$ and

$$\sum_{i=1}^{k} \sum_{j=1}^{m} p_i q_j = \sum_{i=1}^{k} p_i \underbrace{\sum_{j=1}^{m} q_j}_{=1} = \sum_{i=1}^{k} p_i = 1 \,. \tag{3.10}$$

Since $(y^i, z^j) \in X$, this implies $x \in \text{conv}(X)$. $\qquad \square$

**Example 3.34.** Due to Lemma 3.33 the equality $\text{conv}(\{0,1\}^n) = [0,1]^n$ follows directly from the fact that $\text{conv}(\{0,1\}) = [0,1]$.

The following lemmas show that convexity is invariant with respect to linear mappings:[2]

**Lemma 3.35.** Let $X$ be a convex set. Then the set $Z := \{z \mid z = Ax + b, \ x \in X\}$ is convex as well.

*Proof.* Let $z, z' \in Z$. Then they are representable as $z = Ax + b$ and $z' = Ax' + b$ for some $x, x' \in X$ respectively. Consider $p \in (0,1)$ and $pz + (1-p)z' = p(Ax+b) + (1-p)(Ax'+b) = A(px + (1-p)x') + b$. Since $X$ is convex, $px + (1-p)x' \in X$ and, therefore, $pz + (1-p)z' \in Z$. $\quad \square$

**Lemma 3.36.** $\text{conv}\{Ax + b \mid x \in X\} = \{Ax + b \mid x \in \text{conv}(X)\}$, or in other words, convex hull commutes with linear mapping.

*Proof.* $\underline{Z \subseteq \hat{Z}}$: Let

$$Z := \text{conv}\{Ax + b \mid x \in X\} \quad \text{and} \quad \hat{Z} := \{Ax + b \mid x \in \text{conv}(X)\}.$$

Let $z = Ax + b$, $x \in X$. It implies $x \in \text{conv}(X)$ and therefore, $z \in \hat{Z}$.

$\underline{\hat{Z} \subseteq Z}$: Let $z = Ax' + b$, $x' \in \text{conv}(X)$. Therefore, $x = \sum_{i=1}^{N} p_x x$ with $x \in X$ and $p \in \Delta^N$. Hence,

$$z = Ax' + b = \sum_{i=1}^{N} p_x(Ax) + (\sum_{i=1}^{N} p_x)b = \sum_{i=1}^{N} p_x(Ax + b) \in Z \,. \tag{3.11}$$

$\square$

---

[2]In this book we do not distinguish between linear and affine mappings and call both of them linear.

From Definitions 3.24 and 3.28 it follows directly that the convex hull is a convex set. The inverse holds trivially, since each convex set can be seen as its own convex hull. A non-trivial result, stated in the theorem below, describes an important relation between *finite* sets and their convex hulls:

**Theorem 3.1** (Minkowski (1896))**.** The set $P \subset \mathbb{R}^n$ is a polytope if and only if it is representable as the convex hull of a finite set of points.

**Corollary 3.37.** A polytope is the convex hull of its vertices.

Proofs of both statements are beyond the scope of this monograph. In §3.6 we reference the literature, where these proofs are given.

Whereas Corollary 3.37 claims that any polytope is uniquely defined by a finite set (of its vertices) by its convex hull, the following theorem describes vertices of a convex hull of a finite set of points:

**Theorem 3.2.** Let $X \subset \mathbb{R}^n$ be a finite set and $z$ be a vertex of $\mathrm{conv}(X)$. Then $z \in X$.

Theorem 3.2 is illustrated in Figure 3.2. To prove Theorem 3.2 we will need the following simple but very useful Lemma:

**Lemma 3.38.** Let $a \in \mathbb{R}^n$ and $p^* = \arg\min_{p \in \Delta^n} \sum_{i=1}^n p_i a_i$. Then $p_j^* = 0$ for all $j$ such that $a_j > \min_{i=1,\dots,n}\{a_i\}$.

*Proof.* Let $i^* := \arg\min_{i=1,\dots,n}\{a_i\}$ and $a_* := a_{i^*}$. Let also $p_j^* > 0$ for some $a_j > a_*$ and

$$p_i' = \begin{cases} p_i^*, & i \notin \{i^*, j\}, \\ 0, & i = j. \\ p_{i^*}^* + p_j^*, & i = i^*. \end{cases} \tag{3.12}$$

Note that $p' \in \Delta^n$ if $p^* \in \Delta^n$. Then

$$\sum_{i=1}^n p_i' a_i = p_{i^*}^* a_* + p_j^* a_* + \sum_{i \notin \{i^*, j\}} p_i^* a_i$$

$$\overset{a_* < a_j}{<} p_{i^*}^* a_* + p_j^* a_j + \sum_{i \notin \{i^*, j\}} p_i^* a_i = \sum_{i=1}^n p_i^* a_i, \tag{3.13}$$

which contradicts the definition of $p^*$. $\qquad\square$

*Proof of Theorem 3.2.* Let $z \in \mathrm{vrtx}(\mathrm{conv}(X))$. Since $z \in \mathrm{conv}(X)$, it holds that $z = \sum_{x \in X} p_x x$ for some $p \in \Delta^X$. As $z$ is a vertex of $\mathrm{conv}(X)$, there exists $c$ such that

$$z = \underset{x \in \mathrm{conv}(X)}{\arg\min} \ \langle c, x \rangle \tag{3.14}$$

and, moreover, $z$ is a unique solution of this problem. In other words, there is no other vector $z' \in \mathrm{conv}(X)$ such that $\langle c, z \rangle = \langle c, z' \rangle$.

According to Lemma 3.38, it holds that

$$\langle c, z \rangle = \sum_{x \in X} p_x \langle c, x \rangle \geq \min_{x \in X} \langle c, x \rangle \,,$$

and $p_{x'} \neq 0$ only for $x' \in \arg\min_{x \in X} \langle c, x \rangle$ if the equality holds. All such $x'$ are minimizers of the problem (3.14) along with $z$. Therefore, $z$ can be a unique solution only if $z \in \arg\min_{x \in X} \langle c, x \rangle$, which finalizes the proof. $\qquad\square$

**Example 3.39.** According to Example 3.15, $[0, 1]^n$ is a polytope. Let us show that the set of its vertices is precisely the set $\{0, 1\}^n$.

Since $\mathrm{conv}(\{0, 1\}^n) = [0, 1]^n$ (see Example 3.34), Theorem 3.2 implies $\mathrm{vrtx}([0, 1]^n) \subseteq \{0, 1\}^n$. The inverse inclusion $\{0, 1\}^n \subseteq \mathrm{vrtx}([0, 1]^n)$ is proven in Example 3.20, which finalizes the proof.

**Example 3.40.** $\Delta^n \cap \{0, 1\}^n = \mathrm{vrtx}(\Delta^n)$. In other words, the binary vectors of a simplex are all its vertices.

Since $\Delta^n = \mathrm{conv}(\Delta^n \cap \{0, 1\}^n)$ (see Example 3.32), Theorem 3.2 implies $\mathrm{vrtx}(\Delta^n) \subseteq \Delta^n \cap \{0, 1\}^n$. The inverse inclusion $\Delta^n \cap \{0, 1\}^n \subseteq \mathrm{vrtx}(\Delta^n)$ is proven in Example 3.22, which finalizes the proof.

The following simple lemma is useful to find vertices of polytopes obtained as an intersection of polyhedra:

**Proposition 3.41.** Let $X$ be a polytope and $Z$ be a polyhedron in $\mathbb{R}^n$, and $x \in X \cap Z$. If additionally $x$ is a vertex of $X$, then $x$ is a vertex of $X \cap Z$.

The proposition is illustrated in Figure 3.3(b).

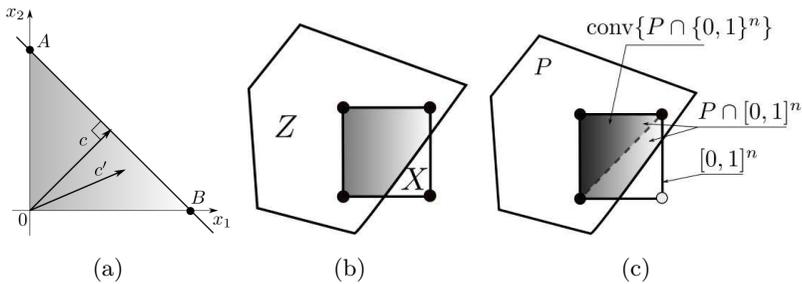

**Figure 3.3: (a)** Illustration of optimal solutions of a linear program. Given the polytope $P$ with vertices $0$, $A$ and $B$, the point $B$ is the unique solution of $\max_{x \in P} \langle c', x \rangle$, whereas for the problem $\max_{x \in P} \langle c, x \rangle$ (with $c$ orthogonal to the line $AB$) any point on the line between points $A$ and $B$ is a solution. **(b)** Illustration of Proposition 3.41. Vertices of a polytope $X$ remain vertices of the intersection $X \cap Z$, if they belong to this intersection. **(c)** Illustration of Corollary 3.56: the set of vertices of the feasible set of the LP relaxation ($P \cap [0,1]^n$) contains the vertices of the integer hull ($\mathrm{conv}(P \cup \{0,1\}^n)$).

*Proof.* Since $x$ is a vertex of $X$, there exists $c \in \mathbb{R}^n$ such that $x$ is the unique solution of $\min_{x' \in X} \langle c, x' \rangle$. Since $x \in Z$, it holds also that it is a unique solution of $\min_{x' \in X \cap Z} \langle c, x' \rangle$ (see Proposition 3.13), and, therefore, $x$ is a vertex of $X \cap Z$. □

## 3.3 Linear programs

Linear programs (LP) are the most important components of integer linear programs – a unified representation for combinatorial problems. We will see below that although linear programs constitute a continuous, non-discrete object, the set of points, which is sufficient to check to obtain their exact solution, is finite. Moreover, it can be shown that this set grows exponentially with the dimension of the problem in general. In this sense linear programs can be seen as combinatorial problems themselves. The main difference, however, is in their solvability: There exist polynomial algorithms for any linear program, whereas the class of combinatorial problems contains $\mathcal{NP}$-hard problems, like, for example, the considered MAP-inference for graphical models. Interestingly, polynomially solvable combinatorial problems often can be equivalently represented in the form of a linear program, whose size grows polynomially with the size of the combinatorial problem.

In this section we will discuss the most important properties of linear programs.

**Definition 3.42.** Let $P$ be a polyhedron in $\mathbb{R}^n$ defined by a set of linear constraints. Optimization problems of the form

$$\min_{x \in P} \langle c, x \rangle \qquad \text{and} \qquad \max_{x \in P} \langle c, x \rangle \qquad (3.15)$$

are called *linear programs (LP)*.

The following lemma is the key to the combinatorial property of linear programs:

**Lemma 3.43.** For any $a \in \mathbb{R}^n$ it holds that

$$\min_{i=1,\dots,n} \{a_i\} = \min_{p \in \Delta^n} \sum_{i=1}^{n} p_i a_i = \min_{\mu \in \text{conv}\{a_i,\ i=1,\dots,n\}} \mu \,. \qquad (3.16)$$

*Proof.* The first equality directly follows from Lemma 3.38 and the second equality holds per definition of the convex hull. □

**Corollary 3.44.** For any $c \in \mathbb{R}^n$ and any finite $X \subset \mathbb{R}^n$ it follows

$$\min_{x \in X} \langle c, x \rangle = \min_{x \in \text{conv}(X)} \langle c, x \rangle \,. \qquad (3.17)$$

*Proof.* The proof is a direct consequence of the transformations:

$$\min_{x \in \text{conv}(X)} \langle c, x \rangle = \min_{p \in \Delta^{|X|}} \left\langle c, \sum_{x \in X} p_x x \right\rangle$$
$$= \min_{p \in \Delta^{|X|}} \sum_{x \in X} p_x \langle c, x \rangle \stackrel{\text{Lem. } 3.38}{=} \min_{x \in X} \langle c, x \rangle \,. \qquad (3.18)$$

□

In general, the claim of Corollary 3.44 holds for any $X \subseteq \mathbb{R}^n$, however, for the purpose of this monograph, considering a finite set is sufficient.

**Corollary 3.45** (Combinatorial property of linear programs)**.** For any $c \in \mathbb{R}^n$ and any polytope $P$ it holds that

$$\min_{x \in P} \langle c, x \rangle = \min_{x \in \mathrm{vrtx}(P)} \langle c, x \rangle \;. \tag{3.19}$$

*Proof.* The statement follows directly from the fact that

$$P = \mathrm{conv}(\mathrm{vrtx}(P))$$

(Corollary 3.37) and Corollary 3.44. □

Corollary 3.45 claims that to solve a linear program over a polytope it is sufficient to evaluate the objective on the finite set of vertices of the polytope. This does not mean, that it is impossible that the minimum is attained in a non-vertex point of the polytope. However, it guarantees, that there will always be a vertex corresponding to the minimal value of the objective, see Figure 3.3(a). Solutions which correspond to vertices of a polytope $P$ are called *basic solutions* and can be found with the famous simplex algorithm (see references in §3.6).

## 3.4 Integer linear programs

The integer linear program (ILP) is a unified format to represent combinatorial problems. The main advantage of formulating combinatorial problems as ILPs, is the unified geometrical (polyhedral) interpretation common to all integer linear programs. This has allowed for significant progress in creating standardized methods to solve them. Below, we will provide a definition of ILPs and show their relation to linear programs.

**Definition 3.46.** A linear program with additional constraints allowing all the variables to take only values 0 or 1, e.g.

$$\min_{x \in P \cap \{0,1\}^n} \langle c, x \rangle \tag{3.20}$$

is called a 0/1 *integer linear program* (further we omit mentioning '0/1' for brevity). Here $P$ is a polyhedron in $\mathbb{R}^n$. Constraints of the form $x \in \{0,1\}^n$ are called *integrality constraints*.

The class of integer linear programs is $\mathcal{NP}$-hard, as is shown by the following two examples, where we represent known $\mathcal{NP}$-hard problem as ILPs:

**Example 3.47** (Binary Knapsack Problem). Let $\{1, \ldots, n\}$ be indexes of $n$ items, each of which has its volume $a_i$ and its cost $c_i$. Therefore, we assume two vectors $a, c \in \mathbb{R}^n$ to be given. The task is to select a subset of items to maximize their total cost, whereas their total volume should not exceed a given (knapsack) volume $b$. This can be equivalently formulated as the integer linear program

$$\max_{x \in \{0,1\}^n} \ \langle c, x \rangle \tag{3.21}$$

$$\text{s.t. } \langle a, x \rangle \le b \,. \tag{3.22}$$

**Example 3.48** (Max-cut problem). Let $\mathcal{G} = (\mathcal{V}, \mathcal{E})$ be an undirected graph and $c_{uv} \in \mathbb{R}$, $uv \in \mathcal{E}$, the costs assigned to its edges. The goal of the max-cut problem is to partition the nodes of the graph into two subsets $\mathcal{V}^0$ and $\mathcal{V}^1$, such that $\mathcal{V}^0 \cup \mathcal{V}^1 = \mathcal{V}$, $\mathcal{V}^0 \cap \mathcal{V}^1 = \emptyset$, and the sum of costs of edges connecting $\mathcal{V}^0$ and $\mathcal{V}^1$, i.e. $\sum_{\substack{uv \in \mathcal{E}: \\ u \in \mathcal{V}^0, \ v \in \mathcal{V}^1}} c_{uv}$, is maximized. Max-cut is known to be $\mathcal{NP}$-hard.

This problem can be shortly written as

$$\max_{x \in \{0,1\}^{\mathcal{V}}} \sum_{uv \in \mathcal{E}} c_{uv} |x_u - x_v| \tag{3.23}$$

by identifying the nodes in $\mathcal{V}^0$ and $\mathcal{V}^1$ with 0 and 1, respectively. To get rid of the absolute value in (3.23) and rewrite the problem as an ILP we introduce three binary variables $z_{uv}$, $z_{uv}^+$ and $z_{uv}^-$ for each edge $uv \in \mathcal{E}$. This is done in such a way that $x_u - x_v = z_{uv}^+$, if $x_u - x_v \ge 0$, and $x_u - x_v = z_{uv}^-$ otherwise. To achieve this, for each edge $uv \in \mathcal{E}$ we will introduce the linear constraint $x_u - x_v = z_{uv}^+ - z_{uv}^-$. Bearing in mind the integrality constraints $z, z^+, z^- \in \{0,1\}^{\mathcal{E}}$ the value $|x_u - x_v|$ becomes equal to $z_{uv} = z_{uv}^+ + z_{uv}^-$.

Putting all constraints together results in the ILP representation of the max-cut problem [3]

---

[3]The standard ILP representation of the max-cut problem differs from (3.24). However, it is more involved and less suitable for the illustration purposes we pursue here.

$$\max \sum_{uv \in \mathcal{E}} c_{uv} z_{uv} \,, \tag{3.24}$$

$$\text{s.t. } x_u - x_v = z_{uv}^+ - z_{uv}^-, \ uv \in \mathcal{E} \,, \tag{3.25}$$

$$z_{uv} = z_{uv}^+ + z_{uv}^- \,, \tag{3.26}$$

$$z, z^+, z^- \in \{0,1\}^{\mathcal{E}} \,, \tag{3.27}$$

$$x \in \{0,1\}^{\mathcal{V}} \,. \tag{3.28}$$

The max-cut problem can also be represented as a MAP-inference problem defined on the same graph $\mathcal{G}$, where each node $u$ is associated with the label set $\mathcal{Y}_u = \{0,1\}$. Unary costs are all equal to zero, and pairwise costs are defined as $\theta_{uv}(s,t) = -c_{uv}[\![s \neq t]\!]$ for all $uv \in \mathcal{E}$ and $(s,t) \in \mathcal{Y}_{uv}$.

In Chapter 4 we will also show how the MAP-inference problem can be represented as an ILP.

**ILP as LP** Somewhat surprisingly, any ILP can be represented as a linear program by a specially constructed polytope. Note that the feasible set of the problem (3.20) is finite, therefore one can use Corollary 3.44, which implies

$$\min_{x \in P \cap \{0,1\}^n} \langle c, x \rangle = \min_{x \in \text{conv}(P \cap \{0,1\}^n)} \langle c, x \rangle \,, \tag{3.29}$$

and, since, according to Theorem 3.1, the convex hull of a finite set is a polytope, the right-hand side of (3.29) constitutes a linear program. This fact does not mean polynomial solvability of the right-hand-side of (3.29) and, therefore, of (3.20). This is because the polytope $\text{conv}(P \cap \{0,1\}^n)$ may require exponentially many linear (in)equalities to be specified explicitly, and computational complexity of the LP optimization methods grows as a polynomial of this amount. For instance, consider the following straightforward (and, therefore, not necessarily shortest) polyhedral representation of $y \in \text{conv}(X)$, where $X = P \cap \{0,1\}^n$:

$$y = \sum_{x \in X} p_x x \,, \quad \sum_{x \in X} p_x = 1, \ p_x \geq 0 \ \forall x \in X \,. \tag{3.30}$$

The length of constraints and their number grows linearly with the size of the set $X = P \cap \{0,1\}^n$, which is exponentially large for typical ILPs (see e.g. Example 3.48). Otherwise such problems would be easily solvable by a direct enumeration of these vectors.

We will call the polytope $\mathrm{conv}(P \cap \{0,1\}^n)$ *the integer hull* of $P$, although this name has a broader meaning as the convex set of all integer points (not only those having coordinates 0 and 1) in $P$. In Chapter 4, we will reformulate the energy minimization problem (1.4) as an ILP and illustrate (3.30) with that example.

**Vertices of the integer hull**

**Lemma 3.49.** Let $X \subseteq \{0,1\}^n$. For any $x \in \mathrm{conv}(X)$ it holds $0 \leq x_i \leq 1$, that is, $\mathrm{conv}(X) \subseteq [0,1]^n$.

*Proof.* Due to Proposition 3.31 and Example 3.34 it holds that $\mathrm{conv}(X) \subseteq \mathrm{conv}(\{0,1\}^n) = [0,1]^n$. $\square$

**Corollary 3.50.** For any $x \in \mathrm{conv}(P \cap \{0,1\}^n)$ it holds that $0 \leq x_i \leq 1$.

Since $P \cap \{0,1\}^n$ is a finite set, it follows from Theorem 3.2 that $\mathrm{vrtx}(\mathrm{conv}(P \cap \{0,1\}^n)) \subseteq P \cap \{0,1\}^n$. The following proposition states that these sets are equal:

**Proposition 3.51.** For any $X \subseteq \{0,1\}^n$ it holds that $\mathrm{vrtx}(\mathrm{conv}(X)) = X$.

*Proof.* Let $z \in X$. To show that it is a vertex of $\mathrm{conv}(X)$, consider the problem $\max\limits_{x \in \mathrm{conv}(X)} \langle c, x \rangle$ with $c_i = 2z_i - 1$. Repeating the proof of Example 3.20 we conclude that

$$\max_{x \in \mathrm{conv}(X)} \langle c, x \rangle = \max_{x \in \mathrm{conv}(X)} \left( \sum_{\substack{i=1,\dots,n: \\ z_i \neq 0}} x_i - \sum_{\substack{i=1,\dots,n: \\ z_i = 0}} x_i \right).$$

Due to Lemma 3.49 it follows that this expression equals to $\sum\limits_{i=1}^{n} z_i$, and is attained in the point $z$ only. This implies that $z$ is the vertex of $\mathrm{conv}(X)$. $\square$

**Corollary 3.52.** $\mathrm{vrtx}(\mathrm{conv}(P \cap \{0,1\}^n)) = P \cap \{0,1\}^n$. In other words, each feasible point of a 0/1 ILP is a vertex of its integer hull [4].

Due to Corollary 3.52 we conclude that for any feasible point of the integer linear program $\min_{x \in P \cap \{0,1\}^n} \langle c, x \rangle$ such a vector $c \in \mathbb{R}^n$ exists, that this point becomes the unique solution of the program, see Figure 3.3(c) for an illustration.

## 3.5 Linear program relaxation

The simplest approximation of an $\mathcal{NP}$-hard integer linear program by a polynomially solvable linear problem can be constructed by omitting integrality constraints. This approximation, known as *linear program relaxation*, or simply *LP relaxation*, is a very important and powerful tool for obtaining approximate and exact solutions of integer linear programs. Below, we investigate some geometric properties of LP relaxations.

**Definition 3.53.** The linear program

$$\min_{x \in P \cap [0,1]^n} \langle c, x \rangle \tag{3.31}$$

is called *linear programming (LP) relaxation* of the integer linear program (3.20).

Note that

- the LP relaxation is constructed by substituting the integrality constraints $x_i \in \{0,1\}$ with *interval* or *box* constraints $x_i \in [0,1]$. Since $x_i \in [0,1]$ is equivalent to $0 \le x_i \le 1$, these constraints define a polyhedron. As the intersection of two polyhedra is a polyhedron (see Proposition 3.18), the problem (3.31) is a linear program;

- the LP relaxation is a relaxation in terms of Definition 3.6, since its feasible set is a superset of the one for the non-relaxed problem and their objectives coincide.

---

[4]This claim does not hold for a more general class of integer linear programs without the 0/1 constraint.

The most important properties of LP relaxations read:

**Proposition 3.54.** Let $f^*$ and $\hat{f}$ be optimal values of the ILP problem (3.20) and its LP relaxation (3.31) respectively. Let also $\hat{x}$ be a solution of the LP relaxation. Then the following assertions hold:

1. $\hat{f} \leq f^*$.

2. If $\hat{x} \in \{0,1\}^n$, then $\hat{x}$ is a solution of the non-relaxed ILP problem (3.20).

3. $P \cap [0,1]^n \supseteq \text{conv}(P \cap \{0,1\}^n)$, that is, the feasible set of the LP relaxation is a superset of the integer hull.

*Proof.* The first two statements follow directly from Proposition 3.9. The third one follows from the fact that $\text{conv}(P \cap [0,1]^n) = P \cap [0,1]^n$ (since $P$ and $[0,1]^n$ are convex, their intersection is a convex set), $P \cap [0,1]^n \supseteq P \cap \{0,1\}^n$, and Proposition 3.31. $\square$

### 3.5.1   Vertices of the feasible set of the LP relaxation

While the statement (2) of Proposition 3.54 claims that integrality of a relaxed solution implies its optimality for the non-relaxed ILP problem, it does not say, whether there has to be a solution to the LP relaxation which is integral. Or in other words, whether any solution of the ILP problem can be obtained as an optimal solution of its LP relaxation. The following proposition fills this gap:

**Proposition 3.55.** For any polyhedron $P$ it holds that $P \cap \{0,1\}^n \subseteq \text{vrtx}(P \cap [0,1]^n)$.

In other words, Proposition 3.55 states that all feasible vectors of an ILP are vertices of the feasible set of its LP relaxation and, therefore, can be reached as a relaxed solution.

*Proof.* Let us consider the intersection of the polyhedron $P$ and the polytope $[0,1]^n$. According to Proposition 3.41, vertices of $[0,1]^n$ belonging to $P$ are also vertices of $P \cap [0,1]^n$. As shown in Example 3.39, any vertex of the $n$-dimensional cube $[0,1]^n$ is a vector from the set $\{0,1\}^n$ and vice versa. Therefore, any point from $P \cap \{0,1\}^n$ is a vertex of $P \cap [0,1]^n$. $\square$

The following corollary additionally relates vertices of the integer hull and of the feasible set of the LP relaxation (see Figure 3.3(c)):

**Corollary 3.56.** $\mathrm{vrtx}(\mathrm{conv}(P \cap \{0,1\}^n)) \subseteq \mathrm{vrtx}(P \cap [0,1]^n)$.

*Proof.* The statement follows from $\mathrm{vrtx}(\mathrm{conv}\{P \cap \{0,1\}^n\}) = P \cap \{0,1\}^n \subseteq \mathrm{vrtx}(P \cap [0,1]^n)$, where the first relation follows from Proposition 3.52 and the second holds due to Proposition 3.55. $\square$

The last proposition of the chapter states that there is no other binary vector in the feasible set of LP relaxation than those, which are feasible for the non-relaxed integer linear program:

**Proposition 3.57.** For any set $P \subset \mathbb{R}^n$ it holds that $\mathrm{vrtx}(P \cap [0,1]^n) \cap \{0,1\}^n = P \cap \{0,1\}^n$.

*Proof.* Let $x \in P \cap \{0,1\}^n$. Trivially $x \in \{0,1\}^n$. According to Proposition 3.55 $x \in \mathrm{vrtx}(P \cap [0,1]^n)$. Therefore, $P \cap \{0,1\}^n \subseteq \mathrm{vrtx}(P \cap [0,1]^n) \cap \{0,1\}^n$.

The other way around, let $x \in \mathrm{vrtx}(P \cap [0,1]^n) \cap \{0,1\}^n$. Trivially, $x \in \mathrm{vrtx}(P \cap [0,1]^n)$ and $x \in \{0,1\}^n$. Since $\mathrm{vrtx}(P \cap [0,1]^n) \subseteq P \cap [0,1]^n$, it holds $x \in P$. Therefore, $x \in P \cap \{0,1\}^n$. $\square$

## 3.6 Bibliography and further reading

One of the best introductory text-books about linear programming and the simplex method is written by the author of the simplex method, Georg Dantzig [23]. For fundamental results on polyhedra, linear inequalities and linear programming the monograph [107] is a classical reference. It includes also fundamental results on the polynomial solvability of linear programs. We refer to the text-book [57] for the Minkowski theorem and its corollary. This text-book as well as [78] can be recommended for learning the basics of combinatorial optimization and its relation to (integer) linear programming.

# 4

## Energy Minimization as Integer Linear Program

In this chapter we will consider the MAP-inference problem in the context of the theoretical background provided in Chapter 3. In particular, we will show how to reformulate this problem as an integer linear program and investigate the properties of its integer hull. Then we will consider its LP relaxation, discuss optimality properties of the relaxed solution and provide a simple rounding technique to obtain approximate solutions of the non-relaxed problem.

### 4.1 MAP-inference as ILP

To formulate MAP-inference as an ILP we first have to represent its objective in the form of a linear function. To do so, recall how the unary and pairwise cost functions have been translated to the cost vector $\theta \in \mathbb{R}^{\mathcal{I}}$ in Chapter 1, where elements of the set $\mathcal{I}$ index all labels in all nodes and all label pairs in each edge.

Consider the set $\{0, 1\}^{\mathcal{I}}$ of binary vectors of the same dimension and let $\mu$ be a vector from this set. This means that each label $s$ in a particular node $u$ is associated with a binary variable $\mu_u(s) \in \{0, 1\}$, and all such binary variables together form a vector $\mu_u = (\mu_u(s),\ s \in \mathcal{Y}_u) \in \{0,1\}^{\mathcal{Y}_u}$. Each pair of labels $(s,t)$ corresponding to an edge $uv$ is associated with $\mu_{uv}(s,t) \in \{0,1\}$. Similarly, the vector $\mu_{uv} = (\mu_{uv}(s,t),\ (s,t) \in \mathcal{Y}_{uv}) \in \{0,1\}^{\mathcal{Y}_{uv}}$ contains all these variables as its coordinates. All vectors $\mu_u$ are combined to a joint unary vector $\mu_{\mathcal{V}} = (\mu_u,\ u \in \mathcal{V})$ and the joint pairwise vector $\mu_{\mathcal{E}}$ is defined as $(\mu_{uv},\ uv \in \mathcal{E})$. Finally, vector $\mu \in \{0,1\}^{\mathcal{I}}$ consists of the joint unary and pairwise vectors, i.e. $\mu = (\mu_{\mathcal{V}}, \mu_{\mathcal{E}})$.

In other words, each elementary cost gets a corresponding 0/1 variable. The values 1 and 0 of this variable "switch" this cost "on" and "off", respectively. Therefore the inner product $\langle \theta, \mu \rangle$ denotes the sum of costs, which correspond to the labels and label pairs with the value 1 in the corresponding coordinate of vector $\mu$, see Figure 4.1 for an illustration.

These labels and label pairs are not necessarily consistent, i.e. they do not necessarily correspond to a labeling, see Figure 4.1(d). To enforce this, let us introduce the function $\delta_{\mathcal{G}} \colon \mathcal{Y}_{\mathcal{V}} \to \{0,1\}^{\mathcal{I}}$, which maps each labeling $y$ onto the corresponding binary vector $\delta_{\mathcal{G}}(y)$ such that $[\delta_{\mathcal{G}}(y)]_u(s) = \llbracket y_u = s \rrbracket$ and $[\delta_{\mathcal{G}}(y)]_{uv}(s,t) = \llbracket (s,t) = (y_u, y_v) \rrbracket$, for all $u \in \mathcal{V}$, $s \in \mathcal{Y}_u$ and $uv \in \mathcal{E}$, $(s,t) \in \mathcal{Y}_{uv}$, respectively. Here, the notation $[\delta_{\mathcal{G}}(y)]_u(s)$ and $[\delta_{\mathcal{G}}(y)]_{uv}(s,t)$ stands for the corresponding coordinates of the vector $\delta_{\mathcal{G}}(y)$. For the sake of notation we will write $\delta(y)$ instead of $\delta_{\mathcal{G}}(y)$ whenever the graph $\mathcal{G}$ is clear from the context.

Using this notation we can write the MAP-inference problem (1.4) as

$$\min_{y \in \mathcal{Y}_{\mathcal{V}}} E(y; \theta) = \min_{y \in \mathcal{Y}_{\mathcal{V}}} \left( \sum_{u \in \mathcal{V}} \theta_u(y_u) + \sum_{uv \in \mathcal{E}} \theta_{uv}(y_u, y_v) \right)$$

$$= \min_{y \in \mathcal{Y}_{\mathcal{V}}} \left( \sum_{u \in \mathcal{V}} \sum_{s \in \mathcal{Y}_u} \theta_u(s)[\delta(y)]_u(s) + \sum_{uv \in \mathcal{E}} \sum_{(s,t) \in \mathcal{Y}_{uv}} \theta_{uv}(s,t)[\delta(y)]_{uv}(s,t) \right)$$

$$= \min_{y \in \mathcal{Y}_{\mathcal{V}}} \langle \theta, \delta(y) \rangle = \min_{\mu \in \delta(\mathcal{Y}_{\mathcal{V}})} \langle \theta, \mu \rangle \ , \quad (4.1)$$

where the set $\delta(\mathcal{Y}_{\mathcal{V}})$ is the set $\{\mu \in \{0,1\}^{\mathcal{I}} \mid \mu = \delta(y),\ y \in \mathcal{Y}_{\mathcal{V}}\}$ of all binary vectors corresponding to labelings.

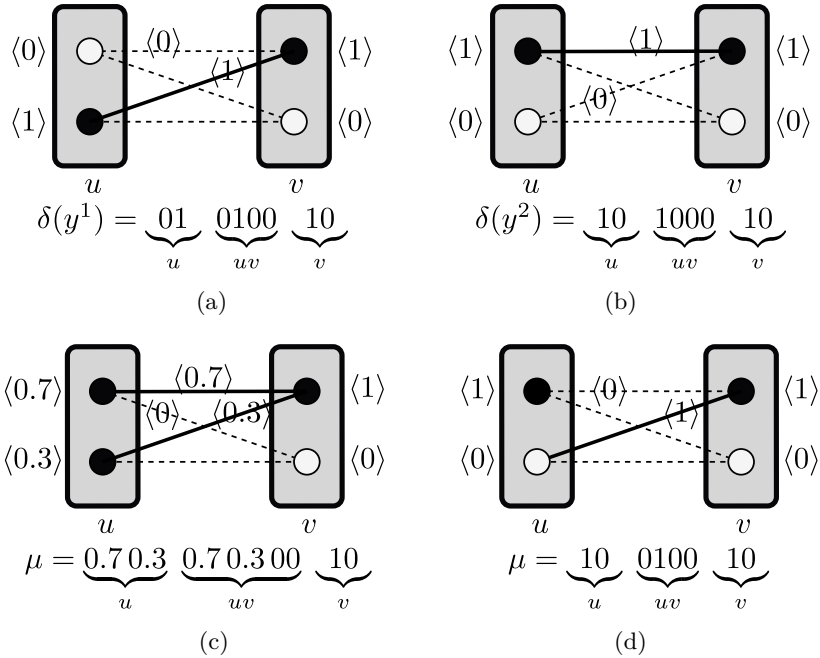

**Figure 4.1: (a,b)** Correspondence between labelings $y^1, y^2 \in \mathcal{Y}_{\mathcal{V}}$ and their representations as binary vectors $\delta(y^1)$ and $\delta(y^2)$. The coordinates of the binary vectors and their convex combinations are given in angular brackets $\langle \rangle$. **(c)** The representation of the relaxed labeling $\mu = 0.3\delta(y^1) + 0.7\delta(y^2)$. **(d)** Binary vector $\mu \in \{0,1\}^{\mathcal{I}}$ not corresponding to any labeling, because the non-zero values assigned to the labels in node $u$ and the label pair on the edge $uv$ are inconsistent.

To turn the last expression of (4.1) into an ILP, we have to give a polyhedral description of the set $\delta(\mathcal{Y}_{\mathcal{V}})$, that is, represent it in the form $P \cap \{0,1\}^{\mathcal{I}}$ for some polyhedron $P$.

### 4.1.1 Local polytope and ILP representation

To find a polyhedral representation of the set $\delta(\mathcal{Y}_{\mathcal{V}})$, let us investigate the linear constraints which are fulfilled for all elements of this set.

**Lemma 4.1.** Let $\mu \in \delta(\mathcal{Y}_{\mathcal{V}})$. Then $\sum_{s \in \mathcal{Y}_w} \mu_w(s) = 1$ for all $w \in \mathcal{V} \cup \mathcal{E}$.

*Proof.* Given that $\delta(\mathcal{Y}_{\mathcal{V}}) \subset \{0,1\}^{\mathcal{I}}$, the statement essentially reduces to the claim that there is a unique label (label pair) $s \in \mathcal{Y}_w$ for each node (edge) $w \in \mathcal{V} \cup \mathcal{E}$ such that $\mu_w(s) = 1$. For other labels (label pairs) $t \in \mathcal{Y}_w \backslash \{s\}$ it holds that $\mu_w(t) = 0$. This claim directly follows from the definition of $\delta(\mathcal{Y}_{\mathcal{V}})$. $\qquad\square$

Lemma 4.1 and the fact that $\delta(\mathcal{Y}_{\mathcal{V}}) \subset \{0,1\}^{\mathcal{I}}$ result in two linear constraints restricting the set $\delta(\mathcal{Y}_{\mathcal{V}})$:

$$\sum_{s \in \mathcal{Y}_w} \mu_w(s) = 1, \quad w \in \mathcal{V} \cup \mathcal{E}, \tag{4.2}$$
$$\mu_w \geq 0.$$

These constraints are insufficient, because they relate only variables within each separate node (edge) and do not provide connections between them.

For example, consider the simple model in Figure 4.1(d). The binary vector $\mu = (10\ 0100\ 10)$ satisfies (4.2), but does not correspond to any labeling.

To change this, we must add connections with the following property:

$$\mu_u(s) = 1 \; \Rightarrow \; \forall v \in \mathcal{N}_b(u) \; \exists t \in \mathcal{Y}_v \colon \mu_{uv}(s,t) = 1. \tag{4.3}$$

In other words, for each label $s$ assigned the value 1 in a node $u$ there must be a label pair $(s,t)$ assigned 1 for each edge incident to $u$. And symmetrically,

$$\mu_{uv}(s,t) = 1 \; \Rightarrow \; \mu_u(s) = 1 \text{ and } \mu_v(t) = 1. \tag{4.4}$$

In other words, if a label pair $(s,t)$ is assigned the value 1 in an edge $uv$, then both the label $s$ in the node $u$ and the label $t$ in the node $v$ must be assigned value 1.

These connections must be represented in the form of linear (in)equalities to define a polyhedron together with the constraints (4.2). One way to do so is as follows:

$$\sum_{s \in \mathcal{Y}_u} \mu_{uv}(s,t) = \mu_v(t), \; \forall \, uv \in \mathcal{E}, \; t \in \mathcal{Y}_v. \tag{4.5}$$

It is easy to see that for $\mu \in \{0,1\}^{\mathcal{I}}$ constraints (4.5) are equivalent to those defined by (4.3) and (4.4).

Together with (4.2) these constraints define the *local marginal polytope*

$$\mathcal{L} := \begin{cases} \sum_{s \in \mathcal{Y}_u} \mu_{uv}(s,t) = \mu_v(t), & uv \in \mathcal{E}, \ t \in \mathcal{Y}_v & \text{(a)} \\ \sum_{t \in \mathcal{Y}_v} \mu_v(t) = 1, & v \in \mathcal{V} & \text{(b)} \\ \mu \geq 0, & & \text{(c)} \\ \sum_{(s,t) \in \mathcal{Y}_{uv}} \mu_{uv}(s,t) = 1, & uv \in \mathcal{E} . & \text{(d)} \end{cases} \tag{4.6}$$

Constraints (4.6) indeed define a polytope, since the lines (b-d) imply that $\mu \in [0,1]^{\mathcal{I}}$. For brevity we will omit the word *marginal* and stick to the term *local polytope*, as it is also often used in the literature.

Constraints (4.6)(a) are called *coupling* or *marginalization* constraints. We will mainly use the first term to underline their role to *couple* variables of individual unary and pairwise factors. The term *marginalization* comes from summing over, or *marginalizing out*, one variable, e.g. $s \in \mathcal{Y}_u$ in (4.6)(a).

Constraints (4.6)(b) and (d) are often referred to as *uniqueness* or *simplex* constraints.

The constraints (4.6)(d) can be obtained by summing up (4.6)(a) over all $t$ and taking into account (4.6)(b). Therefore they can be omitted.

By construction, binary vectors satisfying (4.6) correspond to labelings, i.e. $\mathcal{L} \cap \{0,1\}^{\mathcal{I}} = \delta(\mathcal{Y}_{\mathcal{V}})$. This in its turn results in the ILP representation of the MAP-inference problem, which reads:

$$\min_{y \in \mathcal{Y}_{\mathcal{V}}} E(y; \theta) = \min_{\mu \in \mathcal{L} \cap \{0,1\}^{\mathcal{I}}} \langle \theta, \mu \rangle . \tag{4.7}$$

Naturally, the local polytope contains also non-binary vectors. An example of such a vector is given in Figure 4.1(c).

**Example 4.2** (How large-scale is the ILP representation of the MAP-inference problem?). Consider the depth reconstruction Example 1.2 from Chapter 1. Each element of the nodes set $\mathcal{V}$ corresponds to a pixel in an input image. Let the latter be as small as $100 \times 100 = 10^4$. Let also each node be associated with 10 labels only. It implies, there are in total $10^{10^4}$ labeling.

To estimate the number of variables in the ILP representation (4.7) we count for $10 \times 10^4 = 10^5$ total number of labels (corresponding to the variables $\mu_u(s)$).

The number of edges in the grid graph is $\approx 2 \times 10^4$ and each edge is associated with $10^2$ label pairs (and the corresponding variables $\mu_{uv}(s,t)$). Therefore, the total number of label pairs is $\approx 2 \times 10^4 \times 10^2 = 2 \times 10^6$.

This results in $O(10^6)$ variables in th ILP representation (4.7).

As it follows from the definition of the local polytope (4.6), the number of equality constraints grows linearly with the total number of labels and therefore is of order $O(10^5)$. However, since there is the non-negativity constraint $\mu_w(s) \geq 0$ for each variable, the total number of constraints has the same order $O(10^6)$ as the number of variables.

The millions of constraints and variables turn the MAP-inference (4.7) to the class of *large-scale optimization* problems. Large size of the problem significantly limits the set of methods, which can be applied to tackle it.

## 4.2 ILP properties of the MAP-inference problem

Due to the ILP representation (4.7) of the MAP-inference problem, we can directly make conclusions about its own properties and the properties of its natural LP relaxation:

- The integer hull of the set of feasible binary vectors $\delta(\mathcal{Y}_{\mathcal{V}}) = \mathcal{L} \cap \{0,1\}^{\mathcal{I}}$ is $\mathcal{M} := \mathrm{conv}(\delta(\mathcal{Y}_{\mathcal{V}})) = \mathrm{conv}(\mathcal{L} \cap \{0,1\}^{\mathcal{I}})$ and is called *marginal polytope*. As it follows from Corollary 3.52, $\mathrm{vrtx}(\mathcal{M}) = \delta(\mathcal{Y}_{\mathcal{V}})$. In other words, each vertex of the marginal polytope corresponds to a labeling, and the other way around, each labeling corresponds to a vertex of $\mathcal{M}$.

- Since $\mathcal{L} \cap [0,1]^{\mathcal{I}} = \mathcal{L}$ due to (4.6)(b)-(d), the LP relaxation of the MAP-inference problem (4.7) reads

$$\min_{\mu \in \mathcal{L}} \langle \theta, \mu \rangle \ . \tag{4.8}$$

The local polytope $\mathcal{L}$ gives also the name to the LP relaxation (4.8), which is often referred to as the *local polytope relaxation* of the MAP-inference problem. Like any relaxation, it constitutes a lower bound, i.e. $E(y; \theta) \geq \min_{\mu \in \mathcal{L}} \langle \theta, \mu \rangle$ for any (in particular the optimal) labeling $y$. Cases when this bound is tight, i.e. $\min_{y \in \mathcal{Y}_{\mathcal{V}}} E(y; \theta) = \min_{\mu \in \mathcal{L}} \langle \theta, \mu \rangle$, will be called *LP-tight*.

- The marginal polytope is a subset of the local polytope, i.e. $\mathcal{M} \subseteq \mathcal{L}$. According to Corollary 3.56, moreover, $\mathrm{vrtx}(\mathcal{M}) \subseteq \mathrm{vrtx}(\mathcal{L})$ holds, that is, all vertices of the marginal polytope are also vertices of the local polytope.

## 4.3    Properties of the local polytope relaxation

### 4.3.1    Fractional vertices of the local polytope

So far we have shown that vertices of the local polytope include all binary vectors corresponding to labelings. We will refer to these vertices as *integer*. However, as the following example shows, the local polytope may have *fractional* vertices, which are vectors with non-integer coordinates. These vectors, also called *relaxed labelings*, may correspond to the relaxed solutions of the MAP-inference problem.

**Example 4.3** ($\mathcal{M} \subset \mathcal{L}$). Consider the graphical model with three nodes as in Figure 4.2(a). Assume that all unary costs are zero, pairwise costs corresponding to solid lines are equal to 0, and costs corresponding to dashed lines are $\alpha > 0$.

Consider the vector $\mu'$, whose coordinates corresponding to all possible labels in the unary factors and to solid lines in the pairwise factors are equal to 0.5 and zero otherwise.

The vector $\mu'$ satisfies constraints (4.6) of the local polytope, therefore, $\mu' \in \mathcal{L}$. Since for any $\mu \in [0,1]^{\mathcal{I}}$ it holds that $\langle \theta, \mu \rangle \geq 0$, it holds also for all $\mu \in \mathcal{L} \subset [0,1]^{\mathcal{I}}$. Since $\langle \theta, \mu' \rangle = 0$, this implies that $\mu'$ is a solution of the LP relaxation (4.8).

It is also easy to show that $\mu'$ is the unique solution of the LP relaxation and, therefore, is a vertex of the corresponding local polytope $\mathcal{L}$.

Indeed, for any relaxed solution $\mu$ it must hold $\mu_{uv}(s,t) = 0$ if $\theta_{uv}(s,t) = \alpha$, since otherwise $\langle \theta, \mu \rangle > 0$. Informally speaking, any relaxed solution may "include" only label pairs with zero cost. It remains to show that $\mu'$ is the only point in $\mathcal{L}$, which satisfies this condition. Let $\mu \in \mathcal{L}$ satisfy this condition and let $\mu_1(0) = \gamma$. Sequentially applying the coupling constraints (a) of the local polytope (4.6) one obtains that $\mu_1(0) = \mu_{12}(0,0) = \gamma$, $\mu_{12}(0,0) = \mu_2(0) = \gamma$, $\mu_{23}(0,1) = \mu_3(1) = \gamma$, $\mu_3(1) = \mu_{31}(1,1) = \gamma$ and, finally, $\mu_{31}(1,1) = \mu_1(1) = \gamma$. Due to the simplex constraints (4.6) (b-d) it holds that $\mu_1(1) + \mu_1(0) = 2\gamma = 1$ and, therefore, $\gamma = 0.5$, which finalizes the proof.

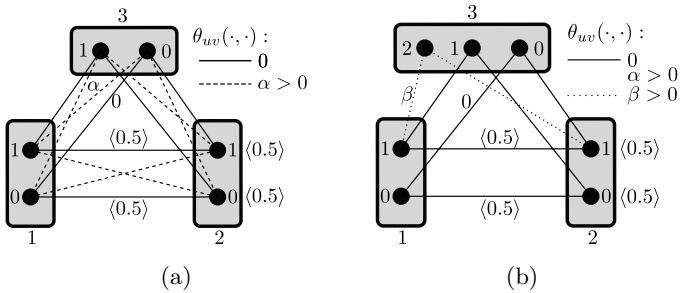

**Figure 4.2: (a)** Illustration for Example 4.3. The fully-connected graph $\mathcal{G}$ contains three nodes $1, 2, 3$. Unary costs as well as the pairwise costs corresponding to solid lines are equal to $0$. Pairwise costs corresponding to dashed lines are equal to $\alpha > 0$. Angular brackets $\langle\rangle$ are used for coordinates of $\mu$. For the unique relaxed solution $\mu$ it holds that $\mu_u(0) = \mu_u(1) = 0.5$ for all $u$. The same value $0.5$ is assigned to the coordinates corresponding to the solid lines, i.e. $\mu_{12}(0,0) = \mu_{12}(1,1) = \mu_{13}(1,1) = \mu_{13}(0,0) = \mu_{32}(1,0) = \mu_{32}(0,1) = 0.5$. The remaining coordinates are zero. **(b)** Illustration for Example 4.5. Similar to (a), but the node $3$ has $3$ labels. All unary costs are zero. All label pairs not connected by a line, have cost $\alpha > 0$, dotted lines correspond to the pairwise cost $\beta > 0$, solid lines to the cost $0$. It holds $\alpha > 2\beta$. The optimal labeling is $(1, 1, 2)$, the optimal relaxed labeling is the same as in (a).

### 4.3.2 Rounding of LP solutions

An LP relaxation of some integer linear program is typically considered to be good, if the lower bound is close to the optimum of the non-relaxed problem and, more importantly, if most of the variables of the relaxed solution have integer values.

Given the relaxed solution one may round the fractional variables to obtain an approximate solution of the non-relaxed problem, as it is shown in Example 4.4.

**Example 4.4** (Deterministic rounding). Let $\mu'$ be a relaxed solution and let some of the nodes be assigned fractional values. For such a fractional-valued node $u$ it holds that $\mu'_u \notin \{0, 1\}^{\mathcal{Y}_u}$, or, in other words, $\mu_u(s) \in (0, 1)$ for at least two labels $s$. A labeling $y'$, which is an approximate non-relaxed solution to the MAP-inference problem can be obtained by rounding of $\mu'$ in the fractional-valued nodes:

$$y'_u := \arg\max_{s \in \mathcal{Y}_u} \mu'_u(s), \ u \in \mathcal{V}. \qquad (4.9)$$

In particular, from equation (4.9) it follows that the label $s$ is assigned to the node $u$ if $\mu_u(s) = 1$, i.e. if the node $u$ is already integer-valued.

Since $y'$ is only an approximate solution of the MAP-inference problem, it holds that $E(y'; \theta) - \min_{y \in \mathcal{Y}_{\mathcal{V}}} E(y; \theta) = d \geq 0$. Since the value of $d$ determines quality of the approximate solution $y'$, the natural question that arises is how big $d$ can be.

Unfortunately, as the following Example 4.5 shows, the value of $d$ can be arbitrarily large.

**Example 4.5** (Deterministic rounding can be arbitrarily bad). Consider the graphical model in Figure 4.2(b). It differs from the one in Example 4.3 (Figure 4.2(a)) by an additional label 2 in node 3. As in Example 4.3 this label has unary cost 0. The pairwise costs are augmented with $\theta_{13}(1, 2) = \theta_{23}(1, 2) = \beta > 0$ and $\theta_{13}(0, 2) = \theta_{23}(0, 2) = \alpha$, where $\alpha > 0$ is the same positive constant as in Example 4.3.

The optimal relaxed labeling for the corresponding MAP-inference problem coincides with the one for Example 4.3. The optimal solution of the non-relaxed problem depends on the ratio between $\alpha$ and $\beta$. If $\alpha - 2\beta > 0$ the optimal labeling is $(1, 1, 2)$ and has energy $2\beta$. The best labeling which can be obtained by deterministic rounding (4.9) has the energy $\alpha$. Therefore, the gap value constitutes $d = \alpha - 2\beta$ and it can be made arbitrarily large.

Unfortunately, as shown by e.g. [70], MAP-inference belongs to the problems for which no polynomial constant factor approximation algorithm is possible. Therefore, one cannot hope that there exists a rounding technique for the (polynomially solvable) LP relaxation which would provide such a bound. This does not imply that the rounding (4.9) does not make sense to use in practice. Indeed, this simple technique often leads to reasonable results in applications.

Although good approximation is impossible for the whole class of MAP-inference problems, some of its subclasses belonging to the so called *metric labeling problem* admit such bounds. Moreover, efficient algorithms exist which are able to attain these bounds. We refer to Chapter 11 for further discussion on this topic.

**Partial optimality of integer-valued coordinates**     Most rounding schemes implicitly assume that integer-valued coordinates of the relaxed solution belong to an optimal solution of the non-relaxed problem. In other words, they assume that these coordinates are *partially optimal*.

This observation indeed often holds in practice, and, moreover, there are exact algorithms for the non-relaxed MAP-inference which efficiently utilize it (see §4.4 for references).

In general, however, this assumption is wrong as shown by Example 4.6.

**Example 4.6** (Integrality does not imply partial optimality). Consider the problem in Figure 4.3. All unary costs are assigned zero values, solid lines correspond to zero pairwise costs, dotted lines correspond to cost $\beta > 0$ and those, which are not shown, to cost $\alpha > 0$. Moreover, $\alpha > 4\beta$.

The graphical model is constructed by extending the one from Example 4.3. In particular, each of the original three nodes is augmented with a third label and an additional fourth node with two labels is added to the model. Costs are selected such that the relaxed solution from Example 4.3 when augmented with the label 0 of the fourth node, becomes the relaxed solution $\mu$ of the problem. A the same time, the optimal labeling $y = (2, 2, 2, 1)$ consists of the labels for which $\mu_u(y_u) = 0$ for all $u$.

This demonstrates that integer-valued coordinates of the relaxed solution need not belong to any optimal non-relaxed one.

Example 4.6 provides the general worst-case scenario. In Chapter 12 we will show that the local polytope relaxation of binary graphical models (where only two labels can be assigned to each graph node) indeed possesses the considered partial optimality property. This is the reason why the graphical model in Example 4.6 is not binary.

**Universality of local polytope**    The relaxed MAP-inference problem $\min_{\mu \in \mathcal{L}} \langle \theta, \mu \rangle$ defines a special class of linear programs, therefore, it can be solved with general-purpose LP solvers in polynomial time. A natural question is whether this class of problems is specific enough to allow for specialized LP solvers able to solve it faster even in the worst-case.

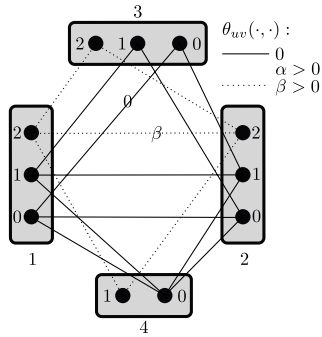

**Figure 4.3:** Illustration of Example 4.6. Notation is similar to Figure 4.2(b): unary costs are 0, as well as pairwise costs associated with solid lines. Dotted lines are associated with pairwise costs $\beta > 0$, not shown lines with the cost $\alpha > 0$. It holds $\alpha > 4\beta$. The optimal labeling reads $(2, 2, 2, 1)$, the non-zero coordinates of the optimal relaxed solution $\mu$ correspond to solid lines and are equal to $0.5$ as well as values assigned to the labels connected by the solid lines in nodes $1, 2, 3$: $\mu_1(0) = \mu_1(1) = \mu_2(0) = \mu_2(1) = \mu_3(0) = \mu_3(1) = 0.5$. Finally, $\mu_4(0) = 1$.

Unfortunately, the answer to this question is negative, as stated by the following theorem:

**Theorem 4.1** (Universality of the local polytope [83])**.** *Any* linear program reduces to the local polytope relaxation of a certain graphical model with 3 labels in linear time.

In other words, there are instances of the problem, which are as difficult as the most difficult linear programs. Therefore, the known complexity estimates for linear programs apply to the local polytope relaxation.

However, problems emerging in numerous applications, typically do not represent the worst-case scenario and can be approximately or even exactly solved with efficient polynomial algorithms. A powerful subclass of such algorithms, which we will consider in further chapters, is based on the local polytope relaxation.

## 4.4 Bibliography and further reading

Two groups of methods try to utilize the observation that integer-valued coordinates of a relaxed solution mostly keep their value in an optimal solution. First of all, these are the so-called *partial optimality* or *persistency* methods [60, 61, 46, 108, 110, 125, 124, 120, 113, 109, 111]. These polynomial methods formulate and verify sufficient optimality conditions for integer-valued coordinates. Those satisfying these conditions can be fixed and further exact combinatorial solvers can be applied [3, 43].

An alternative method [95, 33] builds a series of nested auxiliary subproblems mostly consisting of the fractional-valued coordinates of the relaxed solution. The subproblems grow until a sufficient optimality condition can be proved for a current one. The idea of the method is to restrict application of a combinatorial solver to only a small part of the initial problem related to the fractional-valued coordinates.

Universality of the local polytope was originally shown in [85]. The topic was further developed in [83, 86, 84]. Interestingly, a similar complexity result was obtained also for LP-relaxations of other $\mathcal{NP}$-hard problems [87].

A unified description and comparison of other (convex) relaxations of the MAP-inference problem can be found in the seminal work of [64].

# 5

# Background: Basics of Convex Analysis

This section introduces the notion of a convex optimization problem. This type of problems plays a central role in optimization as a whole, since there are good reasons to treat these problems as efficiently solvable. Most of the methods we consider below are based on approximating a difficult non-convex problem by a much easier convex one and then solving the latter.

One such approximation technique known as *Lagrange duality* and its application to integer linear programs are considered in the last sections of this chapter.

## 5.1 Convex optimization problems

### 5.1.1 Extended value functions

So far, when speaking about the optimization problem

$$\min_{x \in X \subseteq \mathbb{R}^n} f(x)\,, \tag{5.1}$$

we assumed that the value $f(x)$ is finite for all $x \in X$. But this may not always hold. Consider for example, that $f(x)$ is defined implicitly, as the result of another optimization, e.g. $f(x) = \max_{z \in Z} g(x, z)$ with $g \colon X \times$ $Z \to \mathbb{R}$. The maximization problem may turn out to be unbounded for some $x$. Loosely speaking, this would mean that $f(x) = \infty$.

**Example 5.1** (Implicitly defined extended value function). Consider the following extended value function, defined implicitly:

$$f(x) = \max_{z \in \mathbb{R}} (x_1 z + x_2) . \tag{5.2}$$

It is easy to see that $f(x) = \begin{cases} x_2, & x_1 = 0 \\ \infty, & x_1 \neq 0 \end{cases}$ . Therefore, $f$ is extended value, although $x_1 z + x_2$ is defined for all $z, x$ and $z$ is an unconstrained variable.

Let us now take a closer look at the minimization problem

$$\min_{x \in X} f(x)$$

when $f$ is extended value:

- Let now $X'$ be the maximal subset of $X$ such that the value $f(x)$ is finite for all $x \in X'$. Then the optimal value of (5.1) is finite as well and can be attained in some $x \in X'$ only.

- If $f(x)$ is unbounded from above for all $x \in X$, this implies that its minimal value is unbounded as well. In other words, $\min_{x \in X} f(x) = \infty$. This is the same as in the situation, where problem (5.1) is infeasible, since there is no element $x \in X$ such that $f(x)$ is finite.

With these considerations in mind it makes sense to assume the following definition.

**Definition 5.2** (Extended value function). A mapping of the form $f \colon X \to \mathbb{R} \cup \{\infty, -\infty\}$ is called an *extended value function*. The set $\{x \in X \mid -\infty < f(x) < \infty\}$ is called *domain* of $f$ and is denoted by $\mathrm{dom} f$.

Let $f$ be an extended value function. The constrained minimization problem $\min_{x \in X} f(x)$ is defined similarly as in Chapter 3 with the following modifications:

- The feasible set is defined as $\{x \in X \mid f(x) < \infty\}$.

- The problem is called unbounded if $\inf_{x \in X} f(x) = -\infty$. In particular, it is sufficient that there exists an $x \in X$ such that $f(x) = -\infty$ for the problem to be unbounded.

The maximization of extended value functions is treated similarly by substituting min with max, inf with sup, and $-\infty$ with $\infty$ and the other way around.

Most importantly, with the extended-valued functions one can "integrate" constraints into the objective function. Let $X \subset \mathbb{R}^n$ and let $\iota_X \colon \mathbb{R}^n \to \mathbb{R} \cup \{\infty\}$ be the extended value function

$$\iota_X(x) = \begin{cases} 0, & x \in X \\ \infty, & x \notin X \,. \end{cases} \tag{5.3}$$

Unless $f$ is equal to $-\infty$ for some $x \in X$, the optimization problem $\min_{x \in X} f(x)$ can be equivalently written as $\min_{x \in \mathbb{R}^n} (f(x) + \iota_X(x))$. Similarly, $\max_{x \in X} f(x) = \max_{x \in \mathbb{R}^n} (f(x) - \iota_X(x))$.

In the following, when speaking about a function, we will assume an *extended value* function, if the opposite is not stated.

### 5.1.2  Convex functions

**Definition 5.3.** An *epigraph* of a function $f$ is defined as the set $\operatorname{epi} f := \{(x, t) \colon x \in \operatorname{dom} f, \ t \geq f(x)\}$.

Loosely speaking, epigraph is set of points "above" the plot of a function, as illustrated in Figure 5.1(a).

**Example 5.4.** The epigraph of a linear function $f \colon \mathbb{R}^n \to \mathbb{R}$ is a half-space in $\mathbb{R}^{n+1}$.

The following proposition simplifies the analysis of a large class of functions:

**Proposition 5.5.** For any two functions $f$ and $g$ it holds that $\operatorname{epi}(\max\{f, g\}) = \operatorname{epi}(f) \cap \operatorname{epi}(g)$.

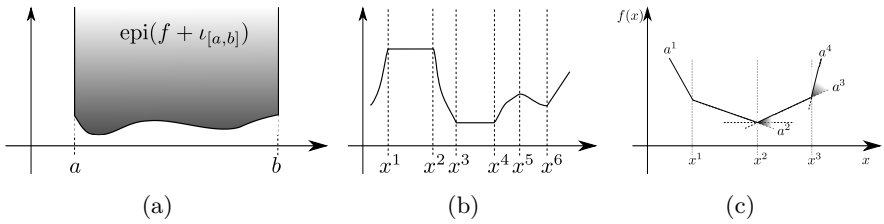

**Figure 5.1: (a)** Epigraph of the function $f + \iota_{[a,b]}$; **(b)** Illustration of local and global minima (maxima) of a function. All points in the closed interval $[x^1, x^2]$ are local maxima, points in the open interval $(x^1, x^2)$ are also local minima. Similarly, points $[x^3, x^4]$ are local and global minima, in $(x^3, x^4)$ are also local maxima. Points $x^5$ and $x^6$ are local maximum and minimum respectively. **(c)** Illustration of the subgradients of a piecewise linear function. The shadowed sectors denote the range of possible subgradients. In the point $x^3$ those are all hyperplanes, which lie (are convex combinations of) between $a^3$ and $a^4$. The point $x^2$ contains a horizontal hyperplane, which corresponds to the zero subgradient. Existence of such a hyperplane is a necessary and sufficient optimality condition for convex functions. Note that along with the zero subgradient there are others, that are convex combinations of $a^2$ and $a^3$.

*Proof.* The proof follows directly from the definition of the epigraph of a function: $t \geq f(x)$ and $t \geq g(x)$ together imply $t \geq \max\{f(x), g(x)\}$, and the other way around. $\qquad\square$

**Definition 5.6.** A function $f\colon \mathbb{R}^n \to \mathbb{R} \cup \{\infty\}$ is called *convex* if epi$f$ is a convex set.

**Definition 5.7.** A function $f\colon \mathbb{R}^n \to \mathbb{R} \cup \{-\infty\}$ is called *concave* if $(-f)$ is convex.

**Exercise 5.8.** Prove that a concave and convex function $f$ is *linear*, that is, it is representable as $\langle c, x \rangle + b$ for some $c \in \mathbb{R}^n$ and $b \in \mathbb{R}$.

**Exercise 5.9.** Prove that for any convex function $f$ the set $X = \{x\colon f(x) \leq 0\}$ is convex.

**Proposition 5.10.** A function $f\colon \mathbb{R}^n \to \mathbb{R} \cup \{\infty\}$ is convex if and only if the inequality:

$$f(px + (1-p)z) \leq pf(x) + (1-p)f(z) \qquad (5.4)$$

holds for any $x, z \in \mathbb{R}^n$ and $p \in (0, 1)$. Here we assume that $p \cdot \infty = \infty$ for any positive $p$, $\infty + x = \infty$, $-\infty + x = -\infty$ for any $x \in \mathbb{R}$ and $-(-\infty) = \infty$.

*Proof.* In case either $f(x)$ or $f(z)$ is infinite, the inequality trivially holds. Otherwise, $(x, f(x)) \in \text{epi} f$ as well as $(z, f(z)) \in \text{epi} f$. Since $\text{epi} f$ is a convex set, if and only if it contains also $(px + (1-p)z, pf(x) + (1-p)f(z))$. By the definition of $\text{epi} f$ this is equivalent to $f(px + (1 - p)z) \leq pf(x) + (1 - p)f(z)$. $\qquad\square$

**Remark 5.11.** The natural addition and multiplication rules for infinite values defined in Proposition 5.10 are considered as defined further in this monograph. Additionally $-\infty \cdot \infty = -\infty$ and the result of $\infty - \infty$ is undefined.

**Exercise 5.12.** Let $f \colon \mathbb{R}^n \to \mathbb{R} \cup \{\infty, -\infty\}$ be an extended value function such that $\text{dom} f \neq \emptyset$, and there is $x \in \mathbb{R}^n$ such that $f(x) = -\infty$. Show that $f$ is not convex.

**Exercise 5.13.** Show that a function $f \colon \mathbb{R}^n \to \mathbb{R} \cup \{\infty\}$ is convex if and only if the inequality

$$f \left( \sum_{i=1}^{N} p_i x^i \right) \leq \sum_{i=1}^{N} p_i f(x^i) \tag{5.5}$$

holds for any natural $N$, any $N$-tuple $x^i \in \mathbb{R}^n$, $i = 1, \ldots, N$, and any $p \in \Delta^N$.

The following proposition allows us to construct convex functions with linear operations applied to other convex functions:

**Proposition 5.14.** Let $f$ and $g$ be convex functions and $\alpha$ and $\beta$ be non-negative numbers. Then $\alpha f + \beta g$ is convex as well.

*Proof.* Let $p \in [0, 1]$. Then convexity of $f$ and $g$ implies the following sequence of inequalities, which proves the statement of the proposition:

$$\begin{aligned}
(\alpha f + \beta g)(px + (1-p)z) &= \alpha f(px + (1-p)z) + \beta g(px + (1-p)z) \\
&\leq \alpha(pf(x) + (1-p)f(z)) + \beta(pg(x) + \beta(1-p)g(z)) \\
&= p(\alpha f + \beta g)(x) + (1-p)(\alpha f + \beta g)(z).
\end{aligned} \tag{5.6}$$

$\qquad\square$

**Definition 5.15.** A function $h\colon \mathbb{R}^k \to \mathbb{R}$ is called *non-decreasing*, if $h(x_1, \ldots, x_k) \geq h(z_1, \ldots, z_k)$ as soon as $x_i \geq z_i$ for all $i = 1, \ldots, k$.

A superposition of convex functions is convex under an additional condition:

**Proposition 5.16.** Let $f(x) = h(f_1(x), \ldots, f_k(x))$ with $f_i\colon \mathbb{R}^n \to \mathbb{R}$ be convex functions, and $h\colon \mathbb{R}^k \to \mathbb{R}$ be convex and non-decreasing. Then $f$ is convex.

*Proof.* The following inequality holds for any $p \in [0,1]$ and arbitrary $x$ and $z$ due to convexity of $f_i$ and the non-decreasing property of $h$:

$$h(f_1(px + (1-p)z), \ldots, f_k(px + (1-p)z))$$
$$\leq h(pf_1(x) + (1-p)f_1(z), \ldots, pf_k(x) + (1-p)f_k(z)) . \quad (5.7)$$

Further, due to convexity of $h$, the right-hand-side of (5.7) does not exceed

$$ph(f_1(x), \ldots, f_k(x)) + (1-p)h(f_1(z), \ldots, f_k(z)) , \quad (5.8)$$

that finalizes the proof. $\qquad \square$

**Proposition 5.17.** Let $f$ be a convex function and $X \subset \mathbb{R}^n$ be a convex set. Then the function $f + \iota_X \colon \mathbb{R}^n \to \mathbb{R}$ is convex, where $\iota_X$ is defined as in (5.3).

*Proof.* Consider the convex function $\hat{\iota}_X(x) := \begin{cases} -\infty, & x \in X \\ \infty, & x \notin X \end{cases}$.

It is straightforward to see that $f + \iota_X = \max\{f, \hat{\iota}_X\}$. Due to Proposition 5.5, it holds that $\mathrm{epi}(f + \iota_X) = \mathrm{epi}(\max\{f, \hat{\iota}_X\}) = \mathrm{epi}(f) \cap \mathrm{epi}(\hat{\iota}_X)$, and, therefore, $\mathrm{epi}(f + \iota_X)$ is convex as it is the intersection of convex sets (see Lemma 3.26). $\qquad \square$

Note that a similar claim holds also for a *concave* function $f$ and a *convex* set $X$. In this case the function $f - \iota_X$ is concave.

### 5.1.3   Local and global minima

Let the set $B_\epsilon(x) = \{x' \in \mathbb{R}^n \mid \|x - x'\| \le \epsilon\}$ denote the ball of radius $\epsilon \ge 0$ around $x \in \mathbb{R}^n$.

**Definition 5.18** (Local and global minima). Let $f\colon X \to \mathbb{R} \cup \{\infty, -\infty\}$ be an extended value function. A point $x$ such that $f(x) < \infty$ is called a *local minimum* of $f$ if there exists $\epsilon > 0$ such that $x \in \arg\min_{x' \in B_\epsilon(x)} \left( f(x') + \iota_X(x') \right)$. The point $x$ is called a *global minimum* of $f$, if $x \in \arg\min_{x' \in X} f(x')$.

Definition 5.18 is illustrated by Figure 5.1(b).

The following proposition states one of the most important properties of convex functions for their minimization:

**Proposition 5.19.** Any local minimum of a convex function is also its global minimum.

*Proof.* Let $x$ be a local, but not a global minimum on $X \subseteq \mathbb{R}^n$. Therefore, there exists $z \in X$ such that $f(z) < f(x)$. Since $f$ is convex it holds that $f(px + (1 - p)z) \le pf(x) + (1 - p)f(z) < f(x)$ for any $p \in (0, 1)$. For any $\epsilon > 0$ there is $p > 0$ such that $(px + (1 - p)z) \in B_\epsilon$, therefore, $f(px + (1 - p)z) < f(x)$ implies that $x$ is not a local minimum.   $\square$

Due to Proposition 5.19, it makes sense to speak about *minima* of a convex function, without subdividing them into local and global ones.

**Proposition 5.20.** The set of minima of a convex function is convex.

*Proof.* Let $x$ and $z$ be two minima and, therefore, $f(x) = f(z)$. Then for all $p \in (0, 1)$ the inequality $f(px + (1 - p)z) \le pf(x) + (1 - p)f(z) = pf(x) + (1 - p)f(x) = f(x)$ implies that $px + (1 - p)z$ is a minimum as well.   $\square$

Note that the set of maxima of a concave function is *convex* as well.

The following important proposition essentially states that a maximum of convex functions is a convex function itself. This observation is crucial for properties of Lagrange relaxations, which will be considered later in this chapter.

**Proposition 5.21.** Let $Z$ be any set and let the function $f\colon \mathbb{R}^n \times Z \to \mathbb{R}$ be convex w.r.t. $x \in \mathbb{R}^n$ for each fixed $z \in Z$. Then the function $g(x) = \sup_{z \in Z} f(x, z)$ is convex.

*Proof.* For each fixed $z$ the epigraph of $f(\cdot, z)$ is a convex set. The epigraph of the function $g$ is the intersection of the epigraphs of the functions $f(\cdot, z)\colon \mathbb{R}^n \to \mathbb{R}$ for all $z \in Z$. Therefore, according to Lemma 3.26 epi$(g)$ is a convex set as an intersection of convex sets. So by Definition 5.6 $g$ is convex. $\qquad\square$

**Corollary 5.22.** The function $g(x) = \inf_{z \in Z} f(x, z)$ is concave if $f\colon \mathbb{R}^n \times Z \to \mathbb{R}$ is concave w.r.t. $x$ for each fixed $z \in Z$.

**Definition 5.23** (Convex optimization problem). Let $f$ be a convex function and $X$ be a convex set. The problems of the form $\min_{x \in X} f(x)$ and $\max_{x \in X}(-f(x))$ are called *convex optimization problems*. In other words, convex optimization is a minimization of a convex or maximization of a concave function on a convex set.

According to Proposition 5.17, for any convex $f$ and $X$ the problem $\min_{x \in \mathbb{R}^n}(f(x) + \iota_X(x))$ is convex and equivalent to $\min_{x \in X} f(x)$. Therefore, we conclude that the set of optimal solutions of this problem is convex.

Linear programs, as well as the optimization problems from Examples 3.1 and 3.3 are convex.

**Definition 5.24.** (Convex relaxation) The problem $\min_{x \in X'} g(x)$ is called a *convex relaxation* of the problem $\min_{x \in X} f(x)$, if it is both, a relaxation (that is, $X' \in X$ and $g(x) \leq f(x)$ for all $x \in X$), and a convex optimization problem.

## 5.2 Subgradient

Most of the functions we will deal with in this monograph will be non-smooth, or, more precisely, non-differentiable. In particular, this implies that they cannot be minimized with such a simple and well-known technique as gradient descent. To circumvent this, we will now introduce a generalization of the notion of gradient to non-differentiable functions. In later chapters we will also consider the corresponding generalizations of the gradient descent algorithm along with other non-smooth optimization methods.

**Definition 5.25.** Let $f$ be convex. A vector $g$ is called *a subgradient* of $f$ in the point $x^0 \in \mathrm{dom}\, f$ if for any $x \in \mathrm{dom}\, f$ it holds that

$$f(x) \geq f(x^0) + \left\langle g, x - x^0 \right\rangle . \tag{5.9}$$

The set $\partial f(x^0)$ of all subgradients of $f$ in $x^0$ is called the *subdifferential* of $f$ in $x^0$.

Definition 5.25 is illustrated in Figure 5.1(c).

For concave functions the inequality (5.9) changes its sign, i.e. if $g$ is a subgradient of $f$ it holds that $f(x) \leq f(x^0) + \langle g, x - x^0 \rangle$. Since $-f$ is concave if $f$ is convex, this implies that if $g$ is a subgradient of $f$, then $-g$ is a subgradient of $-f$. For concave functions the term *supergradient* is also often used instead of subgradient. This reflects the fact that the supergradient defines an upper bound of a concave function contrary to the subgradient defining a lower bound of a convex one.

**Proposition 5.26.** Let $f$ be a convex function. For any $x \in \mathrm{dom}\, f$ the set $\partial f(x)$ is convex.

*Proof.* Expression (5.9) defines a linear inequality w.r.t. $g$ for each $x^0 \in \mathrm{dom}\, f$. Any linear inequality defines a convex set, see Proposition 3.27. The intersection of these convex sets for all $x \in \mathrm{dom}\, f$ is also a convex set according to Lemma 3.26. $\qquad\square$

Proposition 5.26 is illustrated in Figure 5.1(c).

In practice, the following proposition is very important for computing subgradients:

**Proposition 5.27.** If $f$ is convex and differentiable in $x$, then $\partial f(x) = \{\nabla f(x)\}$.

We leave out the proof and here provide a suitable reference in §5.5.

Since subgradient is defined by a linear inequality, its linear properties follow directly from the definition:

**Proposition 5.28** (Linearity of subdifferential)**.** Let $f$ and $\hat{f}$ be convex functions. Then

$$\partial(\alpha f(x) + \beta \hat{f}(x)) = \alpha \partial f(x) + \beta \partial \hat{f}(x) \qquad (5.10)$$

for any $\alpha, \beta \in \mathbb{R}_+$ with $+$ denoting the Minkowski set sum. In particular, if $g \in \partial f(x)$ and $\hat{g} \in \partial \hat{f}(x)$, then $(\alpha g + \beta \hat{g}) \in \partial(\alpha f(x) + \hat{\beta} g(x))$.

Subgradient (but not subdifferential in general!) of function composition can be obtained with the same chain rule as gradient:

**Proposition 5.29** (Subgradient of function composition)**.** Let $f(x) = h(f_1(x), \dots, f_k(x))$, $h \colon \mathbb{R}^k \to \mathbb{R}$ be convex and non-decreasing and $f_i \colon \mathbb{R}^n \to \mathbb{R}$ be convex functions. Let also

$$g_i(x^0) \in \partial f_i(x^0) \quad \text{and} \quad z \in \partial h(f_1(x^0), \dots, f_k(x^0)). \qquad (5.11)$$

Then $g = (\sum_{i=1}^k z_i g_i) \in \partial f(x^0)$.

*Proof.* According to Proposition 5.16 the function $f$ is convex w.r.t. $x$. The proof of the proposition statement follows from the sequence of (in)equalities:

$$\begin{aligned}
f(x) &= h(f_1(x), \dots, f_k(x)) \\
&\geq h\left(f_1(x^0) + \left\langle g_1, x - x^0 \right\rangle, \dots, f_k(x^0) + \left\langle g_k, x - x^0 \right\rangle\right) \\
&\geq h(f_1(x^0), \dots, f_k(x^0)) + \left\langle z, \left(\left\langle g_1, x - x^0 \right\rangle, \dots, \left\langle g_k, x - x^0 \right\rangle\right)^\top \right\rangle \\
&= f(x^0) + \left\langle g, x - x^0 \right\rangle.
\end{aligned} \qquad (5.12)$$

$\square$

Therefore, the subgradient can be seen as a generalization of the gradient for convex functions. An important difference is, however, that for non-differentiable functions the subgradient is not unique, i.e. $\partial f$ may contain multiple elements in the points where the function is non-differentiable. Below, we consider an important class of such functions.

**Definition 5.30.** The function $f \colon \mathbb{R}^n \to \mathbb{R}$ is called *convex piecewise linear*, if for all $x \in \mathbb{R}^n$ it is representable as $f(x) = \max_{i \in I} \left\langle a^i, x \right\rangle + b^i$, with $I$ being a finite set and $a^i \in \mathbb{R}^n$, $b^i \in \mathbb{R}$ for $i = 1, \ldots, |I|$. By exchanging min with max above one obtains *concave piecewise linear* functions.

Note that a convex piecewise linear function is indeed convex as it is a point-wise maximum of convex (here, linear) functions, see Proposition 5.21.

Figure 5.1(c) as well as Example 5.31 below illustrate the notion of piecewise linear functions and their subgradients.

**Example 5.31.** Let $I$ be a finite set and let $f$ be a convex piecewise linear function $f(x) = \max_{i \in I} \langle a^i, x \rangle + b^i$. Then from $\langle a^i, x^0 \rangle + b^i = f(x^0)$ it follows $a^i \in \partial f(x^0)$ . Indeed,

$$
\begin{aligned}
f(x) - f(x^0) &= \max_{j \in I} \left( \left\langle a^j, x \right\rangle + b^j \right) - \left( \left\langle a^i, x^0 \right\rangle + b^i \right) \\
&\geq \left( \left\langle a^i, x \right\rangle + b^i \right) - \left( \left\langle a^i, x^0 \right\rangle + b^i \right) = \left\langle a^i, x - x^0 \right\rangle . \quad (5.13)
\end{aligned}
$$

The following lemma generalizes Example 5.31:

**Lemma 5.32.** Let functions $f_i \colon \mathbb{R}^n \to \mathbb{R}$, $i = 1, \ldots, m$ be convex and differentiable. Then the function $f(x) = \max_{i=1,\ldots,m} f_i(x)$ is convex and $\partial f(x) = \mathrm{conv}\{\nabla f_i(x) \colon i \in I(x)\}$, where $I(x) := \arg\max_{i=1,\ldots,m} f_i(x)$.

*Proof.* We will here concentrate on showing that

$$
\partial f(x) \subseteq \mathrm{conv}\{\nabla f_i(x) \colon i \in I(x)\}.
$$

For the reverse inclusion we refer to the corresponding literature in §5.5 as it requires more advanced analysis.

To prove the above inclusion we first show that $\nabla f_i(x^0) \in \partial f(x^0)$ for all $i \in I(x^0)$, $x^0 \in \mathbb{R}^n$. By Proposition 5.27 we have for all $x \in \mathbb{R}^n$:

$$
f_i(x) \geq f_i(x^0) + \left\langle \nabla f_i(x^0), x - x^0 \right\rangle . \quad (5.14)
$$

Therefore, by definition of $f$, it holds that

$$
f(x) \geq f(x^0) + \left\langle \nabla f_i(x^0), x - x^0 \right\rangle , \quad (5.15)
$$

and, hence, $\nabla f_i(x^0) \in \partial f(x^0)$.

As $\mathrm{conv}\{\nabla f_i(x) \colon i \in I(x)\}$ is the smallest convex set containing all vectors $\nabla f_i(x)$, $i \in I(x)$, and $\partial f(x)$ is convex by Proposition 5.26, this implies $\partial f(x) \subseteq \mathrm{conv}\{\nabla f_i(x) \colon i \in I(x)\}$.                    $\square$

In particular, for a convex piecewise linear function

$$f(x) = \max_{i \in I} \left\langle a^i, x \right\rangle + b^i$$

Lemma 5.32 states that $\partial f(x) = \text{conv}\{a^j \mid f(x) = \langle a^j, x \rangle + b^j, \ j \in J\}$. In other words, the subdifferential is a convex hull of all vectors $a^j$ such that $f(x) = \langle a^j, x \rangle + b^j$, see Figure 5.1(c).

**Example 5.33.** Consider $f(x) = \max\{0, 1 - x\}$.
The subdifferential of $f$ is $\partial f(x) = \begin{cases} -1, & x < 1 \\ [-1, 0], & x = 1 \\ 0, & x > 1 \end{cases}$.

A statement analogous to Lemma 5.32 can be formulated for concave functions by substituting "convex" with "concave" and max with min:

**Corollary 5.34.** Let functions $f_i \colon \mathbb{R}^n \to \mathbb{R}$, $i = 1, \ldots, m$ be concave and differentiable. Then the function $f(x) = \min_{i \in 1, \ldots, m} f_i(x)$ is concave and $\partial f(x) = \text{conv}\{\nabla f_i(x) \colon i \in I(x)\}$, where $I(x) := \arg\min_{i \in 1, \ldots, m} f_i(x)$.

Apart from its use in gradient-based optimization methods, the gradient of a differentiable convex function uniquely characterizes its optima. A similar result holds also for the subgradient:

**Theorem 5.1.** Let $f$ be convex. Then $f(x^*) = \min_{x \in \text{dom} f} f(x)$ if and only if $0 \in \partial f(x^*)$.

*Proof.* Indeed, if $0 \in \partial f(x^*)$, then $f(x) \geq f(x^*) + \langle 0, x - x^* \rangle = f(x^*)$ for all $x \in \text{dom} f$. On the other hand, if $f(x) \geq f(x^*)$ for all $x \in \text{dom} f$, then $f(x) \geq f(x^*) + \langle 0, x - x^* \rangle = f(x^*)$ and $0$ is a subgradient of $f$ in $x^*$. $\qquad\square$

Theorem 5.1 is illustrated in Figure 5.1(c).

## 5.3 Lagrange duality

In Chapter 3 we defined a general notion of relaxation and considered an important example thereof, the linear programming relaxations. The latter are constructed by substituting the integrality constraints with box constraints.

Here we will consider another powerful method for constructing relaxations, where some constraints are substituted with a specially constructed function, which is then added to the original objective.

Although these two types of relaxations seem to have nothing in common at the first glance, they are in fact often equivalent. At the end of the chapter we will consider this point and formulate sufficient conditions for this equivalence.

### 5.3.1   Dual problem

Let us consider an optimization problem of the form

$$\min_{x \in P} f(x) \tag{5.16}$$
$$\text{s.t. } Ax = b,$$

where $P \subseteq \mathbb{R}^n$. We will assume that omitting the equality constraints $Ax = b$ would make the optimization problem easy (or at least *easier*) to solve. Therefore, we write these constraints separately instead of including them into the set $P$. Below we consider a powerful technique, which allows us to partially get rid of these constraints to simplify optimization.

Consider the following function of a vector $\lambda \in \mathbb{R}^n$

$$\min_{x \in P} \left[ f(x) + \langle \lambda, Ax - b \rangle \right], \tag{5.17}$$

which we will call the *Lagrange dual function*. The expression $L(x, \lambda) := f(x) + \langle \lambda, Ax - b \rangle$ is referred to as the *Lagrangean*. One also says that the constraint $Ax = b$ is *dualized* or *relaxed*. Variables $\lambda$ are called *dual variables*. The dimensionality of $\lambda$ equals the number of rows of the matrix $A$. In other words, each elementary constraint gets a corresponding dual variable.

**Example 5.35.** Consider

$$\min_{x \in [0,1]^2} \langle c, x \rangle \tag{5.18}$$
$$\text{s.t. } x_1 = x_2. \tag{5.19}$$

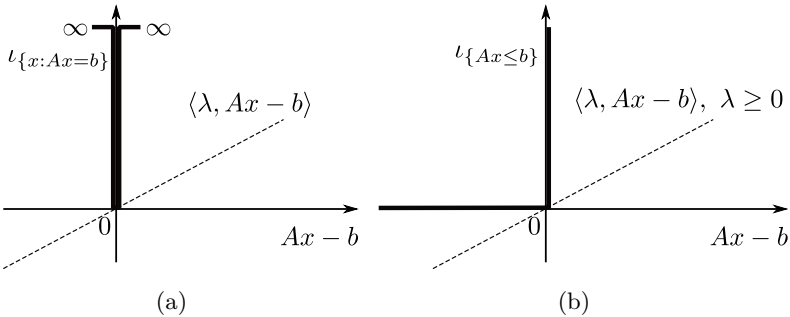

**Figure 5.2:** Illustration of the Lagrange relaxation for **(a)** the equality constraints $Ax = b$. Bold lines denote the function $\iota_{\{x:\, Ax=b\}}(x)$, see (5.3). The dashed line denotes the linear function $\langle \lambda, Ax - b \rangle$, which constitutes a lower bound for $\iota_{\{x:\, Ax=b\}}(x)$. **(b)** The same, but when the inequality $Ax \le b$ is dualized. Due to the, compared to Figure (a), different function $\iota_{\{x:\, Ax \le b\}}(x)$, only non-negative $\lambda$ must be considered to guarantee the lower bound property of $\langle \lambda, Ax - b \rangle$.

Its Lagrange dual function reads

$$\min_{x \in [0,1]^2} \langle c, x \rangle + \lambda(x_1 - x_2) \tag{5.20}$$

with $\lambda$ being a real number.

For any $\lambda$ the problem (5.17) is a relaxation of (5.16), because:

- The feasible set of (5.17) is a superset of the one for (5.16).

- The objective values of both problems are equal for any feasible point of (5.16), since in this case $Ax = b$ holds, and, therefore, $\langle \lambda, Ax - b \rangle = 0$.

These kinds of relaxations are called *Lagrange* relaxations.

The properties above imply:

**Proposition 5.36.** For every $\lambda \in \mathbb{R}$ the optimal value of (5.17) does not exceed the optimal value of (5.16), and they coincide if the constraint $Ax' = b$ holds in an optimum $x'$ of (5.17).

*Proof.* The first claim follows directly from Proposition 3.9. The second one follows from Proposition 3.9 and the fact that the objective functions of (5.16) and (5.17) coincide if $Ax' = b$. $\qquad\square$

**Exercise 5.37.** Show that for every $\lambda \in \mathbb{R}$ the problem

$$g(\lambda) = \min_{x \in P} f(x) + \langle \lambda, Ax - b \rangle \tag{5.21}$$

$$\text{s.t. } \lambda \geq 0 \tag{5.22}$$

is a relaxation of the problem

$$\min_{x \in P} f(x) \tag{5.23}$$

$$\text{s.t. } Ax \leq b. \tag{5.24}$$

In general, Lagrange relaxation is a way to construct a relaxation by dualizing any equality or inequality constraints, not necessarily linear ones. In these cases the constraint is defined as $h(x) = b$ (or $h(x) \leq b$), and is substituted by the term $\lambda(h(x) - b)$ (with $\lambda > 0$), which is then added to the objective. In this monograph we make use of the dualization of linear constraints only, therefore, we restrict our attention to this very important special case.

Figure 5.2 gives an intuitive explanation of the relaxation properties of the Lagrange duals. Indeed, let us move the constraint $Ax = b$ of the problem (5.16) into the objective. Then it becomes equivalent to

$$\min_{x \in P} f(x) + \iota_{\{x: \, Ax = b\}}(x). \tag{5.25}$$

To obtain the dual (5.17) we only substitute the function $\iota_{\{x: \, Ax = b\}}$ with $\langle \lambda, Ax - b \rangle$, and, as shown in Figure 5.2(a), the latter is a linear lower bound for the former one. The same reasoning can be used for the inequality case from Exercise 5.37 and is illustrated in Figure 5.2(b).

Since it is important to have the tightest possible lower bound, one considers the *Lagrange dual problem*

$$\max_{\lambda} \min_{x \in P} f(x) + \langle \lambda, Ax - b \rangle, \tag{5.26}$$

which amounts to maximizing (5.17) w.r.t. $\lambda$. One speaks about *full Lagrange dual* if all constraints are dualized and not only a subset of them.

The most important property of the Lagrange dual problem (5.26) is its concavity w.r.t. $\lambda$, which does not depend on the fact whether the primal problem (5.16) is convex or not:

**Proposition 5.38.** The dual function (5.17) is concave w.r.t. $\lambda$.

*Proof.* The proof follows from Corollary 5.22 and the fact that the objective of problem (5.17) is linear w.r.t. $\lambda$, and, therefore, concave. $\square$

The following example shows that although the dual function is concave, it may not be differentiable:

**Example 5.39** (Piecewise linear Lagrange dual)**.** Consider the Lagrange dual as in (5.16):

$$\min_{x \in P} \langle c, x \rangle + \langle \lambda, Ax - b \rangle \tag{5.27}$$

It is a concave piecewise linear function, when $P$ is

- a finite set, since (5.27) satisfies Definition 5.30 in this case.

- a polytope. In this case

$$\min_{x \in P} \langle c, x \rangle + \langle \lambda, Ax - b \rangle \overset{\text{Cor. } 3.45}{=} \min_{x \in \text{vrtx}(P)} \left\langle c + A^\top \lambda, x \right\rangle - \langle \lambda, b \rangle \,, \tag{5.28}$$

and the right-hand side satisfies Definition 5.30.

**Proposition 5.40.** Let $P$ be a polytope and the set $P \cap \{x \colon Ax = b\}$ be non-empty. Then there exist optimal primal and dual solutions. Moreover, the *strong duality* holds, i.e.

$$\min_{\substack{x \in P \\ Ax = b}} \langle c, x \rangle = \max_{\lambda} \min_{x \in P} \langle c, x \rangle + \langle \lambda, Ax - b \rangle \,. \tag{5.29}$$

*Proof.* Let us first show that there exist optimal primal and dual solutions. From $P$ being a polytope it follows that the primal feasible set $P \cap \{x \colon Ax = b\}$ is a polytope as well (see Proposition 3.18). Therefore, the optimal primal value $p^*$ is finite and attained in one of its vertices, according to Corollary 3.45.

Let $d^*$ denote the optimal dual value. The following sequence of inequalities, which uses Proposition 5.36,

$$\infty > p^* \geq d^* = \max_{\lambda} \min_{x \in P} \langle c, x \rangle + \langle \lambda, Ax - b \rangle \overset{\lambda = 0}{\geq} \min_{x \in P} \langle c, x \rangle > -\infty \tag{5.30}$$

implies that the dual optimal value $d^*$ is also finite.

Consider the dual objective

$$g(\lambda) := \min_{x \in P} \langle c, x \rangle + \langle \lambda, Ax - b \rangle \stackrel{\text{Cor. 3.45}}{=} \min_{x \in \text{vrtx}(P)} \langle c, x \rangle + \langle \lambda, Ax - b \rangle , \tag{5.31}$$

and the set $\hat{P} := \arg\min_{x \in \text{vrtx}(P)} \langle c, x \rangle + \langle \lambda, Ax - b \rangle$. Since $P$ is nonempty, $\hat{P}$ is so as well. The maximal dual value $d^*$ is attained in all vectors $\lambda$ satisfying

$$d^* = \langle c, x' \rangle + \langle \lambda, Ax' - b \rangle , \tag{5.32}$$

where $x'$ is some vector from $\hat{P}$. The set of $\lambda$ satisfying (5.32) is nonempty as soon as $Ax' - b = 0$ implies $d^* = \langle c, x' \rangle$. This implication is given by the following sequence of equalities:

$$d^* = \max_\lambda \min_{x \in \text{vrtx}(P)} \langle c, x \rangle + \langle \lambda, Ax - b \rangle = \max_\lambda \langle c, x' \rangle + \langle \lambda, Ax' - b \rangle$$
$$\stackrel{Ax'=b}{=} \max_\lambda \langle c, x' \rangle = \langle c, x' \rangle . \tag{5.33}$$

Therefore, the dual solution exists. Let $\lambda^*$ be such a solution.

According to Lemma 5.32,

$$\partial g(\lambda^*) = \text{conv}\{Ax^* - b \colon x^* \in \hat{P}\} \stackrel{\text{Lem. 3.36}}{=} \{Ax^* - b \colon x^* \in \text{conv}(\hat{P})\} . \tag{5.34}$$

Since $\lambda^*$ is a dual optimal solution, Theorem 5.1 implies that $\partial g(\lambda^*) \ni 0$. In other words, there exists $x^* \in \text{conv}(\hat{P})$ such that

$$Ax^* - b = 0.$$

Observe that $x^* \in P$, since $\hat{P} \subseteq \text{vrtx}(P)$ and, therefore, $\text{conv}(\hat{P}) \subseteq \text{conv}(\text{vrtx}(P)) = P$. This yields

$$d^* = \langle c, x^* \rangle + \langle \lambda^*, Ax^* - b \rangle \stackrel{Ax^*-b=0}{=} \langle c, x^* \rangle \stackrel{\substack{x^* \in P \\ Ax^*-b=0}}{\geq} \min_{\substack{x \in P \\ Ax=b}} \langle c, x \rangle = p^* . \tag{5.35}$$

Together with (5.30) this implies $p^* = d^*$. $\qquad\square$

**Example 5.41.** Consider Example 5.35, the primal

$$\min_{x \in [0,1]^2} \langle c, x \rangle \tag{5.36}$$

$$\text{s.t. } x_1 = x_2 \,, \tag{5.37}$$

and the dual problem:

$$\max_{\lambda} \min_{x \in [0,1]^2} [\langle c, x \rangle + \lambda(x_1 - x_2)]$$
$$= \max_{\lambda} \min_{x \in [0,1]^2} [(c_1 + \lambda)x_1 + (c_2 - \lambda)x_2] \,. \tag{5.38}$$

The primal problem can be equivalently rewritten as

$$\min_{x_1 \in [0,1]} (c_1 + c_2)x_1 = \begin{cases} 0, & (c_1 + c_2) \geq 0 \\ (c_1 + c_2), & (c_1 + c_2) < 0 \,. \end{cases} \tag{5.39}$$

Consider the dual problem (5.38). It can be equivalently rewritten as

$$\max_{\lambda} \left( \min\{0, (c_1 + \lambda)\} + \min\{0, (c_2 - \lambda)\} \right) \,. \tag{5.40}$$

The necessary and sufficient optimality condition is the existence of the zero subgradient of the function. According to Lemma 5.32, this condition holds when either

$$c_1 + \lambda \geq 0 \text{ and } c_2 - \lambda \geq 0 \,, \tag{5.41}$$

or

$$c_1 + \lambda < 0 \text{ and } c_2 - \lambda < 0 \,. \tag{5.42}$$

hold. One pair of these inequalities is always attained when $c_1 + \lambda = c_2 - \lambda$, therefore, $\lambda = \frac{c_2 - c_1}{2}$. Substituting $\lambda$ with this value in (5.41) we obtain $c_1 + c_2 \geq 0$. Similarly (5.42) turns into $c_1 + c_2 < 0$.

Summing up,

$$\max_{\lambda} \left( \min\{0, (c_1 + \lambda)\} + \min\{0, (c_2 - \lambda)\} \right) \tag{5.43}$$

$$= \begin{cases} 0, & (c_1 + c_2) \geq 0 \\ (c_1 + c_2), & (c_1 + c_2) < 0 \,, \end{cases} \tag{5.44}$$

which coincides with the primal optimal value.

## 5.4   Lagrange relaxation of integer linear programs

Since we are interested in solving integer linear programs, let us consider the Lagrange relaxation of this kind of problem.

In what follows we will deal with integer linear programs of the form

$$\min_{x \in P \cap \{0,1\}^n} \langle c, x \rangle \qquad (5.45)$$

$$\text{s.t. } Ax = b,$$

where $P$ is a polytope. The assumption that $P$ is bounded (i.e. it is a polytope and not a general polyhedron) can be made w.l.o.g., since the problem (5.45) can be turned into an equivalent one with a bounded $P$ by adding the constraints $x \in [0,1]^n$.

As before, we assume that without the constraint $Ax = b$ the problem (5.45) becomes easy or at least easier to solve. Therefore, these constraints are separated out from the polytope $P$. Dualizing these constraints we obtain the Lagrange dual problem to (5.45):

$$\max_{\lambda} \min_{x \in P \cap \{0,1\}^n} \langle c, x \rangle + \langle \lambda, Ax - b \rangle . \qquad (5.46)$$

Note that objective of the dual problem (5.46) is a concave piecewise linear, therefore, non-differentiable, function.

### 5.4.1   Primal of the relaxed problem for ILPs

The Lagrange dual problem (5.26) can be treated as a problem of selecting the best relaxation from a given set of relaxations. This class is parametrized with the dual vector $\lambda$. The solution of the dual problem is a value of this parameter. To obtain the corresponding relaxed solution one has to solve the related relaxed problem $\min_{x \in X} f(x) + \langle \lambda, Ax - b \rangle$ with fixed $\lambda$. In this problem, however, both, the objective and the initial feasible sets are changed compared to the primal problem (5.45). This makes analysis of the relaxed solution much more complicated.

This differs from the case of the LP relaxation for integer linear programs, where only the feasible set is increased, but the objective itself remains unchanged. Interestingly, for integer linear programs, one can define a relaxation, which is equivalent to the Lagrange relaxation, by changing only the feasible set.

To do so, let us consider the Lagrange dual (5.46). Corollary 3.44 implies

$$\min_{x \in P \cap \{0,1\}^n} \langle c, x \rangle + \langle \lambda, Ax - b \rangle = \min_{x \in \mathrm{conv}(P \cap \{0,1\}^n)} \langle c, x \rangle + \langle \lambda, Ax - b \rangle . \tag{5.47}$$

Consider now the following minimization problem

$$\min_{\substack{x \in \mathrm{conv}(P \cap \{0,1\}^n) \\ Ax = b}} \langle c, x \rangle . \tag{5.48}$$

Assuming that the problem (5.45) is feasible implies that (5.48) is so, too. Since the feasible set of (5.48) is a polytope, Proposition 5.40 implies that its Lagrange dual is tight:

$$\min_{\substack{x \in \mathrm{conv}(P \cap \{0,1\}^n) \\ Ax = b}} \langle c, x \rangle = \max_{\lambda} \min_{x \in \mathrm{conv}(P \cap \{0,1\}^n)} \langle c, x \rangle + \langle \lambda, Ax - b \rangle . \tag{5.49}$$

Comparing it to (5.47) we conclude

**Proposition 5.42.**

$$\max_{\lambda} \min_{x \in P \cap \{0,1\}^n} \langle c, x \rangle + \langle \lambda, Ax - b \rangle = \min_{\substack{x \in \mathrm{conv}(P \cap \{0,1\}^n) \\ Ax = b}} \langle c, x \rangle . \tag{5.50}$$

In other words, the Lagrange dual is tight for the problem (5.48). Therefore, we will call the latter the *primal relaxed* problem for the Lagrange relaxation.

Note that (5.48)) is also a convex relaxation of the ILP problem (5.45), since $\mathrm{conv}(P \cap \{0,1\}^n) \supseteq P \cap \{0,1\}^n$ by definition of the convex hull. Contrary to the Lagrange relaxation, the problem (5.48) modifies only the feasible set of the non-relaxed problem (5.45), and keeps its objective unchanged. Note also that (5.48) is a linear program, since its feasible set is a polytope. It differs from the LP relaxation of (5.45)

$$\min_{\substack{x \in P \cap [0,1]^n \\ Ax = b}} \langle c, x \rangle \tag{5.51}$$

by the feasible set only. The following statement compares these two linear programs and establishes conditions under which they coincide.

**Corollary 5.43.** Let $\text{vrtx}(P \cap [0,1]^n) \subseteq \{0,1\}^n$, that is, all vertices of the polytope $P \cap [0,1]^n$ are defined by binary vectors. Then the primal relaxed problem (5.48) is equal to the LP relaxation (5.51). Otherwise, the relaxation (5.48) is tighter.

Since the primal relaxed (5.48) and dual (5.46) problems provide the same bound, one also says that *the Lagrange relaxation is in general tighter than the LP one*.

*Proof.* According to Proposition 3.55, it holds that $\text{vrtx}(P \cap [0,1]^n) \supseteq P \cap \{0,1\}^n$. According to Proposition 3.57, it holds $\text{vrtx}(P \cap [0,1]^n) \cap \{0,1\}^n = P \cap \{0,1\}^n$. Therefore, $\text{vrtx}(P \cap [0,1]^n) \subseteq \{0,1\}^n$ implies $\text{vrtx}(P \cap [0,1]^n) = P \cap \{0,1\}^n$. From Corollary 3.37 it follows that $\text{conv}(P \cap \{0,1\}^n) = P \cap [0,1]^n$, and, therefore, the considered Lagrange and LP relaxations are equivalent.

In general, however, it holds that $\text{conv}(P \cap \{0,1\}^n) \subseteq P \cap [0,1]^n$ (see Proposition 3.54), which implies that relaxation (5.48) is tighter than the LP one (5.51), as their objective functions coincide.    $\square$

### 5.4.2   Reparametrization

For (integer) linear programs the Lagrange relaxation has a simple and intuitive interpretation, expressed in terms of *reparametrization* or *equivalent transformation* of the cost vector. The optimization of the dual problem can be viewed as finding such an equivalent transformation, such that the dualized constraints are fulfilled for the (relaxed) primal problem.

Till the end of this chapter we will assume that the dualized constraints have the form $Ax = 0$. Although the more general case $Ax = b$ can be treated similarly, the resulting expressions for the case $b = 0$ are simpler, and, what is even more important, they correspond to the form of the constraints of the MAP-inference problem in its ILP formulation (4.7).

Let us consider the Lagrange dual problem (5.46). For $b = 0$ its objective is

$$\langle c, x \rangle + \langle \lambda, Ax \rangle = \langle c, x \rangle + \left\langle A^\top \lambda, x \right\rangle = \left\langle c + A^\top \lambda, x \right\rangle . \qquad (5.52)$$

Note that for any feasible $x$, it holds that $Ax = 0$ and therefore

$$\langle c, x \rangle = \left\langle c + A^\top \lambda, x \right\rangle . \tag{5.53}$$

It implies that the problem

$$\min_{x \in P \cap \{0,1\}^n} \left\langle c + A^\top \lambda, x \right\rangle \tag{5.54}$$

$$\text{s.t. } Ax = 0 \tag{5.55}$$

is *equivalent* to (5.45) (with $b = 0$), i.e. for any $\lambda$ and any feasible $x$ the objectives of both problems have equal values. Clearly, the solutions of both problems are therefore also equal. The transformation $c \to c + A^\top \lambda$ is therefore called *an equivalent transformation* or *a reparametrization* of the problem (5.45). The dual problem

$$\max_{\lambda} \left[ D(\lambda) := \min_{x \in P \cap \{0,1\}^n} \left\langle c + A^\top \lambda, x \right\rangle \right] \tag{5.56}$$

therefore consists in finding an *optimal reparametrization* of the problem.

**Example 5.44.** In Example 5.41 the reparametrized cost vector has coordinates $c_1 + \lambda$ and $c_2 - \lambda$.

Note that the equality (5.53) implies also that the *relaxed* primal problem, defined by the right-hand-side of (5.50), is equivalent to its reparametrization for any $\lambda$:

$$\min_{\substack{\text{conv}(x \in P \cap \{0,1\}^n) \\ Ax=0}} \langle c, x \rangle = \min_{\substack{\text{conv}(x \in P \cap \{0,1\}^n) \\ Ax=0}} \left\langle c + A^\top \lambda, x \right\rangle . \tag{5.57}$$

### 5.4.3 Primal-dual optimality conditions

In the following, we formulate necessary and sufficient conditions for dual optimality. These conditions play an important role for constructing a stopping criterion for iterative algorithms as well as for reconstruction of primal solutions from the dual ones, see §6.2.

**Proposition 5.45.** Vector $\lambda^*$ is an optimum of the dual problem (5.56) if and only if there exists

$$x^* \in \underset{x \in \text{conv}(P \cap \{0,1\}^n)}{\arg\min} \left\langle c + A^\top \lambda^*, x \right\rangle \tag{5.58}$$

such that $Ax^* = 0$. In other words, $\lambda^*$ is a dual optimum if and only if there exists $x^*$ satisfying (5.58), which is feasible for the relaxed primal problem defined by the right-hand-side of (5.50). Such an $x^*$ is a solution of the relaxed primal problem.

*Proof.* Due to (3.29), the objective of the dual problem (5.56) can be equivalently expressed as

$$g(\lambda) := \min_{x \in \mathrm{conv}(P \cap \{0,1\}^n)} \left\langle c + A^\top \lambda, x \right\rangle . \qquad (5.59)$$

Since $g$ is concave, its necessary and sufficient optimality condition is the existence of the zero vector in its subdifferential which is the set

$$\partial g(\lambda^*) = \left\{ Ax^* : x^* \in \operatorname*{arg\,min}_{x \in \mathrm{conv}(P \cap \{0,1\}^n)} \left\langle c + A^\top \lambda^*, x \right\rangle \right\}, \qquad (5.60)$$

see the proof of Proposition 5.40 for details.

Therefore, since $Ax^* = 0$, it holds that

$$
\begin{aligned}
x^* \in \arg & \min_{\substack{x \in \mathrm{conv}(P \cap \{0,1\}^n) \\ Ax = 0}} \left\langle c + A^\top \lambda^*, x \right\rangle \\
= \arg & \min_{\substack{x \in \mathrm{conv}(P \cap \{0,1\}^n) \\ Ax = 0}} \langle c, x \rangle + \langle \lambda^*, Ax \rangle \\
= \arg & \min_{\substack{x \in \mathrm{conv}(P \cap \{0,1\}^n) \\ Ax = 0}} \langle c, x \rangle .
\end{aligned}
\qquad (5.61)
$$

Hence, $x^*$ is a solution of the relaxed primal problem, defined by the right-hand-side of (5.50). $\qquad \square$

**Remark 5.46.** Note that if, in addition to (5.58) and $Ax^* = 0$, it is $x^* \in \{0,1\}^n$, then $x^*$ is a solution of the non-relaxed ILP problem (5.45). This directly follows from the fact that problem (5.58) is a relaxation of (5.45), $x^*$ is its minimizer, and $\langle c, x \rangle = \left\langle c + A^\top \lambda^*, x \right\rangle$ holds since $Ax^* = 0$.

## 5.5   Bibliography and further reading

For a more fundamental understanding the notion of convex functions and their subgradients we refer to the textbook [76]. In particular, it contains proofs of Proposition 5.27 and Lemma 5.32. The comprehensive source on the topic is [90].

For an extended introduction to the duality theory we recommend the excellent textbook on convex optimization [15]. It also details the notion of the extended value functions and their relation to convex analysis and Lagrange duality. More general and deep analysis of the mathematical phenomenon of duality is presented in [89].

# 6

---

# Lagrange Duality for MAP-inference

---

This chapter is devoted to the Lagrangean dual of the MAP-inference problem. This dual is obtained by relaxing the coupling constraints of the local polytope. Based on the theoretical background provided in Chapter 5 we will analyze properties of the dual. In particular, we will show that it is equivalent to the local polytope relaxation, that is, it provides the same lower bound.

The second important topic we deal with in this chapter is the optimality conditions for the considered dual problem. We will give exact (necessary and sufficient) as well as approximate (only necessary) optimality conditions. Whereas checking the first ones can be as computationally expensive as optimizing the dual itself, the latter can be verified with a simple and efficient algorithm which we will also describe.

In the last section we return to the acyclic graphical models introduced in Chapter 2, and show that the approximate optimality conditions are exact for such models.

## 6.1 Reparametrization and Lagrange dual

Recall the ILP formulation of the MAP-inference problem

$$\min_{\mu \in \mathcal{L} \cap \{0,1\}^{\mathcal{I}}} \langle \theta, \mu \rangle \tag{6.1}$$

with the local polytope $\mathcal{L}$ defined by (4.6). Note that $\mathcal{L}$ can be rewritten as follows:

$$\mathcal{L} = \begin{cases} \mu_u \in \Delta^{\mathcal{Y}_u}, & \forall u \in \mathcal{V} \\ \mu_{uv} \in \Delta^{\mathcal{Y}_{uv}}, & \forall uv \in \mathcal{E} \\ \sum_{t \in \mathcal{Y}_v} \mu_{uv}(s,t) = \mu_u(s), & \forall u \in \mathcal{V}, \ v \in \mathcal{N}_b(u), \ s \in \mathcal{Y}_u\,. \end{cases} \tag{6.2}$$

As in Chapter 4, the notation $\mu_u$ stands for the vector $(\mu_u(s) \colon s \in \mathcal{Y}_u)$, which encodes the selected label in the node $u$, and $\mu_{uv} = (\mu_{uv}(s,t) \colon (s,t) \in \mathcal{Y}_{uv})$ encodes the selected label pair in the edge $uv$. Let us construct the Lagrange dual to (6.1) by relaxing the coupling constraints (the last line in (6.2)). This is done similarly to the general scheme provided in §5.3.1.

Consider the linear term of the objective, which corresponds to dualizing the coupling constraints. In other words, we will specify the term $\langle \lambda, Ax \rangle$ corresponding to the constraints $Ax = 0$, when the role of the latter is played by the coupling constraints. Since for each node $u \in \mathcal{V}$ there is one constraint for each of its labels $s \in \mathcal{Y}_u$ and each neighboring node $v \in \mathcal{N}_b(u)$, the role of the dual vector $\lambda$ is played by the vector $\phi \in \mathbb{R}^{\mathcal{J}}$ with coordinates $\phi_{u,v}(s)$, where $\mathcal{J} := \{(u,v,s) \mid u \in \mathcal{V}, \ v \in \mathcal{N}_b(u), \ s \in \mathcal{Y}_u\}$. The comma between $u$ and $v$ in the lower index $u,v$ underlines that $\phi_{u,v}(s)$ and $\phi_{v,u}(s)$ are two different vectors, contrary to $\mu_{uv}(s,t)$ and $\mu_{vu}(t,s)$, where $u$ and $v$ are not separated by the comma.

Given the above notation, the considered linear term corresponding to the coupling constraints reads:

$$\sum_{u \in \mathcal{V}} \sum_{v \in \mathcal{N}_b(u)} \sum_{s \in \mathcal{Y}_u} \phi_{u,v}(s) \left( \sum_{t \in \mathcal{Y}_v} \mu_{uv}(s,t) - \mu_u(s) \right)\,. \tag{6.3}$$

By adding it to the objective

$$\langle \theta, \mu \rangle = \sum_{u \in \mathcal{V}} \sum_{s \in \mathcal{Y}_u} \theta_u(s)\mu_u(s) + \sum_{uv \in \mathcal{E}} \sum_{(s,t) \in \mathcal{Y}_{uv}} \theta_{uv}(s,t)\mu_{uv}(s,t) \tag{6.4}$$

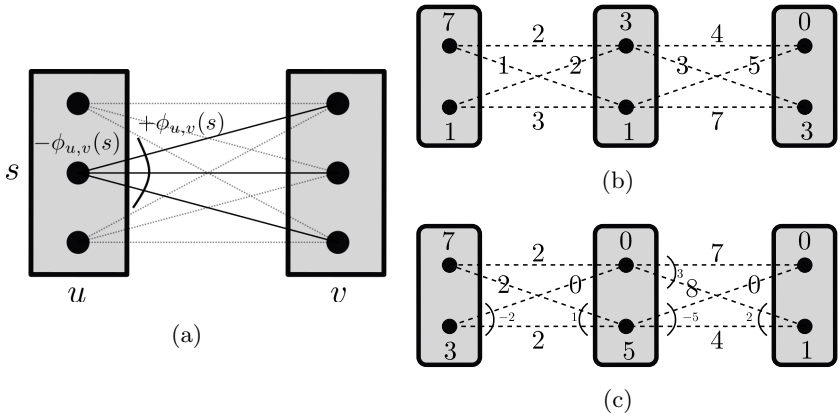

**Figure 6.1: (a)** Illustration of the reparametrization. Energy of all labelings remains the same when the value $\phi_{u,v}(s)$ is subtracted from the cost $\theta_u(s)$ of the label $s$ in the node $u$ and simultaneously added to the costs of all label pairs $(s,t)$, $t \in \mathcal{Y}_v$, associated with one of the edges incident to $u$. **(b-c)** Example of the reparametrization: (b) initial costs, (c) reparametrized costs. Bigger numbers stand for the costs (initial in (b) and reparametrized in (c)), smaller numbers in (c) define values of the reparametrization vector $\phi$ with the same meaning as in (a).

and regrouping terms one obtains

$$\sum_{u\in\mathcal{V}}\sum_{s\in\mathcal{Y}_u}\mu_u(s)\left(\theta_u(s)-\sum_{v\in\mathcal{N}_b(u)}\phi_{u,v}(s)\right)$$

$$+\sum_{uv\in\mathcal{E}}\sum_{(s,t)\in\mathcal{Y}_{uv}}\mu_{uv}(s,t)\left(\theta_{uv}(s,t)+\phi_{u,v}(s)+\phi_{v,u}(t)\right)=\left\langle\theta^\phi,\mu\right\rangle,$$
(6.5)

where we introduced the *reparametrized costs* $\theta^\phi$ defined as

$$\theta_u^\phi(s):=\theta_u(s)-\sum_{v\in\mathcal{N}_b(u)}\phi_{u,v}(s),\ v\in\mathcal{V},\ s\in\mathcal{Y}_u,\qquad(6.6)$$

$$\theta_{uv}^\phi(s,t):=\theta_{uv}(s,t)+\phi_{u,v}(s)+\phi_{v,u}(t),\ uv\in\mathcal{E},\ (s,t)\in\mathcal{Y}_{uv}.$$

See Figure 6.1 for an illustration of the reparametrization.

Note that the reparametrized Lagrangean $\left\langle\theta^\phi,\mu\right\rangle$ is an instance of the reparametrization introduced in §5.4.2 for general (integer) linear programs, see Expression (5.56). In other words, the Lagrangian dual constructed by relaxing the coupling constraints consists in finding an optimal reparametrization:

$$\max_{\phi \in \mathbb{R}^{\mathcal{J}}} \underbrace{\min_{\substack{\mu \in \{0,1\}^{\mathcal{I}} \\ \mu_u \in \Delta^{\mathcal{Y}_u}, u \in \mathcal{V} \\ \mu_{uv} \in \Delta^{\mathcal{Y}_{uv}}, uv \in \mathcal{E}}} \left\langle \theta^\phi, \mu \right\rangle}_{\text{dual function } \mathcal{D}(\phi)} . \tag{6.7}$$

Note that the constraints on $\mu$ decouple into the coordinates indexed by $u$ and $uv$:

$$\min_{\substack{\mu \in \{0,1\}^{\mathcal{I}} \\ \mu_u \in \Delta^{\mathcal{Y}_u}, u \in \mathcal{V} \\ \mu_{uv} \in \Delta^{\mathcal{Y}_{uv}}, uv \in \mathcal{E}}} \left[ \left\langle \theta^\phi, \mu \right\rangle = \sum_{u \in \mathcal{V}} \left\langle \theta_u^\phi, \mu_u \right\rangle + \sum_{uv \in \mathcal{E}} \left\langle \theta_{uv}^\phi, \mu_{uv} \right\rangle \right]$$

$$= \min_{\substack{\mu_u \in \Delta^{\mathcal{Y}_u \cap \{0,1\}^{\mathcal{Y}_u}} \\ u \in \mathcal{V}}} \sum_{u \in \mathcal{V}} \left\langle \theta_u^\phi, \mu_u \right\rangle + \min_{\substack{\mu_{uv} \in \Delta^{\mathcal{Y}_{uv} \cap \{0,1\}^{\mathcal{Y}_{uv}}} \\ uv \in \mathcal{E}}} \sum_{uv \in \mathcal{E}} \left\langle \theta_{uv}^\phi, \mu_{uv} \right\rangle$$

$$= \sum_{u \in \mathcal{V}} \min_{\mu_u \in \Delta^{\mathcal{Y}_u \cap \{0,1\}^{\mathcal{Y}_u}}} \left\langle \theta_u^\phi, \mu_u \right\rangle + \sum_{uv \in \mathcal{E}} \min_{\mu_{uv} \in \Delta^{\mathcal{Y}_{uv} \cap \{0,1\}^{\mathcal{Y}_{uv}}}} \left\langle \theta_{uv}^\phi, \mu_{uv} \right\rangle$$

$$= \sum_{u \in \mathcal{V}} \min_{s \in \mathcal{Y}_u} \theta_u^\phi(s) + \sum_{uv \in \mathcal{E}} \min_{(s,t) \in \mathcal{Y}_{uv}} \theta_{uv}^\phi(s,t) . \tag{6.8}$$

Therefore, the dual problem (6.7) takes the form

$$\max_{\phi \in \mathbb{R}^{\mathcal{J}}} \mathcal{D}(\phi) := \max_{\phi \in \mathbb{R}^{\mathcal{J}}} \left( \sum_{u \in \mathcal{V}} \min_{s \in \mathcal{Y}_u} \theta_u^\phi(s) + \sum_{uv \in \mathcal{E}} \min_{(s,t) \in \mathcal{Y}_{uv}} \theta_{uv}^\phi(s,t) \right) . \tag{6.9}$$

Note that to compute the value of the dual objective in (6.9) for a fixed $\phi$ one has to select *independently* in each node and each edge an optimal label and an optimal label pair, respectively. Such labels and label pairs will be called *locally optimal* in the following. Note also that these locally optimal labels and label pairs need not be consistent. In other words, if the label $s$ is selected in node $u$ and the label pair $(t, t')$ in edge $uv$, this does not generally imply that $s = t$. However, as we will show later in this chapter, some kind of relaxed consistency (known as *arc-consistency*) is enforced in the optimum of the dual function.

**Example 6.1** (How large-scale is the Lagrangean dual of the MAP-inference?)**.** Consider now the Lagrangean dual (6.9). The number of variables is equal to $|\mathcal{J}|$, that grows linearly with the number of labels and the number of edges in a graph. This is in contrast to the primal problem (6.1), where the number of variables grows quadratically with the number of labels.

Consider Example 4.2, where the size of the ILP representation of the MAP-inference problem was estimated for the depth reconstruction problem from Example 1.2. In that case the dual problem has $O(10^5)$ variables, which is an order of magnitude less than the size $O(10^6)$ of the primal problem.

The following proposition summarizes properties of the dual problem, which specialize the general properties of the Lagrangean relaxation from §5.

**Proposition 6.2.** For the dual problem (6.9) it holds that:

1. $\left\langle \theta^\phi, \mu \right\rangle = \langle \theta, \mu \rangle$ holds for any $\phi \in \mathbb{R}^\mathcal{J}$ and $\mu \in \mathcal{L}$ (in particular for $\mu \in \mathcal{L} \cap \{0,1\}^\mathcal{I}$ and $\mu \in \mathcal{M}$). This implies that for any labeling $y \in \mathcal{Y}_\mathcal{V}$ it holds that $E(y; \theta) = E(y; \theta^\phi)$ with $E$ being the energy of $y$ as defined in (1.4).

2. $\mathcal{D}(\phi)$ is a lower bound for the energy minimization, i.e. $\mathcal{D}(\phi) \leq \langle \theta, \mu \rangle$ for any $\mu \in \mathcal{L} \cap \{0,1\}^\mathcal{I}$ and $\phi \in \mathbb{R}^\mathcal{J}$.

3. $\mathcal{D}$ is concave piecewise linear, and, therefore, a non-differentiable function.

4. The primal relaxed problem corresponding to (6.9) is the local polytope relaxation, i.e. $\max_{\phi \in \mathbb{R}^\mathcal{J}} \mathcal{D}(\phi) = \min_{\mu \in \mathcal{L}} \langle \theta, \mu \rangle$.

5. The optimality condition for the dual $\mathcal{D}(\phi)$ reads:

$$\phi \in \arg\max_{\phi' \in \mathbb{R}^\mathcal{J}} \mathcal{D}(\phi')$$

   if and only if the polytope

$$\mathcal{L}(\phi) := \Big\{ \mu \in \mathcal{L} \colon \mu_w(s) = 0 \text{ if } \theta^\phi_w(s) > \min_{s' \in \mathcal{Y}_w} \theta^\phi_w(s'),$$
$$w \in \mathcal{V} \cup \mathcal{E}, \ s \in \mathcal{Y}_w \Big\} \quad (6.10)$$

is non-empty. In other words, there must be a relaxed labeling $\mu \in \mathcal{L}$ with non-zero coordinates assigned only to the locally optimal labels and label pairs. Such a relaxed labeling is also a solution to the local polytope relaxation $\min_{\mu \in \mathcal{L}} \langle \theta, \mu \rangle$.

6. Tightness of the Lagrange dual: Let the set $\mathcal{L}(\phi)$ defined by (6.10) contain an integer labeling, in other words, there is $y \in \mathcal{Y}_{\mathcal{V}}$ such that $\delta(y) \in \mathcal{L}(\phi)$. Then $y$ is the solution of the (non-relaxed) energy minimization problem (1.4). In this case $\mathcal{D}(\phi) = E(y; \theta)$. We will call this case *LP-tight*, since the Lagrange dual is equivalent to the local polytope relaxation.

   This statement holds also in the opposite direction, that is, $\mathcal{D}(\phi) = E(y; \theta)$ if and only if there is $y \in \mathcal{Y}_{\mathcal{V}}$ such that $\delta(y) \in \mathcal{L}(\phi)$.

*Proof.*

1. Equality follows from the fact that for any $\mu \in \mathcal{L}$ the coupling constraints hold, and, therefore, $\left\langle \theta^\phi, \mu \right\rangle$ is a reparametrization of $\langle \theta, \mu \rangle$, as introduced in §5.4.2. Alternatively, this fact can be directly shown by observing that

$$
\begin{aligned}
\left\langle \theta^\phi, \mu \right\rangle &= \langle \theta, \mu \rangle + \sum_{u \in \mathcal{V}} \sum_{s \in \mathcal{Y}_u} \mu_u(s) \left( - \sum_{v \in \mathcal{N}_b(u)} \phi_{u,v}(s) \right) \\
&\quad + \sum_{uv \in \mathcal{E}} \sum_{(s,t) \in \mathcal{Y}_{uv}} \mu_{uv}(s,t) \left( \phi_{u,v}(s) + \phi_{v,u}(t) \right) \\
&= \langle \theta, \mu \rangle - \sum_{u \in \mathcal{V}} \sum_{v \in \mathcal{N}_b(u)} \sum_{s \in \mathcal{Y}_u} \phi_{u,v}(s) \mu_u(s) \\
&\quad + \sum_{u \in \mathcal{V}} \sum_{v \in \mathcal{N}_b(u)} \sum_{s \in \mathcal{Y}_u} \phi_{u,v}(s) \sum_{t \in \mathcal{Y}_v} \mu_{uv}(s,t) \\
&= \langle \theta, \mu \rangle + \sum_{u \in \mathcal{V}} \sum_{v \in \mathcal{N}_b(u)} \sum_{s \in \mathcal{Y}_u} \phi_{u,v}(s) \underbrace{\left( \sum_{t \in \mathcal{Y}_v} \mu_{uv}(s,t) - \mu_u(s) \right)}_{=0} \\
&= \langle \theta, \mu \rangle \quad\quad\quad\quad\quad\quad\quad\quad\quad\quad\quad\quad (6.11)
\end{aligned}
$$

The simpler special case of $\mu$ being a labeling, i.e. $\mu = \delta(y)$, $y \in \mathcal{Y}_{\mathcal{V}}$, is illustrated in Figure 6.1.

2. This follows from Proposition 5.36, since $\mathcal{D}(\phi)$ is a Lagrangian relaxation of the energy minimization problem (6.1). This fact can be also directly proven by observing that for any $\mu \in \mathcal{L} \cup \{0,1\}^{\mathcal{I}}$ and any $\phi \in \mathbb{R}^{\mathcal{J}}$

$$
\begin{aligned}
\langle \theta, \mu \rangle &\geq \min_{y \in \mathcal{Y}_{\mathcal{V}}} \sum_{u \in \mathcal{V}} \left( \theta_u(y_u) + \sum_{uv \in \mathcal{E}} \theta_{uv}(y_u, y_v) \right) \\
&\stackrel{\text{see Item 1}}{=} \min_{y \in \mathcal{Y}_{\mathcal{V}}} \left( \sum_{u \in \mathcal{V}} \theta_u^\phi(y_u) + \sum_{uv \in \mathcal{E}} \theta_{uv}^\phi(y_u, y_v) \right) \\
&\geq \sum_{u \in \mathcal{V}} \min_{s \in \mathcal{Y}_u} \theta_u^\phi(s) + \sum_{uv \in \mathcal{E}} \min_{(s,t) \in \mathcal{Y}_{uv}} \theta_{uv}^\phi(s,t) = \mathcal{D}(\phi) \quad (6.12)
\end{aligned}
$$

3. The dual is concave piecewise linear as any Lagrange dual of an integer linear program is. This also can be seen directly from (6.9): The reparametrized potentials $\theta^\phi$ are linear functions of $\phi$ according to (6.6), and the dual objective has also the representation (6.7) as a minimum over linear functions of $\phi$.

4. The proof is based on Proposition 5.42 applied to problem (6.1). Recall that its dual reads

$$
\max_{\phi \in \mathbb{R}^{\mathcal{J}}} \mathcal{D}(\phi) = \max_{\phi \in \mathbb{R}^{\mathcal{J}}} \min_{\substack{\mu \in \{0,1\}^{\mathcal{I}} \\ \mu_u \in \Delta^{\mathcal{Y}_u}, u \in \mathcal{V} \\ \mu_{uv} \in \Delta^{\mathcal{Y}_{uv}}, uv \in \mathcal{E}}} \left\langle \theta^\phi, \mu \right\rangle . \quad (6.13)
$$

To use it, let us consider the feasible set $\mathrm{conv}(P \cap \{0,1\}^n)$ of the primal relaxed problem in the right-hand-side of (5.50). Applied to the MAP-inference problem (6.1) it corresponds to

$$
\mathrm{conv}(\{\mu \in \{0,1\}^{\mathcal{I}} \colon \mu_u \in \Delta^{\mathcal{Y}_u}, u \in \mathcal{V}, \ \mu_{uv} \in \Delta^{\mathcal{Y}_{uv}}, uv \in \mathcal{E}\}) . \quad (6.14)
$$

This set decomposes into a set of independent constraints for each $\mu_w$, $w \in \mathcal{V} \cup \mathcal{E}$: $\mu_w \in \mathrm{conv}(\Delta^{\mathcal{Y}_w} \cap \{0,1\}^{\mathcal{Y}_w})$. According to Example (3.40) the set $\Delta^{\mathcal{Y}_w} \cap \{0,1\}^{\mathcal{Y}_w}$ is precisely the set of vertices of the simplex $\Delta^{\mathcal{Y}_w}$ and, therefore, its convex hull is the set $\Delta^{\mathcal{Y}_w}$ itself (Proposition 3.37).

Therefore, it holds that

$$\text{conv}(\{\mu \in \{0,1\}^{\mathcal{I}} \mid \mu_u \in \Delta^{\mathcal{Y}_u}, u \in \mathcal{V}, \ \mu_{uv} \in \Delta^{\mathcal{Y}_{uv}}, uv \in \mathcal{E}\})$$
$$= \{\mu \in \mathbb{R}^{\mathcal{I}} \mid \mu_u \in \Delta^{\mathcal{Y}_u}, u \in \mathcal{V}, \ \mu_{uv} \in \Delta^{\mathcal{Y}_{uv}}, uv \in \mathcal{E}\}. \quad (6.15)$$

Combining (6.15) with the coupling constraints one obtains the local polytope $\mathcal{L}$ (6.2).

5. <u>Necessity</u> can be shown using Proposition 5.45. The optimality condition (5.58) translates into

$$\mu \in \underset{\substack{\mu'_u \in \Delta^{\mathcal{Y}_u}, u \in \mathcal{V} \\ \mu'_{uv} \in \Delta^{\mathcal{Y}_{uv}}, uv \in \mathcal{E}}}{\arg\min} \left\langle \theta^{\phi}, \mu' \right\rangle \quad (6.16)$$

due to (6.15). Let $\phi \in \mathbb{R}^{\mathcal{J}}$ be a dual optimum. Then according to Proposition 5.45 there exists $\mu$ satisfying both (6.16) and the coupling constraints. It implies that $\mu \in \mathcal{L}$. To show that $\mu \in \mathcal{L}(\phi)$ it remains to prove that $\mu_w(s) = 0$ if $s \notin \arg\min_{s \in \mathcal{Y}_w} \theta^{\phi}_w(s)$. The latter holds due to

$$\underset{\substack{\mu'_u \in \Delta^{\mathcal{Y}_u}, u \in \mathcal{V} \\ \mu'_{uv} \in \Delta^{\mathcal{Y}_{uv}}, uv \in \mathcal{E}}}{\min} \left\langle \theta^{\phi}, \mu' \right\rangle = \sum_{w \in \mathcal{V} \cup \mathcal{E}} \underset{\mu'_w \in \Delta^{\mathcal{Y}_w}}{\min} \sum_{s \in \mathcal{Y}_w} \mu'_w(s) \theta^{\phi}_w(s) \quad (6.17)$$

and Lemma 3.38.

To prove <u>sufficiency</u> let us denote $\gamma_w := \min_{s \in \mathcal{Y}_w} \theta^{\phi}_w(s)$. If $\mathcal{L}(\phi) \neq \emptyset$ then there exists $\mu \in \mathcal{L}(\phi) \subseteq \mathcal{L}$. Due to Item 1 it holds that

$$\langle \theta, \mu \rangle = \left\langle \theta^{\phi}, \mu \right\rangle = \sum_{w \in \mathcal{V} \cup \mathcal{E}} \sum_{s \in \mathcal{Y}_w} \gamma_w \mu_w(s)$$
$$= \sum_{w \in \mathcal{V} \cup \mathcal{E}} \gamma_w \underbrace{\sum_{s \in \mathcal{Y}_w} \mu_w(s)}_{=1} = \sum_{w \in \mathcal{V} \cup \mathcal{E}} \min_{s \in \mathcal{Y}_w} \theta^{\phi}_w(s) = \mathcal{D}(\phi). \quad (6.18)$$

Since $\langle \theta, \mu \rangle \geq D(\phi)$ for any $\mu \in \mathcal{L}$, this proves that $\phi$ maximizes $\mathcal{D}$ and $\mu$ minimizes $\langle \theta, \mu \rangle$ on $\mathcal{L}$.

6. From Item 5 it follows that $\delta(y)$ is a minimizer of the local polytope relaxation. Since $\delta(y) \in \{0,1\}^{\mathcal{I}}$ it also minimizes the non-relaxed problem (6.1) and, therefore, $y$ is an optimal labeling. The equality $\mathcal{D}(\phi) = E(y; \theta) \equiv \langle \theta, \delta(y) \rangle$ follows from (6.18), which holds for any $\mu \in \mathcal{L}(\phi)$ and, therefore, for $\delta(y)$.

To prove the statement in the opposite direction assume that $\mathcal{L}(\phi) \neq \emptyset$, but that it does not contain any integer labeling, i.e. $\delta(y) \notin \mathcal{L}(\phi)$ for any $y \in \mathcal{Y}$. Due to Item 4 a vector $\mu \in \mathcal{L}$ is an optimum of the local polytope relaxation if and only if $\left\langle \theta^{\phi}, \mu \right\rangle = \mathcal{D}(\phi)$. In its turn, according to Item 5, the equality $\left\langle \theta^{\phi}, \mu \right\rangle = \mathcal{D}(\phi)$ holds only for $\mu \in \mathcal{L}(\phi)$. Since $\delta(y) \notin \mathcal{L}(\phi)$, the vector $\delta(y)$ is not a solution of the relaxed problem. Therefore,

$$\mathcal{D}(\phi) = \left\langle \theta^{\phi}, \mu \right\rangle < \left\langle \theta^{\phi}, \delta(y) \right\rangle = E(y; \theta), \qquad (6.19)$$

which is a contradiction.

$\square$

**Exercise 6.3.** Show that the following linear program is equivalent to the dual problem (6.9), i.e. for any $\phi \in \mathbb{R}^{\mathcal{J}}$ the maximization over the remaining variables gives the value $\mathcal{D}(\phi)$:

$$\max_{\gamma, \phi} \left( \sum_{u \in \mathcal{V}} \gamma_u + \sum_{uv \in \mathcal{E}} \gamma_{uv} \right) \qquad (6.20)$$
$$\gamma_u \leq \theta_u^{\phi}(s), \ u \in \mathcal{V}, \ s \in \mathcal{Y}_u,$$
$$\gamma_{uv} \leq \theta_{uv}^{\phi}(s, t), \ uv \in \mathcal{E}, \ (s, t) \in \mathcal{Y}_{uv}.$$

This linear program can be obtained by dualizing all constraints defining the local polytope. The LP form (6.20) of the Lagrangean dual supports a known fact (see e.g. the text-book [15]), namely, that the dual of a linear program (which is the considered local polytope relaxation) can itself be represented as a linear program.

## 6.2 Primal solutions from the dual problem

### 6.2.1 Naïve dual rounding

Statements (5)-(6) in Proposition (6.2) justify the following heuristics for obtaining approximate solutions of the non-relaxed problem (6.1) from a reparametrization $\phi$: In each $u \in \mathcal{V}$ a locally optimal label is

assigned to $y_u$:

$$y_u := \arg \min_{s \in \mathcal{Y}_u} \theta_u^{\phi}(s) \,. \tag{6.21}$$

Compared to the primal rounding (4.9) the naïve dual rounding (6.21) is weaker in the following sense: from $\mu_u(s) > 0$ it follows that

$$s \in \arg \min_{s' \in \mathcal{Y}_u} \theta_u^{\phi}(s') \,, \tag{6.22}$$

but not the other way around. In other words, a locally optimal label may correspond to a zero coordinate of an optimal relaxed solution. It implies that even in case the Lagrange dual is tight (as in Proposition $(6.2)(6)$), the rounding (6.21) may not allow reconstruction of the optimal solution. The following example illustrates this case:

**Example 6.4** (Naïve dual rounding is weak). Consider the graph $(\mathcal{V}, \mathcal{E})$ with $\mathcal{V} = \{u, v\}$ and $\mathcal{E} = \{uv\}$, see Figure 6.2(a). The label set is the same in both nodes and contains two labels: $\mathcal{Y}_u = \mathcal{Y}_v = \{a, b\}$. All unary costs are equal to 0 except $\theta_u(b) = 1$. Pairwise costs are defined as $\theta_{uv}(s, t) = [\![a \neq b]\!]$.

Since $\theta \geq 0$ and $\theta_u(a) = \theta_v(a) = \theta_{uv}(a, a) = 0$, the labeling $(a, a)$ consists of locally optimal labels and label pairs. Therefore, according to Proposition $(6.2)(6)$ the trivial reparametrization $\phi = 0$ is an optimal one for the dual problem. The optimal primal and dual values are equal to 0. Moreover, the labeling $(a, a)$ is the unique solution of the relaxed as well as the non-relaxed problem.

Consider now the naïve dual rounding (6.21). Since $\theta_v(b) = 0$ this procedure may return the labeling $(a, b)$ with energy 1, which is non-optimal.

Example 6.4 shows that even if the set $\mathcal{L}(\phi)$ contains a single integer labeling, the locally optimal labels and label pairs may not be unique. Therefore, the naive rounding may return sub-optimal results.

Importantly, it can be shown that in the considered LP-tight case $\mathcal{L}(\phi) = \{\delta(y)\}$ there always exists an optimal reparametrization, such that each node contains a *unique* locally optimal label. For such a reparametrization the naïve rounding returns the optimal primal solution. Figure 6.2(b) shows such a reparametrization for the problem from Example 6.4.

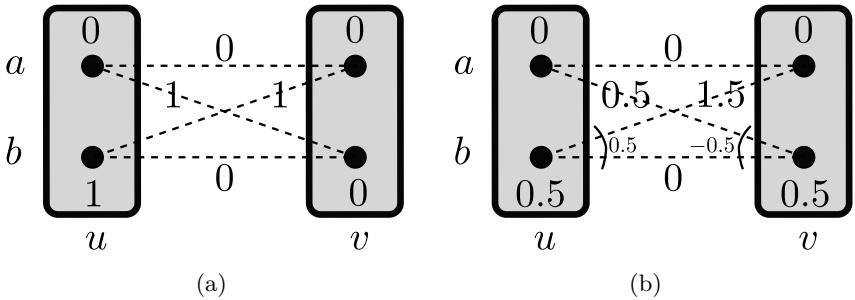

**Figure 6.2:** Illustration of Example 6.4. Characters $u$, $v$ stand for graph nodes and $a, b$ for labels. Big numbers 0, 1, 0.5, 1.5 define costs assigned to labels and label pairs. Small numbers 0.5 in (b) provide values of the dual variables $\phi$, which are required to turn the costs of (a) to the costs of (b). **(a)** In spite of the optimal reparametrization, the naïve dual rounding (6.21) may return a non-optimal labeling $(a, b)$, even in the LP-tight case with a unique integer solution. **(b)** Another optimal reparametrization, where the locally optimal labels are unique and therefore the naïve dual rounding (6.21) returns the correct primal optimal integer solution.

Therefore, there exist three ways to construct reasonably good rounding algorithms for the Lagrange dual:

- Use dual optimization algorithms which return possibly unique locally optimal labels/label pairs in each node/edge. In Chapter 8 we will consider one such algorithm, known as *min-sum diffusion*.

- Create rounding methods which are specialized to a certain type of dual optimization algorithm. In Chapter 10 we will consider such a specialized rounding method related to the TRW-S algorithm.

- Find a rounding method which works reasonably for any optimal reparametrization.

Below we concentrate on the last option. First of all, we will study local properties of optimal dual solutions, which will play a crucial role in constructing dual rounding schemes. Later, in §6.2.4, we consider the relaxation labeling algorithm, which can be used to obtain approximate primal solutions.

### 6.2.2 Arc-consistency

As follows from Proposition 6.2(5), to determine whether a given dual vector $\phi$ is optimal, it suffices to check whether the polytope $\mathcal{L}(\phi)$ defined in (6.10) is non-empty. Such kinds of problems, where it is necessary to verify whether a set is empty, are called *feasibility problems*. In general, feasibility problems for linear programming are as difficult as linear programs themselves. In other words, checking whether a polyhedron is empty can be formulated as a linear program with a different polyhedron.

In practice, a simpler optimality condition is often desirable, which could be checked often enough, for example on each iteration of some iterative algorithm. In this section we formulate such a condition. The price for its simplicity is that it is only necessary for dual optimality and no longer sufficient. Informally, instead of looking for an integer or relaxed labeling consisting of locally optimal labels and label pairs, which would require a *global* reasoning, one may check whether the locally optimal labels and label pairs are *locally* consistent. Below, we give a formal definition of this simpler condition and an algorithm for its verification.

In the following, we will use two functions, which turn real-valued vectors from the primal space $\mathbb{R}^{\mathcal{I}}$ into binary ones:

- For any $\mu \in \mathbb{R}$ let $\mathrm{nz}[\mu] = [\![\mu \neq 0]\!]$ be the indicator function of $\mu$ being **n**on-**z**ero. When applied to a vector $\mu \in \mathbb{R}^n$, it acts coordinate-wise, i.e. $\mathrm{nz}[\mu]_i = \mathrm{nz}[\mu_i]$.

- For $\theta \in \mathbb{R}^{\mathcal{I}}$ let $\mathrm{mi}[\theta]$ be defined such that locally minimal labels and label pairs (w.r.t. $\theta$) obtain the value 1 and others zero, i.e. $\mathrm{mi}[\theta]_w(x_w) := [\![\theta_w(x_w) = \min_{x_w \in \mathcal{Y}_w} \theta_w(x_w)]\!]$ for $w \in \mathcal{V} \cup \mathcal{E}$. Here mi stands for **mi**n.

**Example 6.5.**

$$\mathrm{nz}[\delta(y)] = \delta(y);$$
$$\mathrm{nz}[(0, 0.2, 0.8, 0)] = (0, 1, 1, 0);$$
$$\mathrm{mi}[(0, -2, -1, -2, 3)] = (0, 1, 0, 1, 0);$$
$$\mathrm{mi}[(0, -2), (1, 1, 2, 0), (7, -1)] = (0, 1), (0, 0, 0, 1), (0, 1).$$

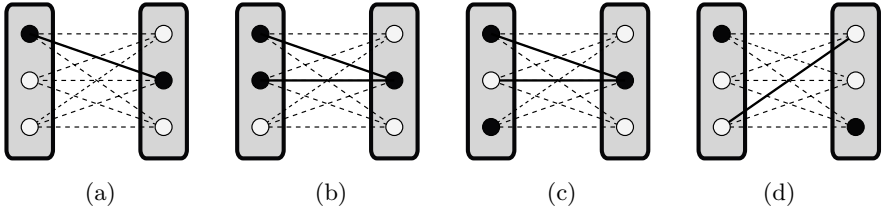

**Figure 6.3:** Illustration to the notions of arc-consistency and node-edge agreement (non-empty closure). Black circles and solid lines correspond to the non-zero coordinates of $\xi \in \{0,1\}^{\mathcal{I}}$. White circles and dashed lines to the zero-valued coordinates. **(a)** Strictly arc-consistent vector. **(b)** Arc-consistent $\xi$. **(c)** Not arc-consistent, however the node-edge agreement holds ($\mathrm{cl}(\xi) \neq \bar{0}$). **(d)** Not arc-consistent, no node-edge agreement ($\mathrm{cl}(\xi) = \bar{0}$).

**Definition 6.6.** A binary vector $\xi \in \{0,1\}^{\mathcal{I}}$ is called *arc-consistent* if

1.  $\xi_{uv}(s,t) = 1$ implies $\xi_u(s) = \xi_v(t) = 1$ for all $uv \in \mathcal{E}$, $(s,t) \in \mathcal{Y}_{uv}$, and

2.  $\xi_u(s) = 1$ implies that for any $v \in \mathcal{N}_b(u)$ there exists $t \in \mathcal{Y}_v$ such that $\xi_{uv}(s,t) = 1$.

The set of arc-consistent vectors from $\{0,1\}^{\mathcal{I}}$ will be denoted as $\mathcal{AC}^{\mathcal{I}}$.

This definition as well as the following one are illustrated in Figure 6.3.

**Definition 6.7.** $\xi \in \{0,1\}^{\mathcal{I}}$ is *strictly arc-consistent* if it is arc-consistent and $\sum_{s \in \mathcal{Y}_u} \xi_u(s) = 1$ for all $u \in \mathcal{V}$ as well as $\sum_{(s,t) \in \mathcal{Y}_{uv}} \xi_{uv}(s,t) = 1$ for all $uv \in \mathcal{E}$.

In other words, strict arc-consistency means that (i) for a single label in each node and a single label pair in each edge the corresponding coordinate of $\xi$ is equal to 1 and (ii) these labels and label pairs are consistent with each other.

The statements of the following proposition follow directly from the definition of the local polytope $\mathcal{L}$. They show that arc consistency is an intrinsic property of all vectors in $\mathcal{L}$, and that strict arc-consistency is equivalent to the notion of integer labeling.

**Proposition 6.8.**

1. $\mu \in \mathcal{L}$ implies that $\mathrm{nz}[\mu]$ is arc-consistent.

2. For any $y \in \mathcal{Y}_\mathcal{V}$, $\delta(y)$ is strictly arc-consistent.

3. For any $\theta \in \mathbb{R}^\mathcal{I}$ the following two statements are equivalent:

   (a) $\mathrm{mi}[\theta]$ is strictly arc-consistent;
   (b) $\mathrm{mi}[\theta] \in \mathcal{L} \cap \{0, 1\}^\mathcal{I}$, i.e. $\mathrm{mi}[\theta]$ is an integer labeling.

**Node-edge agreement is insufficient for dual optimality**

Proposition 6.8(1) defines a *necessary* condition for the dual optimum given in Proposition 6.2(5). Coordinates of $\mu \in \mathcal{L}(\phi)$ may be greater than 0 only for locally optimal labels and label pairs. Therefore, for a dual vector $\phi$ to be optimal it is *necessary* that there exists an arc-consistent subset of locally optimal coordinates (corresponding to labels and label pairs) of $\theta^\phi$:

**Definition 6.9.** One says that nodes and edges *agree* (or there is a *node-edge agreement*) for a cost vector $\theta \in \mathbb{R}^\mathcal{I}$, if for each node (edge) $w \in \mathcal{V} \cup \mathcal{E}$ there is a non-empty subset of locally optimal labels (label pairs) $\mathbb{S}_w \subseteq \arg\min_{s \in \mathcal{Y}_w} \theta_w(s)$ such that the vector $\xi$ with coordinates $\xi_w(s) = [\![s \in \mathbb{S}_w]\!]$ is arc-consistent.

In §6.2.4 we will provide an algorithm able to find the subsets $\mathbb{S}_w$. However, first we show that the arc-consistency of the binary vector $\xi$ defined by these subsets is only necessary, but not sufficient for dual optimality.

Example 6.10 provides a graphical model with costs $\theta$ such that $\mathrm{mi}[\theta]$ is arc-consistent. However, the reparametrization $\phi = 0$ is not optimal.

**Example 6.10** (Arc Consistency $\neq$ Dual Optimality [140])**.** Consider the graphical model in Figure 6.4. Assume all labels are locally optimal and locally optimal label pairs are denoted as solid lines with and without arrows. Let us show that there is no vector $\mu$ in the corresponding local polytope such that from $\mu_{uv}(s,t) > 0$ it follows that the label pair $(s,t)$ is locally optimal.

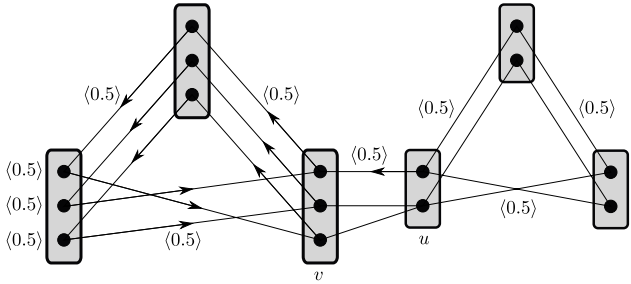

**Figure 6.4:** Illustration corresponding to Example 6.10. Arrows are drawn to rule the reasoning leading to a contradiction. All depicted labels are locally optimal and locally optimal label pairs are denoted as solid lines with and without arrows. Angular brackets $\langle\rangle$ are used for coordinates of $\mu$. It turns out that there is no vector $\mu$ in the corresponding local polytope such that $\mu_{uv}(s,t) > 0$ implies that the label pair $(s,t)$ is locally optimal.

The model consists of two parts each containing three nodes connected by edges. The right "triangle" is the same as in Example 4.3, where we have learned that the only element of the local polytope compatible with the depicted locally optimal edges is the one with the corresponding coordinates equal to 0.5. In other words, coordinates of $\mu$ corresponding to the locally optimal labels and label pairs in the right "triangle" are equal to 0.5.

Let us now show that this implies that all optimal labels and label pairs in the left "triangle" correspond to the value 0.5 of the respective coordinates of $\mu$ as well. This will mean a contradiction to the simplex constraint $\sum_{s \in \mathcal{Y}_u} \mu_u(s) = 1$, since each of the labels in the left "triangle" will be assigned the value 0.5 and $3 \cdot 0.5 = 1.5 \neq 1$, showing that arc consistency is insufficient for optimality.

To show this, it is enough to start with the locally optimal label pair denoted by the arrowed line between the nodes $u$ and $v$. The coordinate of $\mu$ corresponding to this label pair is equal to 0.5 due to the coupling constraints. Going along the arrows and applying the coupling constraint again and again one obtains the required property.

### 6.2.3 Arc-consistency closure

As noted in §6.2.2, the necessary dual optimality condition requires to find an arc-consistent subset of locally optimal labels and label pairs. Here we define this subset formally and provide an efficient algorithm, which either finds such a subset or determines that there is none.

Let us introduce two logical operations: "and", denoted as $\wedge$, and "or", denoted as $\vee$. We assume their coordinate-wise application to binary vectors. In this section we will give an algorithm which, for an input binary vector $\xi \in \{0,1\}^{\mathcal{I}}$, finds out whether an arc-consistent binary vector $\xi'$ exists such that $\xi' \wedge \xi = \xi'$. The vector $\mathrm{mi}[\theta^\phi]$ is used in place of $\xi$, when optimality of $\phi$ must be determined w.r.t. the dual objective $\mathcal{D}$.

**Definition 6.11.** A vector $\xi' \in \mathcal{AC}^{\mathcal{I}}$ is called the *arc-consistency closure* of a vector $\xi \in \{0,1\}^{\mathcal{I}}$ and is denoted as $\mathrm{cl}(\xi)$ if

1. $\xi' \wedge \xi = \xi'$, and

2. for any $\xi'' \in \mathcal{AC}^{\mathcal{I}}$ such that $\xi'' \wedge \xi = \xi''$ it holds that $\xi'' \wedge \xi' = \xi''$.

The closure $\xi'$ is called non-empty, if $\xi' \neq \overline{0}$, where $\overline{0}$ is the vector with all zero coordinates.

The configuration in Figure 6.3(d) has an empty closure, an in Figures 6.3(c), 6.3(a), 6.3(b) and 6.4 a non-empty one. The configurations in latter three figures have non-empty closures, since they are arc-consistent themselves.

The arc-consistency closure is also often called *kernel* in the literature. We will often refer to the arc-consistency closure as *closure*, if it does not lead to misunderstanding.

To analyze properties of the arc-consistency closure we will require the following simple lemma, which directly follows from Definition 6.6.

**Lemma 6.12.** Let $\xi, \xi'$ be two arc-consistent vectors. Then $\xi \vee \xi'$ is arc-consistent as well.

**Proposition 6.13** (Closure properties). For any $\xi \in \{0,1\}^{\mathcal{I}}$ its arc-consistency closure is unique and is equal to $\bigvee_{\substack{\xi' \in \mathcal{AC}^{\mathcal{I}} \\ \xi' \wedge \xi = \xi'}} \xi'$.

In other words, the closure is the arc-consistent vector with the largest number of coordinates equal to 1.

*Proof.* Let us consider the vector $\hat{\xi} := \bigvee_{\substack{\xi' \in \mathcal{AC}^{\mathcal{I}} \\ \xi' \wedge \xi = \xi'}} \xi'$. By construction and due to Lemma 6.12 it satisfies Definition 6.11, therefore, $\hat{\xi} = \mathrm{cl}(\xi)$.

<u>Uniqueness:</u> Let $\xi''$ be another arc-consistency closure such that $\xi'' \neq \hat{\xi}$. Then $\xi'' \in \mathcal{AC}^{\mathcal{I}}$ and $\hat{\xi}$ is a closure, therefore $\xi'' \wedge \hat{\xi} = \xi''$ by Definition 6.11. However, since $\hat{\xi} \in \mathcal{AC}^{\mathcal{I}}$, and $\xi''$ is a closure, it implies $\xi'' \wedge \hat{\xi} = \hat{\xi}$. Therefore, $\xi'' = \hat{\xi}$, which is a contradiction. $\qquad\square$

**Example 6.14.** It is easy to show that arc-consistency of $\xi$ and $\xi'$ does not imply arc-consistency of $\xi`` := \xi \wedge \xi'$. Indeed, consider a binary graphical model defined on the graph $\mathcal{G} = (\{u, v\}, \{uv\})$ with two nodes and an edge between them. Let the non-zero coordinates of $\xi$ and $\xi'$ are $\xi_u(1) = \xi_v(1) = \xi_v(0) = \xi_{uv}(1, 1) = \xi_{uv}(1, 0) = 1$, and $\xi'_u(0) = \xi'_v(0) = \xi'_{uv}(0, 0) = 1$. Then the only non-zero coordinate of $\xi''$ is $\xi''_v(0) = 1$, which implies its inconsistency.

For the sake of notation, in the following corollary we will use the coordinate-wise comparison of binary vectors. For example, $\xi \leq \xi'$ is equivalent to $\xi \wedge \xi' = \xi$.

**Corollary 6.15.** Let $\xi, \xi' \in \{0, 1\}^{\mathcal{I}}$ and $\xi \leq \xi'$. Then $\mathrm{cl}(\xi) \leq \mathrm{cl}(\xi')$.

*Proof.* According to Proposition 6.13 it holds that

$$\mathrm{cl}(\xi') = \bigvee_{\substack{\xi'' \in \mathcal{AC}^{\mathcal{I}} \\ \xi'' \leq \xi'}} \xi''. \qquad (6.23)$$

Consider $\xi'' = \mathrm{cl}(\xi)$. Since $\mathrm{cl}(\xi) \in \mathcal{AC}^{\mathcal{I}}$ and $\mathrm{cl}(\xi) \leq \xi \leq \xi'$, therefore, according to (6.23) $\mathrm{cl}(\xi) \leq \mathrm{cl}(\xi')$. $\qquad\square$

### 6.2.4 Relaxation labeling algorithm

The *relaxation labeling* algorithm turns any $\xi \in \{0, 1\}^{\mathcal{I}}$ into $\mathrm{cl}(\xi)$ iteratively. It repeats the following two operations until $\xi$ does not change

anymore:

$$\forall u \in \mathcal{V}, \ s \in \mathcal{Y}_u, \ v \in \mathcal{N}_b(u)\colon \quad \xi_u(s) := \xi_u(s) \wedge \bigvee_{t \in \mathcal{Y}_v} \xi_{uv}(s,t), \quad (6.24)$$

$$\forall uv \in \mathcal{E}, \ (s,t) \in \mathcal{Y}_{uv}\colon \quad \xi_{uv}(s,t) := \xi_{uv}(s,t) \wedge \xi_u(s) \wedge \xi_v(t). \quad (6.25)$$

The algorithm is illustrated in Figures 6.5 and 6.6.

Operation (6.24) searches for non-zero label pairs $\xi_{uv}(s,t)$ for each label $s$ in each node $u$ and each neighbor $v \in \mathcal{N}_b(u)$. If there is no such non-zero label pair, then $\bigvee_{t \in \mathcal{Y}_v} \xi_{uv}(s,t) = 0$ and, therefore, $\xi_u(s)$ is turned to zero itself. Operation (6.24) need not be performed for those labels $s$, where $\xi_u(s) = 0$, since independently of the value of $\bigvee_{t \in \mathcal{Y}_v} \xi_{uv}(s,t)$, the value of $\xi_u(s)$ will remain zero.

Operation (6.25) checks whether $\xi_u(s) = \xi_v(t)$ for each label pair $(s,t)$ in each edge $uv$. If it is not the case, the label pair $\xi_{uv}(s,t)$ is assigned zero value. As before, this check has to be performed only for those label pairs, where $\xi_{uv}(s,t) = 1$ initially. As soon as $\xi_{uv}(s,t) = 0$ it keeps this value in the future.

This simple analysis implies that the relaxation labeling algorithm (6.24)-(6.25) is a finite step algorithm, because during each iteration non-zero coordinates can potentially be set to zero, but never the other way around. Therefore, the algorithm stops when either $\xi = \bar{0}$ or it has not changed any coordinate of $\xi$ on a current iteration.

Note that the algorithm stops as soon as $\xi$ becomes arc-consistent. In other words, for an input vector $\xi \in \{0,1\}^{\mathcal{I}}$ the relaxation labeling algorithm (6.24)-(6.25) finds the arc-consistent vector $\xi'$ and with the largest number of coordinates equal to 1 such that $\xi' \wedge \xi = \xi'$. Due to Proposition 6.13, it is $\xi' = \mathrm{cl}(\xi)$.

**Output of relaxation labeling algorithm** Let $\theta^\phi$ be a reparametrized cost vector and $\xi = \mathrm{mi}[\theta^\phi]$ be the input to the relaxation labeling algorithm. There are three possible conditions the output vector $\xi'$ may satisfy:

- If $\xi'$ is strictly arc-consistent, then $\phi$ is the optimal reparametrization and the labeling $y$ with coordinates defined as $y_u = s$ if $\xi'_u(s) = 1$, $u \in \mathcal{V}$, is a minimizer of the MAP-inference problem $\min_{y' \in \mathcal{Y}_\mathcal{V}} \langle \theta, \delta(y') \rangle$.

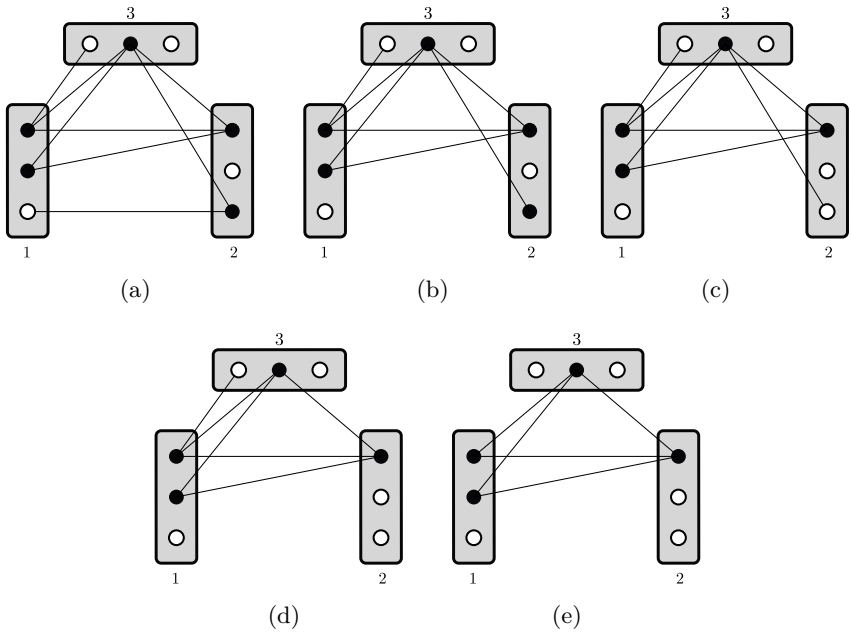

**Figure 6.5:** Illustration of the relaxation labeling algorithm. Black labels and solid lines stand for the value 1, others – for 0. Nodes and edges of the graph are processes in the counter-clock-wise order, i.e. node 1, edge $\{1,2\}$, node 2, edge $\{2,3\}$, node 3 etc. The algorithms stops when no changes have been done on the last iteration over the whole graph. The attained configuration is arc-consistent.

- If $\xi' = \bar{0}$, the closure is empty and therefore, there is no node-edge agreement, which implies that the reparametrization $\phi$ is non-optimal.

- $\xi'$ is neither zero nor strictly arc-consistent. In this case to determine optimality of $\phi$ one has to verify whether there is $\mu \in \mathcal{L}$ such that $\mu \leq \xi'$, which is the feasibility problem of linear programming. This is similar to the optimality condition $\mathcal{L}(\phi) \neq \emptyset$ discussed in § 6.2.2 with the difference that the set of potentially positive coordinates of $\mu$ does not include the locally optimal labels and label pairs set to zero by the relaxation labeling algorithm.

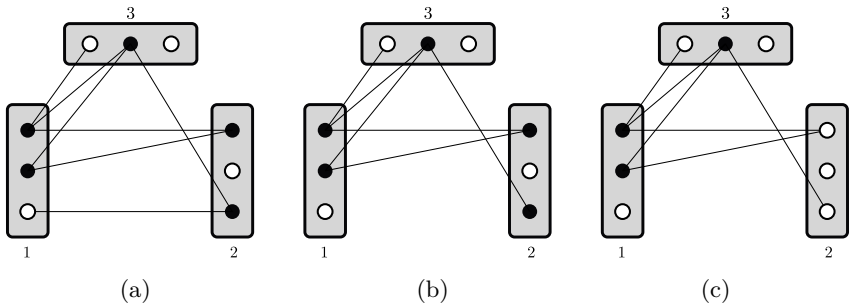

**Figure 6.6:** Illustration of the relaxation labeling algorithm. Note the difference to Figure 6.5, the label pair $\{3, 1\}$ in the edge $\{3, 2\}$ is assigned 0 value. Black labels and solid lines stand for the value 1, others – for 0. Nodes and edges of the graph are processes in the counter-clock-wise order, i.e. node 1, edge $\{1, 2\}$, node 2, edge $\{2, 3\}$, node 3 . . . . The algorithms stops when no changes have been done on the last iteration over the whole graph. The attained closure is empty, since already on when processing the node 2 all its labels where assigned the value 0. It is easy to see that this will inevitably assign values 0 to all other labels and label pairs within a single iteration of the algorithm.

**Relation to Constraint Satisfaction**  For a binary vector $\xi \in \{0, 1\}^{\mathcal{I}}$ the problem to decide whether there exists a strict arc-consistent vector $\xi'$ such that $\xi' \wedge \xi = \xi'$ is $\mathcal{NP}$-hard and constitutes a *constraint satisfaction problem*.

For example, the graphical model obtained by reduction of the Hamiltonian cycle problem in §1.3 can be treated as a constraint satisfaction problem, if zero costs are substituted by cost 1 and infinite costs are substituted by 0.

For $\xi = \text{mi}[\theta^\phi]$ solving the constraint satisfaction problem can be necessary if $\text{cl}(\xi)$ is non-empty and not strictly arc-consistent, and the solution $\mu$ of the feasibility problem "find $\mu \in \mathcal{L}$ such that $\mu \leq \text{cl}(\xi)$" is non-integral. Should the solution of the constraint satisfaction problem exist, then the corresponding labeling is an exact solution of the MAP-inference problem $\min_{y \in \mathcal{Y}_\mathcal{V}} \langle \theta, \delta(y) \rangle$.

**$\epsilon$-Arc-Consistency**  In practice, however, the minimum $\min_{s \in \mathcal{Y}_w} \theta^\phi_w(s)$ is usually attained in only one label due to a finite precision of floating point computations. Therefore, for a potential vector $\theta$, the binary vector $\mathrm{mi}[\theta]$ is either strictly arc-consistent or has an empty arc-consistency closure. Both conditions in this case can be checked directly and without the relaxation labeling algorithm.

To also take labels into account which correspond to almost minimal costs, one has to resort to an approximate arc-consistency, so called $\epsilon$-arc-consistency. We denote as $\mathrm{mi}_\epsilon[\theta]$ the vector with coordinates

$$\mathrm{mi}_\epsilon[\theta]_w(s) = \begin{cases} 1, & \theta_w(s) \le \min_{s \in \mathcal{Y}_w} \theta_w(s) + \epsilon \\ 0, & \text{otherwise} \end{cases} , \ w \in \mathcal{V} \cup \mathcal{E} . \quad (6.26)$$

**Definition 6.16.** We say that for the costs $\theta$ the $\epsilon$-*arc-consistency* property holds for $\epsilon \ge 0$ if $\mathrm{mi}_\epsilon[\theta]$ is arc-consistent. We also say that there is the $\epsilon$-*node-edge agreement*, if $\mathrm{cl}(\mathrm{mi}_\epsilon[\theta]) \ne \bar{0}$.

There always exists a large enough $\epsilon$ such that $\mathrm{mi}_\epsilon[\theta]$ has a non-empty arc-consistency closure. This obviously holds for

$$\epsilon = \max_{w \in \mathcal{V} \cup \mathcal{E}} \left( \max_{s \in \mathcal{Y}_w} \theta_w(s) - \min_{s \in \mathcal{Y}_w} \theta_w(s) \right) , \quad (6.27)$$

since in this case all coordinates of the vector $\mathrm{mi}_\epsilon[\theta]$ are equal to 1.

Moreover, one can compute a *minimal* $\epsilon$, which guarantees $\epsilon$-arc-consistency, by a min / max generalization of the relaxation labeling algorithm. This algorithm resembles the relaxation labeling algorithm (6.24)-(6.25) with the only difference that the operation $\vee$ is substituted with min and $\wedge$ is turned into max. The algorithm starts with $\xi := \theta$ and proceeds as follows:

$$\forall u \in \mathcal{V}, \ s \in \mathcal{Y}_u, \ v \in \mathcal{N}_b(u): \quad \xi_u(s) := \max \left\{ \xi_u(s), \min_{t \in \mathcal{Y}_v} \xi_{uv}(s,t) \right\} \quad (6.28)$$

$$\forall uv \in \mathcal{E}, \ (s,t) \in \mathcal{Y}_{uv}: \quad \xi_{uv}(s,t) := \max \{ \xi_{uv}(s,t), \xi_u(s), \xi_v(t) \} .$$

Algorithm stops when no changes to vector $\xi$ are made on the last iteration run over the whole graph. Let

$$m = \min_{\substack{w \in \mathcal{V} \cup \mathcal{E} \\ s \in \mathcal{Y}_w}} \xi_w(s) \quad (6.29)$$

the minimal coordinate of the vector $\xi$ after the algorithm stops. Independently of the order in which all nodes as edges of the graph are processed, the set of coordinates of $\xi$ with the value $m$ are arc-consistent. In other words, the vector $\xi'$ such that $\xi'_w(s) = [\![\xi_w(s) = m]\!]$ for $w \in \mathcal{V} \cup \mathcal{E}$ and $s \in \mathcal{Y}_w$, is arc-consistent. Figure 6.7 illustrates the work of the algorithm.

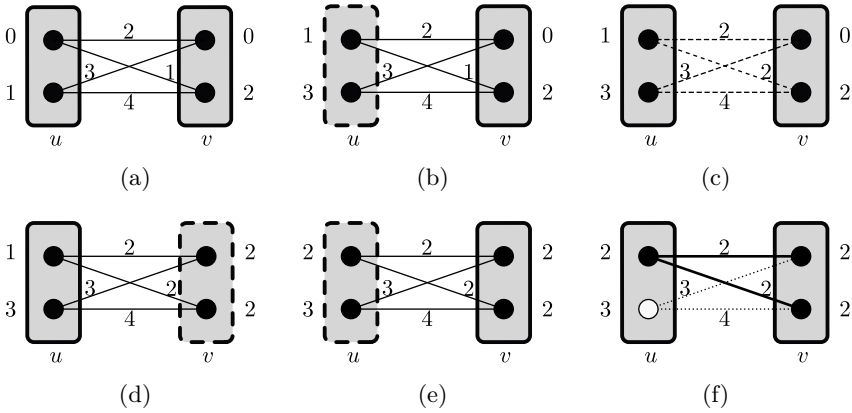

**Figure 6.7:** Illustration to the min-max variant of the relaxation labeling algorithm (6.28). **(a)** Initial graphical model, numbers denote costs, both unary and pairwise. **(a-e)** Algorithm progress, dashed node is the one processed at the current step of the algorithm. Dashed lines in the pairwise factor have the similar meaning of the currently processed factor. **(f)** Final configuration obtained by the algorithm. The labels marked with black circles and the label pairs corresponding to solid lines are those having the smallest cost 2. Note that they build an arc-consistent configuration. The value $\epsilon$ computed as in (6.30) is equal to 2.

The resulting value $\epsilon$ is obtained as

$$\epsilon := \max_{\substack{w \in \mathcal{V} \cup \mathcal{E} \\ s \in \mathcal{Y}_w}} \left( [\![\xi_w(s) = m]\!] \cdot (m - \theta_w(s)) \right) . \tag{6.30}$$

Note that the value $[\![\xi_w(s) = m]\!]$ differs from 0 and is equal to 1 only for those coordinates of $\xi$, which are equal to $m$.

Analysis of the algorithm is similar to the one for the relaxation labeling algorithm with $\lor$ and $\land$ operations. Therefore, we omit it here and provide the corresponding reference in §6.4.

## 6.3   Dual optimality for acyclic graphical models

As shown in §6.2.2 node-edge agreement is only necessary, but not sufficient dual optimality condition in general. However, there are a few cases, where this condition is sufficient as well. This holds in particular for binary graphical models, i.e. when only two labels are associated with each node (see Chapter 12 for details), or when the graphical models are acyclic. The latter case is considered below.

**Proposition 6.17.** Let $(\mathcal{G}, \mathcal{Y}_{\mathcal{V}}, \theta)$ be a graphical model with the graph $\mathcal{G}$ being acyclic. Then $\mathrm{cl}(\mathrm{mi}[\theta^\phi]) \neq \bar{0}$ implies $\phi$ is the dual optimum. Moreover, there is an optimal integer labeling $y \in \mathcal{Y}_{\mathcal{V}}$ such that $\delta(y) \leq \mathrm{cl}(\mathrm{mi}[\theta^\phi])$.

Furthermore, for any $(w, s) \in (\mathcal{V} \cup \mathcal{E}) \times \mathcal{Y}_w$ such that $\mathrm{cl}(\mathrm{mi}[\theta^\phi])_w(s) \neq 0$ there is an optimal integer labeling $y^* \in \mathcal{Y}_{\mathcal{V}}$ such that $y_w^* = s$.

*Proof.* Without loss of generality we assume that $\mathcal{G}$ is connected, i.e. is a tree. Otherwise the proof can be done for each connected component separately.

Let $\xi := \mathrm{cl}(\mathrm{mi}[\theta^\phi]) \neq \bar{0}$. The proof will be done by constructing a labeling $y$ such that $\delta(y) \leq \xi$.

We will iteratively build a *connected* subgraph $(\mathcal{V}', \mathcal{E}')$ and assign labels to the nodes of $\mathcal{V}'$ such that if label pair $(y_u, y_v)$ is assigned to the edge $(u, v) \in \mathcal{E}'$ it holds that $\xi_u(y_u) = 1$, $\xi_v(y_v) = 1$ and $\xi_{uv}(y_u, y_v) = 1$. Our procedure finalizes when $\mathcal{V}' = \mathcal{V}$ and $\mathcal{E}' = \mathcal{E}$.

The sets $\mathcal{V}'$ and $\mathcal{E}'$ are empty at the beginning of the procedure. Let $u \in \mathcal{V}$ be any vertex. Consider a label $s \in \mathcal{Y}_u$ for which $\xi_u(s) = 1$. Assign $y_u := s$ and $\mathcal{V}' := \mathcal{V}' \cup \{u\}$.

On each iteration a node $v \in \mathcal{V} \backslash \mathcal{V}'$ is considered, which is incident to some node of $\mathcal{V}'$. Note that any node $v \in \mathcal{V} \backslash \mathcal{V}'$ is incident to *at most* one node in $\mathcal{V}'$. Indeed, if $u'$ and $u''$ would be two nodes of $\mathcal{V}'$ incident to $v$, there would be a path $p$ between them in the graph $(\mathcal{V}', \mathcal{E}')$, since the latter is connected. Therefore, there would a cycle $u', v, u'', p, u'$, which would mean a contradiction, as the initial graph is acyclic.

Let therefore $v \in \mathcal{V} \backslash \mathcal{V}'$ be any node connected to some node $u \in \mathcal{V}'$. By definition of closure there is a $t \in \mathcal{Y}_v$ such that $\xi_{uv}(y_u, t) = \xi_v(t) = 1$. Assign $y_v := t$ and $\mathcal{V}' := \mathcal{V}' \cup \{v\}$, $\mathcal{E}' := \mathcal{E}' \cup \{uv\}$.

Repeat the process until $\mathcal{V}' = \mathcal{V}$ and therefore $\mathcal{E} = \mathcal{E}'$.

By construction it holds that $\delta(y) \leq \xi$ for the labeling $y$, which finalizes the proof of the first two statements.

For the last statement, the labeling $y^*$ can be constructed with $(u, s)$ being selected at the first step of the procedure. □

**Corollary 6.18.** Let $(\mathcal{G}, \mathcal{Y}_\mathcal{V}, \theta)$ be a graphical model with the graph $\mathcal{G}$ being acyclic. Let $\mathcal{L}$ and $\mathcal{M}$ be the corresponding local and marginal polytopes. Then $\mathcal{L} = \mathcal{M}$.

*Proof.* Since $\mathcal{M} \subseteq \mathcal{L}$ if suffices to show that $\mathcal{L} \subseteq \mathcal{M}$. Consider a maximizer $\phi$ of the Lagrange dual $\mathcal{D}$ defined by (6.9) for an acyclic problem with the cost vector $\theta$. Since $\mathrm{cl}(\mathrm{mi}[\theta^\phi]) \neq \bar{0}$ is necessary for optimality of $\phi$ (see §6.2.2), Proposition 6.17 implies that $\mathcal{L}(\phi)$ contains an integer labeling. In its turns, it implies (see statement 5 of Proposition 6.2) that this labeling is a solution to the local polytope relaxation of the MAP-inference problem.

In other words, Proposition 6.17 implies that for any cost vector $\theta$ there is always an integer solution of the local polytope relaxation of the MAP-inference problem for acyclic graphs. It implies that only for vectors $\delta(y)$, $y \in \mathcal{Y}_\mathcal{V}$, corresponding to integer labelings, a cost vector $\theta$ may exist such that $\delta(y)$ is the unique solution of the relaxed problem $\min_{\mu \in \mathcal{L}} \langle \theta, \mu \rangle$.

Due to Definition 3.19 this implies $\mathrm{vrtx}(\mathcal{L}) \subseteq \{\delta(y) \mid y \in \mathcal{Y}_\mathcal{V}\}$. Therefore,

$$\mathcal{L} \subseteq \mathrm{conv}\{\delta(y) \mid y \in \mathcal{Y}_\mathcal{V}\} = \mathcal{M} \subseteq \mathcal{L}, \tag{6.31}$$

which finalizes the proof. □

**Corollary 6.19.** Let $(\mathcal{G}, \mathcal{Y}_\mathcal{V}, \theta)$ be a graphical model with the graph $\mathcal{G}$ being acyclic. Let also $\mu^* \in \mathcal{L}$ be a solution of the local polytope relaxation of the MAP-inference problem. Then from $\mu_u^*(s) > 0$ for some $u \in \mathcal{V}$, $s \in \mathcal{Y}_v$ it follows that there is an optimal integer labeling $y^*$ such that $y_u = s$.

*Proof.* The labeling $y^*$ can be constructed as in the proof of Proposition 6.17 with $(u, s)$ being selected as the first node and label. □

## 6.4   Bibliography and further reading

The local polytope relaxation for graphical models, its dual and its analysis were first provided by Schlesinger (see [114, 59]). An algorithm equivalent to the relaxation labeling was independently proposed by [136], Schlesinger (see Shlezinger [114]) and revisited by [91].

The comprehensive overview [140] is a standard reference for the relation of the primal and dual formulations of the energy minimization problem. It includes definitions of arc-consistency, the relaxation labeling algorithm and a lot of related notions and facts.

How to estimate a primal relaxed solution during dual optimization of the MAP-inference problem is described in [97].

We refer to [92] and references therein for a definition of constraint satisfaction problems in general. The description of its solvable subclasses was first given in [17], see also [20]. An analysis of solvability of constraint satisfaction problems, where min/max operations instead of $\vee/\wedge$ is given in [133]. A simple generalization of the relaxation labeling algorithm for such problems can be found in [104].

# 7

# Background: Basics of Non-Smooth Convex Optimization

The goal of this chapter is to give a brief overview of the simplest optimization methods applicable to large-scale convex programs, such as the Lagrange dual of the MAP-inference problem. To keep our exposition short, we restrict our overview to the three most basic techniques: gradient descent, the subgradient method and block-coordinate descent. The latter two will be used for MAP-inference, whereas the first one plays the role of a baseline for comparisons. These methods are applicable to smooth and non-smooth problems with different convergence guarantees. Since the convergence proofs and the derivation of the convergence rates are often quite involved and the techniques used in the proofs do not play any significant role for further chapters of the monograph, we mostly omit the proofs in this chapter and refer to the corresponding literature instead.

Throughout the chapter we will use the notation $\|\cdot\|$ for the Euclidean norm in $\mathbb{R}^n$.

## 7.1 Gradient descent

We start the chapter with differentiable functions and the most basic algorithm for their optimization. Although for MAP-inference we will deal mostly with non-smooth optimization, we provide the basic results about gradient descent as a baseline for comparison to other methods.

**Lipschitz-continuity**  We will be interested in differentiable functions with continuous gradient. Its "degree of continuity" is determined by the following definition:

**Definition 7.1.** A mapping $f\colon \mathbb{R}^n \to \mathbb{R}^m$ is called *Lipschitz-continuous* on $X \subseteq \mathbb{R}^n$ if there is a constant $L \geq 0$ such that for any $x, z \in X$ it holds that

$$\|f(x) - f(z)\| \leq L\|x - z\|. \tag{7.1}$$

The smallest value $L$ satisfying the condition above is called the *Lipschitz constant* for $f$ on $X$.

Figure 7.1 illustrates Definition 7.1. Informally speaking, the Lipschitz constant is an upper bound for the speed of change of the value of the mapping $f$. The following statement considers the limit case of this inequality and allows us to estimate the value of the Lipschitz constant in simple cases:

**Proposition 7.2.** If $f\colon \mathbb{R}^n \to \mathbb{R}^m$ is differentiable, then

$$L = \sup_{x \in X} \|\nabla f(x)\|_2,$$

where $\|\cdot\|_2$ is a spectral norm, i.e. the largest eigenvalue of the matrix $\nabla f$. In a special case, when $m = 1$, it coincides with the Euclidean norm.

**Example 7.3.** Proposition 7.2 implies that:

- $f(x) = x$ is Lipschitz-continuous on $\mathbb{R}$;

- $f(x) = x^n$, $n = 2, 3, \ldots$ is Lipschitz-continuous on any bounded subset of $\mathbb{R}$, e.g. on $[0, 1]$, but not on $\mathbb{R}$ itself.

We will use the notation $C_L^{1,1}(X)$ for functions $f\colon X \to \mathbb{R}$, which are differentiable on $X$ and whose gradient is Lipschitz-continuous on $X$ with the Lipschitz constant $L$. In general $C_L^{k,m}(X)$ is a standard notation for the set of $k$-times differentiable functions such that their $m$th derivative is Lipschitz-continuous on $X$ with the constant $L$.

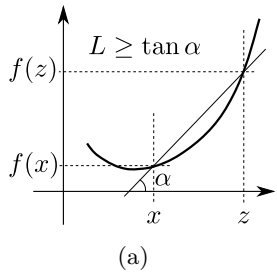

(a)

**Figure 7.1:** Illustration of the Lipschitz continuity and Lipschitz constant $L$.

**Example 7.4.** $C_L^{1,1}(\mathbb{R})$ includes such functions as $x$ and $x^2$. The functions $x^n$ for $n = 3, 4, \ldots$ belong to $C_L^{1,1}(X)$ for any bounded $X \subset \mathbb{R}$, but not for $X = \mathbb{R}$.

**Theorem 7.1.** Let $f \colon \mathbb{R}^n \to \mathbb{R}$ be convex and differentiable. Then $\nabla f(x) = \bar{0}$ is the necessary and sufficient condition for $x \in \mathbb{R}^n$ to be a minimum of $f$.

*Proof.* This follows directly from the similar property of the subgradient expressed by Theorem 5.1 together with Proposition 5.27. □

**Gradient descent algorithm**   Let $f \colon \mathbb{R}^n \to \mathbb{R}$ and $x^0 \in \mathbb{R}^n$ be an initial point. The iterative process of the form

$$x^{t+1} = x^t - \alpha^t \nabla f(x^t) \tag{7.2}$$

for $t = 0, \ldots, \infty$ and some $\alpha^t > 0$ is called *gradient descent*. The value $\alpha^t$ is called the *step-size* of the algorithm.

Let $x^* \in \arg\min_{x \in \mathrm{dom} f} f(x)$ and $f^* = f(x^*)$. The following statement specifies convergence properties of the gradient descent:

**Theorem 7.2** ([81, 76]). Let $f \in C_L^{1,1}(\mathbb{R}^n)$ be convex. Then with $\alpha^t = \frac{1}{L}$ for the sequence $x^t$ defined by (7.2) it holds that:

$$\|x^{t+1} - x^*\|^2 \le \|x^t - x^*\|^2 - \frac{1}{L^2}\|\nabla f(x^t)\|^2, \tag{7.3}$$

$$f(x^{t+1}) \le f(x^t) - \frac{1}{2L}\|\nabla f(x^t)\|^2, \tag{7.4}$$

$$f(x^t) - f^* \leq \frac{2L\|x^0 - x^*\|}{t + 4} \, . \tag{7.5}$$

Let us consider the statement of the theorem in detail. First of all, it considers a *constant* step-size $\alpha := \alpha^t = \frac{1}{L}$, which is inverse-proportional to the speed of change of the gradient $\nabla f$ expressed by its Lipschitz constant. In other words, the faster the gradient changes, the smaller step-size must be selected. Therefore, the Lipschitz constant defines the vicinity in which the considered function has "approximately the same" gradient. The larger the constant the smaller is the vicinity and, hence, the smaller the step-size.

Expression (7.3) states that each step of the algorithm produces a solution estimate $x^{t+1}$ which is closer to the optimum $x^*$ than the previous estimate $x^t$.

Expression (7.4) shows that the algorithm is strictly monotonous, i.e. unless $\nabla f = \bar{0}$, the objective value strictly decreases on each iteration.

Expression (7.5) provides the *convergence rate* of the gradient descent algorithm. The accuracy $\epsilon = f(x^t) - f^*$ is attained after at most $O(\frac{L}{t})$ iterations. This complexity is often formulated in a "dual" fashion, when one says that the algorithm requires $O(\frac{L}{\epsilon})$ iterations to attain a given accuracy $\epsilon$.

**Step-size selection for gradient descent**     One may object that a single Lipschitz constant $L$ can be a too rough estimation for a behavior of a function on its whole domain and, therefore, the constant step-size policy $\alpha = \frac{1}{L}$ may not work well in practice. This is indeed often the case. Consider the function $x^3$ on the interval $[0, b]$. The Lipschitz constant of its derivative $3x^2$ on this interval is bounded by the value $6b$ according to Proposition 7.2. If $b = 1$ this results in $L = 6$ and if $b = 100$ Proposition 7.2 gives $L = 600$. The step-size $\frac{1}{600}$ would lead to a much slower convergence in the vicinity of $x = 1$ than the step-size $\frac{1}{6}$. In other words, the step-size and the convergence speed depend on the domain on which the estimation of the Lipschitz constant was done. Therefore, a variable step-size quite often works better in practice.

There are three typical ways to define $\alpha^t$ in (7.2):

- $\alpha^t = \frac{1}{L}$ - constant step size. Requires knowing the Lipschitz constant or a good estimate of it. May be very inefficient if $\|\nabla f\|$ significantly varies over the domain of $f$.

- $\alpha^t = \arg\min_{\alpha \in \mathbb{R}_+} f(x^t - \alpha \nabla f(x^t))$, which is the step-size that minimizes $f$ in the negative gradient direction. For efficiency of the gradient descent algorithm as a whole, a closed form solution for $\alpha^t$ is typically required, which is often not available.

- $\alpha^t$ is selected to satisfy

$$f(x^{t+1}) \le f(x^t) - \frac{\alpha^t}{2} \|\nabla f(x^t)\|^2. \tag{7.6}$$

Here $\frac{1}{\alpha^t}$ plays the role of a local estimate for $L$. To this end $\alpha^t$ is searched in the form $\alpha^t = \frac{2\alpha^{t-1}}{2^n}$ for $n = 0, \ldots, \infty$ until (7.6) is satisfied. In other words, one tries to double the step size, checks condition (7.6) and divides the step-size by two until the condition is not satisfied anymore. The initial value $\alpha^0$ can be selected arbitrary, as soon as more than just few iterations of the algorithm are used. See the *Goldstein-Armijo* rule in [76, 15] for a substantiation.

Though the practical behavior can differ for the above three cases, their worst-case analysis is similar and for any of these strategies it results in the convergence rate given by Theorem 7.2.

## 7.2 Sub-gradient method

However, we can not use the gradient descent for our dual MAP-inference problem (6.9) directly, because it is not differentiable.

**Definition 7.5.** Let $x^0 \in \mathbb{R}^n$ be a starting point and $f \colon \mathbb{R}^n \to \mathbb{R}$ be a convex function. The iterative process of the form

$$x^{t+1} = x^t - \alpha^t g(x^t), \tag{7.7}$$

for $t = 0, \ldots, \infty$, some $\alpha^t > 0$ and $g(x^t) \in \partial f(x^t)$, is called *a subgradient method*.

The update rule (7.7) is basically the same as (7.2). However, when the function $f$ is non-smooth, the step-size $\alpha^t$ must be selected differently than in the smooth case:

**Theorem 7.3** ([76, 8])**.** Let $f$ be convex and Lipschitz continuous in a ball $B_R(x^*) = \{x \in \mathbb{R}^n \colon \|x^* - x\| \leq R\}$ with Lipschitz constant $M$. Let $x^0 \in B_R(x^*)$ be the initial point and let the step-size $\alpha^t$ satisfy

$$\alpha_t > 0, \ \alpha^t \xrightarrow{t \to \infty} 0, \ \text{and} \ \sum_{t=1}^{\infty} \alpha_t = \infty. \tag{7.8}$$

Then the following holds:

- $f(x^t) - f^* \xrightarrow{t \to \infty} 0$;

- [8] $0 < \alpha^t < 2\frac{f(x^t) - f^*}{\|g^t\|}$ implies

$$\|x^{t+1} - x^*\| \leq \|x^t - x^*\|; \tag{7.9}$$

- [76] In particular, $\alpha^t = \frac{R}{\|g^t\|\sqrt{t+1}}$ for all $t = 1, 2 \ldots$ implies

$$f(x^t) - f^* \leq \frac{MR}{\sqrt{t+1}}. \tag{7.10}$$

Note the following important properties of the sub-gradient method stated by Theorem 7.3:

- For a non-differentiable function $f$ the convergence is guaranteed only if the step-size $\alpha^t$ satisfies the *diminishing step-size rule* (7.8). This is due to the fact that the subgradient need not vanish in the vicinity of the minimum, contrary to the gradient, see Figure 7.2 for illustration. Since neither $f(x^t) - f^*$ nor $R$ is typically known, the claims of Theorem 7.3 can be simplified as *there exists $\alpha^t$* such that (7.9) holds, and *there exist $\{\alpha^t\}$* such that (7.10) holds. Step-sizes which decrease as $O(\frac{1}{t})$ or $O(\frac{1}{\sqrt{t}})$ are typical choices that satisfy condition (7.8).

- The update rule (7.7) does not guarantee the monotonic improvement of the function value if $f$ is non-differentiable. In other words,

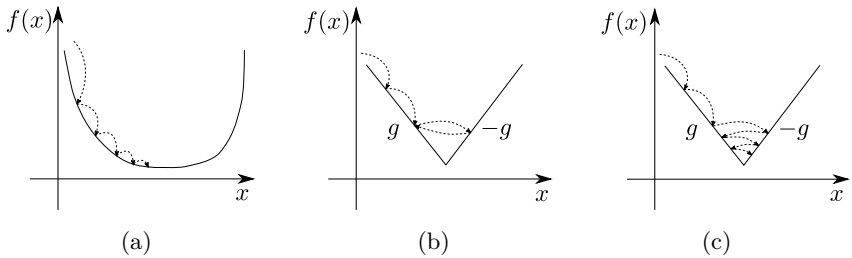

**Figure 7.2:** Illustration of the step-size rules for minimizing smooth and non-smooth functions. **(a)** The function $f$ is smooth, therefore its gradient vanishes as one gets closer to the optimum. Therefore, the size of the update steps $\alpha^t \nabla f$ will vanish as well, even for a constant step-size $\alpha^t = \alpha$. **(b)** The function $f$ is non-smooth, hence its subgradient does not vanish in general. Therefore, a constant step-size leads to oscillations around the optimum. **(c)** A vanishing step-size allows for convergence to the optimum even for non-smooth functions.

it may happen that $f(x^{t+1}) \geq f(x^t)$, even for an arbitrary small $\alpha^t$, see Figure 7.3 for illustration. It implies that minimization of $f$ in the direction of a negative subgradient does not make sense, contrary to the case when $f$ is smooth. On the other hand, the distance to the optimum $x^*$ never grows during iterations, according to (7.9).

- The convergence rate defined by (7.10) can be written as $O(\frac{1}{\sqrt{t}})$, or, alternatively, as $O(\frac{1}{\epsilon^2})$. This is significantly slower than the convergence rate $O(\frac{L}{t})$ in the smooth case. For example, to attain the precision $\epsilon = 0.1$ the subgradient algorithm must perform 100 times more iterations than to obtain the precision $\epsilon = 1$. Note that the gradient descent would require only 10 times more iterations given that the function to be optimized is smooth.

Table 7.1 summarizes the difference between gradient- and subgradient descent algorithms, or, more precisely, between the cases where the function $f$ to be optimized is smooth (is in $C_L^{1,1}$) and non-smooth.

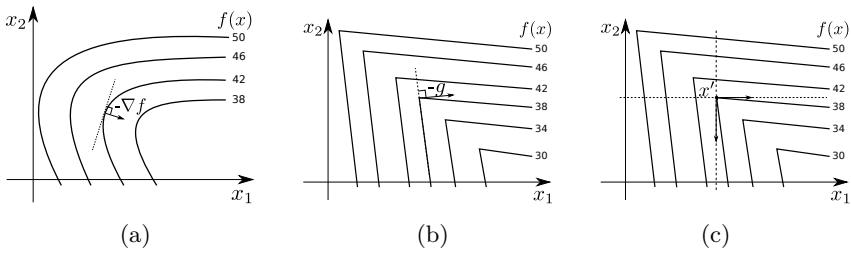

(a)             (b)             (c)

**Figure 7.3:** Level sets, i.e. lines of a constant function value, of smooth (a) and non-smooth (b-c) convex functions are depicted. **(a)** For differentiable functions the negated gradient always points to the direction of the function decrease, therefore for small enough step-sizes the gradient descent (7.2) is monotonous. **(b)** A negated subgradient of a non-differentiable function may point to the direction in which the function increases. Therefore, the subgradient algorithm (7.7) does not guarantee the monotonous decrease of the function on each iteration. **(c)** The point $x'$ is the directional minimum of the function $f$, i.e. when moving along the axes $x_1$ and $x_2$ from the point $x'$, the value of $f$ increases compared to $f(x')$.

**Table 7.1:** Difference of (sub)gradient descent algorithms for smooth and non-smooth optimization.

| | Smooth $f$ | Non-smooth $f$ |
|---|---|---|
| Update rule | $x^{t+1} = x^t - \alpha^t \nabla f(x^t)$ | $x^{t+1} = x^t - \alpha^t g^t, \ g^t \in \partial f(x^t)$ |
| Step-size $\alpha_t > 0$ | $\alpha^t = \frac{1}{L}$ | $\alpha^t \xrightarrow{t \to \infty} 0, \ \sum_{t=1}^{\infty} \alpha_t = \infty$ |
| Distance to optimum | $\|x^{t+1} - x^*\|^2 \le \|x^t - x^*\|^2 - \frac{1}{L^2}\|\nabla f(x^t)\|^2$ | $\|x^{t+1} - x^*\| \le \|x^t - x^*\|$ |
| Monotonicity | $f(x^{t+1}) \le f(x^t) - \frac{1}{2L}\|\nabla f(x^t)\|^2$ | — |
| Convergence rate | $f(x^t) - f^* \le \frac{2L\|x^0 - x^*\|}{t+4}$ | $f^t - f^* \le \frac{MR}{\sqrt{t+1}}$ |

## 7.3   Coordinate descent

Coordinate or block-coordinate descent is the last optimization method in our short overview. Like the subgradient method it is well-defined for both smooth and non-smooth functions, although the corresponding convergence guarantees differ significantly.

**Definition 7.6.** For a function $f \colon \prod_{j=1}^{n} \mathbb{R}^{n_j} \to \mathbb{R}$ the iterative process

$$x_i^{t+1} = \arg\min_{\xi \in \mathbb{R}^{n_i}} f(x_1^{t+1}, \ldots, x_{i-1}^{t+1}, \xi, x_{i+1}^t, \ldots, x_n^t) \tag{7.11}$$

$$i = (i+1) \mod n \tag{7.12}$$

$$t = 1, 2, \ldots \tag{7.13}$$

is called *cyclic* (or *Gauss-Seidel*) *coordinate descent*. Here $i \mod n$ denotes the remainder from the integer division of $i$ to $n$.

The method is also known as *alternating minimization*, since one alternates the optimization w.r.t. the different variables. When the dimensionality $n_j$ of individual variables is greater than one, one often speaks also about *block*-coordinate descent.

Note that by definition the coordinate descent is monotonous, i.e. it holds that

$$\underbrace{f(x_1^{t+1}, \ldots, x_{i-1}^{t+1}, x_i^t, x_{i+1}^t, \ldots, x_n^t)}_{f_{i-1}^t}$$

$$\geq \underbrace{f(x_1^{t+1}, \ldots, x_{i-1}^{t+1}, x_i^{t+1}, x_{i+1}^t, \ldots, x_n^t)}_{f_{i-1}^{t+1}} . \tag{7.14}$$

and therefore, the sequence $f_i^t$ is monotonously non-increasing w.r.t. $k = n \cdot t + i \to \infty$. Assuming that the minimal value of $f$ exists, i.e. $f^* > -\infty$, this implies convergence of the sequence $f_i^t$. To analyze this convergence and the convergence of the argument sequences $x_i^t$ for each $i$ as $t \to \infty$, we will require the following simple lemma:

**Lemma 7.7.** Let $f \colon \mathbb{R}^n \times \mathbb{R}^m \to \mathbb{R}$ be a convex function of two (vector-)variables. Then for any fixed $z' \in \mathbb{R}^m$ the function $f_{z'}(x) := f(x, z')$ is convex w.r.t. $x$ and its optimum is attained in $x \in \mathbb{R}^n$ if and only if $\bar{0} \in \partial f_{z'}(x)$.

*Proof.* The epigraph of $f_{z'}$ is an intersection of the epigraph of $f$ and the linear subspace $z = z'$. Therefore, it is convex as an intersection of two convex sets. This proves convexity of $f_{z'}$. The claim about optimality of $x$ follows from Theorem 5.1. $\qquad\square$

Let us now consider a function $f\colon \prod_{j=1}^{n} \mathbb{R}^{n_j} \to \mathbb{R}$ of $n$ variables. Let $f_i(\xi; x') := f(x_1', \ldots, x_{i-1}', \xi, x_{i+1}', \ldots, x_n')$ be the function of the $i$-th variable given all others are fixed to the corresponding coordinates of $x'$. Assume $x'$ is a fixed point of the algorithm (7.11), i.e. $x_i' \in \arg\min_{\xi \in \mathbb{R}^{n_i}} f_i(\xi; x')$. Due to Lemma 7.7 it holds that

$$\overline{0} \in \partial f_i(x_i'; x') \text{ for all } i = 1, \ldots, n. \tag{7.15}$$

If $f$ is convex differentiable, then condition (7.15) is equivalent to $\nabla f(x') = \overline{0}$ and, therefore, a fixed point is an optimum.

However, if $f$ is non-differentiable, this implication does not work, as shown in Figure 7.3(c).

This is due to the non-uniqueness of the subgradient, which causes this behavior of the fixed point for non-differentiable functions. In other words, there might exist a zero subgradient for each coordinate, but not for all of them together.

Moreover, in general, the algorithm (7.11) does not even guarantee convergence to such a point where the conditions (7.15) are fulfilled. This negative statement holds even if $f$ is differentiable, as shown by [82]. Indeed, to guarantee convergence, an additional *uniqueness* property is required:

**Theorem 7.4** ([8]). Let $f$ be continuously differentiable on $\mathbb{R}^n$. Furthermore, let for each $i$ and $x \in \mathbb{R}^n$ the minimum in (7.11) be *uniquely* attained. Then for every limit point $x^*$ of the sequence $x^t$, $t = 1, 2, \ldots$ defined by (7.11) it holds that $\nabla f(x^*) = 0$.

If $f$ is convex differentiable the optimality of $x$ follows from the condition $\nabla f(x) = \overline{0}$. Theorem 7.4 does not claim neither existence nor uniqueness of the limit points. It only proofs their properties if they exist. In particular, if the set $F(x^0) := \{x\colon f(x) \leq f(x^0)\}$ is bounded, then $x^t \in F(x^0)$ due to the monotonicity of the coordinate descent, i.e. $f(x^t) \leq f(x^0)$. Hence $\{x^t\}$ is bounded and has limit points, which all correspond to the optimum of $f$ according to Theorem 7.4.

**Convergence rate** A little is known about the convergence rate of cyclic coordinate descent (7.11) for general convex differentiable functions. Recently, the convergence rate $O(\frac{1}{t})$ was proven by [6] for the case $n = 2$, i.e. when only two blocks of variables are considered in (7.11).

## 7.4 Bibliography and further reading

In this chapter we have reviewed classical results which can be found in the text-books [8, 76, 81].

The subgradient method as well as a number of its variants and applications to the non-smooth minimization have been proposed by N. Z. Shor in 1960's, see [115] and references therein.

Although coordinate descent belongs to one of the earliest algorithms, the results about its convergence and the corresponding convergence rates in general are still missing. Earlier works on this topic mostly considered specializations to different function subclasses e.g. [71, 129, 131, 130]. The most general convergence results for variants of the coordinate descent were recently given in [6].

# 8

# Subgradient and Coordinate Descent for MAP-Inference

In this chapter we apply the optimization methods from Chapter 7 to the MAP-inference problem.

First of all, we consider one of the simplest optimization methods for this problem, known as *Iterated Conditional Modes (ICM)*. This algorithm is a direct implementation of the coordinate descent technique for discrete energy functions.

Afterwards, we consider more powerful techniques addressing the Lagrangean relaxation and the respective dual problem. As we learned in Chapter 6, the latter can be seen as an unconstrained concave non-smooth maximization problem. Therefore, the subgradient method is applicable directly and we provide the corresponding analysis.

However, the most efficient existing methods for the Lagrangean dual are based on the block-coordinate ascent technique. Although these methods do not attain the dual optimum in general, they are often able to attain practically good approximate solutions in a moderate number of iterations.

Moreover, we show that most of the existing block-coordinate ascent methods converge to node-edge agreement, as introduced in Chapter 6. To this end we introduce a broad class of block-coordinate ascent algorithms called *anisotropic diffusion* and provide a general convergence proof scheme for these and similar methods.

We learned in Chapter 6 that the LP relaxation is tight for acyclic graphical models. Here, we additionally show that exact inference for acyclic models, based on the dynamic programming algorithm, can be interpreted as a dual block-coordinate ascent, specifically, as the anisotropic diffusion algorithm, which we introduce below.

We conclude this chapter with an empirical comparison of the considered methods on several popular benchmark MAP-inference problem instances.

## 8.1 Primal coordinate descent

### 8.1.1 Iterated Conditional Modes

Iterated Conditional Modes (ICM) is one of the first and simplest algorithms for energy minimization. As we show below, it is an implementation of the coordinate descent (7.11) for the primal non-relaxed energy minimization objective $E(y; \theta)$.

The algorithm starts with some labeling $y \in \mathcal{Y}_\mathcal{V}$ and iteratively tries to improve it. In its elementary step it minimizes the objective w.r.t. the label of a node $u \in \mathcal{V}$, whereas the labels of other nodes are kept fixed.

Let the notation

$$l(y, u, s) := y' \in \mathcal{Y}_\mathcal{V} \quad \text{with} \quad y'_v = \begin{cases} y_v, & v \neq u \\ s, & v = u \end{cases} \qquad (8.1)$$

define a labeling where all coordinates but the one corresponding to node $u$ coincide with the labeling $y$ and the coordinate $u$ is assigned the label $s$.

Then the elementary step of the algorithm can be written as

$$y_u := \arg\min_{s \in \mathcal{Y}_u} E(l(y, u, s); \theta). \qquad (8.2)$$

Algorithm 4 applies this rule to each node sequentially and iterates this procedure until the labeling $y$ does not change anymore.

---

**Algorithm 4** ICM: Coordinate Descent for $E(y; \theta)$

---
1: **Init:** $y \in \mathcal{Y}_\mathcal{V}$
2: **repeat**
3:     **for** $u \in \mathcal{V}$ **do**
4:         $y_u := \arg\min_{s \in \mathcal{Y}_u} E(l(y, u, s); \theta)$
5:     **end for**
6: **until** the labeling $y$ does not change anymore

---

As any coordinate descent, Algorithm 4 guarantees that the objective value of $E$ is monotonically nonincreasing. However, since $E$ as a function of $y$ is neither convex nor differentiable (moreover, it is defined on a discrete set $\mathcal{Y}_\mathcal{V}$ only), Algorithm 4 does not guarantee attainment of the optimal value of $E$ and its output significantly depends on the initial labeling. In practice Algorithm 4 typically returns labelings with significantly higher energies than many other techniques.

The following example demonstrates that Algorithm 4 may fail to attain the optimum even in pretty simple cases:

**Example 8.1.** Consider a simple two-node model with two labels in each node as shown in Figure 8.1. Suppose that the unary costs are all zero and the pairwise are $\theta_{uv}(0,0) = 0$, $\theta_{uv}(1,1) = 1$, $\theta_{uv}(0,1) = \theta_{uv}(1,0) = 2$. If the initial labeling is $(1,1)$, the ICM updates cannot improve the labeling, since they can switch only to the labelings $(1,0)$ or $(0,1)$, each with total cost 2, while the labeling $(1,1)$ has total cost 1. However, the optimum is attained in labeling $(0,0)$, which has total cost 0

## 8.1.2   Block-ICM algorithm

The idea of the ICM algorithm goes beyond Algorithm 4. In general, one has to fix labels in a subset of nodes and optimize with respect to the labels in the remaining nodes. Should this optimization be tractable then the whole algorithm is as well. Below we describe one such tractable case.

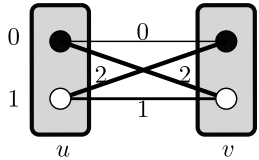

**Figure 8.1:** The ICM algorithm does not guarantee attainment of the optimum even for a simple binary 2-node model. Suppose that the unary costs are all zero and the pairwise costs are $0, 1, 2$, the thickness of the lines grows with the cost: $\theta_{uv}(0,0) = 0$, $\theta_{uv}(1,1) = 1$, $\theta_{uv}(0,1) = \theta_{uv}(1,0) = 2$. If the initial labeling is $(1,1)$, the ICM updates cannot improve the labeling, see Example 8.1 for a detailed explanation.

**Definition 8.2.** A subgraph $(\mathcal{V}', \mathcal{E}')$ of a graph $(\mathcal{V}, \mathcal{E})$ is called *induced* by its set of nodes $\mathcal{V}'$, if $\mathcal{E}' = \{uv \in \mathcal{E} : u, v \in \mathcal{V}'\}$, i.e. its set of edges contains exactly those edges from $\mathcal{E}$, that connect nodes from $\mathcal{V}'$.

Let $\mathcal{G} = (\mathcal{V}, \mathcal{E})$ be a graph defining a graphical model and $\mathcal{V}^f \subset \mathcal{V}$ be a set of nodes with fixed labels. The corresponding *partial labeling* will be denoted as $y^f \in \mathcal{Y}_{\mathcal{V}^f}$.

Introducing $\mathcal{V}' = \mathcal{V} \backslash \mathcal{V}^f$ as complement of $\mathcal{V}^f$ we can generalize the notation defined in (8.1) such that

$$l(y^f, \mathcal{V}^f, y') := \tilde{y} \in \mathcal{Y}_{\mathcal{V}} \quad \text{with} \quad \tilde{y}_u = \begin{cases} y_u^f, & u \in \mathcal{V}^f \\ y_u', & u \in \mathcal{V}' \end{cases} \quad (8.3)$$

defines a labeling, where nodes from the set $\mathcal{V}^f$ are labeled with $y^f$, otherwise with $y'$.

Let $(\mathcal{V}', \mathcal{E}')$ be the subgraph of $\mathcal{G}$ induced by $\mathcal{V}'$. It turns out that if it is acyclic, one can optimize w.r.t. the labeling on $\mathcal{V}'$,

$$y := \arg\min_{y' \in \mathcal{Y}_{\mathcal{V}'}} E(l(y^f, \mathcal{V}^f, y'); \theta), \quad (8.4)$$

using dynamic programming as shown below.

As Example 8.1 in particular shows, Algorithm 4 does not solve the MAP-inference on acyclic models in general. Therefore, considering blocks of variables associated with such acyclic models for a block-coordinate descent may result in better approximate solutions.

Below we derive the expression for the unary and pairwise costs of the acyclic model to be optimized over on each iteration of the considered block-ICM algorithm. Let $(\mathcal{V}^f, \mathcal{E}^f)$ be the subgraph of $\mathcal{G}$ induced by $\mathcal{V}^f$ and let $\mathcal{E}'^f = \{uv \in \mathcal{E} : u \in \mathcal{V}', \ v \in \mathcal{V}^f\}$ be the set of edges of $\mathcal{G}$ connecting the nodes from $\mathcal{V}'$ and $\mathcal{V}^f$ (see Figure 8.2 for illustration). Consider the energy $E(l(y^f, \mathcal{V}^f, y'); \theta)$, defined by (8.3):

$$
\begin{aligned}
E(&l(y^f, \mathcal{V}^f, y'); \theta) \\
&= \sum_{u \in \mathcal{V}'} \theta_u(y'_u) + \sum_{uv \in \mathcal{E}'} \theta_{uv}(y'_u, y'_v) + \sum_{uv \in \mathcal{E}'^f} \theta_{uv}(y'_u, y^f_v) \\
&\quad + \sum_{u \in \mathcal{V}^f} \theta_u(y^f_u) + \sum_{uv \in \mathcal{E}^f} \theta_{uv}(y^f_u, y^f_v) \\
&= \sum_{u \in \mathcal{V}'} \left( \theta_u(y'_u) + \sum_{v \in \mathcal{N}_b(u) \cap \mathcal{V}^f} \theta_{uv}(y'_u, y^f_v) \right) \\
&\quad + \sum_{uv \in \mathcal{E}'} \theta_{uv}(y'_u, y'_v) + \sum_{u \in \mathcal{V}^f} \theta_u(y^f_u) + \sum_{uv \in \mathcal{E}^f} \theta_{uv}(y^f_u, y^f_v) . \quad (8.5)
\end{aligned}
$$

Recall that $y^f$ is fixed. Therefore, the last expression represents the energy of a problem on the graph $(\mathcal{V}', \mathcal{E}')$ with modified unary costs, plus a constant represented by the last two terms of the expression. It implies that minimization of $E(l(y^f, \mathcal{V}^f, y'); \theta)$ w.r.t. $y'$ can be done efficiently by dynamic programming if the graph $(\mathcal{V}', \mathcal{E}')$ is acyclic.

Algorithm 5 summarizes these observations. At each iteration it generates a subset $\mathcal{V}'$ of nodes over which it minimizes the corresponding auxiliary energy (8.5) with dynamic programming, as described in Chapter 2, until the current labeling does not change for a certain number of iterations or the total iteration limit is reached.

The way the sets $\mathcal{V}'$ are generated can be either quite generic or depend on a particular graph structure. For grid-graphs a simple and well-parallelizable approach is to consider four graphs $(\mathcal{V}', \mathcal{E}')$ in a cyclic order: columns with even indexes, columns with odd indexes, rows with even indexes and rows with odd indexes. Note that minimization in line 5 of Algorithm 5 can be done for all columns (rows) of each of these subgraphs in parallel, see Figures 8.2(b) and 8.2(c).

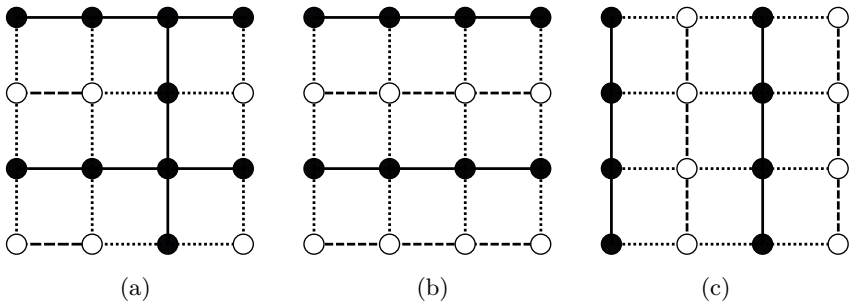

(a)                          (b)                          (c)

**Figure 8.2:** A grid graphical model and different acyclic induced subgraphs. Black and white circles stand for graph nodes, solid, dashed and dotted lines for graph edges. White circles and dashed lines denote $\mathcal{V}^f$ and $\mathcal{E}^f$ respectively, i.e. the part of the model with fixed labels. Dotted lines stand for the edges from the set $\mathcal{E}'^f$ connecting nodes from $\mathcal{V}'$ and $\mathcal{V}^f$. The block-ICM algorithm optimizes over the labels in the nodes $\mathcal{V}'$ denoted as black circles and connected by solid lines denoting the set of edges $\mathcal{E}'$. **(b,c)** Note that the subgraphs $(\mathcal{V}', \mathcal{E}')$ consist of several connected components, which can be processed independently in parallel.

---

**Algorithm 5** Block-ICM: Block-Coordinate Descent for $E(y; \theta)$

---

1: **Init:** $y \in \mathcal{Y}_\mathcal{V}$
2: **repeat**
3:     Generate $\mathcal{V}' \subseteq \mathcal{V}$ such that $(\mathcal{V}', \mathcal{E}')$ induced by $\mathcal{V}'$ is acyclic
4:     Define $\mathcal{V}^f := \mathcal{V} \backslash \mathcal{V}'$ and $y^f := y|_{\mathcal{V}^f}$
5:     Compute $y^* := \arg\min_{y' \in \mathcal{Y}_{\mathcal{V}'}} E(l(y^f, \mathcal{V}', y'); \theta)$
6:     Re-assign the values of $y$ on the coordinates of $\mathcal{V}'$: $y|_{\mathcal{V}'} := y^*$
7: **until** some stopping condition holds

---

## 8.2 Dual sub-gradient method

Now we turn to the maximization of the Lagrange dual. Recall, it has the form of an unconstrained piecewise linear concave function

$$\mathcal{D}(\phi) = \sum_{u \in \mathcal{V}} \min_{s \in \mathcal{Y}_u} \theta_u^\phi(s) + \sum_{uv \in \mathcal{E}} \min_{(s,t) \in \mathcal{Y}_{uv}} \theta_{uv}^\phi(s, t), \qquad (8.6)$$

with the reparametrized costs $\theta^\phi$ defined by (6.6).

The subgradient method (7.7) is the simplest way to optimize concave non-smooth functions in general. To be able to use it we have to compute a subgradient of $\mathcal{D}$.

Consider its single coordinate $\frac{\partial \mathcal{D}}{\partial \phi_{u,v}(s)}$ and the function

$$\mathcal{D}_{u,v}(\phi_{u,v}) := \min_{s' \in \mathcal{Y}_u} \theta_u^\phi(s') + \min_{(s',t') \in \mathcal{Y}_{uv}} \theta_{uv}^\phi(s',t'), \qquad (8.7)$$

containing exactly those two terms of $\mathcal{D}$, which depend on $\phi_{u,v}$. Other terms do not contribute to $\frac{\partial \mathcal{D}}{\partial \phi_{u,v}(s)}$, due to the linear properties of the subgradient (Proposition 5.28) and due to the fact that the subgradient of a constant function is identical to zero (since any constant function is differentiable and has gradient zero, which is also its subgradient, see Proposition 5.27). In other words, $\frac{\partial \mathcal{D}}{\partial \phi_{u,v}(s)} = \frac{\partial \mathcal{D}_{u,v}}{\partial \phi_{u,v}(s)}$.

Due to the linearity of the subgradient we can compute it for each term of $\mathcal{D}_{u,v}$ separately. Recall also how $\theta^\phi$ depends on $\phi$:

$$\theta_u^\phi(s) := \theta_u(s) - \sum_{v \in \mathcal{N}_b(u)} \phi_{u,v}(s), \ u \in \mathcal{V}, \ s \in \mathcal{Y}_u, \qquad (8.8)$$

$$\theta_{uv}^\phi(s,t) := \theta_{uv}(s,t) + \phi_{u,v}(s) + \phi_{v,u}(t), \ uv \in \mathcal{E}, \ (s,t) \in \mathcal{Y}_{uv}.$$

Consider the first term of $\mathcal{D}_{u,v}$,

$$\min_{s' \in \mathcal{Y}_u} \theta_u^\phi(s'). \qquad (8.9)$$

Its subgradients are by virtue of Lemma 5.32 equal to the convex hull of the vectors $\frac{\partial \theta_u^\phi(s')}{\partial \phi_{u,v}(s)}$ for all minimizers $s'$ of (8.9). Since the numerator is a linear and, therefore, differentiable function of $\phi$, we can apply the standard differentiation rules, yielding $\frac{\partial \theta_u^\phi(s')}{\partial \phi_{u,v}(s)} = -[\![s = s']\!]$.

Similarly, subgradients of the second term in (8.7) can be computed as $\frac{\partial \theta_{uv}^\phi(s'',t'')}{\partial \phi_{u,v}(s)} = [\![s = s'']\!]$ for each minimizer $(s'',t'')$ of the second term. Note that if $s' = s'' = s$ the subgradients of both terms cancel out. We now combine our observations to construct a subgradient of $\mathcal{D}$. Let

$$y_u' \in \operatorname*{arg\,min}_{s \in \mathcal{Y}_u} \theta_u^\phi(s) \qquad (8.10)$$

and

$$(y_u'', y_v'') \in \operatorname*{arg\,min}_{(s,t) \in \mathcal{Y}_{uv}} \theta_{uv}^\phi(s,t) \qquad (8.11)$$

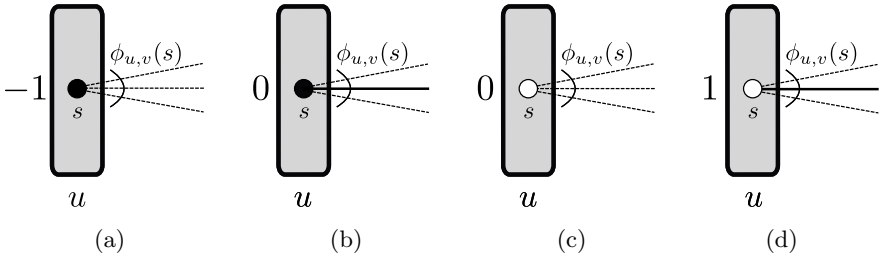

**Figure 8.3:** Illustration to Equation (8.12). Black circles correspond to $y'_u$, solid lines to $(y''_u, y''_v)$. White circles and dashed lines are labels and label pairs, which are not equal to $y'_u$ and $(y''_u, y''_v)$ respectively. The number to the left of the rectangles denoting the node $u$, is the value of the corresponding subgradient coordinate.

be defined for all $u \in \mathcal{V}$ and all $uv \in \mathcal{E}$. The vector $\frac{\partial \mathcal{D}}{\partial \phi}$ with coordinates

$$\frac{\partial \mathcal{D}}{\partial \phi_{u,v}(s)} = \begin{cases} 0, & s \neq y'_u \quad \text{and} \quad s \neq y''_u \\ 0, & s = y'_u \quad \text{and} \quad s = y''_u \\ -1, & s = y'_u \quad \text{and} \quad s \neq y''_u \\ 1, & s \neq y'_u \quad \text{and} \quad s = y''_u, \end{cases} \tag{8.12}$$

is a subgradient of $\mathcal{D}$. Figure 8.3 illustrates all possible combinations occurring in (8.12).

Note that (8.12) defines only one possible subgradient. It is unique only if $y'_u$ and $(y''_u, y''_v)$ are *unique* solutions of the respective local minimization problems (8.10) and (8.11). Multiple solutions imply multiple subgradients of the form (8.12) and the convex combinations thereof.

However, any of these subgradients can be used in the subgradient method (7.7), which is now defined up to the step-size $\alpha^t$:

$$\phi^{t+1} = \phi^t + \alpha^t \frac{\partial \mathcal{D}}{\partial \phi}. \tag{8.13}$$

As noted in Chapter 7, the step-size $\alpha^t$ must be diminishing, i.e. satisfy (7.8). There is no single best performing rule, its choice can significantly depend on the specific problem type or even problem instance. One practical way to deal with this is to search for a suitable step-size in parametric form, e.g.

$$\alpha^t = \beta(1+t)^\gamma, \tag{8.14}$$

and tune the parameters, here $\beta > 0$ and $0 > \gamma \geq -1$, on similar problem instances.

**Remark 8.3.** Note that subgradients computed according to (8.12) are quite sparse. For a given node $u$ and edge $uv$ there are $|\mathcal{Y}_{uv}|$ entries of the subgradient, corresponding to the $|\mathcal{Y}_{uv}|$ dual variables $\phi_{u,v}(s)$, $s \in \mathcal{Y}_u$. According to (8.12) at most two coordinates of the subgradient have non-zero entries. The rest $|\mathcal{Y}_u| - 2$ are zero. This, in turn, implies that only a small number of coordinates of the dual vector $\phi$ are changed on each iteration of the subgradient method. This is one of the reasons, why this algorithm converges quite slowly in practice, especially for problems with a large number of labels $\mathcal{Y}_u$.

**Sufficient dual optimality condition**    Note that Theorem 5.1 provides a necessary and sufficient optimality condition for the Lagrange dual $\mathcal{D}$. It is the existence of a zero element in $\partial \mathcal{D}(\phi)$. However, checking this condition is in general computationally costly, as it would require to check whether the convex polyhedral set $\partial \mathcal{D}$ is non-empty. This is in general as difficult as solving a linear program of the size comparable to the size of the relaxed MAP-inference problem. Therefore, a simplified *sufficient* optimality condition is typically used, where a zero subgradient is searched only in the form given by (8.12). To this end it must be verified, whether $y'_u = y''_u$ for all $u \in \mathcal{V}$, or, in other words, if there is an integer labeling consisting of locally optimal labels and label pairs. Indeed, Item 5 of Proposition 6.2 also implies sufficiency of this condition.

**Exercise 8.4.** Consider Example 4.3. Show that there is a zero subgradient in the point $\phi = 0$, although it does not have the form (8.12).

**Exercise 8.5.** Consider Example 6.10. Show that there is no zero subgradient in the point $\phi = 0$.

## 8.3   Min-sum diffusion as block-coordinate ascent

Till the end of this chapter we will concentrate on block-coordinate ascent methods for the Lagrange dual (8.6). We start with one of the first such methods proposed for the MAP-inference problem, *min-sum diffusion*, see §8.3.1. Although it is not the most efficient algorithm of the considered type, its simplicity allows us to illustrate the properties typical for most block-coordinate ascent algorithms for the Lagrange dual of the MAP-inference problem. An additional advantage of min-sum diffusion as a starting point of our consideration is that it can be turned into the state-of-the-art block-coordinate ascent method by a small modification. This modification will be considered in §8.3.2.

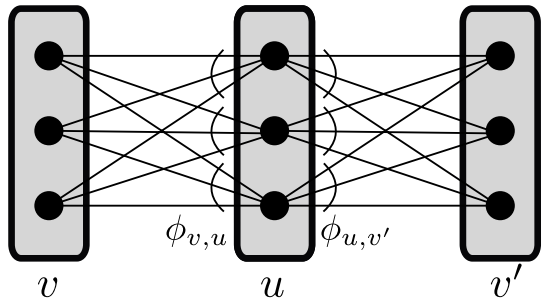

**Figure 8.4:** One elementary step of the min-sum diffusion algorithm optimizes over dual variables $\phi_{u,v}(s)$, $v \in \mathcal{N}_b(u)$, $s \in \mathcal{Y}_u$, "attached" to the current node $u$. The variables are "sitting" on the pencils, associated with each label of the node $u$ and denoted by small arcs.

### 8.3.1 Min-sum diffusion

To construct a block-coordinate ascent algorithm one has to first partition the coordinates into blocks. To get an efficient algorithm, the minimization with respect to each block should be efficient. Ideally, the minimum should have a closed form or be attained by a finite-step algorithm with linear complexity.

One such possibility is to assume that one block consists of variables $(\phi_{u,v}(s)\colon v \in \mathcal{N}_b(u), s \in \mathcal{Y}_u)$ "attached" to one node of the graph, see Figure 8.4. For a given edge $uv$ and label $s \in \mathcal{Y}_u$ we will refer to the set of label pairs $\{(s,l),\ l \in \mathcal{Y}_v\}$ as a *pencil*. Each pencil is defined by the triple $(s, u, v)$, where $u \in \mathcal{V}$, $v \in \mathcal{N}_b(u)$ and $s \in \mathcal{Y}_u$. Pencils $(s, u, v)$, $v \in \mathcal{N}_b(u)$, will be called *associated* with the label $s$ in node $u$.

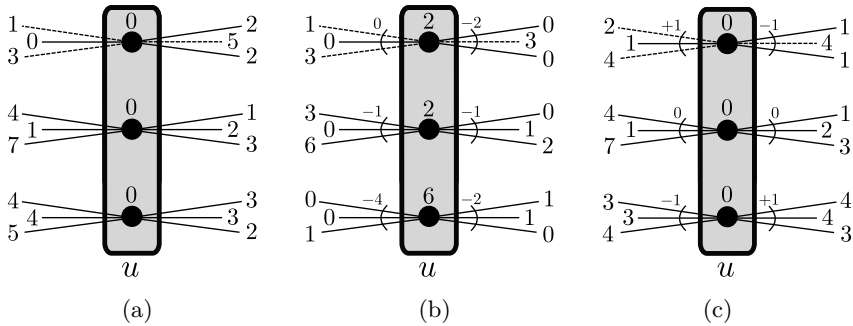

**Figure 8.5:** Example of an elementary update step of the diffusion algorithm. Larger numbers stand for (reparametrized) costs, smaller numbers for the corresponding coordinates of the dual vector $\phi$, assuming that initially $\phi = 0$. Solid and dashed lines denote minimal and non-minimal label pairs in each pencil respectively. **(a)** Initial costs. **(b)** Costs and dual variables after (8.15). **(c)** Costs and dual variables after (8.16). Note that unary costs are zero and minimal label pairs in each pencil get equal weights, as stated in (8.17) and (8.18) respectively.

**The elementary update step**   of the diffusion algorithm consists of the following two operations,

$$\forall v \in \mathcal{N}_b(u) \; \forall s \in \mathcal{Y}_u \quad \phi_{u,v}^{t+1}(s) := \phi_{u,v}^t(s) - \min_{l \in \mathcal{Y}_v} \theta_{uv}^{\phi^t}(s, l)\,, \qquad (8.15)$$

$$\forall v \in \mathcal{N}_b(u) \; \forall s \in \mathcal{Y}_u \quad \phi_{u,v}^{t+2}(s) := \phi_{u,v}^{t+1}(s) + \frac{\theta_u^{\phi^{t+1}}(s)}{|\mathcal{N}_b(u)|}\,, \qquad (8.16)$$

and is illustrated in Figure 8.5. Loosely speaking, executing one elementary update step for a node $u$ redistributes the minimum costs evenly between all pencils of a given label while moving all costs away from the node itself.

For fixed $s$ and $v$ the first operation (8.15) computes the minimal cost in the pencil $(s, u, v)$, subtracts it from the costs of all label pairs of the pencil and adds it to the costs of the label $s$. As a result, the minimal cost in each pencil associated with the label $s$ becomes zero.

The second operation then redistributes the reparametrized unary costs $\theta_u^\phi(s)$ now corresponding to the label $s$ equally between all pencils $(s, u, v)$, $v \in \mathcal{N}_b(u)$. As a result, the unary cost of the label $s$ becomes zero and the minimal pairwise costs in all pencils associated with $s$ become equal to each other.

In other words, after operations (8.15) and (8.16) have been executed for node $u \in \mathcal{V}$ the following holds:

$$\theta_u^\phi(s) = 0, \ \forall s \in \mathcal{Y}_u \,, \tag{8.17}$$

$$\min_{t \in \mathcal{Y}_v} \theta_{uv}^\phi(s,t) = \min_{t \in \mathcal{Y}_{v'}} \theta_{uv'}^\phi(s,t), \ \forall v,v' \in \mathcal{N}_b(u), \ s \in \mathcal{Y}_u \,. \tag{8.18}$$

**Proposition 8.6.** Operations (8.15)-(8.16) maximize the Lagrange dual $\mathcal{D}$ w.r.t. the block of variables $(\phi_{u,v}(s) \colon v \in \mathcal{N}_b(u), s \in \mathcal{Y}_u)$.

*Proof.* By virtue of Lemma 7.7 it is sufficient to prove the existence of a zero subgradient in a subdifferential of $\mathcal{D}$ when the latter is treated as a function of only the considered block of variables.

We will search for a zero subgradient in the form (8.12). Since only the costs in the node $u$ and its incident edges $uv$, $v \in \mathcal{N}_b(u)$, are dependent on the considered variable block, it is sufficient to show that

$$\exists y_u \in \mathcal{Y}_u \ \forall v \in \mathcal{N}_b(u) \ \exists y_v \in \mathcal{Y}_v, :$$
$$y_u \in \operatorname*{arg\,min}_{s \in \mathcal{Y}_u} \theta_u^\phi(s) \text{ and } (y_u, y_v) \in \operatorname*{arg\,min}_{(s,l) \in \mathcal{Y}_{uv}} \theta_{uv}^\phi(s,l) \,. \tag{8.19}$$

See Figure 8.3(b) for illustration.

Note that $y_u \in \arg\min_{s \in \mathcal{Y}_u} \theta_u^\phi(s)$ holds for any $y_u \in \mathcal{Y}_u$ since all labels in node $u$ have the same cost after the elementary diffusion step due to (8.17), namely cost 0. Let us, therefore, assign $y_u$ such that $(y_u, y_{v'})$ is a minimizer of $\theta_{uv'}^\phi$ for some $v' \in \mathcal{N}_b(u)$. Condition (8.18) implies that it then also minimizes $\theta_{uv}^\phi$ for all $v \in \mathcal{N}_b(u)$. This finalizes the proof. $\qquad\square$

Note that the sufficient block-optimality condition (8.19) can be equivalently expressed as

$$\bigvee_{s \in \mathcal{Y}_u} \left( \xi_u(s) \wedge \bigwedge_{v \in \mathcal{N}_b(u)} \bigvee_{t \in \mathcal{Y}_v} \xi_{uv}(s,t) \right) = 1 \,. \tag{8.20}$$

for the binary vector $\xi := \mathrm{mi}[\theta^\phi]$. Note that (8.20) is precisely the value computed by the node-wise operation (6.24) of the relaxation labeling algorithm.

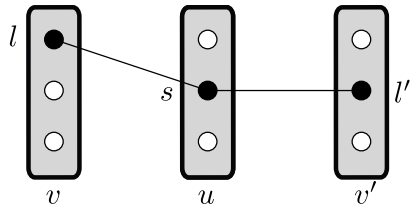

**Figure 8.6:** Illustration to the proof of Lemma 8.8. Black labels $s, l, l'$ and the solid lines connecting them correspond to locally optimal label pairs $(s, l)$ and $(s, l')$.

**Remark 8.7.** If the block-optimality condition (8.20) holds, the elementary step (8.15)-(8.16) of the diffusion algorithm applied to node $u$ does not change the dual value, since the dual vector is already optimal with respect to the considered block of variables.

**The min-sum diffusion algorithm**   consists in a cyclic repetition of the operations (8.15)-(8.16) for all $u \in \mathcal{V}$, which corresponds to the cyclic block-coordinate ascent method (7.11), up to substitution of the convex function $f$ by the concave function $\mathcal{D}$ and minimization w.r.t. each coordinate block by maximization.

Since $\mathcal{D}$ is neither differentiable nor does it have a unique optimum, the min-sum diffusion is not guaranteed to optimize the dual $\mathcal{D}(\phi)$. Its convergence properties are analyzed below.

First, we will show that node-edge agreement is a fixed point of min-sum diffusion, i.e. once this condition is satisfied at some iteration of the algorithm, it is satisfied at all further iterations as well. Afterwards we show that min-sum diffusion actually converges to node-edge agreement.

Let the function $F_u \colon \mathbb{R}^{\mathcal{I}} \to \mathbb{R}^{\mathcal{I}}$ define the transformation of costs performed by a single elementary update step (8.15)-(8.16) applied to the node $u$. In other words, $F_u$ transforms a cost vector $\theta \in \mathbb{R}^{\mathcal{I}}$ to another cost vector $F_u(\theta) \in \mathbb{R}^{\mathcal{I}}$.

Then the following holds:

**Lemma 8.8.** $\mathrm{cl}(\mathrm{mi}[\theta]) \leq \mathrm{cl}(\mathrm{mi}[F_u(\theta)])$.

*Proof.* Since $\mathrm{cl}(\mathrm{mi}[F_u(\theta)])$ is the maximal arc-consistent vector contained in $\mathrm{mi}[F_u(\theta)]$, it suffices to show that $\mathrm{cl}(\mathrm{mi}[\theta]) \leq \mathrm{mi}[F_u(\theta)]$ due to Proposition 6.13.

The application of a single elementary update step (8.15)-(8.16) only affects the entries at $u$ and its incident edges. For any other $w \in \mathcal{V} \backslash \{u\} \cup \mathcal{E} \backslash \{uv \colon v \in \mathcal{N}_b(u)\}$, $s \in \mathcal{Y}_w$, it is $F_u(\theta)_w(s) = \theta_w(s)$, and, therefore, $\mathrm{mi}[F_u(\theta)]_w(s) = \mathrm{mi}[\theta]_w(s)$. This in turn implies $\mathrm{cl}(\mathrm{mi}[\theta])_w(s) \leq \mathrm{mi}[F_u(\theta)]_w(s)$. It remains to prove this property for $u$ and its incident edges.

Let $\theta_u^*$ denote the cost of a locally minimal label on $u$. Analogously, let $\theta_{uv}^*$ denote the cost of a locally minimal pair on the edge $uv$ for any $v \in \mathcal{N}_b(u)$, $(s,t) \in \mathcal{Y}_{uv}$,

$$F_u(\theta)_{uv}(s,t) = \theta_{uv}(s,t) - \min_{t' \in \mathcal{Y}_v} \theta_{uv}(s,t')$$

$$+ \frac{1}{|\mathcal{N}_b(u)|} \left( \theta_u(s) + \sum_{v' \in \mathcal{N}_b(u)} \min_{t' \in \mathcal{Y}_{v'}} \theta_{uv'}(s,t') \right)$$

$$\geq 0 + \frac{1}{|\mathcal{N}_b(u)|} \left( \theta_u^* + \sum_{v' \in \mathcal{N}_b(u)} \theta_{uv'}^* \right). \qquad (8.21)$$

Consider now a node $v \in \mathcal{N}_b(u)$ and a label pair $(s,t) \in \mathcal{Y}_{uv}$ such that $\mathrm{cl}(\mathrm{mi}[\theta])_{uv}(s,t) = 1$. Since $\mathrm{cl}(\mathrm{mi}[\theta])$ is arc-consistent this implies

$$\theta_{uv}(s,t) = \theta_{uv}^*,$$
$$\theta_u(s) = \theta_u^*, \qquad\qquad (8.22)$$
$$\min_{t' \in \mathcal{Y}_{v'}} \theta_{uv'}(s,t') = \theta_{uv'}^* \text{ for all } v' \in \mathcal{N}_b(u).$$

$$(8.23)$$

Thus, equality holds in (8.21), which implies

$$F_u(\theta)_{uv}(s,t) = \min_{(s',t') \in \mathcal{Y}_{uv}} F_u(\theta)_{uv}(s',t').$$

Therefore, $\mathrm{mi}[F_u(\theta)]_{uv}(s,t) = 1$. So the above statement holds for all edges incident to $u$.

Finally, note that by (8.17) all label costs in $u$ are minimal after an elementary update step in $u$, which means that $\mathrm{mi}[F_u(\theta)]_u(s) = 1$ holds for all $s \in \mathcal{Y}_u$. Therefore, $\mathrm{cl}(\mathrm{mi}[\theta])_u(s) \leq \mathrm{mi}[F_u(\theta)]_u(s)$ for all $s \in \mathcal{Y}_u$ which finalizes the proof. □

**Corollary 8.9.** The statement of Lemma 8.8 holds for an arbitrary number of elementary step applied to any nodes of the graph.

*Proof.* We can prove this by induction. By Lemma 8.8 the statement holds for the first transformation. Let $F\colon \mathbb{R}^{\mathcal{I}} \to \mathbb{R}^{\mathcal{I}}$ be the mapping describing some number of update steps and $\mathrm{cl}(\mathrm{mi}[\theta]) \leq \mathrm{cl}(F(\mathrm{mi}[\theta]))$. By Lemma 8.8 it holds that $\mathrm{cl}(F(\mathrm{mi}[\theta])) \leq \mathrm{cl}(\mathrm{mi}[F_u(F(\theta))])$. Combining these two inequalities implies $\mathrm{cl}(\mathrm{mi}[\theta]) \leq \mathrm{cl}(\mathrm{mi}[F_u(F(\theta))])$, which finalizes the proof. □

Let $F\colon \mathbb{R}^{\mathcal{I}} \to \mathbb{R}^{\mathcal{I}}$ stand for the transformation of the cost vector after some number of elementary steps (8.15)-(8.16) has been applied to arbitrary nodes. Corollary 8.9 together with the Remark 8.7 implies the following statement:

**Theorem 8.1.** Let $\theta^\phi \in \mathbb{R}^{\mathcal{I}}$ satisfy node-edge agreement. Then $F(\theta^\phi)$ does so as well. Moreover, in this case the transformation $F$ does not change the dual value, i.e. $\mathcal{D}(\phi') = \mathcal{D}(\phi)$, where $\theta^{\phi'} = F(\theta^\phi)$.

*Proof.* The first part of the statement follows directly from Corollary 8.9.

To prove the second one, note that according to the definition, node-edge agreement implies the block-optimality condition (8.20). Therefore, as noted in Remark 8.7, $\mathcal{D}(\phi') = \mathcal{D}(\phi)$. □

Theorem 8.1 guarantees that as soon as at some iteration of the min-sum diffusion algorithm the node-edge agreement holds, it holds for further iterations as well. The following theorem states that min-sum diffusion converges to node-edge agreement.

**Theorem 8.2** ([102]). Let $\theta^t$ be the cost vector produced on iteration $t$ of the min-sum diffusion algorithm. Let also $\epsilon^t \geq 0$ be the smallest value such that $\mathrm{cl}(\mathrm{mi}_{\epsilon^t}[\theta^t]) \neq \bar{0}$. Then $\epsilon^t \xrightarrow{t \to \infty} 0$.

In Section 8.4 we will prove a more general statement, which includes Theorem 8.2 as a special case.

Since min-sum diffusion converges to node-edge agreement and not necessarily to the optimum, its result in general depends on the initial point (reparametrization). However, compared to the ICM method, this dependence is usually quite weak for diffusion. In other words, typically

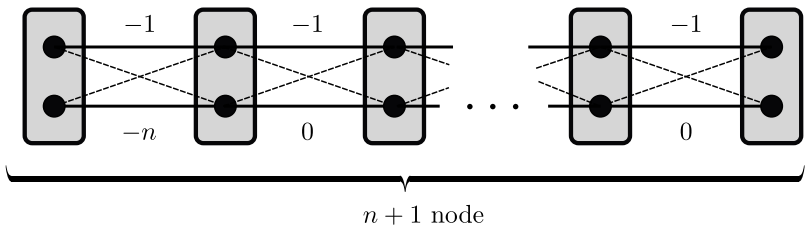

$n + 1$ node

**Figure 8.7:** The figure shows a chain-structured graph, where the diffusion algorithm converges only in the limit. The number of nodes in the chain is $n + 1$, which implies that the number of edges is $n$. Dashed label pairs have very high (infinite) costs, which does not influence the performance of the algorithm, while thick lines indicate label pairs assigned the given initial costs. Basically, there are only two labelings allowed: the "upper" one $(0, 0, \ldots, 0)$ and the "lower" one $(1, 1, \ldots, 1)$. All unary costs are equal to zero, therefore, the energies of both labelings are equal to $-n$. However, while the energy of the upper labeling is equally distributed between edges, the energy of the lower one is concentrated in a single edge. Updates of the diffusion algorithm do not change the upper labeling and redistribute the weights between label pairs of the lower one. However, the arc-consistent configuration is attained only when all lower edges have the same weight equal to $-1$, but the min-sum diffusion algorithm can not attain this point in a finite number of iterations, although it converges to it.

quite similar (or even the same) results are obtained when the algorithm is initialized differently.

**Min-sum diffusion for acyclic graphs**  According to Theorem 8.2 and Proposition 6.17 the diffusion algorithm converges to the dual optimum for acyclic graphs, where node-edge agreement implies dual optimality (see §6.3). However, contrary to dynamic programming, diffusion may require an infinite number of iterations to attain node-edge agreement even in this case. Example 8.10 illustrates this observation.

**Example 8.10.** Consider the chain-structured model as in Figure 8.7. The min-sum diffusion algorithm attains an arc-consistent configuration, which for acyclic models coincides with the dual optimum, only in the limit. Note also that on this problem instance the subgradient algorithm exhibits the same convergence only in the limit.

**Rounding for min-sum diffusion**  To obtain an integer labeling one has to apply some rounding method. Unfortunately, the naïve dual rounding procedure (6.21) can not be used directly with the min-sum diffusion, since it takes into account only the reparametrized unary costs and those are all equal to zero in the diffusion algorithm. One way to deal with that is to remember the locally optimal label after each "collection step" (8.15) of the diffusion algorithm, when the pairwise costs are shifted to the unary costs.

Another one, which typically gives similar result when the algorithm converged, is to "equally" split the pairwise costs between corresponding unary factors, i.e. to perform

$$\forall v \in \mathcal{N}_b(u) \ \forall s \in \mathcal{Y}_u \quad \phi_{u,v}(s) := \phi_{u,v}(s) - \frac{1}{2} \min_{l \in \mathcal{Y}_v} \theta_{uv}^\phi(s, l) \,, \qquad (8.24)$$

for all nodes $u \in \mathcal{V}$. Afterwards the naïve dual rounding (6.21) can be applied.

### 8.3.2 Anisotropic diffusion

Note that the proof of the block-coordinate ascent property (Proposition 8.6) for min-sum diffusion is based on fulfilling condition (8.20) after the corresponding update step (8.15)-(8.16). Let $\omega_{u,v}$ for all $v \in \mathcal{N}_b(u)$ be non-negative numbers such that $\sum_{v \in \mathcal{N}_b(u)} \omega_{u,v} \leq 1$. As the Proposition 8.11 below states, the condition (8.20) is fulfilled also after an update defined as

$$\forall v \in \mathcal{N}_b(u) \ \forall s \in \mathcal{Y}_u \quad \phi_{u,v}^{t+1}(s) := \phi_{u,v}^t(s) - \min_{t \in \mathcal{Y}_v} \theta_{uv}^{\phi^t}(s, t) \,, \qquad (8.25)$$

$$\forall v \in \mathcal{N}_b(u) \ \forall s \in \mathcal{Y}_u \quad \phi_{u,v}^{t+2}(s) := \phi_{u,v}^{t+1}(s) + \omega_{u,v} \, \theta_u^{\phi^{t+1}}(s) \,. \qquad (8.26)$$

Updates (8.25)-(8.26) differ from (8.15)-(8.16) by the coefficients $\omega_{u,v}$, which weight different edges of the graph and may allow for some of the costs to remain in the node labels after reparametrization. For the min-sum diffusion these coefficients correspond to a uniform distribution over the edges, i.e. $\omega_{u,v} = \frac{1}{|\mathcal{N}_b(u)|}$.

**Proposition 8.11.** For any $\omega_{u,v} \geq 0$, $uv \in \mathcal{E}$, such that $\sum_{v \in \mathcal{N}_b(u)} \omega_{u,v} \leq 1$, the operations (8.25)-(8.26) maximize the Lagrange dual $\mathcal{D}$ w.r.t. the block of variables $(\phi_{u,v}(s) \colon v \in \mathcal{N}_b(u), s \in \mathcal{Y}_u)$.

The proof of Proposition 8.11 is similar to the proof of Proposition 8.6 for the elementary updates of the min-sum diffusion algorithm and is omitted.

We will refer to algorithms based on the updates (8.25)-(8.26) as *anisotropic diffusion* algorithms. Particular settings of $\omega_v$ correspond to different block-coordinate ascent algorithms. The setting $\omega_{u,v} = \frac{1}{|\mathcal{N}_b(u)|+1}$ is known as the *convex message passing* algorithm [35]. Another setting, which relates the anisotropic diffusion and the dynamic programming and leads to a state-of-the-art dual optimization method is considered below.

**Dynamic programming as anisotropic diffusion**

Recall Algorithm 1 for computing forward min-marginals on a chain graph using dynamic programming. Assume as in Chapter 2 the set of nodes to be fully ordered, i.e. $\mathcal{V} = \{1, \ldots, n\}$. Algorithm 1 can be expressed in terms of reparametrizations as shown in Algorithm 6.

---

**Algorithm 6** Dynamic programming as anisotropic diffusion

---

1: Given: $(\mathcal{G}, \mathcal{Y}_\mathcal{V}, \theta)$ - a chain-structured graphical model, $\phi = 0$
2: **for** $i = 1$ **to** $n-1$ **do**
3:    $\phi_{i,i+1}(s) := \theta_i^\phi(s)$ for all $s \in \mathcal{Y}_i$
4:    $\phi_{i+1,i}(s) := -\min_{t \in \mathcal{Y}_i} \theta_{i,i+1}^\phi(t,s)$ for all $s \in \mathcal{Y}_{i+1}$
5:    Set $r_{i+1}(s)$ to arg min of the computation in line 4.
6: **end for**
7: **return** $E^* = \min_{s \in \mathcal{Y}_n} \theta_n^\phi(s)$

---

Step 3 of Algorithm 6 adds the reparametrized unary costs of node $i$ to the pairwise costs of the "next" edge $(i, i+1)$. This turns the reparametrized unary costs of node $i$ to 0. Step 4 selects an optimal "preceding" pairwise cost and adds it to the unary cost of the node $i+1$. In total, the dual variables $\phi_{i+1,i}(s)$ become equal to the negated forward min-marginals $(-F_{i+1}(s))$.

In other words, to turn the forward min-marginals $F_i(s)$ into the reparametrization returned by Algorithm 6 it is sufficient to assign

$$\phi_{i,i-1}(s) := -F_i(s), \ i = 2, \ldots, n \,, \tag{8.27}$$
$$\phi_{i,i+1}(s) := \theta_i(s) + F_i(s), \ i = 1, \ldots, n-1 \,. \tag{8.28}$$

By construction, the reparametrized costs $\theta^\phi$ computed by Algorithm 6 possess the following properties:

- $\theta_i^\phi(s) = 0$ for all $i \in \mathcal{V} \backslash \{n\}$, $s \in \mathcal{Y}_i$;

- $\min_{s \in \mathcal{Y}_{i-1}} \theta_{i-1,i}^\phi(s,t) = 0$ for all $i \in \mathcal{V} \backslash \{1\}$ and $t \in \mathcal{Y}_i$ and, therefore,

- $\theta_{i,i+1}^\phi(s,t) \geq 0$ for all $\{i, i+1\} \in \mathcal{E}$, $(s,t) \in \mathcal{Y}_{i,i+1}$.

These properties imply:

**Proposition 8.12.** The reparametrization $\theta^\phi$ obtained by Algorithm 6 is optimal, i.e. it maximizes the dual $\mathcal{D}(\phi)$.

*Proof.* According to Proposition 6.17 it is sufficient to show that there is node-edge agreement. For that it is enough to find a labeling consisting of locally minimal labels and label pairs. Let us show that the optimal labeling $y$ reconstructed by Algorithm 2 from pointers $r_i(s)$ possesses this property. Indeed, by construction $\theta_{i,i+1}^\phi(y_i, y_{i+1}) = \theta_i^\phi(y_i) = 0$ for all $i \in \mathcal{V} \backslash \{n\}$. Additionally, $\theta_n^\phi(y_n) = \min_{s \in \mathcal{Y}_n} \theta_n^\phi(s) = \min_{s \in \mathcal{Y}_n} F_n(s) + \theta_n(s)$ due to line 2 of Algorithm 2.     □

Note that Algorithm 6 can be seen as anisotropic diffusion with $\omega_{i,i+1} = 1$ (and, therefore, $\omega_{i,i-1} = 0$) applied to the nodes in the order of the nodes in the chain.

**Anisotropic diffusion and state-of-the-art methods for general graphs.** In Chapter 10 we will show another important property of anisotropic min-sum diffusion. Namely, we will show that it can be seen as the main building block for the state-of-the-art Sequential Tree-Reweighted Message Passing algorithm (TRW-S) [49, 51].

**Convergence of anisotropic diffusion**    Note that although all anisotropic diffusion algorithms perform block-coordinate ascent, not all of them converge to a non-empty arc-consistent closure in the sense of Theorem 8.2.

**Exercise 8.13.** Show that the anisotropic diffusion algorithm with $\omega_{u,v} = 0$ for all $v \in \mathcal{N}_b(u)$ gets to its fix-point after a single iteration. Moreover, the arc-consistency closure of the fix point potentials is in general empty.

## 8.4    Convergence of dual block-coordinate ascent methods

The following section is quite technical and has a single goal: provide a general scheme for the convergence proof for dual coordinate ascent algorithms. In particular, this can allow us to show convergence of the diffusion and TRW-S algorithms to the node-edge agreement. This proof largely follows and generalizes the proof for the diffusion algorithm given by [102].

  This section consists of two parts. First, we formulate and prove the convergence theorem in general, assuming that an optimization algorithm has certain specified properties. Afterwards, we show that these properties hold for a subclass of anisotropic diffusion. This will imply convergence to node-edge agreement for this kind of algorithm.

## 8.5    Abstract convergence theorem

Let $f \colon \mathbb{R}^n \to \mathbb{R}$ be a *continuous* function which has to be maximized. The function $f(x)$ corresponds to the Lagrangean dual $\mathcal{D}(\phi)$ introduced in Chapter 6. Let also the mapping $F \colon \mathbb{R}^n \to \mathbb{R}^n$ stands for an iterative algorithm such that $F$ maps the current iterate $x^t$ to the next one $x^{t+1}$, i.e. $x^{t+1} = F(x^t)$. Let $F^i(\theta) := \underbrace{F(F(\ldots F(\theta)))}_{i \text{ times}}$ denote the $i$th composition of $F$. In our case the mapping $F$ corresponds to a single iteration of a block-coordinate ascent algorithm, such as (anisotropic) diffusion (8.25)-(8.26).

  Let further $\varepsilon \colon \mathbb{R}^n \to \mathbb{R}_+$ be a *continuous* non-negative function such that its value $\varepsilon(x^t)$ specifies how well a certain property of $x^t$ is satisfied.

Equality $\varepsilon(x) = 0$ means that the property is fully satisfied for $x$. For the dual block-coordinate ascent algorithm the value $\epsilon = \varepsilon(x^t)$ defines the smallest value $\epsilon$ such that $\epsilon$-node-edge agreement holds. This $\epsilon$ can be computed by the corresponding modification (6.28) of the relaxation labeling algorithm. Since all steps of this algorithm represent continuous mappings and a superposition of a finite number of continuous mappings is a continuous mapping itself, the function $\varepsilon$ is indeed continuous for that case.

As before, let $\|\cdot\|$ stand for the $\ell_2$-norm.

**Definition 8.14.** A mapping $F\colon \mathbb{R}^n \to \mathbb{R}^n$, which satisfies the conditions below is called *consistency-enforcing*:

(1) $F$ is continuous;

(2) $F$ is monotone w.r.t. $f$: For all $x \in \mathbb{R}^n$ it holds that $f(F(x)) \geq f(x)$;

(3) $F$ is bounded w.r.t. $f$: For any $x \in \mathbb{R}^n$ there exists $M \in \mathbb{R}$ such that $f(F^i(x)) < M$ for any $i$;

(4) For all $x \in \mathbb{R}^n$ there exists $C_x > 0$ such that $\|F^i(x)\| \leq C_x$ for any $i$.

(5) $F$ is bounded w.r.t. $\varepsilon$: For all $x \in \mathbb{R}^n$ there exists $M_x \in \mathbb{R}$ such that $\varepsilon(F^i(x)) < M_x$ for any $i$;

(6) There is a natural number $N > 0$ such that $f(F^N(x)) = f(x)$ implies $\varepsilon(x) = 0$.

The following main convergence theorem generalizes Theorem 8.2 formulated for the diffusion algorithm in §8.3.1:

**Theorem 8.3.** For any $x \in \mathbb{R}^n$ and any consistency-enforcing $F$ it holds that $\lim_{i \to \infty} \varepsilon(F^i(x)) = 0$.

*Proof.* Consider some arbitrary fixed $x \in \mathbb{R}^n$. The proof consists of four steps, enumerated respectively:

1. By virtue of (4) the sequence $x^i = F^i(x)$ is bounded and, therefore, contains a converging subsequence $x^{i(j)}$, $j = 1, 2, \ldots$, where $j > j'$ implies $i(j) > i(j')$, i.e. the limit $x^* := \lim_{j \to \infty} x^{i(j)}$ exists. Let us show that

$$\varepsilon(x^*) = 0 \tag{8.29}$$

holds for any converging subsequence of $x^i$.

Due to (2) and (3) the sequence $f(x^i)$ is non-decreasing and bounded from above. Therefore, it converges and there exists $f^* := \lim_{i \to \infty} f(x^i)$. Hence, it also holds

$$f^* = \lim_{j \to \infty} f(x^{i(j)}) = \lim_{j \to \infty} f(x^{i(j)+N}) \tag{8.30}$$

for any natural number $N$. This implies

$$0 = \lim_{j \to \infty} \left( f(x^{i(j)}) - f(x^{i(j)+N}) \right)$$
$$= \lim_{j \to \infty} \left( f(x^{i(j)}) - f(F^N(x^{i(j)})) \right) . \tag{8.31}$$

Since $f$ and $F$ are continuous, it holds

$$0 = \lim_{j \to \infty} \left( f(x^{i(j)}) - f(F^N(x^{i(j)})) \right) = f(x^*) - f(F^N(x^*)) \tag{8.32}$$

and, therefore, (8.29) holds by virtue of (6).

2. Since $\varepsilon$ is a continuous function, (8.29) implies

$$\lim_{j \to \infty} \varepsilon(x^{i(j)}) = 0 \tag{8.33}$$

for any converging subsequence $x^{i(j)}$ of $x^i$.

3. Let us now consider the sequence $\varepsilon(x^i)$. By virtue of (5) it is bounded, and, therefore, there exists $\epsilon^i := \sup_{j \geq i} \varepsilon(x^j)$. The sequence $\epsilon^i$ is a monotonically non-increasing sequence of non-negative numbers and therefore it has a limit $\epsilon^* = \lim_{i \to \infty} \epsilon^i$.

4. According to the "Theorem about superior and inferior limits" from Analysis there exists a subsequence $\varepsilon(x^{i'(j)})$ of $\varepsilon(x^i)$ such that

$$\lim_{j \to \infty} \varepsilon(x^{i'(j)}) = \lim_{i \to \infty} \epsilon^i = \epsilon^* .$$

The sequence $x^{i'(j)}$ is bounded by virtue of (4) and, therefore, contains a converging subsequence $x^{i'(j(k))}$. For this subsequence it also holds

$$\lim_{k \to \infty} \varepsilon(x^{i'(j(k))}) = \epsilon^*.$$

At the same time, as proved in item 2, for any converging subsequence it holds that

$$0 = \lim_{k \to \infty} \varepsilon(x^{i'(j(k))}) = \epsilon^* = \lim_{i \to \infty} \sup_{k \geq i} \varepsilon(x^k).$$

Finally, $0 \leq \varepsilon(x^i) \leq \sup_{k \geq i} \varepsilon(x^k)$ implies $\lim_{i \to \infty} \varepsilon(x^i) = 0$.

<div style="text-align: right;">□</div>

## 8.6 Convergence of diffusion algorithms

To prove convergence of the diffusion algorithms it remains to show that the iterations of the algorithms are consistency-enforcing, that is, they comply with Definition 8.14.

To simplify the notation, we will associate the reparametrized costs $\theta^\phi$ with the variable $x$ from Definition 8.14, instead of the dual variables $\phi$ themselves. Therefore, to denote the reparametrized costs on the $t$-th iteration we will write $\theta^t$ instead of $\theta^{\phi^t}$.

First, we will show that properties (1-5) from Definition 8.14 hold for the elementary step (8.25)-(8.26) of the anisotropic diffusion algorithm with arbitrary weights $\omega_{u,v}$. It trivially implies that they also hold for any finite number of such steps. This, in particular, includes a sequence of updates performed for all nodes of the graph in a prespecified order, or in other words, for an iteration of the algorithm.

Afterwards, we will show that the most complicated condition (6) holds for min-sum diffusion (8.15)-(8.16). The proof of condition (6) for a wide subclass of the anisotropic diffusion methods is provided in [51].

**Proof of properties (1-5) for the anisotropic diffusion updates** (8.25)-(8.26)

(1) follows from the fact that operations (8.25)-(8.26) are defined as a superposition of continuous mappings, which is a continuous mapping itself.

(2) follows from Proposition 8.11, which states the ascent property of anisotropic diffusion.

(3) follows from the fact that the maximal dual value is finite and bounded from the top by the minimal energy value.

(4-5) are closely related to each other. To prove them it is sufficient to show that there is $M(\theta)$ such that $\max_{w \in \mathcal{V} \cup \mathcal{E}} \max_{s \in \mathcal{Y}_w} |\theta_w^t(s)| \leq M(\theta)$ for any iteration $t$. Property (4) would follow directly since

$$\|F^t(\theta)\| = \|\theta^t\| = \sqrt{\sum_{w \in \mathcal{V} \cup \mathcal{E}} \sum_{s \in \mathcal{Y}_w} (\theta_w^t(s))^2}$$

$$\leq M(\theta) \sqrt{\sum_{w \in \mathcal{V} \cup \mathcal{E}} |\mathcal{Y}_w|}. \quad (8.34)$$

Property (5) follows from the definition of $\epsilon$-node-edge agreement. As noted in §6.2.4, see (6.27), for

$$\epsilon = \max_{w \in \mathcal{V} \cup \mathcal{E}} \left( \max_{s \in \mathcal{Y}_w} \theta_w^t(s) - \min_{s \in \mathcal{Y}_w} \theta_w^t(s) \right)$$

$$\leq 2 \max_{w \in \mathcal{V} \cup \mathcal{E}} \max_{s \in \mathcal{Y}_w} |\theta_w^t(s)| \leq 2M(\theta) \quad (8.35)$$

the $\epsilon$-node-edge agreement always holds, which implies

$$\varepsilon(\theta^t) \leq \epsilon \leq 2M(\theta). \quad (8.36)$$

Therefore, the following statement finalizes the proof of properties (4) and (5):

**Lemma 8.15.** Let $F$ denote one step of the anisotropic diffusion algorithm (8.25)-(8.26) applied to some node $u$. Let also $F^t$ stands for $t$ such updates applied to a specified sequence of nodes. Then for any $\theta$ there is $M(\theta)$ such that

$$\max_{w \in \mathcal{V} \cup \mathcal{E}} \max_{s \in \mathcal{Y}_w} |F^t(\theta_w(s))| \leq M(\theta) \,. \qquad (8.37)$$

*Proof.* The proof consists of two part. First we study behavior of the updates (8.25)-(8.26) if the cost vector has some special properties like e.g. non-negativity of its components. Afterwards we consider a general cost vector $\theta$.

*Non-negative, non-positive and constant cost vectors.*

- Equations (8.25)-(8.26) imply that if all components of $\theta$ are non-negative, so are the components of $F^t(\theta)$ as well.

- Let $\theta$ be a cost vector such that the cost function is constant for each node and edge, i.e. :

$$\theta_w(s) = \theta_w(s') \text{ for all } w \in \mathcal{V} \cup \mathcal{E} \text{ and } s, s' \in \mathcal{Y}_w \,. \qquad (8.38)$$

Applying equations (8.25)-(8.26) one obtains that the same holds for $F^t(\theta)$ as well:

$$F^t(\theta)_w(s) = F^t(\theta)_w(s'), \ w \in \mathcal{V} \cup \mathcal{E}, \ s, s' \in \mathcal{Y}_w \,. \qquad (8.39)$$

- Let $\theta' = \theta + \hat{\theta}$, where $\theta$ satisfies (8.38). Then (8.25)-(8.26) imply

$$F^t(\theta') = F^t(\theta) + F^t(\hat{\theta}) \,. \qquad (8.40)$$

- Equations (8.25)-(8.26) imply that if $\theta$ is constant for each node and edge (satisfies (8.38)) and all its components are non-positive, so are the components of $F^t(\theta)$.

- Let us first assume that all components of the initial cost vector $\theta$ are non-negative. Since the anisotropic diffusion performs a reparametrization, it holds that $E(y; F^t(\theta)) = E(y; \theta)$ for all labelings $y$. Due to the non-negativity $F^t(\theta)_w(s) \leq E(y; F^t(\theta)) = E(y; \theta)$ for any $y$ such that $y_w = s$. It implies

$$F^t(\theta)_w(s) \leq E(y; \theta)) \leq \sum_{w \in \mathcal{V} \cup \mathcal{E}} \sum_{s \in \mathcal{Y}_w} |\theta_w(s)| \,. \qquad (8.41)$$

Therefore,

$$|F^t(\theta)_w(s)| \leq \sum_{w \in \mathcal{V} \cup \mathcal{E}} \sum_{s \in \mathcal{Y}_w} |\theta_w(s)| \,. \qquad (8.42)$$

*General cost vector $\theta$.* Let $\theta$ not satisfy the non-negativity constraint, which is, $\min_{\substack{w \in \mathcal{V} \cup \mathcal{E} \\ s \in \mathcal{Y}_w}} \theta_w(s) < 0$. Consider the non-negative cost vector $\theta'$ with coordinates $\theta'_w(s) := \theta_w(s) - \hat{\theta}_w(s)$, where $\hat{\theta}$ is a cost vector with only negative coordinates equal to each other and defined by

$$\hat{\theta}_w(s) := \min_{\substack{w' \in \mathcal{V} \cup \mathcal{E} \\ s' \in \mathcal{Y}_w}} \theta_{w'}(s') < 0, \ w \in \mathcal{V} \cup \mathcal{E}, \ s \in \mathcal{Y}_w. \quad (8.43)$$

Then due to (8.40) $F^t(\theta) = F^t(\theta') + F^t(\hat{\theta})$ and, therefore,

$$|F^t(\theta)| \leq |F^t(\theta')| + |F^t(\hat{\theta})|. \quad (8.44)$$

Due to non-positivity of $\hat{\theta}$, for any $y$ such that $y_w = s$ it holds that

$$F^t(\hat{\theta})_w(s) \geq E(y; F^t(\hat{\theta})) = E(y; \hat{\theta}) \geq (|\mathcal{V}| + |\mathcal{E}|)\hat{\theta}_w(s), \quad (8.45)$$

and, therefore,

$$|F^t(\hat{\theta})_w(s)| \leq (|\mathcal{V}| + |\mathcal{E}|)|\hat{\theta}_w(s)|. \quad (8.46)$$

Due to non-negativity of $\theta'$, according to (8.42),

$$|F^t(\theta')_w(s)| \leq \sum_{w \in \mathcal{V} \cup \mathcal{E}} \sum_{s \in \mathcal{Y}_w} |\theta'_w(s)|. \quad (8.47)$$

Adding the last two inequalities and taking into account (8.44) one obtains

$$|F^t(\theta)| \leq (|\mathcal{V}| + |\mathcal{E}|)|\hat{\theta}_w(s)| + \sum_{w \in \mathcal{V} \cup \mathcal{E}} \sum_{s \in \mathcal{Y}_w} |\theta'_w(s)|, \quad (8.48)$$

that finalizes the proof. $\qquad \square$

**Proof of property (6) for min-sum diffusion** Let, as before, $F_u$ be the mapping $\mathbb{R}^{\mathcal{I}} \to \mathbb{R}^{\mathcal{I}}$ corresponding to the elementary update step of the min-sum diffusion algorithm (8.15)-(8.16) applied to node $u$. Let also for two vectors $\xi, \xi' \in \mathbb{R}^n$ the relation $\xi < \xi'$ ($\xi \leq \xi'$) be applied coordinate-wise, i.e. for each coordinate $i$ it holds that $\xi_i < \xi'_i$ ($\xi_i \leq \xi'_i$). Without loss of generality we assume that all occurring unary costs are zero, as holds already after the first iteration of the min-sum diffusion algorithm. Abusing notation we will assume the Lagrange dual $\mathcal{D}$ to be a function of the (reparametrized) costs $\theta^\phi$, i.e. the notation $\mathcal{D}(\theta^\phi)$ till the end of this section corresponds to the previously used one $\mathcal{D}(\phi)$ and stands for $\sum_{uv \in \mathcal{E}} \min_{(s,l) \in \mathcal{Y}_{uv}} \theta^\phi_{uv}(s, l)$.

**Lemma 8.16.** Let $\mathcal{D}(\theta) = \mathcal{D}(F_u(\theta))$. Then $\mathrm{mi}[F_u(\theta))] \leq \mathrm{mi}[\theta]$, i.e. the set of locally optimal label pairs may increase only when the dual value increases.

*Proof.* Let $\theta' := F_u(\theta)$. Since the mapping $F_u$ does not change the values of $\theta_w$ for $w \in \mathcal{V} \cup \mathcal{E} \setminus \{uv \colon v \in \mathcal{N}_b(u)\}$, the lemma states that the equality

$$\theta'_{uv}(s, l) = \min_{(s',l') \in \mathcal{Y}_{uv}} \theta'_{uv}(s', l') \tag{8.49}$$

implies the equality

$$\theta_{uv}(s, l) = \min_{(s',l') \in \mathcal{Y}_{uv}} \theta_{uv}(s', l') \tag{8.50}$$

for all $v \in \mathcal{N}_b(u)$, as soon as $\mathcal{D}(\theta) = \mathcal{D}(\theta')$.

In particular, (8.49) implies

$$\min_{l' \in \mathcal{Y}_v} \theta'_{uv}(s, l') = \min_{(s',l') \in \mathcal{Y}_{uv}} \theta'_{uv}(s', l') \tag{8.51}$$

Let $l^*(v) \in \arg\min_{l \in \mathcal{Y}_v} \theta'_{uv}(s, l)$ define a minimal label pair $(s, l^*(v))$ in the pencil $(s, u, v)$. By virtue of (8.18), equality (8.51) then holds for all $v \in \mathcal{N}_b(u)$ and fixed $s \in \mathcal{Y}_u$. In the introduced notation it reads

$$\theta'_{uv}(s, l^*(v)) = \min_{(s',l') \in \mathcal{Y}_{uv}} \theta'_{uv}(s', l') . \tag{8.52}$$

Summing up over all $v \in \mathcal{N}_b(u)$ we obtain

$$\sum_{v \in \mathcal{N}_b(u)} \theta'_{uv}(s, l^*(v)) = \sum_{v \in \mathcal{N}_b(u)} \min_{(s',l') \in \mathcal{Y}_{uv}} \theta'_{uv}(s', l') . \tag{8.53}$$

Since $\mathcal{D}(\theta) = \mathcal{D}(\theta')$, it holds that

$$\sum_{v \in \mathcal{N}_b(u)} \min_{(s',l') \in \mathcal{Y}_{uv}} \theta_{uv}(s', l') = \sum_{v \in \mathcal{N}_b(u)} \min_{(s',l') \in \mathcal{Y}_{uv}} \theta'_{uv}(s', l') . \tag{8.54}$$

On the other hand, since $F_u$ performs a reparametrization, it holds that

$$\sum_{v \in \mathcal{N}_b(u)} \theta_{uv}(s, l^*(v)) = \sum_{v \in \mathcal{N}_b(u)} \theta'_{uv}(s, l^*(v)) . \tag{8.55}$$

Combining this with (8.53) and (8.54) we obtain

$$\sum_{v \in \mathcal{N}_b(u)} \theta_{uv}(s, l^*(v)) = \sum_{v \in \mathcal{N}_b(u)} \min_{(s', l') \in \mathcal{Y}_{uv}} \theta_{uv}(s', l') , \tag{8.56}$$

which is fulfilled if and only if the corresponding summands on the right-
and left-hand sides are equal:

$$\theta_{uv}(s, l^*(v)) = \min_{(s', l') \in \mathcal{Y}_{uv}} \theta_{uv}(s', l') , \tag{8.57}$$

To finalize the proof note that $l^*(v)$ can always be selected such that
$l = l^*(v)$ and, therefore, (8.50) follows from (8.49). □

The condition $\mathrm{cl}(\mathrm{mi}[\theta^\phi]) = \bar{0}$ implies that there are two incident
edges $uv$ and $uv'$ such that the locally optimal label pair in $uv$ is
"incident" only to locally non-optimal ones in $uv'$. Formally, it holds
that

$$\theta_{uv}(s, l) = \min_{(s', l') \in \mathcal{Y}_{uv}} \theta_{uv}(s', l') \tag{8.58}$$

and

$$\min_{l' \in \mathcal{Y}_{v'}} \theta_{uv'}(s, l') > \min_{(s', l') \in \mathcal{Y}_{uv}} \theta_{uv'}(s', l') . \tag{8.59}$$

In this case we will call the node $u$ *locally inconsistent* with respect to $\theta$.

**Lemma 8.17.** Let $u$ be locally inconsistent with respect to $\theta$ and let
$\mathcal{D}(\theta) = \mathcal{D}(F_u(\theta))$. Then $\mathrm{mi}[F_u(\theta))] < \mathrm{mi}[\theta]$.

*Proof.* Let $l^*(v) \in \arg\min_{l \in \mathcal{Y}_v} \theta_{uv}(s, l)$ defines a minimal label pair
$(s, l^*(v))$ in the pencil $(s, u, v)$. Let also $\theta' := F_u(\theta)$. With (8.59) in
mind summing up over all $v \in \mathcal{N}_b(v)$ we obtain

$$\sum_{v \in \mathcal{N}_b(u)} \theta_{uv}(s, l^*(v)) > \sum_{v \in \mathcal{N}_b(u)} \min_{(s', l') \in \mathcal{Y}_{uv}} \theta_{uv}(s', l') . \tag{8.60}$$

Since $\mathcal{D}(\theta) = \mathcal{D}(\theta')$ and all unary costs are zero, it holds that

$$\sum_{v \in \mathcal{N}_b(u)} \min_{(s', l') \in \mathcal{Y}_{uv}} \theta_{uv}(s', l') = \sum_{v \in \mathcal{N}_b(u)} \min_{(s', l') \in \mathcal{Y}_{uv}} \theta'_{uv}(s', l') . \tag{8.61}$$

Since $F_u$ performs a reparametrization, it holds that

$$\sum_{v \in \mathcal{N}_b(u)} \theta_{uv}(s, l^*(v)) = \sum_{v \in \mathcal{N}_b(u)} \theta'_{uv}(s, l^*(v)) \,. \tag{8.62}$$

Combining this with (8.60) and (8.61) we obtain

$$\sum_{v \in \mathcal{N}_b(u)} \theta'_{uv}(l, s^*(v)) > \sum_{v \in \mathcal{N}_b(u)} \min_{\substack{l' \in \mathcal{Y}_u \\ s' \in \mathcal{Y}_v}} \theta'_{uv}(l', s') \,. \tag{8.63}$$

By virtue of (8.18) all summands in the right-hand side are equal to each other, as well as all summands in the left-hand side. This immediately implies

$$\theta'_{uv'}(l, s^*(v)) > \min_{\substack{l' \in \mathcal{Y}_u \\ s' \in \mathcal{Y}_v}} \theta'_{uv'}(l', s') \,, \tag{8.64}$$

for all nodes $v \in \mathcal{N}_b(u)$ including the one where equality (8.58) holds for $\theta$. Together with Lemma 8.16 this finalizes the proof. $\square$

Let $F$ be the mapping $\mathbb{R}^{\mathcal{I}} \to \mathbb{R}^{\mathcal{I}}$ corresponding to one iteration of the min-sum diffusion algorithm, that is, $F$ corresponds to the cost transformation (8.15)-(8.16) applied to all nodes $u \in \mathcal{V}$ in a predefined order. Property (6) in Definition 8.14 is equivalent to the following statement:

**Theorem 8.4.** $\mathrm{cl}(\mathrm{mi}[\theta]) \neq \bar{0}$ if and only if $\mathcal{D}(\theta) = \mathcal{D}(F^N(\theta))$, where $N = \sum_{uv \in \mathcal{E}} |\mathcal{Y}_{uv}|$.

*Proof.* Necessity: The condition $\mathrm{cl}(\mathrm{mi}[\theta]) = \bar{0}$ implies that there is at least one inconsistent node on each iteration. By Lemma 8.17 if $\mathcal{D}(\theta) = \mathcal{D}(F^N(\theta))$ the number of locally minimal label pairs strictly decreases on each iteration and becomes zero after at most $N$ iterations, which is impossible by definition. Therefore, it must hold $\mathcal{D}(\theta) < \mathcal{D}(F^N(\theta))$.

Sufficiency: If $\mathrm{cl}(\mathrm{mi}[\theta]) \neq \bar{0}$, then $\mathcal{D}(\theta) = \mathcal{D}(F^i(\theta))$ for any natural $i$ according to Theorem 8.1. $\square$

## 8.7 Empirical comparison of algorithms

To get a feeling for how the algorithms considered in this chapter perform and compare to each other, we run them on four different problem instances. Although in general results obtained on so few test instances do not generalize well, this is not the case here. We intentionally selected some typical problem instances to get the results, which are well correlated with those obtained in studies of much larger scale, such as [126] and [41].

**Algorithms** For comparison we selected the following six inference algorithms and their variants, where the names of the algorithms used in the description further on are emphasized with the `typewriter` font:

- `Naïve` node-wise rounding: This method simply selects the label with the lowest unary cost in each node independently:

$$y_u := \arg\min_{s \in \mathcal{Y}_u} \theta_u(s) \, . \tag{8.65}$$

- The `ICM` method as described in §8.1.1, with the initial labeling selected with the `naïve` method.

- The min-sum `diffusion` algorithm with the initial reparametrization $\phi = 0$. Two different methods were used to obtain an integer labeling from the current reparametrization:

  - `diffusion+naïve`, the labeling was obtained by applying the preprocessing (8.24) followed by the naïve rounding (8.65) applied to the reparametrized costs.

  - `diffusion+ICM`, one iteration of the `ICM` algorithm was run to obtain a labeling with the starting point being obtained as in `diffusion+naïve`.

- Finally, the `subgradient` method with the step-size selected as $\alpha^t = \frac{0.1}{1+t}$ was used to optimize the dual. The constant 0.1 was selected empirically. To obtain an integer labeling the same two techniques were used as in the case of `diffusion` with the difference that the preprocessing (8.24) was not applied before the `naïve` rounding. The corresponding algorithms are referred to as `subgradient+naïve` and `subgradient+ICM`.

**Datasets** We run all algorithms listed above on the following four problem instances taken from the OpenGM2 [41] and MRF Middlebury [126] benchmarks. Below we provide only a short description and refer to these benchmarks for a more detailed one:

- The `color segmentation` N4 dataset [41]: the instance `pfau -small` is represented by a 4-connected grid graph of the size $320 \times 240$ nodes corresponding to image pixels. The considered applied problem is to approximate an input image with an image where only 12 fixed colors are used. Therefore, the label set $\mathcal{Y}_u$ in all pixels $u \in \mathcal{V}$ is the same and contains $|\mathcal{Y}_u| = 12$ labels, one for each of the fixed colors. The unary cost $\theta_u(s)$ is determined by the difference between the actual image color in the pixel $u$ and the predefined color corresponding to the label $s \in \mathcal{Y}_u$. The pairwise costs represent the *compactness assumption* (colors in the neighboring pixels are most likely the same) discussed in §1.1 and are given by the expression:

$$\theta_{uv}(s, t) = \begin{cases} 0, & s = t \\ \alpha_{uv} > 0, & s \neq t. \end{cases} \tag{8.66}$$

- The MRF-`stereo` dataset [41, 126], the instance `tsu-gm`, is represented by a 4-connected grid graph of the size $384 \times 288$ nodes corresponding to image pixels. The considered applied problem consists in recovering a depth value in each pixel. The label set $\mathcal{Y}_u$ in all pixels $u \in \mathcal{V}$ is the same and contains $|\mathcal{Y}_u| = 12$ labels, one for each of the considered discretized depth values. The unary costs are computed based on the color difference between the corresponding pixels of the stereo pair of images, where correspondence is determined by the depth value. The pairwise costs represent the assumption that the 3D surface to be reconstructed is piecewise smooth. Therefore, the cost is proportional to the difference of the depth values within each piece and one pays a fixed cost (although dependent on the color information taken from the image, which is expressed by a constant $\alpha_{uv}$ below) if two neighboring nodes belong to different pieces (in other words, if the jump between the depth values exceeds a certain threshold $T$):

$$\theta_{uv}(s,t) = \alpha_{uv} \min\{|s-t|, T\}. \tag{8.67}$$

- The `worms` dataset [39], the considered instance
  `pha4A7L1_1213062-lowThresh-hyp.surf`
  is represented by a non-regular densely connected graph with 556
  nodes, with the number of labels varying from 4 to 47 and the
  number of edges being about 10% of all possible $556 \times 555$ node
  pairs. Most of the edges, however, represent so-called *uniqueness
  constraints*, similar to the ones in §1.3, where they were used to
  guarantee a one-to-one mapping between nodes of the graphical
  model and the graph where a Hamiltonian cycle had to be found:

$$\theta_{uv}(s,t) = \begin{cases} 0, & s \neq t \\ \infty, & s = t \end{cases}. \tag{8.68}$$

- The `protein`-folding dataset [41], the instance `2BBN` is represented
  by a fully-connected graph with 37 nodes and the number of labels
  varies from 2 to 404. The pairwise costs vary in a very large range,
  similar to the one from the `worms` dataset and dissimilar to the
  relatively small range of pairwise costs in the `stereo` and `color
  segmentation` datasets.

**Primal-dual plots** The performance of algorithms is represented by
the *primal-dual* convergence plots in Figures 8.8 and 8.9. These plots
show both the obtained dual value on each iteration of the `diffusion`
and `subgradient` methods as well as energies of the reconstructed
labelings. In addition to the six methods listed above we plot also the
value of the optimal solution of the non-relaxed MAP-inference problem
obtained with a modern exact inference technique [33].

It can be seen that the `naïve` method performs worst, which is not
surprising taking into account that it completely ignores the pairwise
costs. In contrast, applied after several iterations of the `diffusion` or
`subgradient`, it delivers much better results, as these dual methods
reparametrize the unary and pairwise costs to make them possibly
consistent.

Still, the naïve rounding does not perform well even when coupled with these dual algorithms for the `protein` and `worms` datasets. The reason for that is the existence of the very large (infinite) pairwise costs. If there is at least one pair of labels (selected naïvely by the dual rounding in each node independently) corresponding to such a very large (infinite) pairwise cost, the resulting energy becomes very large (infinite) as well.

A simple and efficient way to cope with the described problem is to run the `ICM` algorithm after obtaining the naïve rounding. Even one iteration of this method allows us to obtain a labeling without very large (infinite) pairwise costs. Such `ICM` "postprocessing" improves the approximate primal solution also for the `color segmentation` and `stereo` datasets, where all pairwise costs are quite modest. Therefore, the obtained improvement is less pronounced in this case.

When applied as a standalone method, the `ICM` algorithm converges after a couple of iterations, but typically obtains only approximate solutions, which are very far from the global optimum.

Note also, that the subgradient method is always outperformed by the min-sum diffusion after at most several dozens of iterations. We attribute this fact mainly to the sparsity of the subgradient updates, as pointed out in Remark 8.3. Contrary to the subgradient method, the min-sum diffusion algorithm updates *all* dual variables on each of its iterations.

## 8.8   Bibliography and further reading

The ICM algorithm [9] is among the first methods addressing the MAP-inference problem. There are several publications suggesting different variants of the block-ICM method with tree-structured subproblems. Among the recent ones are [44] and [19, 127], where its parallelization properties were used to achieve faster convergence.

Another efficient block-ICM scheme for optimizing over small, but arbitrary structured subproblems is proposed in [4] and known as *Lazy Flipper*.

The dual subgradient algorithm for the MAP-inference problem in the described form was proposed in [101]. Independently, the subgradient method was evaluated for a slightly different dual problem (considered in Chapter 9) in [119, 55, 56].

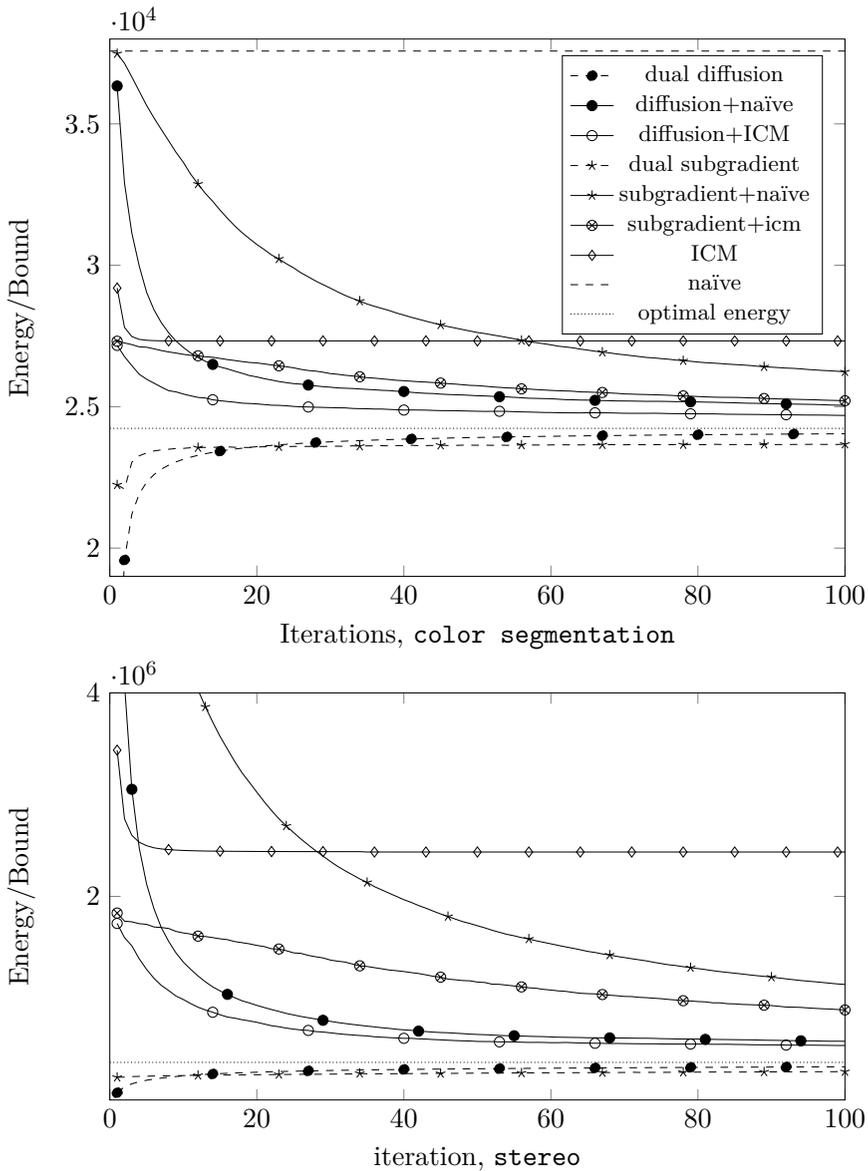

**Figure 8.8:** Primal-dual plots for **(top)** `color segmentation` and **(bottom)** `stereo` problem instances. The legend for the bottom plot is the same as for the top one. The `naïve` method is not shown in the bottom plot, since it lies beyond its field of view. Note how a small difference in the dual domain between `diffusion` and `subgradient` turns into a much larger one between energies of the primal solutions obtained by the corresponding `diffusion+naïve` and `subgradient+naïve` methods.

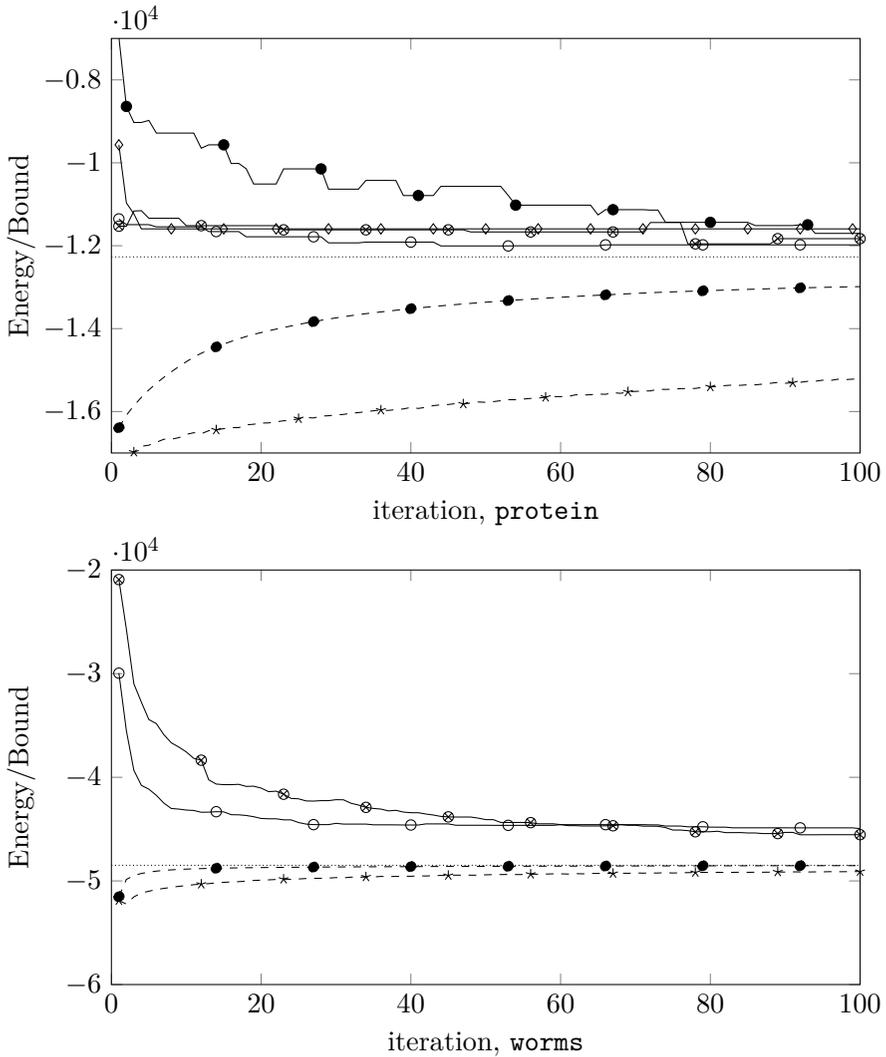

**Figure 8.9:** Primal-dual plots for **(top)** `protein` and **(bottom)** `worms` problem instances. See also the legend in the top plot in Figure 8.8. The `naïve`, `diffusion+naïve` and `subgradient+naïve` methods, as well as `ICM` in the bottom plot, are not shown, since they lie beyond the field of view of the plot.

Although the min-sum diffusion algorithm is known at least since the 70's, its detailed convergence analysis was given only recently in [102]. The anisotropic diffusion was introduced in the work [51], along with a unified analysis and evaluation of different dual coordinate ascent algorithms including the CMP (Convex Message Passing) algorithm [35] and MPLP (Message Passing for Linear Programming) [28]. The latter was recently significantly improved in [128].

A min-sum diffusion-like algorithm for general convex piecewise linear functions was suggested in [139]. A unified analysis of the convergence properties of block-coordinate descent algorithms for non-smooth convex functions was recently given in [141].

There is an alternative method that guarantees the monotonic improvement of the dual objective as well as the convergence to the node-edge agreement. Its representative is the *Augmenting DAG (directed acyclic graph)* algorithm, proposed by [59] and recently reviewed in [140]. A generalization of this algorithm to higher order models was suggested in [21]. This algorithm can be seen as an instance of the combinatorial primal-dual method [79]. The latter covers in particular such famous techniques as the algorithm of Ford-Fulkerson for max-flow and the Hungarian method for the linear assignment problem.

# 9

# Lagrange (Dual) Decomposition

In the previous chapters we learned about two different types of optimization methods with (different) optimality guarantees. On one side we have dynamic programming and on the other side the dual subgradient- and min-sum diffusion methods.

Dynamic programming:

- is a finite step algorithm and has linear complexity with respect to the problem size;

- solves the energy minimization problem exactly;

- is limited to acyclic graphs.

The subgradient and min-sum diffusion methods:

- are iterative algorithms with no finite-step convergence, even for acyclic graphs (see Example 8.10);

- address the dual LP and, therefore, provide a lower bound for the optimal energy value;

- are applicable to graphical models with arbitrary structure, but, in general solve even acyclic problems only in the limit.

Our goal is to construct a method, which combines the advantages of dynamic programming as well as the dual iterative methods. This method should

- be applicable to graphical models with an arbitrary structure;

- be as efficient as dynamic programming for acyclic graphs;

- provide a lower bound and guarantee convergence either to the dual optimum or at least to node-edge agreement for general graphs.

In §8.3.2 we have seen that anisotropic diffusion updates with specially constructed weights are equivalent to dynamic programming. Indeed, in Chapter 10 we will construct an algorithm based on anisotropic diffusion which satisfies all the properties above. Although one could prove these properties directly, such a proof would lack interpretation power. In other words, with a direct proof it would be difficult to explain *why* the proposed algorithm is very efficient. Additionally, it would be difficult to explain how the correct weights $\omega_i$ for the anisotropic diffusion steps must be selected to result in an efficient algorithm.

Therefore, we will make a detour and consider another convex relaxation technique known as *Lagrange (or dual) decomposition*. This technique is widely used for constructing convex relaxations of combinatorial problems. An advantage of this detour is not only another interpretation of anisotropic diffusion, but also another dimension for efficiency analysis applicable to most of the existing dual iterative inference algorithms for graphical models.

## 9.1 Lagrange decomposition in a nutshell

Let $f\colon X \to \mathbb{R}$ be the function which has to be minimized. Assume that the minimization is difficult, however, $f$ is representable as a sum of two other functions on the same domain, i.e. $f(x) = f^1(x) + f^2(x)$. Furthermore, assume that the optimization of $f^1(x) + \langle \lambda, x \rangle$ and $f^2(x) + \langle \lambda, x \rangle$ independently is easy for any vector $\lambda$. *Lagrange decomposition* is the method for constructing a relaxation of the problem $\min_{x \in X} f(x)$

as follows:

$$
\begin{aligned}
\min_{x \in X} f(x) = \min_{x \in X}(f^1(x) + f^2(x)) &= \min_{\substack{x^1, x^2 \in X \\ x^1 = x^2}} (f^1(x^1) + f^2(x^2)) \\
&\geq \min_{x^1, x^2 \in X} \left( f^1(x^1) + f^2(x^2) + \left\langle \lambda, (x^1 - x^2) \right\rangle \right) \\
&= \min_{x^1 \in X} \left( f^1(x^1) + \left\langle \lambda, x^1 \right\rangle \right) + \min_{x^2 \in X} \left( f^2(x^2) - \left\langle \lambda, x^2 \right\rangle \right) .
\end{aligned}
\tag{9.1}
$$

In a nutshell, the Lagrange decomposition is a Lagrange relaxation, where the equality constraint $x^1 = x^2$ is dualized, as is done in the second line of (9.1). This fact also proves the respective inequality in (9.1), see Chapter 5.3.1. Alternatively, this inequality can be proven directly by noting that

$$
\min_{\substack{x^1, x^2 \in X \\ x^1 = x^2}} (f^1(x^1) + f^2(x^2)) = \min_{x^1, x^2 \in X} (f^1(x^1) + f^2(x^2) + \iota_{\{x^1 = x^2\}}(x^1, x^2))
\tag{9.2}
$$

with $\iota_{\{x^1 = x^2\}}$ being a function such that $\iota_{\{x^1 = x^2\}} = 0$ if $x^1 = x^2$ and $\infty$ otherwise. Since $\iota_{\{x^1 = x^2\}}(x^1, x^2) \geq \langle \lambda, (x^1 - x^2) \rangle$, the inequality (9.1) holds.

The problem of finding the best lower bound of the form (9.1) is the (Lagrange) *decomposition-based* dual:

$$
\max_{\lambda} \left( \min_{x^1 \in X} \left( f^1(x^1) + \left\langle \lambda, x^1 \right\rangle \right) + \min_{x^2 \in X} \left( f^2(x^2) - \left\langle \lambda, x^2 \right\rangle \right) \right) .
\tag{9.3}
$$

As a Lagrange dual it is a concave (see Proposition 5.38), and (in general, see Example 5.39) non-smooth problem.

## 9.2 Lagrange decomposition for grid graphs

In this section we will consider applying of the Lagrange decomposition technique to energy minimization in graphical models. We start with the special case of grid graphs to avoid cluttering the main idea with technical details.

Let $\mathcal{G}$ be a grid graph with $m^r$ rows and $m^c$ columns and let $E(y; \theta)$ be the energy of a labeling $y \in \mathcal{Y}_\mathcal{V}$. Let us split the graph $\mathcal{G}$ into chain subgraphs $\mathcal{G}^r$ and $\mathcal{G}^c$, corresponding to the rows and columns of the graph $\mathcal{G}$ respectively, see Figure 9.1. In other words, if $(i,j)$, $i = 1,\dots,m^r$, $j = 1,\dots,m^c$, stands for the $j$-th node in the $i$-th row, then $\mathcal{G}^r = (\mathcal{V}, \mathcal{E}^r)$, with $\mathcal{E}^r = \{\{(i,j),(i,j+1)\}\colon i = 1,\dots,m^r,\ j = 1,\dots,m^c - 1\}$ and, similarly, $\mathcal{G}^c = (\mathcal{V}, \mathcal{E}^c)$, with $\mathcal{E}^c = \{\{(i,j),(i+1,j)\}\colon i = 1,\dots,m^r - 1,\ j = 1,\dots,m^c\}$. We will call the subgraphs $\mathcal{G}^r$ and $\mathcal{G}^c$ *slave* in contrast to the *master* graph $\mathcal{G}$.

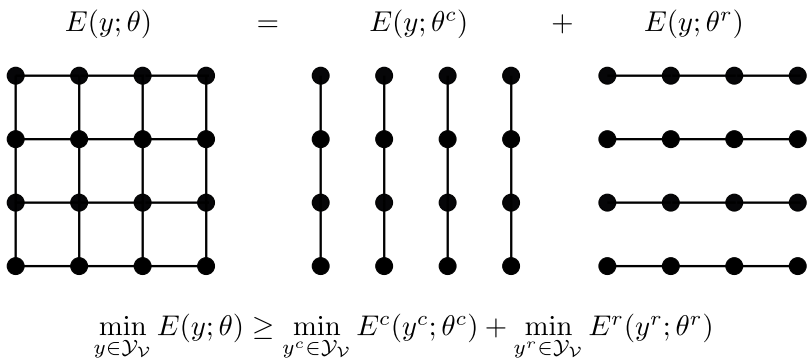

**Figure 9.1:** An illustration of the Lagrange decomposition-based dual construction for grid graphs. The master graph is split into two acyclic slave subgraphs corresponding to the rows and columns of the graph. The total energy of the corresponding slave subproblems is equal to the energy of the master problem. This is achieved by splitting the unary costs between slave subproblems. The decomposition-based dual problem consists in finding such a splitting that maximizes the corresponding lower bound.

The total energy of all columns will be denoted as $E^c$ and the energy of all rows as $E^r$. We will split costs between these two subgraphs in a way that the sum of energies of the slave subgraphs is equal to the energy of the master graph for any labeling $y \in \mathcal{Y}_\mathcal{V}$:

$$E(y; \theta) = E^c(y; \theta^c) + E^r(y; \theta^r). \tag{9.4}$$

Energy minimizations for columns ($E^c$) and rows ($E^r$) are called the *slave (sub)problems*, in contrast to the energy minimization over the whole graph, called the *master problem*.

Let $\mathcal{I}^c$ and $\mathcal{I}^r$ denote the number of primal variables in the ILP representations of the corresponding subproblems, similar to $\mathcal{I}$ defining it for the *master* problem.

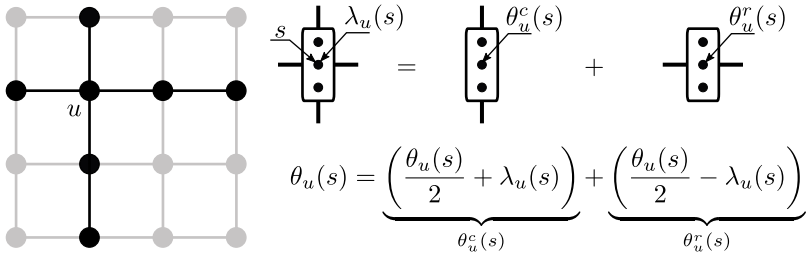

**Figure 9.2:** An illustration of the Lagrange decomposition-based dual construction for grid graphs. The master graph is split into two acyclic slave subgraphs corresponding to the rows and columns of the graph. The unary costs split between slave subproblems. This can be achieved by considering an unconstrained dual variable $\lambda_u(s)$ for each node $u$ and each label $s$. See also comments to (9.5).

To ensure fulfillment of equality (9.4) we will split the unary costs between $\mathcal{G}^r$ and $\mathcal{G}^c$ and assign the pairwise ones to the corresponding subgraphs without changing them. In other words, the costs $\theta^c \in \mathbb{R}^{\mathcal{I}^c}$ and $\theta^r \in \mathbb{R}^{\mathcal{I}^r}$ are defined such that $\theta_u^c(s) + \theta_u^r(s) = \theta_u(s)$ for any node $u \in \mathcal{V}$ and label $s \in \mathcal{Y}_u$. This is equivalent to the following representation (see Figure 9.2):

$$\theta_u^c(s) = \theta_u^c[\lambda](s) := \frac{\theta_u(s)}{2} + \lambda_u(s)\,, \quad \theta_u^r(s) = \theta_u^r[\lambda](s) := \frac{\theta_u(s)}{2} - \lambda_u(s)\,, \tag{9.5}$$

where the values $\lambda := \{\lambda_u(s) \colon u \in \mathcal{V},\ s \in \mathcal{Y}_u\} \in \mathbb{R}^{\sum_{u \in \mathcal{V}} |\mathcal{Y}_u|}$ can be arbitrary. Therefore, further in the formulas we will omit the vector space $\mathbb{R}^{\sum_{u \in \mathcal{V}} |\mathcal{Y}_u|}$ they belong to for the sake of notation. We will write $\theta[\lambda]$ where it is important to emphasize that costs $\theta^c$ and $\theta^r$ depend on $\lambda$.

Since $\mathcal{G}^r$ and $\mathcal{G}^c$ have no common edges, there is no need to change pairwise costs to enforce the equality (9.4). Therefore, they read:

$$\theta_{uv}^c = \theta_{uv}^c[\lambda] := \theta_{uv}\,,\ uv \in \mathcal{E}^c\,, \qquad \theta_{uv}^r = \theta_{uv}^r[\lambda] := \theta_{uv}\,,\ uv \in \mathcal{E}^r\,. \tag{9.6}$$

It follows from (9.4) that

$$\min_{y \in \mathcal{Y}_{\mathcal{V}}} E(y; \theta) \geq \underbrace{\min_{y^c \in \mathcal{Y}_{\mathcal{V}}} E^c(y^c; \theta^c[\lambda]) + \min_{y^r \in \mathcal{Y}_{\mathcal{V}}} E^r(y^r; \theta^r[\lambda])}_{\mathcal{U}(\lambda)} \tag{9.7}$$

for any $\lambda$.

Note also that the right-hand-side of (9.7), denoted as $\mathcal{U}(\lambda)$, can be easily evaluated by dynamic programming due to the fact that both $\mathcal{G}^c$ and $\mathcal{G}^r$ are acyclic. Inequality (9.7) implies that for any $\lambda$ the value $\mathcal{U}(\lambda)$ is a lower bound for the optimal value of the energy $E$.

Since $\mathcal{U}(\lambda)$ and $\min_{y \in \mathcal{Y}_\mathcal{V}} E(y; \theta)$ can be rewritten as

$$\mathcal{U}(\lambda) = \min_{y^c, y^r \in \mathcal{Y}_\mathcal{V}} \left( E^c(y^c; \theta^c[\lambda]) + E^r(y^r; \theta^r[\lambda]) \right) , \qquad (9.8)$$

and

$$\min_{y \in \mathcal{Y}_\mathcal{V}} E(y; \theta) = \min_{\substack{y^c, y^r \in \mathcal{Y}_\mathcal{V} \\ y^c = y^r}} \left( E^c(y^c; \theta^c[\lambda]) + E^r(y^r; \theta^r[\lambda]) \right) , \qquad (9.9)$$

respectively, the former can be seen as a relaxation of the latter one in the sense of Definition 3.6. Let us show that it is a Lagrange relaxation, and, therefore, the problem $\max_\lambda \mathcal{U}(\lambda)$ can be seen as a decomposition-based dual of the MAP-inference problem, as in (9.3).

To this end, we will introduce two cost vectors $\hat{\theta}^r$ and $\hat{\theta}^c$ for the master graph $\mathcal{G}$. Coordinates of these vectors represent costs of the slave subproblems associated with the subgraphs $\mathcal{G}^r$ and $\mathcal{G}^c$ and otherwise are equal to zero:

$$\forall uv \in \mathcal{E}: \quad \hat{\theta}^c_{uv} = \begin{cases} \theta_{uv}, & uv \in \mathcal{E}^c, \\ \bar{0}, & uv \in \mathcal{E}^r \end{cases} , \quad \hat{\theta}^r_{uv} = \begin{cases} \theta_{uv}, & uv \in \mathcal{E}^r \\ \bar{0}, & uv \in \mathcal{E}^c \end{cases} . \quad (9.10)$$

We will also assume that the unary costs are equally split between column- and row-wise subproblems:

$$\forall u \in \mathcal{V}: \quad \hat{\theta}^c_u = \hat{\theta}^r_u = \frac{1}{2}\theta_u . \qquad (9.11)$$

Together with (9.10) this implies $\hat{\theta}^c + \hat{\theta}^r = \theta$.

Similarly, we will introduce a vector $\hat{\lambda} \in \mathbb{R}^\mathcal{I}$ with zero coordinates corresponding to the edges of the master graph $\mathcal{G}$:

$$\forall uv \in \mathcal{E}: \quad \hat{\lambda}_{uv} = \bar{0} . \qquad (9.12)$$

In what follows we will also distinguish between notations $\delta_\mathcal{G}$ for the master graph and $\delta_{\mathcal{G}^r}$ and $\delta_{\mathcal{G}^c}$ for the slave subgraphs. To give a meaning to the inner product in $(9.3)$ we will write $E$ as a linear function of the cost vector $\theta$:

$$
\begin{aligned}
\min_{y \in \mathcal{Y}_\mathcal{V}} E(y; \theta) &= \min_{y \in \mathcal{Y}_\mathcal{V}} \langle \theta, \delta_\mathcal{G}(y) \rangle \\
&= \min_{y \in \mathcal{Y}_\mathcal{V}} \left( \left\langle \hat{\theta}^c, \delta_\mathcal{G}(y) \right\rangle + \left\langle \hat{\theta}^r, \delta_\mathcal{G}(y) \right\rangle \right) \\
&= \min_{\substack{y^c, y^r \in \mathcal{Y}_\mathcal{V} \\ y^c = y^r}} \left( \left\langle \hat{\theta}^c, \delta_\mathcal{G}(y^c) \right\rangle + \left\langle \hat{\theta}^r, \delta_\mathcal{G}(y^r) \right\rangle \right) \\
&= \min_{\substack{\mu^c, \mu^r \in \delta_\mathcal{G}(\mathcal{Y}_\mathcal{V}) \\ \mu^c_\mathcal{V} = \mu^r_\mathcal{V}}} \left( \left\langle \hat{\theta}^c, \mu^c \right\rangle + \left\langle \hat{\theta}^r, \mu^r \right\rangle \right) \\
&\geq \min_{\mu^c, \mu^r \in \delta_\mathcal{G}(\mathcal{Y}_\mathcal{V})} \left( \left\langle \hat{\theta}^c, \mu^c \right\rangle + \left\langle \hat{\theta}^r, \mu^r \right\rangle + \left\langle \hat{\lambda}_\mathcal{V}, \mu^c_\mathcal{V} - \mu^r_\mathcal{V} \right\rangle \right) \\
&\stackrel{(9.12)}{=} \min_{\mu^c \in \delta_\mathcal{G}(\mathcal{Y}_\mathcal{V})} \left\langle \hat{\theta}^c + \hat{\lambda}, \mu^c \right\rangle + \min_{\mu^r \in \delta_\mathcal{G}(\mathcal{Y}_\mathcal{V})} \left\langle \hat{\theta}^r + \hat{\lambda}, \mu^r \right\rangle \\
&= \min_{\mu^c \in \delta_{\mathcal{G}^c}(\mathcal{Y}_\mathcal{V})} \left\langle \theta^c[\lambda], \mu^c \right\rangle + \min_{\mu^r \in \delta_{\mathcal{G}^r}(\mathcal{Y}_\mathcal{V})} \left\langle \hat{\theta}^r[\lambda], \mu^r \right\rangle \\
&= \min_{y^c \in \mathcal{Y}_\mathcal{V}} \left\langle \theta^c[\lambda], \delta_{\mathcal{G}^c}(y^c) \right\rangle + \min_{y^r \in \mathcal{Y}_\mathcal{V}} \left\langle \theta^r[\lambda], \delta_{\mathcal{G}^r}(y^r) \right\rangle \\
&= \min_{y^c \in \mathcal{Y}_\mathcal{V}} E^c(y^c; \theta^c[\lambda]) + \min_{y^r \in \mathcal{Y}_\mathcal{V}} E^r(y^r; \theta^r[\lambda]) = \mathcal{U}(\lambda). \quad (9.13)
\end{aligned}
$$

Note that the derivation $(9.13)$ is a specialization of $(9.1)$ and therefore $\mathcal{U}(\lambda)$ is a *Lagrange decomposition-based dual* as the objective in $(9.3)$. Since the slave subgraphs $\mathcal{G}^c$ and $\mathcal{G}^r$ are acyclic, we will also call it *tree-reweighted-* or *tree-decomposition-based* dual to distinguish it from the Lagrange dual $D(\phi)$ introduced in Chapter 6 (see $(6.9)$).

Note that according to $(9.5)$ the costs $\theta^c[\lambda]$ and $\theta^r[\lambda]$ linearly depend on $\lambda$. Therefore, the function $\mathcal{U}$ is piecewise linear and concave, but non-smooth. As we learned in Chapter 5 (see Example 5.39), these are the general properties of the Lagrange duals for discrete minimization problems when linear constraints are dualized.

## 9.3 General graph decomposition scheme

### 9.3.1 Dual problem

In this section we generalize the results of the previous section for general graphs and their general decomposition. Let $(\mathcal{G}, \mathcal{Y}_\mathcal{V}, \theta)$ be a graphical model, which will be referred to as a *master* model. Let $\mathcal{T}$ be a finite set and $\mathcal{G}^\tau$, $\tau \in \mathcal{T}$, be a collection of subgraphs of a graph $\mathcal{G}$ such that each edge and each node belongs to at least one subgraph of the collection. We will call $\mathcal{G}$ *master* graph and $\mathcal{G}^\tau$, $\tau \in \mathcal{T}$, its *slave* subgraphs. Several exemplary decompositions are shown in Figure 9.3.

For each slave subgraph $\mathcal{G}^\tau = (\mathcal{V}^\tau, \mathcal{E}^\tau)$, $\tau \in \mathcal{T}$, we define the set

$$\mathcal{I}^\tau = \{(w,s) | w \in \mathcal{V}^\tau \cup \mathcal{E}^\tau, \ s \in \mathcal{Y}_w\} \tag{9.14}$$

of index variables in the ILP representation of the MAP-inference problem associated with the subgraph $\mathcal{G}^\tau$. The set of labelings of the subgraph $\mathcal{G}^\tau$ will be denoted as $\mathcal{Y}_{\mathcal{V}^\tau} := \prod_{u \in \mathcal{V}^\tau} \mathcal{Y}_u$. We also introduce the set $\mathcal{T}_w := \{\tau \in \mathcal{T} \colon w \in \mathcal{V}^\tau \cup \mathcal{E}^\tau\}$ as the set of all slave subgraphs containing the node or edge $w \in \mathcal{V} \cup \mathcal{E}$.

With each slave sub-graph $\mathcal{G}^\tau$ we associate the set $\Theta^\tau := \mathbb{R}^{\mathcal{I}^\tau}$ of its possible cost vectors. We also define the set

$$\Theta^{\mathcal{T}} = \left\{ \theta^{\mathcal{T}} \in \prod_{\tau \in \mathcal{T}} \Theta^\tau | \sum_{\tau \in \mathcal{T}_w} \theta^\tau_w = \theta_w, \ w \in \mathcal{V} \cup \mathcal{E} \right\} \tag{9.15}$$

of all *feasible* costs for all subgraphs, i,e. costs which sum up to the cost of the master model, see Figure 9.4 for illustration. Let also $\mathcal{M}^\tau := \mathrm{conv}(\delta_{\mathcal{G}^\tau}(\mathcal{Y}_{\mathcal{V}^\tau}))$ be the marginal polytope corresponding to the slave graphical model $(\mathcal{G}^\tau, \theta^\tau)$.

Then for any $\theta^{\mathcal{T}} \in \Theta^{\mathcal{T}}$ the following sequence of (in)equalities holds:

$$\min_{y \in \mathcal{Y}_{\mathcal{V}}} E(y; \theta) = \min_{y \in \mathcal{Y}_{\mathcal{V}}} \langle \theta, \delta_{\mathcal{G}}(y) \rangle = \min_{y \in \mathcal{Y}_{\mathcal{V}}} \sum_{w \in \mathcal{V} \cup \mathcal{E}} \langle \theta_w, \delta_{\mathcal{G}}(y)_w \rangle$$

$$= \min_{y \in \mathcal{Y}_{\mathcal{V}}} \sum_{w \in \mathcal{V} \cup \mathcal{E}} \left\langle \sum_{\tau \in \mathcal{T}_w} \theta^\tau_w, \delta_{\mathcal{G}}(y)_w \right\rangle = \min_{y \in \mathcal{Y}_{\mathcal{V}}} \sum_{w \in \mathcal{V} \cup \mathcal{E}} \sum_{\tau \in \mathcal{T}_w} \langle \theta^\tau_w, \delta_{\mathcal{G}}(y)_w \rangle$$

$$= \min_{y \in \mathcal{Y}_{\mathcal{V}}} \sum_{\tau \in \mathcal{T}} \sum_{w \in \mathcal{V}^\tau \cup \mathcal{E}^\tau} \langle \theta^\tau_w, \delta_{\mathcal{G}}(y)_w \rangle = \min_{y \in \mathcal{Y}_{\mathcal{V}}} \sum_{\tau \in \mathcal{T}} \langle \theta^\tau, \delta_{\mathcal{G}^\tau}(y_{\mathcal{V}^\tau}) \rangle$$

$$\geq \sum_{\tau \in \mathcal{T}} \min_{y^\tau \in \mathcal{Y}_{\mathcal{V}^\tau}} \langle \theta^\tau, \delta_{\mathcal{G}^\tau}(y^\tau) \rangle = \sum_{\tau \in \mathcal{T}} \min_{\mu^\tau \in \mathcal{M}^\tau} \langle \theta^\tau, \mu^\tau \rangle =: \mathcal{U}(\theta^{\mathcal{T}}) . \tag{9.16}$$

We will show below (see Remark 9.1 below) that the function $\mathcal{U}$ is the Lagrange decomposition-based dual. It is a direct generalization of the Lagrange decomposition-based dual constructed in (9.7) for

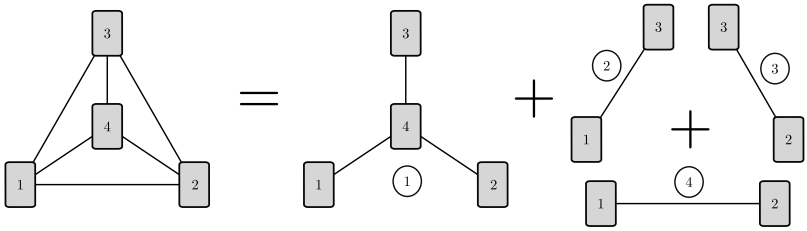

(a) An acyclic decomposition. The local sets of subproblem indexes are
$\mathcal{T}_1 = \{1, 2, 4\}$, $\mathcal{T}_2 = \{1, 3, 4\}$, $\mathcal{T}_3 = \{1, 2, 3\}$, $\mathcal{T}_4 = \{1\}$; $\mathcal{T}_{14} = \mathcal{T}_{24} = \mathcal{T}_{34} = \{1\}$,
$\mathcal{T}_{12} = \{4\}$, $\mathcal{T}_{13} = \{2\}$, $\mathcal{T}_{23} = \{3\}$.

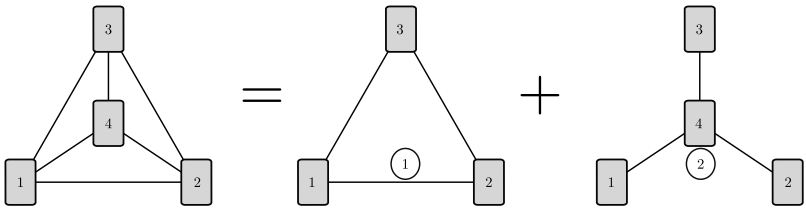

(b) A decomposition consisting of a cyclic and an acyclic graphs. The local
sets of subproblem indexes are $\mathcal{T}_1 = \mathcal{T}_2 = \mathcal{T}_3 = \{1, 2\}$, $\mathcal{T}_4 = \{1\}$; $\mathcal{T}_{12} = \mathcal{T}_{23} =$
$\mathcal{T}_{31} = \{1\}$, $\mathcal{T}_{14} = \mathcal{T}_{24} = \mathcal{T}_{34} = \{2\}$.

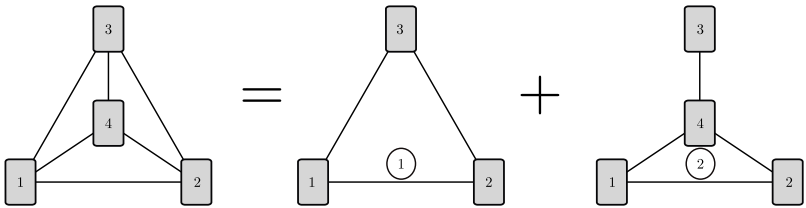

(c) A decomposition consisting of two cyclic graphs. The local sets of sub-
problem indexes are $\mathcal{T}_1 = \mathcal{T}_2 = \mathcal{T}_3 = \{1, 2\}$, $\mathcal{T}_4 = \{1\}$; $\mathcal{T}_{12} = \{1, 2\}$,
$\mathcal{T}_{12} = \mathcal{T}_{23} = \mathcal{T}_{31} = \{1\}$, $\mathcal{T}_{14} = \mathcal{T}_{24} = \mathcal{T}_{34} = \{2\}$.

**Figure 9.3:** Examples of different decompositions. Numbers inside graph nodes
denote their indexes in the set of nodes $\mathcal{V}$ of the master graph. Numbers in circles
are indexes $\tau \in \mathcal{T}$ of the slave subgraphs.

grid graphs and their decomposition into columns and rows. Likewise,
the subproblems $\min_{y^\tau \in \mathcal{Y}_{\mathcal{V}^\tau}} \langle \theta^\tau, \delta_{\mathcal{G}^\tau}(y^\tau) \rangle$ are called *slave* subproblems
in relation to the *master* problem $\min_{y \in \mathcal{Y}_\mathcal{V}} \langle \theta, \delta_{\mathcal{G}}(y) \rangle$. Note also that
evaluating $\mathcal{U}(\theta^\mathcal{T})$ is as difficult as minimizing the energy for each slave
model. For example if all graphs $\mathcal{G}^\tau$ are acyclic, the dual $\mathcal{U}$ can be evaluated by applying dynamic programming to each slave subproblem independently.

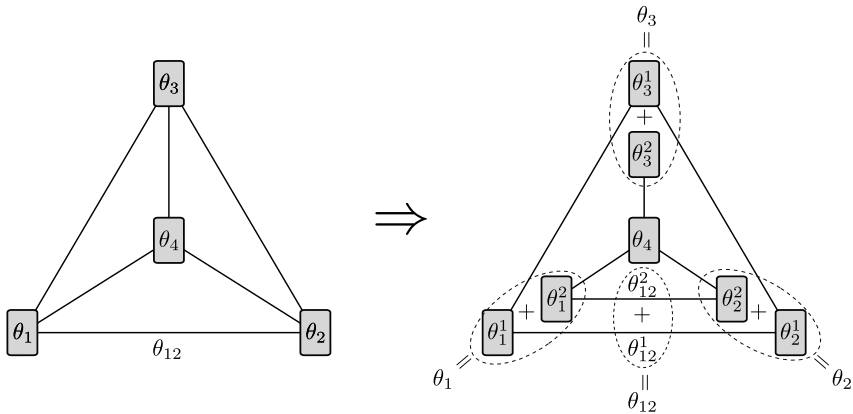

**Figure 9.4:** Illustration of cost splitting for a general Lagrange decomposition. The vectors of unary costs $\theta_i$ are sitting inside the corresponding graph nodes. A single pairwise cost vector $\theta_{12}$ is shown as it is the only one that is shared by several slave subproblems. The dashed ellipses contain nodes that are shared by several slave subproblems. The total cost inside each ellipse is equal to the corresponding cost of the master problem.

The feasible set $\Theta^{\mathcal{T}}$ is obviously convex, the dual objective $\mathcal{U}$ is concave as a sum of minima of linear functions. Therefore, finding its maximal value, which constitutes the best lower bound of the form (9.16) is a concave (non-smooth) maximization problem:

$$\max_{\theta^{\mathcal{T}} \in \Theta^{\mathcal{T}}} \mathcal{U}(\theta^{\mathcal{T}}) . \tag{9.17}$$

**Unconstrained form of decomposition-based dual**  The constrained problem (9.17) can be transformed into an unconstrained one having the form specified in (9.3). To this end let $\lambda_w^\tau(s) \in \mathbb{R}$, $w \in \mathcal{V} \cup \mathcal{E}$, $s \in \mathcal{Y}_w$, be arbitrary numbers. For any $w \in \mathcal{V} \cup \mathcal{E}$ let $\tau_w$ be some fixed index in $\mathcal{T}_w$. Consider the following representation of $\theta_w^\tau(s)$:

$$\theta_w^\tau[\lambda](s) := \begin{cases} \frac{\theta_w(s)}{|\mathcal{T}_w|} + \lambda_w^\tau(s), & \tau \neq \tau_w \\ \frac{\theta_w(s)}{|\mathcal{T}_w|} - \sum_{\sigma \in \mathcal{T}_w \setminus \{\tau_w\}} \lambda_w^\sigma(s), & \tau = \tau_w . \end{cases} \tag{9.18}$$

Note that $\theta^\tau[\lambda]$ defined as in (9.18) automatically satisfies the constraint $\sum_{\tau \in \mathcal{T}_w} \theta_w^\tau[\lambda] = \theta_w$ and therefore $\theta^{\mathcal{T}} \in \Theta^{\mathcal{T}}$. On the other hand, any element of $\Theta^{\mathcal{T}}$ can be represented in the form (9.18), since numbers $\lambda_w^\tau(s)$ are arbitrary. For the sake of notation we will use $\mathcal{J}^{\mathcal{T}}$ to denote the number of dimensions of the vector $\lambda$, that is $\mathcal{J}^{\mathcal{T}} = \sum_{w \in \mathcal{V} \cup \mathcal{E}}(|\mathcal{T}_w| - 1)|\mathcal{Y}_w|$, and, therefore, $\lambda \in \mathbb{R}^{\mathcal{J}^{\mathcal{T}}}$.

**Remark 9.1.** We can use the Lagrange decomposition method as in (9.1) to obtain the decomposition-based dual problem (9.17) with an explicit representation of the costs as in (9.18). The derivation is similar to the one for grid graphs covered by (9.13). To this end let $\theta_w^\tau[\bar{0}] \in \Theta^\tau$, $\tau \in \mathcal{T}$, be the costs defined by (9.18) with $\lambda = \bar{0}$. For each $\tau \in \mathcal{T}$ consider the cost vector of the master graph $\hat{\theta}^\tau \in \mathbb{R}^{\mathcal{I}}$, with the coordinates assigned zero values everywhere except those corresponding to nodes and edges of $\mathcal{G}^\tau$:

$$\forall w \in \mathcal{V} \cup \mathcal{E}: \quad \hat{\theta}_w^\tau = \begin{cases} \theta_w^\tau[\bar{0}], & w \in \mathcal{V}^\tau \cup \mathcal{E}^\tau, \\ \bar{0}, & \text{otherwise} \end{cases}. \qquad (9.19)$$

It is easy to see that the above definition implies $\sum_{\tau \in \mathcal{T}} \hat{\theta}^\tau = \theta$.

Let also $\lambda^\tau \in \mathbb{R}^{\sum_{w \in \mathcal{V} \cup \mathcal{E}} |\mathcal{Y}_w|}$, $\tau \in \mathcal{T}$, be arbitrary vectors with coordinates indexed as $\lambda_w^\tau(s) \in \mathbb{R}$, $w \in \mathcal{V} \cup \mathcal{E}$, $s \in \mathcal{Y}_w$, like in (9.18).

For each $\tau \in \mathcal{T}$ consider the vector $\hat{\lambda}^\tau \in \mathbb{R}^{\mathcal{I}}$ with the coordinates assigned zero values everywhere except those corresponding to nodes and edges of $\mathcal{G}^\tau$:

$$\forall w \in \mathcal{V} \cup \mathcal{E}, \ \tau \in \mathcal{T}: \quad \hat{\lambda}_w^\tau = \begin{cases} \lambda_w^\tau, & w \in \mathcal{V}^\tau \cup \mathcal{E}^\tau, \ \tau \in \mathcal{T}_w \backslash \{\tau_w\} \\ \bar{0}, & \text{otherwise} \end{cases}.$$
$$(9.20)$$

Let $\tau_0 \in \mathcal{T}$ be arbitrary. Applying the transformations similar to (9.13) for grid graphs, one obtains:

$$\min_{y \in \mathcal{Y}_\mathcal{V}} E(y; \theta) = \min_{y \in \mathcal{Y}_\mathcal{V}} \langle \theta, \delta_\mathcal{G}(y) \rangle$$

$$= \min_{y \in \mathcal{Y}_\mathcal{V}} \left\langle \sum_{\tau \in \mathcal{T}} \hat{\theta}^\tau, \delta_\mathcal{G}(y) \right\rangle = \min_{y \in \mathcal{Y}_\mathcal{V}} \sum_{\tau \in \mathcal{T}} \left\langle \hat{\theta}^\tau, \delta_\mathcal{G}(y) \right\rangle$$

$$
\begin{aligned}
&= \min_{(y^\tau \in \mathcal{Y}_\mathcal{V})_{\tau \in \mathcal{T}}} \sum_{\tau \in \mathcal{T}} \left\langle \hat{\theta}^\tau, \delta_\mathcal{G}(y^\tau) \right\rangle \quad \text{s.t. } y^\tau = y^{\tau_0} \text{ for all } \tau \in \mathcal{T} \backslash \{\tau_0\} \\
&= \min_{(\mu^\tau \in \delta_\mathcal{G}(\mathcal{Y}_\mathcal{V}))_{\tau \in \mathcal{T}}} \sum_{\tau \in \mathcal{T}} \left\langle \hat{\theta}^\tau, \mu^\tau \right\rangle \\
&\qquad\qquad\qquad \text{s.t. } \mu_w^\tau = \mu_w^{\tau_w} \text{ for all } w \in \mathcal{V} \cup \mathcal{E}, \tau \in \mathcal{T}_w \backslash \{\tau_w\} \\
&\geq \min_{(\mu^\tau \in \delta_\mathcal{G}(\mathcal{Y}_\mathcal{V}))_{\tau \in \mathcal{T}}} \sum_{\tau \in \mathcal{T}} \left( \left\langle \hat{\theta}^\tau, \mu^\tau \right\rangle + \left\langle \hat{\lambda}^\tau, \mu^\tau - \mu^{\tau_w} \right\rangle \right) \\
&= \sum_{\tau \in \mathcal{T} \backslash \{\tau_w\}} \min_{\mu^\tau \in \delta_\mathcal{G}(\mathcal{Y}_\mathcal{V})} \left\langle \hat{\theta}^\tau + \hat{\lambda}^\tau, \mu^\tau \right\rangle \\
&\quad + \min_{\mu^{\tau_w} \in \delta_\mathcal{G}(\mathcal{Y}_\mathcal{V})} \left\langle \hat{\theta}^{\tau_w} - \sum_{\sigma \in \mathcal{T} \backslash \{\tau_w\}} \hat{\lambda}^\sigma, \mu^{\tau_w} \right\rangle \\
&= \sum_{\tau \in \mathcal{T}} \min_{\mu^\tau \in \delta_{\mathcal{G}^\tau}(\mathcal{Y}_{\mathcal{V}^\tau})} \left\langle \theta^\tau[\lambda], \mu^\tau \right\rangle = \sum_{\tau \in \mathcal{T}} \min_{\mu^\tau \in \mathcal{M}^\tau} \left\langle \theta^\tau[\lambda], \mu^\tau \right\rangle = \mathcal{U}(\theta^\mathcal{T}) \,.
\end{aligned}
\tag{9.21}
$$

Note that when constructing the decomposition-based dual $\mathcal{U}$ the constraints $\mu_w^\tau = \mu_w^{\tau_w}$, $w \in \mathcal{V} \cup \mathcal{E}$, $\tau \in \mathcal{T}_w$, were dualized. Each such constraint corresponds to the dual vector $\lambda_w^\tau \in \mathbb{R}^{\mathcal{Y}_w}$.

### 9.3.2 Relation to Lagrange dual

Let $\theta^\phi \in \mathbb{R}^\mathcal{I}$ be a vector of reparametrized costs, and the Lagrange dual objective $D(\phi)$ be defined as in Chapter 6.1:

$$
\mathcal{D}(\phi) = \sum_{w \in \mathcal{V} \cup \mathcal{E}} \min_{s \in \mathcal{Y}_w} \theta_w^\phi(s) \,. \tag{9.22}
$$

Let us consider numbers $\rho_w^\tau \geq 0$ such that $\sum_{\tau \in \mathcal{T}} \rho_w^\tau = 1$ for all $w \in \mathcal{V} \cup \mathcal{E}$ and $\tau \in \mathcal{T}_w$. These conditions are satisfied for example for
$\rho_w^\tau = \begin{cases} \frac{1}{|\mathcal{T}_w|}, & \tau \in \mathcal{T}_w \,, \\ 0, & \text{otherwise} \end{cases}$.

The following proposition claims that any decomposition-based dual is at least as good as the Lagrange dual $\mathcal{D}(\phi)$:

**Proposition 9.2.** Let $(\mathcal{G}^\mathcal{T}, \theta^\mathcal{T})$ be a decomposition with $\theta^\tau$, $\tau \in \mathcal{T}$, defined as $\theta_w^\tau = \rho_w^\tau \theta_w^\phi$ for all $w \in \mathcal{V}^\tau \cup \mathcal{E}^\tau$. Then $\mathcal{U}(\theta^\tau) \geq \mathcal{D}(\phi)$.

*Proof.* We can rewrite the decomposition-based dual as:

$$
\begin{aligned}
\mathcal{U}(\theta^{\mathcal{T}}) &= \sum_{\tau \in \mathcal{T}} \min_{y^\tau \in \mathcal{Y}_{\mathcal{V}^\tau}} \sum_{w \in \mathcal{V}^\tau \cup \mathcal{E}^\tau} \rho_w^\tau \theta_w^\phi(y_w^\tau) \geq \sum_{\tau \in \mathcal{T}} \sum_{w \in \mathcal{V}^\tau \cup \mathcal{E}^\tau} \min_{s \in \mathcal{Y}_w^\tau} \rho_w^\tau \theta_w^\phi(s) \\
&= \sum_{w \in \mathcal{V} \cup \mathcal{E}} \sum_{\tau \in \mathcal{T}_w} \min_{s \in \mathcal{Y}_w^\tau} \rho_w^\tau \theta_w^\phi(s) = \sum_{w \in \mathcal{V} \cup \mathcal{E}} \sum_{\tau \in \mathcal{T}_w} \rho_w^\tau \min_{s \in \mathcal{Y}_w^\tau} \theta_w^\phi(s) \\
&= \sum_{w \in \mathcal{V} \cup \mathcal{E}} \left( \min_{s \in \mathcal{Y}_w} \theta_w^\phi(s) \cdot \sum_{\tau \in \mathcal{T}_w} \rho_w^\tau \right) \\
&= \sum_{w \in \mathcal{V} \cup \mathcal{E}} \min_{s \in \mathcal{Y}_w} \theta_w^\phi(s) = \mathcal{D}(\phi) \,.
\end{aligned}
\tag{9.23}
$$

$\square$

### 9.3.3   Primal problem

Tightness of different dual bounds is often easier to estimate by comparing the *corresponding* primal problems. By "corresponding" we mean the primal problem having the same optimal value as the considered dual. To find the primal for the dual (9.17) let us recall the general method of constructing primal problems for Lagrange duals, as given in §5.4.1. First, we will transform the decomposition-based dual $\mathcal{U}$ to the form of (5.50).

Switching the order of the summation and minimization operations one obtains:

$$
\begin{aligned}
\mathcal{U}(\theta^{\mathcal{T}}) &= \sum_{\tau \in \mathcal{T}} \min_{\mu^\tau \in \mathcal{M}^\tau} \langle \theta^\tau, \mu^\tau \rangle \\
&= \min_{(\mu^\tau \in \mathcal{M}^\tau)_{\tau \in \mathcal{T}}} \sum_{\tau \in \mathcal{T}} \langle \theta^\tau, \mu^\tau \rangle = \min_{\mu \in \prod_{\tau \in \mathcal{T}} \mathcal{M}^\tau} \sum_{\tau \in \mathcal{T}} \langle \theta^\tau, \mu|_{\mathcal{I}^\tau} \rangle \,, \quad (9.24)
\end{aligned}
$$

where $\prod$ stands for the Cartesian product, $\mu$ is the vector stacked from $\mu^\tau \in \mathcal{M}^\tau$, and $\mu|_{\mathcal{I}^\tau}$ is the restriction of $\mu$ to the coordinates $\mu^\tau$ corresponding to the subgraph $\tau$. For the sake of notation we will denote it by $\mu^\tau$ in the following.

Let us now recall that $\theta^\tau$ can be represented in the unconstrained form (9.18) and therefore the dual maximization problem reads

$$
\max_{\lambda \in \mathbb{R}^{\mathcal{J}^{\mathcal{T}}}} \mathcal{U}(\theta^{\mathcal{T}}[\lambda]) = \max_{\lambda \in \mathbb{R}^{\mathcal{J}^{\mathcal{T}}}} \min_{\mu \in \prod_\tau \mathcal{M}^\tau} \sum_{\tau \in \mathcal{T}} \langle \theta^\tau[\lambda], \mu^\tau \rangle \,.
\tag{9.25}
$$

Note that when constructing the dual $\mathcal{U}$, we dualized the equality constraints $\mu_w^\tau = \mu_w^{\tau_w}$, $\tau \in \mathcal{T}_w$, $w \in \mathcal{V} \cup \mathcal{E}$, see Remark 9.1 for details.

Therefore, according to (5.50) the primal problem corresponding to (9.25) is

$$\min_{\mu \in \mathrm{conv}\left(\mathcal{M}^\mathcal{T} \cap \{0,1\}^{\sum_\tau |\mathcal{I}^\tau|}\right)} \sum_{\tau \in \mathcal{T}} \langle \theta^\tau, \mu^\tau \rangle \tag{9.26}$$

$$\text{s.t. } \mu_w^\tau = \mu_w^{\tau_w}, \text{ for all } w \in \mathcal{V} \cup \mathcal{E}, \tau \in \mathcal{T}_w, \tag{9.27}$$

where $\mathcal{M}^\mathcal{T}$ stands for $\prod_{\tau \in \mathcal{T}} \mathcal{M}^\tau$. We now analyze the expression

$$\mathrm{conv}\left(\mathcal{M}^\mathcal{T} \cap \{0,1\}^{\sum_\tau |\mathcal{I}^\tau|}\right)$$

and show that it is equal to $\mathcal{M}^\mathcal{T}$. Indeed,

$$\mathrm{conv}\left(\left(\prod_{\tau \in \mathcal{T}} \mathcal{M}^\tau\right) \cap \{0,1\}^{\sum_\tau |\mathcal{I}^\tau|}\right) = \mathrm{conv}\left(\prod_{\tau \in \mathcal{T}} \mathcal{M}^\tau \cap \{0,1\}^{|\mathcal{I}^\tau|}\right)$$

$$\overset{\text{Lem. } 3.33}{=} \prod_{\tau \in \mathcal{T}} \mathrm{conv}\left(\mathcal{M}^\tau \cap \{0,1\}^{|\mathcal{I}^\tau|}\right) = \prod_{\tau \in \mathcal{T}} \mathcal{M}^\tau = \mathcal{M}^\mathcal{T}. \tag{9.28}$$

Therefore, the primal problem (9.26) can be rewritten as

$$\min_{(\mu^\tau \in \mathcal{M}^\tau)_{\tau \in \mathcal{T}}} \sum_{\tau \in \mathcal{T}} \langle \theta^\tau, \mu^\tau \rangle \quad \text{s.t. } \mu_w^\tau = \mu_w^{\tau_w}, \text{ for all } w \in \mathcal{V} \cup \mathcal{E}, \tau \in \mathcal{T}_w. \tag{9.29}$$

Since $\mathcal{M}^\tau$ are polytopes, the expression (9.29) represents a linear problem. However, recall that in general, the definition of the polytopes $\mathcal{M}^\tau$ may require an exponential number of inequalities.

### 9.3.4 Primal of acyclic decomposition, equivalence to Lagrange dual

To give an example of the primal problem (9.29) let us consider an *acyclic* decomposition, i.e. a decomposition where all slave subgraphs $\mathcal{G}^\tau$ are acyclic. An important special case of such a decomposition was the Lagrange decomposition for grid graphs considered in §9.2.

In the case of an acyclic decomposition each marginal polytope $\mathcal{M}^\tau$ is equivalent to the corresponding local polytope $\mathcal{L}^\tau$. Therefore, the primal problem (9.29) reads

$$\min_{(\mu^\tau \in \mathcal{L}^\tau)_{\tau \in \mathcal{T}}} \sum_{\tau \in \mathcal{T}} \langle \theta^\tau, \mu^\tau \rangle, \;\; \text{s.t. } \mu_w^\tau = \mu_w^{\tau_w}, \text{ for all } w \in \mathcal{V} \cup \mathcal{E}, \tau \in \mathcal{T}_w. \tag{9.30}$$

Let us rewrite the objective and all constraints of the problem (9.30) in more detail:

$$\min_{\left(\mu^\tau \in \mathbb{R}_+^{\mathcal{I}^\tau}\right)_{\tau \in \mathcal{T}}} \sum_{w \in \mathcal{V} \cup \mathcal{E}} \sum_{\tau \in \mathcal{T}_w} \langle \theta_w^\tau, \mu_w^\tau \rangle \tag{9.31}$$

$$\text{s.t. } \sum_{s \in \mathcal{Y}_w} \mu_w^\tau(s) = 1, \quad w \in \mathcal{V} \cup \mathcal{E}, \tau \in \mathcal{T}_w \tag{9.32}$$

$$\begin{cases} \sum_{s \in \mathcal{Y}_u} \mu_{uv}^\tau(s,t) = \mu_v^\tau(t), & t \in \mathcal{Y}_v \\ \sum_{t \in \mathcal{Y}_v} \mu_{uv}^\tau(s,t) = \mu_u^\tau(s), & s \in \mathcal{Y}_u \end{cases}, \;\; uv \in \mathcal{E}, \tau \in \mathcal{T}_{uv} \tag{9.33}$$

$$\mu_w^\tau = \mu_w^{\tau_w}, \quad \tau \in \mathcal{T}_w, w \in \mathcal{V} \cup \mathcal{E} \tag{9.34}$$

First, note that for each node and edge $w \in \mathcal{V} \cup \mathcal{E}$ there are $|\mathcal{T}_w|$ copies of the uniqueness (9.32) and the coupling constraints (9.33). Note also that,

(i) due to (9.34), the variables $\mu_w^\tau$ occurring in these constrains, have the same values for all $\tau \in \mathcal{T}_w$;

(ii) for each $w \in \mathcal{V} \cup \mathcal{E}$ it holds that $\sum_{\tau \in \mathcal{T}_w} \theta_w^\tau = \theta_w$.

Therefore, when substituting all $\mu_w^\tau$, $\tau \in \mathcal{T}_w$, with a single $\mu_w$ the problem (9.31)-(9.34) reduces to the local polytope relaxation over the master graph:

$$\min_{\mu \in \mathcal{L}} \sum_{w \in \mathcal{V} \cup \mathcal{E}} \sum_{\tau \in \mathcal{T}_w} \langle \theta_w^\tau, \mu_w \rangle = \min_{\mu \in \mathcal{L}} \sum_{w \in \mathcal{V} \cup \mathcal{E}} \left\langle \sum_{\tau \in \mathcal{T}_w} \theta_w^\tau, \mu_w \right\rangle$$

$$= \min_{\mu \in \mathcal{L}} \sum_{w \in \mathcal{V} \cup \mathcal{E}} \langle \theta_w, \mu_w \rangle = \min_{\mu \in \mathcal{L}} \langle \theta, \mu \rangle. \tag{9.35}$$

In other words, since each marginal polytope $\mathcal{M}^\tau$ can be described with constraints of a local polytope, the resulting relaxation coincides with the local polytope one. Note that this fact does not depend on a particular decomposition as long as it contains acyclic graphs only. This also means that all such decomposition-based duals have the same maximum. Moreover, since the primal problems for both the Lagrange and decomposition-based duals coincide, this proves the following fact:

**Proposition 9.3.** For any acyclic decomposition the optimal value of the decomposition-based dual $\mathcal{U}$ defined in (9.16) coincides with the optimal value of the Lagrange dual $D$ defined in (6.9), i.e.:

$$\max_{\lambda \in \mathbb{R}^{\mathcal{J}^{\mathcal{T}}}} \mathcal{U}(\theta^{\mathcal{T}}[\lambda]) = \max_{\phi \in \mathbb{R}^{\mathcal{J}}} D(\phi) \,. \tag{9.36}$$

Proposition 9.3 can also be illustrated by the following example:

**Example 9.4** (Complete decomposition)**.** Consider the following acyclic decomposition of the master graph $\mathcal{G}$ (see Figure 9.5):

- $\mathcal{T} = \mathcal{V} \cup \mathcal{E}$, i.e. each subgraph corresponds either to a single node or to a single edge.

- Node-subgraphs have the form $(\{u\}, \emptyset)$ for all $u \in \mathcal{V}$.

- Edge-subgraphs are defined as $(\{u, v\}, \{uv\})$ for all $uv \in \mathcal{E}$.

The corresponding decomposition-based dual coincides with the Lagrange dual defined by (6.9) and the explicit representation of the cost functions (9.18) turns into the reparametrization (6.6). This shows equivalence of this particular acyclic decomposition and the Lagrange dual.

The acyclic decomposition defined as in Example 9.4 will be called *complete*, since the master graph is *completely* decomposed into its nodes and edges.

We learned that all acyclic decompositions are as tight as the Lagrange one. Combining this knowledge with Proposition 9.2 claiming that any decomposition-based dual is at least as tight as the Lagrange one, we conclude:

**Proposition 9.5.** Any decomposition-based dual is at least as tight as an acyclic one.

Therefore, to construct tighter relaxations one has to consider non-acyclic graphs.

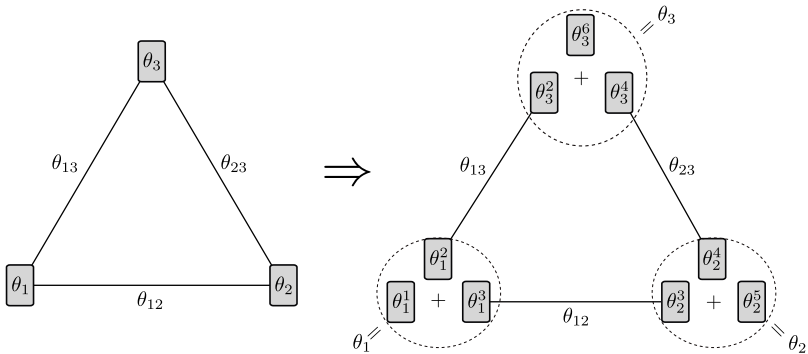

**Figure 9.5:** Example of a complete decomposition, the same notation as in Figure 9.4 is used. Note that the constraints on the costs of the slave subproblems can be automatically fulfilled, if (i) $\theta_1^2 = \phi_{1,3}$, $\theta_1^3 = \phi_{1,2}$ and $\theta_1^1 = \theta_1 - \phi_{1,2} - \phi_{1,3}$; (ii) $\theta_2^3 = \phi_{2,1}$, $\theta_2^4 = \phi_{2,3}$ and $\theta_2^5 = \theta_2 - \phi_{2,1} - \phi_{2,3}$; (iii) $\theta_3^2 = \phi_{3,1}$, $\theta_3^4 = \phi_{3,2}$ and $\theta_3^6 = \theta_3 - \phi_{3,1} - \phi_{3,2}$ with $\phi_{i,j}$ being arbitrary numbers. Note that the vector $\phi$ can be seen as a reparametrization of the master problem with $\theta_1^\phi = \theta_1^1$, $\theta_2^\phi = \theta_2^5$, $\theta_3^\phi = \theta_3^6$, $\theta_{1,2}^\phi = \theta_{1,2} + \theta_1^3 + \theta_2^3 = \theta_{1,2} + \phi_{1,2} + \phi_{2,1}$, $\theta_{1,3}^\phi = \theta_{1,3} + \theta_1^2 + \theta_3^2 = \theta_{1,3} + \phi_{1,3} + \phi_{3,1}$, and $\theta_{2,3}^\phi = \theta_{2,3} + \theta_2^4 + \theta_3^4 = \theta_{2,3} + \phi_{2,3} + \phi_{3,2}$.

## 9.3.5 Optimality conditions

Let us recall Proposition 5.45 from §5.4.3, which provides an optimality condition for Lagrange duals. Applying it to the dual problem in the form (9.25) and taking into account that $\text{conv}(\mathcal{M}^\tau \cap \{0,1\}^{|\mathcal{I}^\tau|}) = \mathcal{M}^\tau$, we obtain the following necessary and sufficient dual optimality condition:

**Proposition 9.6.** The set of costs $\theta^\mathcal{T} \in \Theta^\mathcal{T}$ is a solution to the dual problem (9.17) if and only if for all $\tau \in \mathcal{T}$ there exist

$$\mu^\tau \in \arg\min_{\mu \in \mathcal{M}^\tau} \langle \theta^\tau, \mu \rangle$$

such that $\mu_w^\tau(s) = \mu_w^\sigma(s)$ for all $\tau, \sigma \in \mathcal{T}$ and $(w, s) \in \mathcal{I}^\tau \cap \mathcal{I}^\sigma$. Moreover, for the vector $\mu^* \in \mathbb{R}^\mathcal{I}$ such that $\mu_w^* = \mu_w^\tau$ for all $w \in \mathcal{V} \cup \mathcal{E}$, it holds that $\langle \theta, \mu^* \rangle = \mathcal{U}(\theta^\mathcal{T})$, i.e. $\mu^*$ is a solution of the corresponding primal relaxed problem. If additionally $\mu^* \in \{0,1\}^\mathcal{I}$, then $\mu^*$ is a solution of the initial non-relaxed energy minimization problem.

The last two statements of Proposition 9.6 follow directly from Remark 5.46 in §5.4.3.

**Subproblem (tree) agreement**  Example 9.4 shows that the Lagrange dual is a special case of the decomposition-based one. Therefore, checking for optimality of the dual problem (9.17) is at least as difficult as checking for optimality of the Lagrange dual (see §6.2.2 for a detailed discussion). Similar to the latter case, we will develop necessary local optimality conditions for decomposition-based duals. For the Lagrange dual, the role of such conditions was taken up by the node-edge agreement.

To construct such necessary optimality conditions we will generalize Definition 6.9 of the node-edge agreement. Let a subset of optimal labelings

$$\mathbb{S}^\tau \subseteq \argmin_{y \in \mathcal{Y}_{\mathcal{V}^\tau}} \langle \theta^\tau, \delta_{\mathcal{G}}(y) \rangle \,, \tau \in \mathcal{T}\,, \tag{9.37}$$

be given for each subproblem. For any such subset we will denote by $\mathbb{S}^\tau_w = \{s \in \mathcal{Y}_w | \exists y \in \mathbb{S}^\tau : y_w = s\}$ the set of coordinates which are taken by optimal labelings in the node or edge $w \in \mathcal{V}^\tau \cup \mathcal{E}^\tau$.

**Definition 9.7.** For a given cost vector $\theta^\mathcal{T} \in \Theta^\mathcal{T}$ we will say that there is a *weak subproblem agreement*, or *subproblem agreement* for short, if there are *non-empty* subsets of optimal labelings of all subproblems $\mathbb{S}^\tau$ defined as in (9.37) such that $\mathbb{S}^\tau_w = \mathbb{S}^\sigma_w$ holds for all $\tau, \sigma \in \mathcal{T}_w$.

For acyclic decompositions a special case of the subproblem agreement is known as *tree agreement*. For the complete decomposition, tree agreement is equivalent to the node-edge agreement condition. However, the fact that node-edge agreement holds (i.e. arc-consistency closure is non-empty) can be efficiently verified using the relaxation labeling algorithm (6.24)-(6.25), whereas Definition 9.7 is non-constructive and can not be verified directly. To this end a generalization of the relaxation labeling algorithm can be constructed, which is applicable to arbitrary decompositions and allows us to check whether the subproblem agreement condition is fulfilled. This generalization consists in iterating the following operation

$$\forall w \in \mathcal{V} \cup \mathcal{E}, \tau \in \mathcal{T}: \quad \mathbb{S}^\tau := \mathbb{S}^\tau \cap \{y \in \mathcal{Y}_{\mathcal{V}^\tau} : y_w \in \bigcap_{\tau' \in \mathcal{T}_w} \mathbb{S}^{\tau'}_w\} \tag{9.38}$$

until the sets $\mathbb{S}^\tau$, $\tau \in \mathcal{T}$, do not change anymore. Similarly to the relaxation labeling, this algorithm returns non-empty sets $\mathbb{S}_w$ for all $w \in$ $\mathcal{V} \cup \mathcal{E}$ if and only if there is subproblem agreement for the decomposition costs $\theta^{\mathcal{T}}$. Moreover, the returned sets $\mathbb{S}^\tau$, $\tau \in \mathcal{T}$, are those defining the subproblem agreement as in Definition 9.7. Although this algorithm is well-defined, its implementation can be costly and non-trivial, due to the necessity to operate with the sets of *all* solutions of the slave problems. In case of the complete representation these solutions can be explicitly enumerated, however, this will in general be nontrivial and computationally costly. This is due to the number of optimal solutions, which may grow exponentially with the size of the (sub)problem.

**Strong subproblem agreement**   One can also speak about *strong* subproblem agreement, that is, subproblem agreement, where each slave problem has a unique optimum. In this case the conditions of Proposition 9.6 hold for $\mu^* = \delta(y)$, where the restriction of $y$ to the subgraph $\mathcal{G}^\tau$ is the optimal labeling for the slave problem $\tau$. Similar to strict arc consistency, in this case the lower bound defined by the decomposition-based dual becomes equal to the value of the (non-relaxed) energy function and, therefore, the optimum is reached for the dual, primal relaxed and non-relaxed problems.

In general, however, (weak) subproblem agreement is only necessary and not a sufficient condition. This follows from the fact that it coincides with node-edge agreement in the special case of the complete decomposition.

## 9.4   Equivalence of all acyclic decompositions

**Transformation between acyclic duals**   In this section we will concentrate on acyclic decompositions only and show an even closer relation between the Lagrange dual and the decomposition-based one. In particular, we will show how the decomposition-based dual variable values $\lambda$ can be transformed into the Lagrange dual variables $\phi$ while keeping the value of the dual function unchanged. As a side result, we will show equivalence of node-edge agreement and tree agreement for an arbitrary acyclic decomposition and not only for the complete one.

Proposition 9.2 shows that to guarantee $\mathcal{U}(\theta^{\mathcal{T}}) \geq \mathcal{D}(\phi)$ it is sufficient that $\theta^\tau = \rho_w^\tau \theta_w^\phi$ is satisfied for some $\rho_w^\tau \geq 0$ with $\sum_{\tau \in \mathcal{T}} \rho_w^\tau = 1$ for all $w \in \mathcal{V} \cup \mathcal{E}$. Let us now show that there exists an inverse transformation, $\theta^\tau \to \phi$ satisfying $\mathcal{D}(\phi) = \mathcal{U}(\theta^\mathcal{T})$.

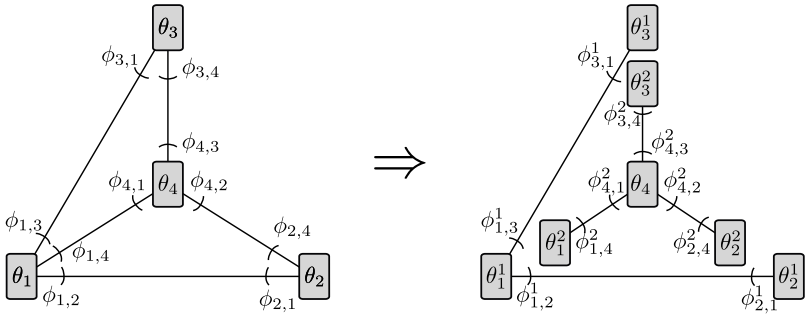

**Figure 9.6:** An example of a canonical decomposition. Small arcs on graph edges and associated vectors $\phi_{u,v}$ and $\phi_{u,v}^\tau$ stand for reparametrizations of the master problem and the slave subproblems respectively. Note that the reparametrizations of all slave subproblems in total have the same dimensionality as the reparametrization of the master problem. Therefore, one may build a one-to-one mapping between them as $\phi_{u,v} = \phi_{u,v}^\tau$ for the (unique) $\tau$ being the index of the subproblem containing the edge $uv$. In the figure these equalities are $\phi_{3,1} = \phi_{3,1}^1$, $\phi_{3,4} = \phi_{3,4}^2$, $\phi_{1,4} = \phi_{1,4}^2$ etc.

Let $(\mathcal{G}^\mathcal{T}, \theta^\mathcal{T})$ be an acyclic decomposition. Let us also assume that (i) each slave subgraph is connected and contains at least one edge; (ii) the sets of edges of the slave subgraphs are non-intersecting, i.e.

$$\mathcal{E}^\tau \cap \mathcal{E}^\sigma = \emptyset, \ \tau \neq \sigma. \tag{9.39}$$

Such decompositions will be called *canonical*.

Although the definition of acyclic canonical decompositions may look restrictive, most of the used decompositions in practice fall into this class (see, e.g. the example of a grid graph decomposition in §9.2). This is because (i) splitting the pairwise costs between different subproblems is computationally expensive, since the number of corresponding dual variables grows quadratically with the number of labels; and (ii) splitting the pairwise costs does not lead to a tighter relaxation with an acyclic decomposition as it follows from Proposition 9.3.

Note that due to (9.39) if one considers a reparametrization $[\theta^\tau]^\phi$ of each slave subproblem, there will be precisely $|\mathcal{Y}_u| + |\mathcal{Y}_v|$ dual variables assigned to each edge $uv \in \mathcal{E}$, since there is exactly one slave subproblem containing the edge $uv$ for each $uv$, see Figure 9.6 for illustration. Therefore, such a reparametrization can be considered as a reparametrization of the initial master problem, that is, for any $uv \in \mathcal{E}$ we can assume $\phi_{u,v} := \phi_{u,v}^\tau$, where $\tau$ is the index of the only subproblem containing the edge $uv$.

Note also, that for a canonical decomposition for any reparametrization $\phi$ and all $w \in \mathcal{V} \cup \mathcal{E}$ and $s \in \mathcal{Y}_w$ it holds that

$$\sum_{\tau \in \mathcal{T}_w} [\theta^\tau]_w^\phi(s) = \theta_w^\phi(s). \tag{9.40}$$

Indeed, when $w$ is an edge, since $|\mathcal{T}_{uv}| = 1$ for all $uv \in \mathcal{E}$, equality (9.40) transforms into $[\theta^{uv}]_{uv}^\phi(s) = \theta_w^\phi(s)$, which holds by definition of a canonical decomposition. When $w$ is a node $u \in \mathcal{V}$, the equality can be proved straightforwardly:

$$\begin{aligned}
\sum_{\tau \in \mathcal{T}_u} [\theta^\tau]_u^\phi(s) &= \sum_{\tau \in \mathcal{T}_u} \left( \theta_u^\tau(s) - \sum_{v \in \mathcal{V}^\tau \cap \mathcal{N}_b(u)} \phi_{u,v}(s) \right) \\
&= \sum_{\tau \in \mathcal{T}_u} \theta_u^\tau(s) - \sum_{\tau \in \mathcal{T}_u} \sum_{v \in \mathcal{V}^\tau \cap \mathcal{N}_b(u)} \phi_{u,v}(s) \\
&= \theta_u(s) - \sum_{v \in \mathcal{N}_b(u)} \phi_{u,v}(s) = \theta_u^\phi(s). \tag{9.41}
\end{aligned}$$

**Lemma 9.8.** Let $(\mathcal{G}^\mathcal{T}, \theta^\mathcal{T})$ be a canonical acyclic decomposition. Let also $[\theta^\tau]^\phi$, $\tau \in \mathcal{T}$, be optimal reparametrizations of the slave problems such that $[\theta_u^\tau]^\phi = [\theta_u^\sigma]^\phi$ for all $u \in \mathcal{V}$, $\tau, \sigma \in \mathcal{T}_u$. Then $\mathcal{D}(\phi) = \mathcal{U}(\theta^\mathcal{T})$.

*Proof.* Note that due to (9.39) it holds (i) $\theta_{uv}^\phi = [\theta_{uv}^\tau]^\phi$ for $\tau \in \mathcal{T}_{uv}$ and, since $[\theta_u^\tau]^\phi = [\theta_u^\sigma]^\phi$, it holds also that (ii) $[\theta_u^\tau]^\phi = \frac{\theta_u^\phi}{|\mathcal{T}_u|}$. Therefore

$$\begin{aligned}
\min_{s \in \mathcal{Y}_u} \theta_u^\phi(s) &= \sum_{\tau \in \mathcal{T}_u} \frac{1}{|\mathcal{T}_u|} \min_{s \in \mathcal{Y}_u} \theta_u^\phi(s) = \sum_{\tau \in \mathcal{T}_u} \min_{s \in \mathcal{Y}_u} \frac{1}{|\mathcal{T}_u|} \theta_u^\phi(s) \\
&= \sum_{\tau \in \mathcal{T}_u} \min_{s \in \mathcal{Y}_u} [\theta_u^\tau]^\phi(s). \tag{9.42}
\end{aligned}$$

we can use this in the following sequence of equalities:

$$\begin{aligned}
\mathcal{D}(\phi) &= \sum_{u \in \mathcal{V}} \min_{s \in \mathcal{Y}_u} \theta_u^\phi(s) + \sum_{uv \in \mathcal{E}} \min_{s \in \mathcal{Y}_{uv}} \theta_{uv}^\phi(s) \\
&= \sum_{u \in \mathcal{V}} \sum_{\tau \in \mathcal{T}_u} \min_{s \in \mathcal{Y}_u} [\theta_u^\tau]^\phi(s) + \sum_{\tau \in \mathcal{T}} \sum_{uv \in \mathcal{E}^\tau} \min_{s \in \mathcal{Y}_{uv}} [\theta_{uv}^\tau]^\phi(s)
\end{aligned}$$

$$
\begin{aligned}
&= \sum_{\tau \in \mathcal{T}} \sum_{u \in \mathcal{V}^\tau} \min_{s \in \mathcal{Y}_u} [\theta_u^\tau]^\phi(s) + \sum_{\tau \in \mathcal{T}} \sum_{uv \in \mathcal{E}^\tau} \min_{s \in \mathcal{Y}_{uv}} [\theta_{uv}^\tau]^\phi(s) \\
&= \sum_{\tau \in \mathcal{T}} \left( \sum_{u \in \mathcal{V}^\tau} \min_{s \in \mathcal{Y}_u} [\theta_u^\tau]^\phi(s) + \sum_{uv \in \mathcal{E}^\tau} \min_{s \in \mathcal{Y}_{uv}} [\theta_{uv}^\tau]^\phi(s) \right) \\
&\overset{(*)}{=} \sum_{\tau \in \mathcal{T}} \min_{y \in \mathcal{Y}_{\mathcal{V}^\tau}} \left( \sum_{u \in \mathcal{V}^\tau} \theta_u^\tau(y_u) + \sum_{uv \in \mathcal{E}^\tau} \theta_{uv}^\tau(y_{uv}) \right) = \mathcal{U}(\theta^\mathcal{T}).
\end{aligned}
\tag{9.43}
$$

The equality marked with $(*)$ follows from the fact that $[\theta^\tau]^\phi$ is an optimal reparametrization and the Lagrange relaxation is tight for acyclic problems (Corollary 6.19). $\qquad\square$

The equality $[\theta_u^\tau]^\phi = [\theta_u^\sigma]^\phi$ required by Lemma 9.8 can be satisfied e.g. by setting the unary costs to zero. This can always be done in practice, as follows from the following lemma:

**Lemma 9.9.** Let each node of a graphical model have at least one incident edge. Then there is always a dual optimal reparametrization $\phi$ of the costs $\theta$ such that $\theta_u^\phi = 0$ for all nodes $u$.

*Proof.* Let $\theta$ be an arbitrary dual optimal reparametrization and $p^u \in \Delta^{\mathcal{N}_b(u)}$, i.e. coordinate of $p^u$ are non-negative and satisfy $\sum_{v \in \mathcal{N}_b(u)} p_v^u = 1$. The transformation

$$
\phi_{u,v}(s) := p_v^u \theta_u(s), \ \forall v \in \mathcal{N}_b(u)
\tag{9.44}
$$

applied to each node $u \in \mathcal{V}$ results in $\theta_u^\phi = 0$. In other words, the unary costs are redistributed between pairwise costs of incident edges.

It remains to show that $\theta^\phi$ is optimal as well. This would follow from the fact that the dual values $\mathcal{D}(0)$ and $\mathcal{D}(\phi)$ corresponding to the initial reparametrization and the new one, respectively, coincide.

Indeed, since $\theta$ is dual optimal, it fulfills the node-edge agreement criterion (see Definition 6.9 in Chapter 6). This implies that for each edge $uv \in \mathcal{E}$ there exists a label pair $(s, l) \in \mathcal{Y}_{uv}$ such that for all $s' \in \mathcal{Y}_u$, $l' \in \mathcal{Y}_v$:

$$
\begin{aligned}
\theta_{uv}(s, l) &\leq \theta_{uv}(s', l') \\
\theta_u(s) &\leq \theta_u(s') \\
\theta_v(l) &\leq \theta_v(l').
\end{aligned}
\tag{9.45}
$$

Since $\theta_{uv}^{\phi}(s,l)$ is obtained by adding to the first line in (9.45) the second and the third ones multiplied by non-negative numbers, for all $(s',l') \in \mathcal{Y}_{uv}$ it holds also that:

$$\theta_{uv}^{\phi}(s,l) \leq \theta_{uv}^{\phi}(s',l').\tag{9.46}$$

Moreover, the condition (9.46) holds for all label pairs satisfying (9.45).

Therefore, taking into account that $\theta_u^{\phi} = 0$ for all $u$, it holds

$$
\begin{aligned}
\mathcal{D}(\phi) &= \sum_{uv \in \mathcal{E}} \min_{(s,t) \in \mathcal{Y}_{uv}} \theta_{uv}^{\phi}(s,t) \\
&\overset{(9.44)}{=} \sum_{uv \in \mathcal{E}} \min_{(s,t) \in \mathcal{Y}_{uv}} \left( \theta_{uv}(s,t) + p_v^u \theta_u(s) + p_u^v \theta_v(l) \right) \\
&\overset{(9.45)}{=} \sum_{uv \in \mathcal{E}} \left( \min_{(s,t) \in \mathcal{Y}_{uv}} \theta_{uv}(s,t) + p_v^u \min_{s \in \mathcal{Y}_u} \theta_u(s) + p_u^v \min_{l \in \mathcal{Y}_v} \theta_v(l) \right) \\
&= \sum_{uv \in \mathcal{E}} \min_{(s,t) \in \mathcal{Y}_{uv}} \theta_{uv}(s,t) + \sum_{u \in \mathcal{V}} \sum_{v \in \mathcal{N}_b(u)} p_v^u \min_{s \in \mathcal{Y}_u} \theta_u(s) \\
&= \sum_{uv \in \mathcal{E}} \min_{(s,t) \in \mathcal{Y}_{uv}} \theta_{uv}(s,t) + \sum_{u \in \mathcal{V}} \left( \min_{s \in \mathcal{Y}_u} \theta_u(s) \cdot \sum_{v \in \mathcal{N}_b(u)} p_v^u \right) \\
&= \sum_{uv \in \mathcal{E}} \min_{(s,t) \in \mathcal{Y}_{uv}} \theta_{uv}(s,t) + \sum_{u \in \mathcal{V}} \min_{s \in \mathcal{Y}_u} \theta_u(s) = \mathcal{D}(0).\tag{9.47}
\end{aligned}
$$

$\square$

Lemma 9.8, therefore, implies the following practically applicable statement:

**Proposition 9.10.** Let $(\mathcal{G}^{\mathcal{T}}, \theta^{\mathcal{T}})$ be a canonical acyclic decomposition. Let $[\theta^{\tau}]^{\phi}$, $\tau \in \mathcal{T}$, be optimal reparametrizations of the slave problems such that $[\theta_u^{\tau}]^{\phi} = 0$ for all $u \in \mathcal{V}^{\tau}$. Then $\mathcal{D}(\phi) = \mathcal{U}(\theta^{\mathcal{T}})$.

Proposition 9.10 together with Proposition 9.2 give us a tool for switching between the Lagrange dual $\mathcal{D}$ and canonical acyclic decomposition-based duals $\mathcal{U}$:

$$\phi \to \theta^{\mathcal{T}}: \quad \theta_w^{\tau} := \rho_w^{\tau} \theta_w, \ \rho \geq 0, \ \sum_{\tau \in \mathcal{T}_w} \rho_w^{\tau} = 1, \ w \in \mathcal{V} \cup \mathcal{E}\tag{9.48}$$

$$\theta^{\mathcal{T}} \to \phi: \quad [\theta_u^\tau]^\phi = 0,$$

$$\sum_{uv \in \mathcal{E}^\tau} \min_{s \in \mathcal{Y}_{uv}} [\theta_{uv}^\tau(s)]^\phi = \min_{y \in \mathcal{Y}_{\mathcal{V}^\tau}} \sum_{w \in \mathcal{V}^\tau \cup \mathcal{E}^\tau} \theta_w^\tau(y_w), \; \tau \in \mathcal{T}. \tag{9.49}$$

In particular, one can change decompositions during optimization without decreasing the corresponding dual value.

**Tree agreement** for acyclic decompositions is preserved by the transformations (9.48)-(9.49), as the propositions below suggest.

The following technical lemma is used in the proof of all propositions below:

**Lemma 9.11.** Let $\rho_w \geq 0$ and $\theta'$ be defined as $\theta'_w = \rho_w \theta_w$ for all $w \in \mathcal{V} \cup \mathcal{E}$. Then $\mathrm{mi}[\theta'] \geq \mathrm{mi}[\theta]$. If $\rho > 0$, then $\mathrm{mi}[\theta'] = \mathrm{mi}[\theta]$.

The proof follows directly from the definition of mi as the vector indicating all locally minimal labels and label pairs with respect to $\theta$.

**Proposition 9.12.** Let $(\mathcal{G}^{\mathcal{T}}, \theta^{\mathcal{T}})$ be an acyclic decomposition. Let $\theta_w^\tau = \rho_w^\tau \theta_w$ for all $w \in \mathcal{V} \cup \mathcal{E}$ and some $\rho_w^\tau > 0$ such that $\sum_{\tau \in \mathcal{T}_w} \rho_w^\tau = 1$ for all $w \in \mathcal{V} \cup \mathcal{E}$. Then node-edge agreement for $\theta$ implies tree agreement for $\theta^{\mathcal{T}}$.

*Proof.* Consider $\mathbb{S}_w = \{s \in \mathcal{Y}_w : \mathrm{cl}(\mathrm{mi}[\theta])_w(s) = 1\}$ for all $w \in \mathcal{V} \cup \mathcal{E}$. Observe that $\mathbb{S}_w$ is not empty due to node-edge agreement. Due to Lemma 9.11 the restriction $\xi^\tau$ of the binary vector $\mathrm{cl}(\mathrm{mi}[\theta])$ to any subgraph $\tau$ is arc-consistent and $\xi^\tau \neq \bar{0}$. To prove the tree agreement it is sufficient to show that for each $w \in \mathcal{V} \cup \mathcal{E}$, each $s \in \mathbb{S}_w$ and each $\tau \in \mathcal{T}_w$ there exists $y^\tau \in \mathcal{Y}^\tau$ such that $y_w^\tau = s$ and $y^\tau \in \arg\min_{y \in \mathcal{Y}^\tau} \langle \theta^\tau, \delta(y) \rangle$. Since each subgraph in the decomposition is acyclic, the latter condition follows directly from Proposition 6.17. □

Example 9.13 below shows that the inverse claim does not hold: In general tree agreement does not imply node-edge agreement. However, it holds if an additional condition is imposed, as stated by Proposition 9.14.

**Example 9.13.** Consider the decomposed problem in Figure 9.7. Each of the slave subproblems has a unique optimal labeling consisting of the upper labels assigned to each node. These labelings coincide in the split node $u$, therefore, tree agreement holds. However, there is no node-edge agreement in the master model.

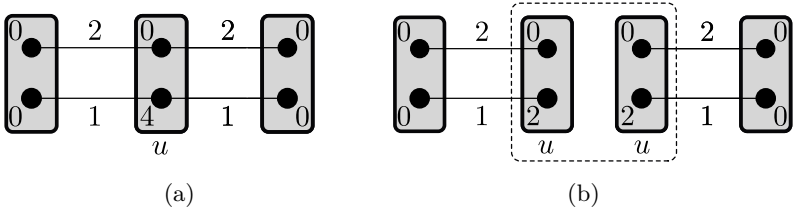

**Figure 9.7:** Illustration to Example 9.13. Numbers denote costs of labels and label pairs. Not shown label pairs are assigned a very high cost and, therefore, are irrelevant. **(a)** Master model. **(b)** Decomposition of the master model. The two copies of the decomposed node $u$ are encircled by a dashed line. Each of the slave subproblems has a unique optimal labeling consisting of upper labels assigned to each node. These labelings coincide in the split node $u$, therefore, the tree agreement holds. However, there is no node-edge agreement in the master model.

**Proposition 9.14.** Let $\theta^{\mathcal{T}}$ be a canonical decomposition and $[\theta^\tau]^\phi$, $\tau \in \mathcal{T}$, be optimal reparametrizations of the slave problems. Then if the tree agreement holds for $\theta^{\mathcal{T}}$ the node-edge agreement holds for $\theta^\phi$.

*Proof.* From optimality of $[\theta^\tau]^\phi$ the node-edge agreement within the subproblem $\tau$ follows for any $w \in \mathcal{V}^\tau \cup \mathcal{E}^\tau$. The node-edge agreement of the master problem follows from (9.40), which is, for any $s \in \mathcal{Y}_w$ it holds that

$$\sum_{\tau \in \mathcal{T}_w} [\theta^\tau]_w^\phi(s) = \theta_w^\phi(s) \,. \tag{9.50}$$

Therefore, local optimality of the label $s$ for each subproblem $\tau$ implies its local optimality in the master model. See Figure 9.8 for illustration. $\qquad\square$

## 9.5   Bibliography and further reading

The classical reference to the Lagrange decomposition technique is the work of [32]. The follow-up paper [31] gives a number of examples of the Lagrange decomposition for different integer linear programming problems. The survey [30] summarizes the experience gained from using the technique over 15 years and gives a number of mathematical and applied insights.

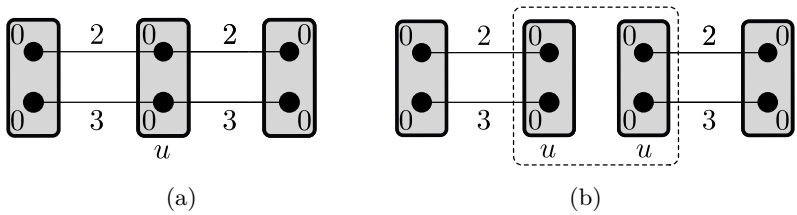

**Figure 9.8:** The same decomposition as in Figure 9.7, but with optimal reparametrization of each slave subproblem. This reparametrization implies node-edge agreement for the respectively reparametrized master problem.

For MAP-inference in graphical models the Lagrange decomposition was proposed in [119]. Later on, and independently, similar work was done in [101, 37] and [55], followed by [56].

The notion of tree agreement was originally proposed in [134], and its current, generalized definition is given in [49].

Decompositions into subgraphs with cycles were used by multiple researchers [58, 116, 118, 42], certain optimality guarantees for selecting relevant cycles in a sequential (cutting plane) fashion were given in [116].

# 10

## Maximization of the Decomposition-Based Dual

This chapter is devoted to optimization of the decomposition-based dual introduced in Chapter 9.

As in the case of Lagrange decomposition we will concentrate on two methods: the subgradient method and block-coordinate ascent. The former possesses convergence guarantees, although it is typically quite slow in practice. The latter is usually much faster, however, it does not guarantee convergence to the dual optimum.

The general block-coordinate ascent method that we describe below, is a generalization of the min-sum diffusion algorithm. However, its straightforward implementation turns out to be very inefficient. Therefore, we will consider a special case of the acyclic decomposition and modify accordingly the method's implementation until it becomes the state-of-the-art Sequential Tree-Reweighted Message Passing algorithm (TRW-S). The latter also can be seen as a variant of the anisotropic diffusion method, introduced in Chapter 8. This variant is also known as a Sequential Reweighted Message Passing algorithm (SRMP), that for the considered pairwise graphical models is equivalent to TRW-S.

At the end of the chapter we provide an empirical comparison of the sub-gradient and block-coordinate ascent methods for different acyclic decompositions.

## 10.1 Subgradient method

We will assume the decomposition $\mathcal{G}^{\mathcal{T}}$ of the master graph $\mathcal{G}$ to be given. The task is to maximize the decomposition-based dual

$$\mathcal{U}(\theta^{\mathcal{T}}[\lambda]) = \sum_{\tau \in \mathcal{T}} \min_{y \in \mathcal{Y}_{\mathcal{V}^\tau}} \langle \theta^\tau[\lambda], \delta(y) \rangle \tag{10.1}$$

with respect to the subproblem cost distribution $\theta^{\mathcal{T}} \in \Theta^{\mathcal{T}}$ defined by (9.15).

**Subgradient $\frac{\partial \mathcal{U}}{\partial \lambda}$ of decomposition-based dual $\mathcal{U}$** w.r.t. $\lambda$ can be obtained as a subgradient of a composite function (Proposition 5.29):

$$\frac{\partial \mathcal{U}}{\partial \lambda} = \frac{\partial \mathcal{U}}{\partial \theta^{\mathcal{T}}} \frac{\partial \theta^{\mathcal{T}}}{\partial \lambda} \,, \tag{10.2}$$

where the two terms on the right-hand-side are combined via matrix multiplication. Note that $\theta^{\mathcal{T}}$ is a linear function of $\lambda$, given by (9.18). As a result, the corresponding rows $\frac{\partial \theta_w^\tau[\lambda](s)}{\partial \lambda_{w'}^\sigma(t)}$ of the matrix $\frac{\partial \theta^{\mathcal{T}}}{\partial \lambda}$ have non-zero entries only for $s = t$, $w = w'$ and $\tau = \sigma$. These entries are equal to 1 for $\tau \neq \tau_w$ and $-1$ for $\tau = \tau_w$.

To compute $\frac{\partial \mathcal{U}}{\partial \theta^{\mathcal{T}}}$ note that $\mathcal{U}$ is sum of piecewise linear concave functions and a subgradient of each such functions can be computed according to Corollary 5.34. Since only a single summand $\min_{y \in \mathcal{Y}_{\mathcal{V}^\tau}} \langle \theta^\tau[\lambda], \delta(y) \rangle$ in (10.1) is dependent on the costs $\theta^\tau$, there is a subgradient of $\mathcal{U}$ with coordinates

$$\frac{\partial \mathcal{U}}{\partial \theta_w^\tau} = \delta(y^\tau)_w \,, \tag{10.3}$$

where $y^\tau \in \arg\min_{y \in \mathcal{Y}_{\mathcal{V}^\tau}} \langle \theta^\tau[\lambda], \delta(y) \rangle$ is an optimal labeling of the subproblem $\tau$.

Putting (10.3) and (10.2) together gives the following formula for computing subgradient coordinates:

$$\frac{\partial \mathcal{U}}{\partial \lambda_w^\tau(s)} = \left\{ \begin{array}{rl} 0, & y_w^\tau = s \text{ and } y_w^{\tau_w} = s\,, \\ 0, & y_w^\tau \neq s \text{ and } y_w^{\tau_w} \neq s\,, \\ 1, & y_w^\tau = s \text{ and } y_w^{\tau_w} \neq s\,, \\ -1, & y_w^\tau \neq s \text{ and } y_w^{\tau_w} = s\,. \end{array} \right. \tag{10.4}$$

Note that computing the subgradient does not require any additional computations compared to those needed to evaluate $\mathcal{U}$. This is an advantage of the subgradient method. The subgradient algorithm (7.7) therefore specializes to Algorithm 7.

---

**Algorithm 7** Subgradient method for $\mathcal{U}(\theta^{\mathcal{T}}[\lambda])$

---

1: **Init:** $\lambda^0$ - starting point, $N$ - number of iterations
2: **for** $t = 1$ **to** $N$ **do**
3:     For all $\tau \in \mathcal{T}$ compute $y^\tau = \arg\min_{y \in \mathcal{Y}_{\mathcal{V}^\tau}} \langle \theta^\tau[\lambda^{t-1}], \delta(y) \rangle$
4:     Compute $\frac{\partial \mathcal{U}}{\partial \lambda}$ according to (10.4)
5:     $\lambda^t := \lambda^{t-1} + \alpha^t \frac{\partial \mathcal{U}}{\partial \lambda}$ with $\alpha^t$ defined as in §7.2
6: **end for**
7: **return** $\lambda^N$

---

The only prerequisite of Algorithm 7 is that all slave subproblems must be efficiently solvable, since their solutions $y^\tau$ are required on each iteration for computing the subgradient.

**Example 10.1.** Consider the grid-graph example and the corresponding decomposition from §9.2, see Figures 9.1, 9.2 and 10.1. The subgradient $\frac{\partial \mathcal{U}}{\partial \lambda}$ in this case reads

$$\frac{\partial \mathcal{U}}{\partial \lambda} = \delta(y^c) - \delta(y^r)\,, \tag{10.5}$$

with $y^c \in \arg\min_{y \in \mathcal{Y}_\mathcal{V}} \langle \theta^c[\lambda], \delta(y) \rangle$ and $y^r \in \arg\min_{y \in \mathcal{Y}_\mathcal{V}} \langle \theta^r[\lambda], \delta(y) \rangle$ being solutions to the column- and row subproblems and $\theta^c[\lambda]$, $\theta^r[\lambda]$ defined by (9.5).

Therefore, Algorithm 7 specializes to Algorithm 8. Let $s$ be a label in node $u$ such that $y_u^c = s \neq y_u^r$. Then the subgradient will increase the cost of this label in the column subproblem and (therefore) decrease it in the row subproblem. This eventually will either make this label non-optimal for the column subproblem or make it optimal for the row subproblem. Therefore, while getting closer to the dual optimum, more and more subproblem agreement is attained. The latter supports the fact that subproblem agreement is a necessary optimality condition.

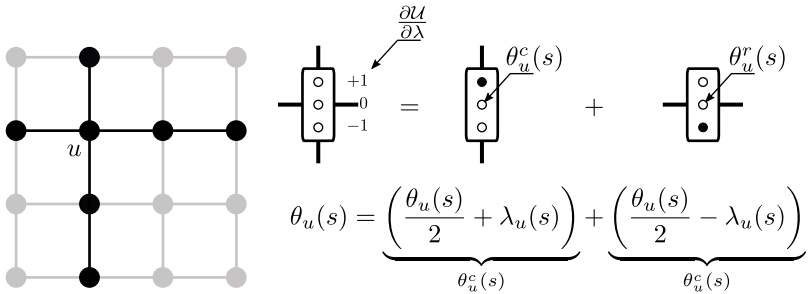

**Figure 10.1:** An illustration of the subgradient computation for the row-column decomposition of grid graphs. Filled circles in the right part of the figure denote the labels belonging to the optimal corresponding row- and column labelings. If they do not match, the corresponding subgradient coordinates obtain values $+1$ and $-1$.

---

**Algorithm 8** Subgradient method for column-row decomposition of grid graphs

---

1: **Init:** $\lambda^0$ - start point, $N$ - number of iterations
2: **for** $t = 1$ **to** $N$ **do**
3:    $y^c \in \arg\min_{y \in \mathcal{Y}_\mathcal{V}} \langle \theta^c[\lambda], \delta(y) \rangle, \quad y^r \in \arg\min_{y \in \mathcal{Y}_\mathcal{V}} \langle \theta^r[\lambda], \delta(y) \rangle$
4:    Compute $g_u(s) = \begin{cases} 0, & y_u^c = y_u^r = s \\ 0, & y_u^c \neq s, \ y_u^r \neq s \\ 1, & y_u^c = s \neq y_u^r \\ -1, & y_u^c \neq s = y_u^r \end{cases}$
5:    $\lambda_u^t(s) = \lambda_u^{t-1}(s) + \alpha^t g_u(s), \ u \in \mathcal{V}, \ s \in \mathcal{Y}_u$
6: **end for**
7: **return** $\lambda^N$

---

**Exercise 10.2.** Consider the complete decomposition and check that the corresponding subgradient coincides with the one provided by (8.12) for the Lagrange dual in §8.2.

## 10.2   Tree Reweighted Message Passing (TRW-S)

### 10.2.1   Block-coordinate ascent: Subproblems averaging

**Edge decomposition**   To construct a block-coordinate ascent method for the decomposition-based dual we will generalize the min-sum diffusion algorithm (8.15)-(8.16) introduced in §8.3.1 for the Lagrange dual. To this end, let us consider the *edge decomposition* of the master graph $\mathcal{G}$ given as follows:

- Reformulate the labeling problem $(\mathcal{G}, \mathcal{Y}_{\mathcal{V}}, \theta)$ in a way that all unary costs $\theta_u$, $u \in \mathcal{V}$, are assigned zero values. This can be done for example by the reparametrization

$$\phi_{u,v}(s) := \frac{\theta_u(s)}{|\mathcal{N}_b(u)|}, \ v \in \mathcal{N}_b(u), \ s \in \mathcal{Y}_u \,, \qquad (10.6)$$

  where the unary costs are uniformly distributed between adjacent pairwise costs.

- Let $\mathcal{T} = \mathcal{E}$, i.e. each subproblem is associated with an edge and has the form $\mathcal{G}^\tau = (\{u, v\}, \{uv\})$ for $\tau = uv \in \mathcal{E}$.

In contrast to the complete decomposition, the edge decomposition does not contain graph nodes as separate subproblems.

   In the edge decomposition each edge belongs to a single subproblem only and the unary costs are zero. Therefore, abusing notation we write $\theta_{uv}$ for the cost vector of the subproblem $uv \in \mathcal{T}$, instead of using the notation $\theta^\tau \equiv \theta^{uv}$ above.

**Diffusion as Subproblem Averaging**   Given the edge decomposition, one can rewrite the diffusion algorithm in the form of *subproblem averaging*, where for each node $u$ the costs $\theta_{uv}^t$, $v \in \mathcal{N}_b(u)$, on iteration $t$ are computed in a way that the following equality marked with $(*)$ holds:

$$\min_{l \in \mathcal{Y}_v} \theta_{uv}^t(s, l) \stackrel{(*)}{=} \frac{\sum\limits_{u'v' \in \mathcal{T}_u} \min\limits_{l' \in \mathcal{Y}_{v'}} \theta_{u'v'}^t(s, l')}{|\mathcal{T}_u|}$$

$$\equiv \frac{\sum\limits_{v' \in \mathcal{N}_b(u)} \min\limits_{l' \in \mathcal{Y}_{v'}} \theta_{uv'}^t(s, l')}{|\mathcal{N}_b(u)|}, \ s \in \mathcal{Y}_u \,. \qquad (10.7)$$

In other words, for each label $s$ in the current node $u$, all optimal labelings containing $s$ obtain the same cost in all subproblems intersecting in $u$. In the considered case of the edge decomposition, these optimal labelings consist of label pairs $(s, l^*)$ with $l^* = \arg\min_{l \in \mathcal{Y}_v} \theta_{uv}^t(s, l)$. Let us denote the total cost of the labeling $(s, l^*)$ in the subproblem $uv$ as $E_u^\tau(s)$, where the upper index $\tau$ indexes the subproblem and the lower one the node the label $s$ belongs to.

Omitting the iteration index $t$, one can rewrite Equation (10.7) as follows:

$$E_u^\tau(s) = \frac{\sum_{\tau \in \mathcal{T}_u} E_u^\tau(s)}{|\mathcal{T}_u|}, \ s \in \mathcal{Y}_u. \tag{10.8}$$

**General subproblem averaging** Equation (10.8) is not specific to the edge decomposition anymore. Therefore, it can be applied to a general decomposition $\mathcal{G}^\mathcal{T}$ as an equation which must be fulfilled after a step of the diffusion-like algorithm applied to node $u$ and label $s$. In this case $E_u^\tau(s)$ is defined as the total cost of the best labeling of the subproblem $\tau$, if label $s$ is assigned in node $u$:

$$E_u^\tau(s) := \min_{y \in \mathcal{Y}_{\mathcal{V}^\tau}} \langle \theta^\tau, \delta(y) \rangle, \ \text{s.t.} \ y_u = s. \tag{10.9}$$

Note that to compute $E_u^\tau(s)$ one may simply reduce the label set $\mathcal{Y}_u$ to $\{s\}$ and find an optimal labeling in the resulting problem.

**Proposition 10.3.** Let $(\mathcal{G}^\mathcal{T}, \theta^\mathcal{T}[\lambda])$ be a graph decomposition, $u \in \mathcal{V}$ and let condition (10.8) hold for all $s \in \mathcal{Y}_u$ with $E_u^\tau(s)$ defined as in (10.9). Then $\theta^\mathcal{T}[\lambda]$ is an optimum of $\mathcal{U}$ with respect to the block of variables $(\lambda_u^\tau(s), s \in \mathcal{Y}_u, \tau \in \mathcal{T}_u \backslash \{\tau_u\})$.

The index $\tau_u$ is excluded in the last expression, since there is no variable $\lambda_u^{\tau_u}$ and the costs $\theta_u^{\tau_u}$ are expressed through $\lambda_u^\tau$, $\tau \neq \tau_u$, see (9.18).

*Proof.* According to Lemma 7.7 the coordinate optimum condition is equivalent to the existence of the zero subgradient with respect to the variables one minimizes over, i.e. $\bar{0} \in \frac{\partial \mathcal{U}}{\partial \lambda_u^\tau(s)}$ for all $s \in \mathcal{Y}_u$ and $\tau \in \mathcal{T}_u \backslash \{\tau_u\}$.

According to (10.4), this condition holds, if there are

$$y^\tau \in \arg\min_{y \in \mathcal{Y}_{\mathcal{V}^\tau}} \langle \theta^\tau, \delta(y) \rangle$$

such that $y_u^\tau = y_u^{\tau_u}$ for all $\tau \in \mathcal{T}_u$. Indeed, the condition (10.8) means that $E_u^\tau(s) = E_u^{\tau_u}(s)$ for all $\tau \in \mathcal{T}_u$ and all $s \in \mathcal{Y}_u$. For $s^* \in \arg\min_{s \in \mathcal{Y}_u} E_u^\tau(s)$ it implies that $y_u^\tau = s^*$ for all $\tau \in \mathcal{T}_u$, which finalizes the proof. $\qquad\square$

**Remark 10.4.** An analogous statement can be proved similarly for the coordinate optimum w.r.t. the pairwise costs $\theta_{uv}^\tau$. First, we define

$$E_{uv}^\tau(s,l) := \min_{y \in \mathcal{Y}_{\mathcal{V}^\tau}} \langle \theta^\tau, \delta(y) \rangle, \text{ s.t. } y_u = s, \ y_v = l. \tag{10.10}$$

The coordinate optimum condition then reads

$$E_{uv}^\tau(s,l) = \frac{\sum_{\tau \in \mathcal{T}_{uv}} E_{uv}^\tau(s,l)}{|\mathcal{T}_{uv}|}. \tag{10.11}$$

Proposition 10.3 shows a possible way in which coordinate ascent algorithms can be constructed for the decomposition dual $\mathcal{U}$. To this end let us write an explicit update formula, which guarantees fulfillment of the coordinate optimality conditions (10.8). Let $[\cdot]^t$ denote the value of $\cdot$ on the iteration $t$ of an iterative algorithm. Then the following property holds:

**Proposition 10.5.** Let $[\theta^\mathcal{T}]^{t-1} \in \Theta^\mathcal{T}$ be the decomposition costs and let $w \in \mathcal{V} \cup \mathcal{E}$ be the factor processed at iteration $t$. Let $[\theta_w^\tau(s)]^t$ be updated such that for all $s \in \mathcal{Y}_w$ and $\tau \in \mathcal{T}_w$ it holds that

$$[\theta_w^\tau(s)]^t = [\theta_w^\tau(s)]^{t-1} - [E_w^\tau(s)]^{t-1} + \frac{1}{|\mathcal{T}_w|} \sum_{\tau' \in \mathcal{T}_w} [E_w^{\tau'}(s)]^{t-1} \tag{10.12}$$

and $[\theta_{w'}^\tau(s)]^t = [\theta_{w'}^\tau(s)]^{t-1}$ for all $w' \in \mathcal{V} \cup \mathcal{E} \backslash \{w\}$, $s \in \mathcal{Y}_{w'}$, $\tau \in \mathcal{T}_{w'}$.

Then $[\theta^\mathcal{T}]^t \in \Theta^\mathcal{T}$ and $[E_w^\tau(s)]^t = \frac{\sum_{\tau' \in \mathcal{T}_w} [E_w^{\tau'}(s)]^t}{|\mathcal{T}_w|}$.

*Proof.* The first statement follows from the definition of $\Theta^\mathcal{T}$ (see (9.15)), the assumption $[\theta^\mathcal{T}]^{t-1} \in \Theta^\mathcal{T}$ and the following sequence of equations:

$$\begin{aligned}
\sum_{\tau \in \mathcal{T}_w} [\theta_w^\tau(s)]^t &= \sum_{\tau \in \mathcal{T}_w} \left( [\theta_w^\tau(s)]^{t-1} - [E_w^\tau(s)]^{t-1} + \frac{1}{|\mathcal{T}_w|} \sum_{\tau' \in \mathcal{T}_w} [E_w^{\tau'}(s)]^{t-1} \right) \\
&= \sum_{\tau \in \mathcal{T}_w} [\theta_w^\tau(s)]^{t-1} - \sum_{\tau \in \mathcal{T}_w} [E_w^\tau(s)]^{t-1} + \sum_{\tau \in \mathcal{T}_w} [E_w^\tau(s)]^{t-1} \\
&= \sum_{\tau \in \mathcal{T}_w} [\theta_w^\tau(s)]^{t-1} = \theta_w(s). \tag{10.13}
\end{aligned}$$

To prove the second statement note that from $[\theta_w^\tau(s)]^t = [\theta_w^\tau(s)]^{t-1} + \alpha$ it follows $[E_w^\tau(s)]^t = [E_w^\tau(s)]^{t-1} + \alpha$. Applying this to (10.12) one obtains

$$[E_w^\tau(s)]^t = [E_w^\tau(s)]^{t-1} - [E_w^\tau(s)]^{t-1} + \frac{1}{|\mathcal{T}_w|} \sum_{\tau' \in \mathcal{T}_w} [E_w^{\tau'}(s)]^{t-1}$$

$$= \frac{1}{|\mathcal{T}_w|} \sum_{\tau' \in \mathcal{T}_w} [E_w^{\tau'}(s)]^{t-1} = \frac{1}{|\mathcal{T}_w|} \sum_{\tau' \in \mathcal{T}_w} [E_w^{\tau'}(s)]^t, \quad (10.14)$$

which finalizes the proof. $\square$

Proposition 10.5 implies that updates of the form (10.12) can be applied directly to the costs $\theta^{\mathcal{T}}$ without considering their unconstrained representation (9.18) explicitly. However, condition (10.12) together with the representation (9.18) implies that those coordinates of the dual vector $\lambda \in \mathbb{R}^{\mathcal{J}^{\mathcal{T}}}$ which correspond to $\tau \in \mathcal{T}_w \backslash \{\tau_w\}$, must be updated similarly:

$$[\lambda_w^\tau(s)]^t = [\lambda_w^\tau(s)]^{t-1} - [E_w^\tau(s)]^{t-1} + \frac{1}{|\mathcal{T}_w|} \sum_{\tau' \in \mathcal{T}} [E_w^{\tau'}(s)]^{t-1}. \quad (10.15)$$

Due to Proposition 10.5, a block-coordinate ascent algorithm for a general decomposition-based dual may look like Algorithm 9.

---

**Algorithm 9** General block-coordinate ascent algorithm for decomposition-based dual

---

1: **Init:** $[\theta^{\mathcal{T}}]^0$ - starting point, $N$ - number of iterations
2: **for** $t = 1$ **to** $N$ **do**
3:    **for** $\{w \in \mathcal{V} \cup \mathcal{E} : |\mathcal{T}_w| > 1\}$ **do**
4:       **for** $s \in \mathcal{Y}_w$ **do**
5:          Compute $[E_w^\tau(s)]^{t-1} := \min_{y \in \mathcal{Y}_{\mathcal{V}^\tau}} \{\langle [\theta^\tau]^{t-1}, \delta(y) \rangle \mid y_w = s\}$

6:          Compute $[\theta^\tau]^t$ as in (10.12)
7:       **end for**
8:    **end for**
9: **end for**
10: **return** $[\theta^{\mathcal{T}}]^N$

---

In its straightforward implementation, Algorithm 9 is very inefficient, since it requires each subproblem to be solved multiple times – as many as $\sum_{w \in \mathcal{V}^\tau \cup \mathcal{E}^\tau \,:\, |\mathcal{T}_w| > 1} |\mathcal{Y}_w|$ times, to be precise. Compare this to Algorithm 7 where each slave subproblem has to be solved only once per iteration. However, in some cases one can reuse computations in Algorithm 9 to make its implementation efficient. We will consider one such case in the following section.

### 10.2.2   Tree-reweighted sequential message passing

Let $(\mathcal{G}^\mathcal{T}, \theta^\mathcal{T})$ be an acyclic canonical decomposition. Let us consider Algorithm 9 and find an efficient implementation for such a decomposition.

**Complexity of Algorithm 9**  Note that the value $E_u^\tau(s)$ is nothing else but the min-marginal of label $s$ in node $u$, introduced in Chapter 2. The complexity of its computation is $O(\sum_{uv \in \mathcal{E}^\tau} |\mathcal{Y}_{uv}|)$ or simply $O(|\mathcal{E}^\tau| L^2)$ if we assume the number of labels in all nodes to be equal to $L$. Step 5 of Algorithm 9 must be performed $L$ times for a given node and $O(L|\mathcal{V}|)$ times in total. Therefore, a naive implementation of Algorithm 9 would require $O(L^3 \sum_{u \in \mathcal{V}} \sum_{\tau \in \mathcal{T}_u} |\mathcal{E}^\tau|)$ operations. Since for any canonical decomposition $\sum_{u \in \mathcal{V}} \sum_{\tau \in \mathcal{T}_u} 1 \leq \sum_{u \in \mathcal{V}} \sum_{v \in \mathcal{N}_b(u)} 1 = 2|\mathcal{E}|$, the total iteration complexity of Algorithm 9 is approximately $O(L^3 |\mathcal{E}||\bar{\mathcal{E}}|)$ with $|\bar{\mathcal{E}}|$ being the average number of edges in a slave subgraph.

Our goal in this section is to modify the algorithm in such a way that we attain an iteration complexity of $O(|\mathcal{E}| L^2)$. The latter is basically the size of the problem, since we have $L^2$ pairwise costs per edge.

**Joint computation of node min-marginals**  First note that for acyclic subgraphs min-marginals $E_u^\tau(s)$ can be computed with dynamic programming for all $s \in \mathcal{Y}_u$ simultaneously, see §2.2 for details. This allows for a reduction of the iteration complexity to $O(L^2 |\mathcal{E}||\bar{\mathcal{E}}|)$.

**Monotonic chain subproblems**  For further complexity reduction we must reuse computations performed multiple times to obtain the min-marginals $E_u^\tau(s)$. To this end we will assume the following:

- There is a complete order on the set of graph nodes, e.g. $\mathcal{V} = \{1, \ldots, |\mathcal{V}|\}$. We will assume expressions $u > v$ (and $u < v$) to be well-defined for any nodes $u \neq v$. With respect to this operation we will speak about *monotonic* node sequences, i.e. a sequence $u^1, u^2, \ldots, u^n$ is *monotonic*, if $u^i < u^{i+1}$ for $i = 1, \ldots, n-1$.

- All slave subgraphs $\mathcal{G}^\tau$, $\tau \in \mathcal{T}$, have a chain structure, i.e. there is a natural node order within each subgraph.

**Definition 10.6.** Let $\mathcal{G}^\tau$ be a chain subgraph of $\mathcal{G}$. We will call the chain $\mathcal{G}^\tau$ *monotonic*, if the sequence of consecutive nodes in the chain is monotonic.

**Example 10.7** (Row-column decomposition)**.** For a row- or column-wise ordering of nodes in a grid graph its row-column decomposition, as introduced in §9.2, consists of monotonic chains.

**Example 10.8** (Edge decomposition)**.** For any ordering of nodes in a graph, its complete and edge decompositions consist of monotonic chains.

Let $F_u^\tau(s)$ be forward min-marginals as introduced in §2.2. Since we are interested in min-marginals computed for a subgraph of $\mathcal{G}$, index $u$ will refer to node $u$ in the master-graph $\mathcal{G}$. Therefore, despite the complete order on $\mathcal{V}$ the node following $u$ in the subgraph may not have index $u + 1$, but can also exceed it. Taking this into account, the dynamic programming Algorithm 1 specializes to Algorithm 10 for the computation of forward min-marginals for the subgraph $\mathcal{G}^\tau$.

Backward marginals $B_i^\tau$ are computed using the same algorithm with the inverse order of nodes. Taking this into account, Algorithm 9 specializes to Algorithm 11.

Further speed-up of Algorithm 11 is based on the sequential character of the algorithm, see Figure 10.2 for the illustration. Each $F_u^\tau$ depends only on $\theta_v^\tau$ with $v < u$. In other words, if the value $\theta_v^\tau$ changes, values $F_u^\tau$ with $u \leq v$ will nevertheless remain correct forward min-marginals. The same holds for $B_u^\tau$ for $u \geq v$.

Therefore, if all chains $\mathcal{G}^\tau$ are monotonic w.r.t. the order in which nodes in Algorithm 11 are processed, there is no need to recompute $F_u^\tau$ and $B_u^\tau$. In this case Algorithm 11 turns into Algorithm 12.

---

**Algorithm 10** Forward marginals computation for a chain subgraph

---

1: **Init:** $\mathcal{V}^\tau$ - ordered set

   $v \in \mathcal{V}^\tau$ - node with the smallest index in $\mathcal{V}^\tau$;

   $F_v^\tau(s) = 0,\ s \in \mathcal{Y}_v$

2: **while** $v$ - has not the largest index in $\mathcal{V}^\tau$ **do**

3:    Find $u \in \mathcal{N}_b(v)\ :\ uv \in \mathcal{V}^\tau, u > v$ - next node in $\mathcal{V}^\tau$ after $v$

$$F_u^\tau(s) = \min_{l \in \mathcal{Y}_v} \left( F_v^\tau(l) + \theta_v^\tau(l) + \theta_{uv}^\tau(s, l) \right),\ s \in \mathcal{Y}_u \qquad (10.16)$$

4: **end while**

5: **return** $F_u^\tau(s), u \in \mathcal{V}^\tau, s \in \mathcal{Y}_u$

---

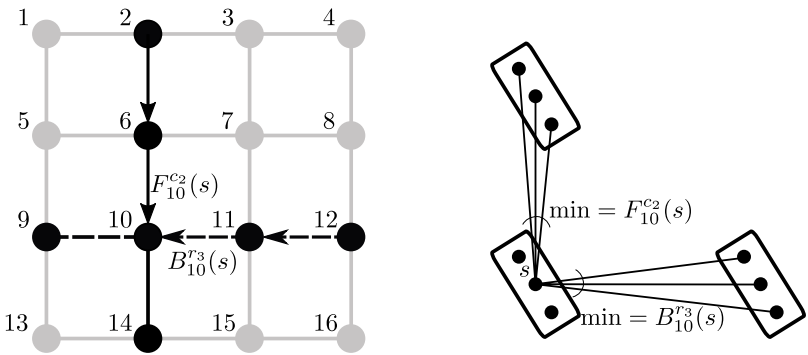

**Figure 10.2:** An illustration of the fact that $F_u^\tau$ depends only on $\theta_v^\tau$ with $v < u$ as well as $B_u^\tau$ for $v > u$. Consider a grid graph with the row-column decomposition. The nodes of the graph are enumerated row-wise from top to bottom. Rows and columns of the graph are enumerated as $r_1$ to $r_4$ from top to bottom and as $c_1$ to $c_4$ from left to right respectively. To compute the value $F_{10}^{c_2}$ the costs associated with the nodes 2 and 6 as well as with the edges $\{2, 6\}$, $\{6, 10\}$ are required. The costs associated with nodes 10 and 14 do not influence $F_{10}^{c_2}$. Similarly the costs associated with the nodes having indexes 10 and below do not influence the value of $B_{10}^{r_3}$.

Note that the very first iteration of Algorithm 12 does not correspond to the considered block-coordinate ascent scheme, since all values $B_u^\tau(s)$ are zero on this iteration instead of the backward min-marginals. This may lead to a decrease in the dual value on the first iteration. To avoid this, one could compute $B_u^\tau(s)$ at the start of the algorithm. However, in practice this does not pay off, because: (i) it requires almost as much

---

**Algorithm 11** Block-coordinate ascent algorithm for a chain decomposition

---

1: **Init:** $\mathcal{G}^{\mathcal{T}}$-a chain decomposition,
   $\theta^{\mathcal{T}}$ - initial cost distribution,
   $N$ - number of iterations
2: **for** $t = 1$ **to** $N$ **do**
3:    **for** $u \in \mathcal{V}$ **do**
4:       Compute $F_u^{\tau}(s)$, $B_u^{\tau}(s)$ for all $s \in \mathcal{Y}_u$, $\tau \in \mathcal{T}_u$;
5:       $E_u^{\tau}(s) := F_u^{\tau}(s) + B_u^{\tau}(s) + \theta_u^{\tau}(s)$, $s \in \mathcal{Y}_u$, $\tau \in \mathcal{T}_u$;
6:       $\theta_u^{\tau}(s) := \frac{1}{|\mathcal{T}_u|} \sum_{\sigma \in \mathcal{T}_u} E_u^{\sigma}(s) - F_u^{\tau}(s) - B_u^{\tau}(s)$, $s \in \mathcal{Y}_u$, $\tau \in \mathcal{T}_u$;
7:    **end for**
8: **end for**
9: **return** $\theta^{\mathcal{T}}$

---

---

**Algorithm 12** Block-coordinate ascent for a monotonic chain decomposition

---

1: **Init:** $\mathcal{G}^{\mathcal{T}}$-a monotonic chain decomposition,
   $\theta^{\mathcal{T}}$ - initial cost distribution,
   $N$ - number of iterations,
   $B_u^{\tau}(s) = F_u^{\tau}(s) := 0$ for $u \in \mathcal{V}$, $s \in \mathcal{Y}_u$, $\tau \in \mathcal{T}$
2: **for** $t = 1$ **to** $N$ **do**
3:    **for** $u = 1$ **to** $|\mathcal{V}|$ **do**
4:       **for** $\tau \in \mathcal{T}_u$ **do**
5:          **if** exists $v \in \mathcal{N}_b(u) : uv \in \mathcal{V}^{\tau}, u > v$ - previous node for $u$
             in $\mathcal{V}^{\tau}$ **then**
6:             Compute $F_u^{\tau}(s) = \min_{l \in \mathcal{Y}_v} \left( F_v^{\tau}(l) + \theta_v^{\tau}(l) + \theta_{uv}^{\tau}(s, l) \right)$ for
                all $s \in \mathcal{Y}_u$
7:             Compute $E_u^{\tau}(s) := F_u^{\tau}(s) + B_u^{\tau}(s) + \theta_u^{\tau}(s)$ for all $s \in \mathcal{Y}_i$
8:          **end if**
9:       **end for**
10:      $\theta_u^{\tau}(s) := \frac{1}{|\mathcal{T}_u|} \sum_{\sigma \in \mathcal{T}_u} E_u^{\sigma}(s) - F_u^{\tau}(s) - B_u^{\tau}(s)$, $s \in \mathcal{Y}_u$, $\tau \in \mathcal{T}_u$
11:   **end for**
12:   Revert node order, and swap $F_u^{\tau} \leftrightarrow B_u^{\tau}$
13: **end for**
14: **return** $\theta^{\mathcal{T}}$

---

time as one iteration of Algorithm 12 itself; (ii) the first iteration usually *does increase* the dual objective, in spite of the zero initialization of the backward min-marginals.

Algorithm 12 already has the required iteration complexity $O(|\mathcal{E}|L^2)$: this is the complexity of computing forward $F_u^\tau(s)$ and backward $B_u^\tau(s)$ marginals performed in line 4. These values are computed for each node $u$ and all slave subproblems contained this node. Since each edge belong to exactly one slave subproblem, the number of subproblems containing a node $u$ is at most the number of the edges incident to $u$. Other operations have complexity $O(|\mathcal{E}|L)$ at most.

Its actual performance, i.e. speed of increasing the dual bound, may significantly depend on the ordering of the nodes and a monotonic decomposition used. However, little is known about the optimal selection of those parameters. We will assume the ordering to be given and will not discuss its selection.

In the following section, we will consider a particular type of decompositions, which have shown good performance in a number of applications. Additionally, we will show that for these decompositions Algorithm 12 becomes particularly simple and specializes to a particularly simple form of the anisotropic min-sum diffusion algorithm introduced in §8.3.2.

### 10.2.3  Maximal monotonic chains, TRW-S and SRMP algorithms

Note that Algorithm 12 requires an explicit definition of the graph decomposition $\mathcal{G}^\mathcal{T}$ to compute corresponding min-marginals $F_u^\tau(s)$. Below, we will show that an explicit definition of the decomposition may not actually be required. We will see that the only numbers, which are important for Algorithm 12, are $|\mathcal{T}_u|$, defining the number of chains containing node $u$. This, in particular, allows for a much simpler implementation of the algorithm and shows equivalence of all decompositions having equal numbers $|\mathcal{T}_u|$.

Let us introduce the notation

$$\hat{\theta}_u(s) := \sum_{\tau \in \mathcal{T}_u} E_u^\tau(s) = \sum_{\tau \in \mathcal{T}_u} \left( F_u^\tau(s) + B_u^\tau(s) + \theta_u^\tau(s) \right)$$
$$= \theta_u(s) + \sum_{\tau \in \mathcal{T}_u} \left( F_u^\tau(s) + B_u^\tau(s) \right) \quad (10.17)$$

Now we can rewrite the averaging step 10 in Algorithm 12 as

$$\theta_u^\tau(s) := \frac{1}{|\mathcal{T}_u|} \hat{\theta}_u(s) - F_u^\tau(s) - B_u^\tau(s) . \quad (10.18)$$

Let us now express $F_v^\tau(l)$ as a function of $F_u^\tau(s)$ for $v$ being the next node after $u$ (which means, in particular, that $v > u$) in the chain indexed by $\tau$:

$$F_v^\tau(l) = \min_{s \in \mathcal{Y}_u} \left( F_u^\tau(s) + \theta_u^\tau(s) + \theta_{uv}^\tau(s, l) \right)$$
$$\overset{(10.18)}{=} \min_{s \in \mathcal{Y}_u} \left( F_u^\tau(s) + \frac{1}{|\mathcal{T}_u|} \hat{\theta}_u(s) - F_u^\tau(s) - B_u^\tau(s) + \theta_{uv}^\tau(s, l) \right)$$
$$= \min_{s \in \mathcal{Y}_u} \left( \frac{1}{|\mathcal{T}_u|} \hat{\theta}_u(s) - B_u^\tau(s) + \theta_{uv}^\tau(s, l) \right) \quad (10.19)$$

Note that, since we consider canonical decompositions, each edge $uv$ has exactly one forward- and one backward-marginal vector, associated with it, namely $F_u^\tau$ and $B_v^\tau$ such that $uv \in \mathcal{E}^\tau$. Therefore, we introduce an alternative notation for these vectors independent of $\tau$: $F_{v,u}$ and $B_{u,v}$, respectively. Here comma between $v$ and $u$ shows the "direction" of computation: For forward marginals the notation $F_{v,u}$ implies that $u > v$ and the other way around, for backward marginals $B_{v,u}$ it implies that $u < v$.

Applying this new notation to (10.17) we obtain:

$$\hat{\theta}_u(s) = \theta_u(s) + \sum_{\substack{v \in \mathcal{N}_b(u) \\ v < u}} F_{v,u}(s) + \sum_{\substack{v \in \mathcal{N}_b(u) \\ v > u}} B_{v,u}(s) . \quad (10.20)$$

Expression (10.19) transforms to

$$F_{u,v}(l) = \min_{s \in \mathcal{Y}_u} \left( \frac{1}{|\mathcal{T}_u|} \hat{\theta}_u(s) - B_{v,u}(s) + \theta_{uv}(s, l) \right) , \quad (10.21)$$

where we omit the index $\tau$ in pairwise costs $\theta_{uv}^\tau(s, l)$, since $\theta_{uv}^\tau(s, l) = \theta_{uv}(s, l)$ as we consider only canonical decompositions.

Now we can rewrite Algorithm 12 as Algorithm 13 by using expressions (10.20) and (10.21).

---

**Algorithm 13** Modified block-coordinate ascent for a monotonic chain decomposition

---

 1: **Init:** Numbers $|\mathcal{T}_u|$ for $u \in \mathcal{V}$,
    $F_{u,v}(l) = B_{v,u}(s) := 0$ for $uv \in \mathcal{E}$, $l \in \mathcal{Y}_v$, $s \in \mathcal{Y}_u$
    $N$ - number of iterations
 2: **for** $t = 1$ **to** $N$ **do**
 3: 　　**for** $u = 1$ **to** $|\mathcal{V}|$ **do**
 4: 　　　　Compute $\hat{\theta}_u(s) = \theta_u(s) + \sum_{\substack{v \in \mathcal{N}_b(u) \\ v<u}} F_{v,u}(s) + \sum_{\substack{v \in \mathcal{N}_b(u) \\ v>u}} B_{v,u}(s)$
    　　　　for all $s \in \mathcal{Y}_u$
 5: 　　　　Compute $F_{u,v}(l) = \min_{s \in \mathcal{Y}_u} \left( \frac{1}{|\mathcal{T}_u|} \hat{\theta}_u(s) - B_{v,u}(s) + \theta_{uv}(s,l) \right)$
    　　　　for all $v \in \mathcal{N}_b(u)$ such that $v > u$ and $l \in \mathcal{Y}_v$
 6: 　　**end for**
 7: 　　Revert node order, and swap $F_{u,v} \leftrightarrow B_{u,v}$ for all $uv \in \mathcal{E}$
 8: **end for**
 9: **return** $F_{u,v}(l)$, $B_{v,u}(s)$ for all $uv \in \mathcal{E}$, $l \in \mathcal{Y}_v$, $s \in \mathcal{Y}_u$

---

Note that Algorithm 13 depends on a particular decomposition only through the numbers $|\mathcal{T}_u|$, i.e. all canonical monotonic decompositions with the same $|\mathcal{T}_u|$ result in exactly the same algorithm.

**Example 10.9** (Edge decomposition). For the edge decomposition the value $|\mathcal{T}_u|$ is equal to the number of incident edges for a node $u$.

**Example 10.10** (Maximal monotonic chains). Maximal monotonic chains is a straightforward way to construct a canonical decomposition consisting of monotonic chains that are as long as possible. Given a graph, the chains are constructed in a greedy manner as follows. First, starting with the first node w.r.t. the introduced order, a monotonic chain is greedily generated step-by-step until it is impossible to extend it, i.e. for the current node $u$ there is no neighbor $v$ such that $v > u$. After that, the edges of the constructed chain are removed from the graph $\mathcal{G}$ and the procedure repeats from the node with the smallest index that still has incident edges.

Although the decomposition into maximal monotonic chains is not unique, all such decompositions are characterized by the same numbers

$$|\mathcal{T}_u| = \max\{|v \in \mathcal{N}_b(u)\colon v < u|, |v \in \mathcal{N}_b(u)\colon v > u|\}, \qquad (10.22)$$

which are equal to the maximum of the number of incoming and outgoing edges for node $u$, respectively. Here 'incoming' and 'outgoing' are considered with respect to the node ordering of the graph.

A variant of Algorithm 13 with $|\mathcal{T}_u|$ corresponding to maximal monotonic chains is known as *Sequential Tree-Reweighted Message Passing (TRW-S)*. Nowadays it is one of the most efficient MAP-inference algorithms for graphical models.

**Computation of dual value. Relation to anisotropic diffusion**  It remains to show how the dual variables $\phi$ can be computed from the forward and backward min-marginals returned by Algorithm 13. The answer is based on generalizing the result of §8.3.2, where we have shown that dynamic programming can be seen as anisotropic diffusion. Recall Algorithm 6 and the comments on it. The expression (8.27) translates to

$$\text{for } v > u: \qquad\qquad\qquad\qquad\qquad\qquad (10.23)$$
$$\phi_{v,u}(l) = -F_{u,v}(l), \ l \in \mathcal{Y}_v, \qquad\qquad\qquad (10.24)$$
$$F_{u,v}(l) = \min_{s \in \mathcal{Y}_u} (\theta_{uv}(s,l) + \phi_{u,v}(s)), \ s \in \mathcal{Y}_u. \qquad (10.25)$$

Comparing the last formula with line 5 of Algorithm 13 we conclude:

$$\text{for } v > u:$$
$$\phi_{v,u}(l) = -F_{u,v}(l), \ l \in \mathcal{Y}_v, \qquad\qquad (10.26)$$
$$\phi_{u,v}(s) = \frac{1}{|\mathcal{T}_u|}\hat{\theta}_u(s) - B_{v,u}(s), \ s \in \mathcal{Y}_u.$$

The values $\hat{\theta}_u(s)$ can be computed from forward and backward min-marginals as in line 4 of Algorithm 13.

Taking (10.26) into account, Algorithm 13 can be solely expressed in terms of the reparametrization. In this form it transforms to Algorithm 14.

---

**Algorithm 14** Sequential reweighted message passing (SRMP)

---

1: **Init:** Numbers $|\mathcal{T}_u|$ for $u \in \mathcal{V}$,
   $\phi$ - starting reparametrization (typically $\phi = 0$)
   $N$ - number of iterations
2: **for** $t = 1$ **to** $N$ **do**
3:    **for** $u = 1$ **to** $|\mathcal{V}|$ **do**
4:       Move pairwise costs to unary costs:
         $\phi_{u,v}(s) := \phi_{u,v}(s) - \min_{l \in \mathcal{Y}_v} \theta_{uv}^{\phi}(s, l)$ for $s \in \mathcal{Y}_u$ and $v \in \mathcal{N}_b(u)$

5:       Redistribute unary costs to the outgoing edges:
         $\phi_{u,v}(s) := \phi_{u,v}(s) + \frac{\theta_u^{\phi}(s)}{|\mathcal{T}_u|}$ for $s \in \mathcal{Y}_u$ and $v > u$
6:    **end for**
7:    Revert node order.
8: **end for**
9: **return** $\phi$

---

It is easy to see that Algorithm 14 is a special instance of anisotropic diffusion. Its most efficient variant corresponds to the decomposition into maximal monotonic chains, i.e. $|\mathcal{T}_u| = \max\{|v \in \mathcal{N}_b(u) \colon v < u|, |v \in \mathcal{N}_b(u) \colon v > u|\}$, see (10.22). This variant of the algorithm is known as *Sequential Reweighted Message Passing (SRMP)*. For the pairwise models considered in this book, SRMP is equivalent to TRW-S and differs only in the numbers, which are computed during its run. In TRW-S the forward and backward min-marginals are computed, whereas in SRMP the dual variables are computed directly. The equivalence of the Algorithms follows from the fact that the computed values are related by (10.26).

Similarly to min-sum diffusion, the TRW-S/SRMP algorithm converges to node-edge agreement. This implies also convergence to tree agreement due to Proposition 9.12. We omit the proof here and refer to the corresponding literature [49, 51].

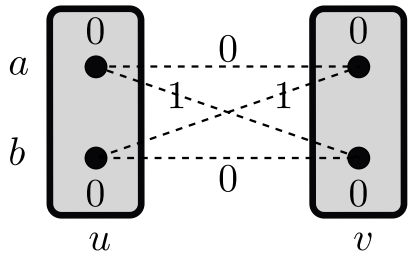

**Figure 10.3:** Depicted is a graphical model with two nodes and two labels, where the numbers stand for unary and pairwise costs. This model is a counterexample for the rounding technique for the TRW-S/SRMP algorithm where the minimal label $y_u$ with respect to the reparametrized cost $\theta_u^\phi$ is selected in line 4 of Algorithm 14. The reparametrized costs in this case are the same as the original costs. This implies that either the label $a$ or the label $b$ can be selected for each node. Since these decisions are taken independently, it may happen that the labeling $(a, b)$ is selected as the rounding result, which is not optimal. Using the relaxation labeling algorithm would not help to resolve the problem, as the locally optimal labels and label pairs are already arc-consistent and form the arc-consistency closure.

### 10.2.4 Rounding technique for the TRW-S/SRMP algorithm

Similarly to the case of the min-sum diffusion algorithm, naïve rounding may not perform well if applied to the result of the TRW-S/SRMP algorithm directly. The reason is the same – for all nodes, where the numbers of incoming and outgoing edges coincide, the unary cost vector contains only zero values. For example, in a grid-graph with column- or row-wise ordered nodes and a monotonic chain decomposition all non-border nodes, i.e. nodes $(i, j)$ with $1 < i < h$ and $1 < j < w$ (with $h$ and $w$ the total number of rows and columns in the graph).

The simplest rounding next to naïve rounding is to remember during each iteration of the algorithm the locally optimal label in each node $u$ after the cost vector $\theta_u^\phi$ is computed in line 4 of Algorithm 14 (or, analogously, the cost vector $\hat{\theta}_u(s)$ of Algorithm 13).

However, this method does not guarantee obtaining the optimal solution even in case the graph $\mathcal{G}$ is acyclic. Indeed, Figure 10.3 shows this problem appearing already for a graph with two nodes and two labels.

Therefore, another approach was proposed in the original work of [49]. It generalizes the way how an optimal solution is obtained with dynamic programming and, therefore, guarantees obtaining an optimal labeling for acyclic graphs after the first iteration of the algorithm.

Consider the dynamic programming Algorithm 1. The pointers $r_i(s)$ are computed for each node $i$ and each label $s$ of the graph to be able to reconstruct an optimal labeling. However, this reconstruction can also be performed without saving the pointers $r_i(s)$, although with slightly higher computational cost.

Indeed, according to Algorithm 1

$$r_i(s) := \arg\min_{t \in \mathcal{Y}_{i-1}} \left( F_{i-1}(t) + \theta_{i-1}(t) + \theta_{i-1,i}(t,s) \right) . \tag{10.27}$$

This computation can be done on the backward move of the dynamic programming algorithm when label $s$ is already known. This is summarized in Algorithm 15, which is a straightforward modification of Algorithm 2.

---

**Algorithm 15** Reconstructing an optimal labeling for a chain graph

---

1: $y_n = \arg\min_{s \in \mathcal{Y}_n}(F_n(s) + \theta_n(s))$
2: **for** $i = n-1$ **to** 1 **do**
3:     $y_i = \arg\min_{t \in \mathcal{Y}_i} (F_i(t) + \theta_i(t) + \theta_{i,i+1}(t, y_{i+1}))$
4: **end for**
5: **return** $y$

---

The primal rounding for TRW-S/SRMP given by Algorithm 16 is a straightforward generalization (with $\phi_{u,v}(s) = -F_u^{uv}(s)$) of Algorithm 15 to arbitrary graphs with the introduced ordering of nodes. Note that in line 4 of Algorithm 16 the original, non-reparametrized costs are used.

## 10.3   Empirical comparison of algorithms

**Algorithms**   Similarly to §8.7, we provide an experimental evaluation of the algorithms described in this chapter. Additionally, we compare them to the methods considered in Chapter 8. Altogether, we consider the following methods:

- We consider two variants of the min-sum `diffusion` algorithm, `diffusion+naïve` and `diffusion+ICM`, for generating primal solutions as described in §8.7;

---

**Algorithm 16** Reconstructing an optimal labeling for the SRMP algorithm

---

1: $u := |\mathcal{V}|$;
2: $y_u = \arg\min_{s \in \mathcal{Y}_{|\mathcal{V}|}} (\theta_u(s) - \sum_{v \in \mathcal{N}_b(u)} \phi_{u,v}(s))$
3: **for** $u = |\mathcal{V}| - 1$ **to** $1$ **do**
4: $\quad y_u = \arg\min_{s \in \mathcal{Y}_i} \left( \theta_u(s) + \sum_{\substack{v \in \mathcal{N}_b(u) \\ v > u}} \theta_{u,v}(s, y_v) - \sum_{\substack{v \in \mathcal{N}_b(u) \\ v < u}} \phi_{u,v}(s) \right)$
5: **end for**
6: **return** $y$

---

- We also look at the subgradient algorithm represented by two variants. The first one is `subgradient`, which is the same Lagrange dual-based algorithm as considered in §8.7. As follows from Example 9.4, it can be seen as a subgradient method based on the complete decomposition. The second variant, `max.chain subgradient` is based on the decomposition into maximal monotonic chains.

- The decomposition-based block coordinate ascent is represented by two variants of Algorithm 13. The first one, termed as `TRWS-edge`, is the specialization of the algorithm for the edge decomposition, considered in §10.2.1. The second one is the `TRWS` algorithm with the decomposition into maximal monotonic chains. By default, this algorithm reconstructs the primal integer solution as described in §10.2.4. Its variant `TRWS+ICM` additionally runs one iteration of the ICM algorithm to improve the primal solution.

Additionally, we plot also the value of the optimal solution of the non-relaxed MAP-inference problem obtained with the modern exact inference technique [33].

**Datasets** The problem instances used for comparison are the same as in §8.7.

**Dual plots and their analysis** Figures 10.4 and 10.5 compare values of the dual bounds obtained by the considered algorithms. We provide only dual bounds, since any particular rounding technique typically yields the lower energies for the estimated primal solutions the higher the lower bound.

TRWS consistently outperforms all competitors on all problem instances and is typically followed by TRWS-edge and diffusion. A notable exception is the stereo problem, where the second best algorithm turns out to be max.chain subgradient. Note that only TRWS and max.chain subgradient are based on a maximal monotonic chains decomposition, all other methods are based on the complete decomposition or the edge decomposition. We explain the high performance of the maximal chain decomposition-based methods by their more efficient propagation of the information across the graph. This is particularly important for sparse graphical models, which virtually contain "long-range" dependencies between distant nodes of the graph due to strong pairwise costs. Among our problem instances only stereo has such properties, since protein and worms are densely-connected and color segmentation has very weak pairwise costs.

Interestingly, diffusion is always outperformed by the nearly identical TRWS-edge. The difference between both algorithms is the order in which edges of the graph are processed. Diffusion always starts from the first node and proceeds to the last one, whereas TRWS-edge inverts the order on each iteration. Apparently, changing the order allows for better information propagation in both directions.

The difference between block-coordinate-ascent and subgradient algorithms is more pronounced for models with a larger label set, namely protein and worms, which is explained by the fact that the subgradient is more sparse for such models.

**Primal-dual plots and their analysis**    Figures 10.6 and 10.7 show both primal and dual bounds obtained by the considered algorithms.

As in the comparison in Chapter 8, even one ICM iteration may significantly improve the primal estimates. It can also be seen that the ICM rounding brings more improvement to simple rounding techniques like those used in the diffusion algorithm. At the same time, it only slightly improves the results of the more sophisticated and theory-based rounding used in the TRWS algorithm. Noticeably, the ICM iteration in combination with the slower diffusion algorithm leads to better primal estimates than the `TRWS+ICM` rounding for the `protein` dataset. We attribute this to the fact that even the stand-alone `ICM` algorithm performs pretty well for this dataset - see Figure 8.9. This may happen because of the fact that `TRWS` primal estimates are close to the "local optimum" of the `ICM` algorithm, whereas the `diffusion` ones are not.

**Conclusions**   The above comparisons show that a decomposition into larger subproblems (e.g. a maximal chain decomposition), as well as the block-coordinate updates are important components of a well performing dual solver. The importance of the former increases as the graph gets more sparse. The latter becomes especially pronounced for large label sets.

## 10.4   Bibliography and further reading

The subgradient technique for the Lagrange decomposition of the MAP-inference problem in graphical models was first applied in [119]. Later on, independently, similar work was done in [101] and [55], followed by [56].

   A number of non-smooth convex optimization techniques have been developed that generalize the subgradient method. Many of them have been tested on the MAP-inference problem for graphical models. These include, but not limited to the mirror descent [72], bundle methods [40], proximal algorithms [88, 105, 75, 73, 26] and the smoothing technique [38, 22, 96, 94]. Although most of these methods guarantee convergence to the dual optimum, in application to the MAP-inference problem they are typically inferior to the best block-coordinate ascent techniques (see the comparison [41]).

   The original TRW-S algorithm was proposed by [48], followed by a journal publication [49]. The generalization of the TRW-S algorithm for higher-order graphical models was proposed in [106]. Although the resulting algorithm was reported to be quite efficient, it is much more complicated than the original TRW-S. Later on, a novel view on TRW-S as anisotropic diffusion was proposed in [51]. This allowed another, elegant and efficient generalization of the method to higher-order models.

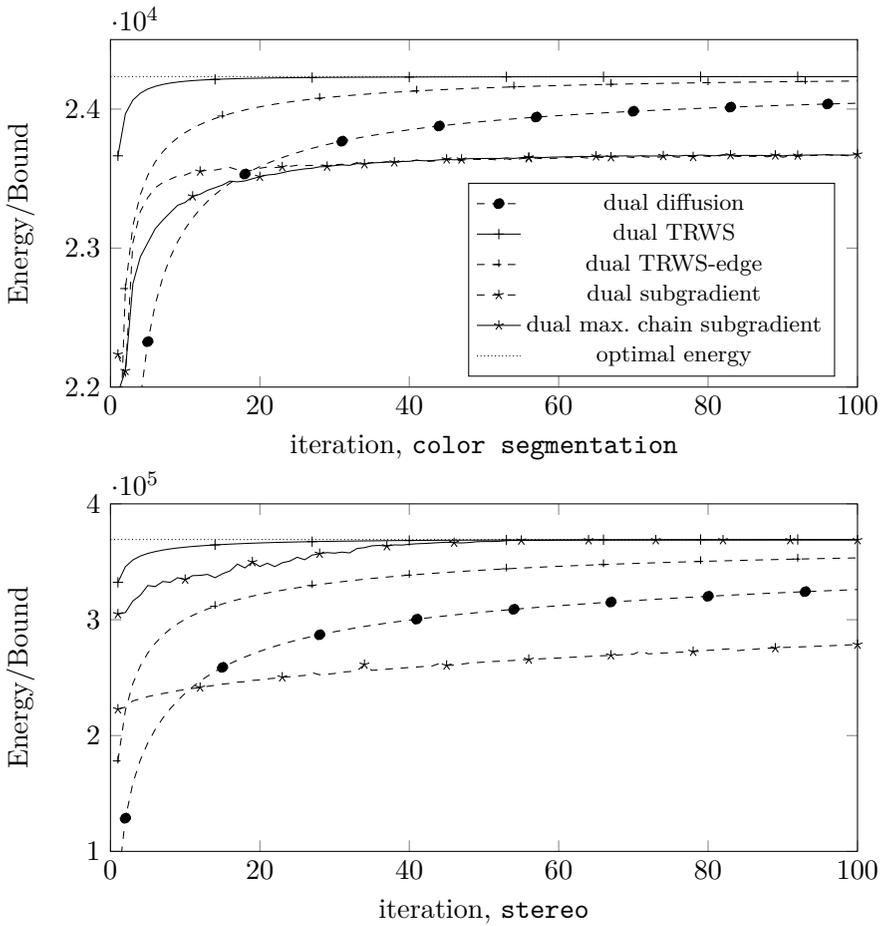

**Figure 10.4:** Dual plots for the **(top)** `color segmentation` and **(bottom)** `stereo` problem instances. The legend for the bottom plot is the same as for the top one. Note that on the `stereo` problem instance `max. chain subgradient` outperforms all other methods except `TRWS`. See main text for the analysis.

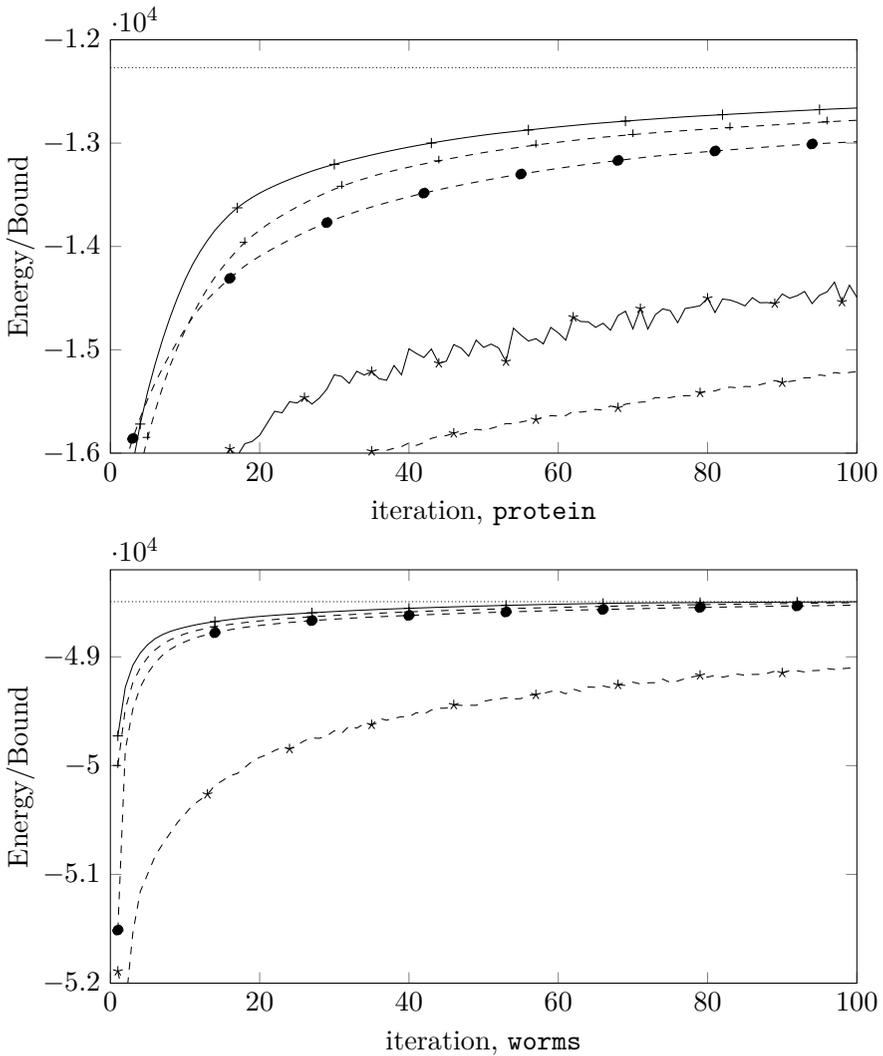

**Figure 10.5:** Dual plots for the **(top)** protein and **(bottom)** worms problem instances. The legend for the plots is the same as in the top plot of Figure 10.4. See main text for the analysis.

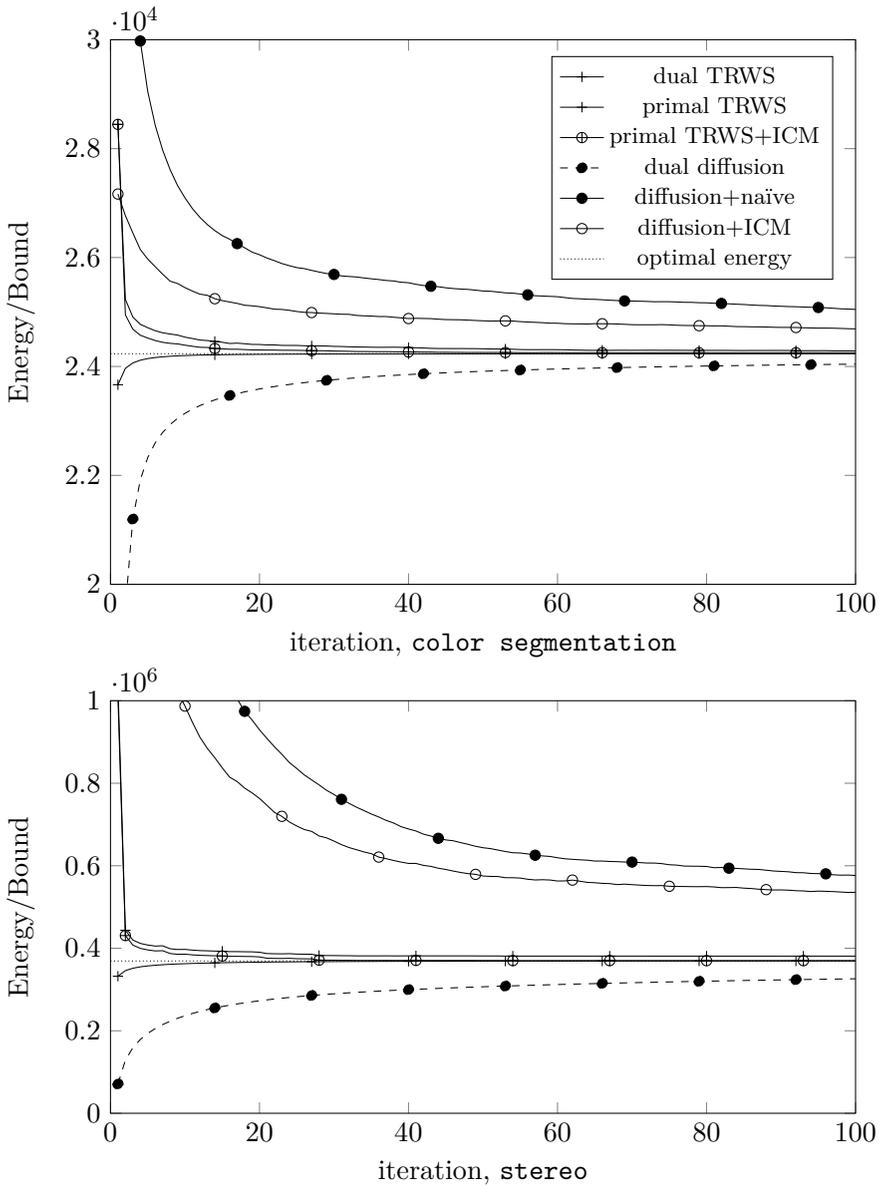

**Figure 10.6:** Primal-dual plots for the **(top)** `color segmentation` and **(bottom)** `stereo` problem instances. The legend for the bottom plot is the same as for the top one. Note that adding just one `ICM` iteration significantly improves the primal estimates of the `diffusion`, whereas the more sophisticated rounding of the `TRWS` algorithm can be only marginally improved. See the main text for the analysis.

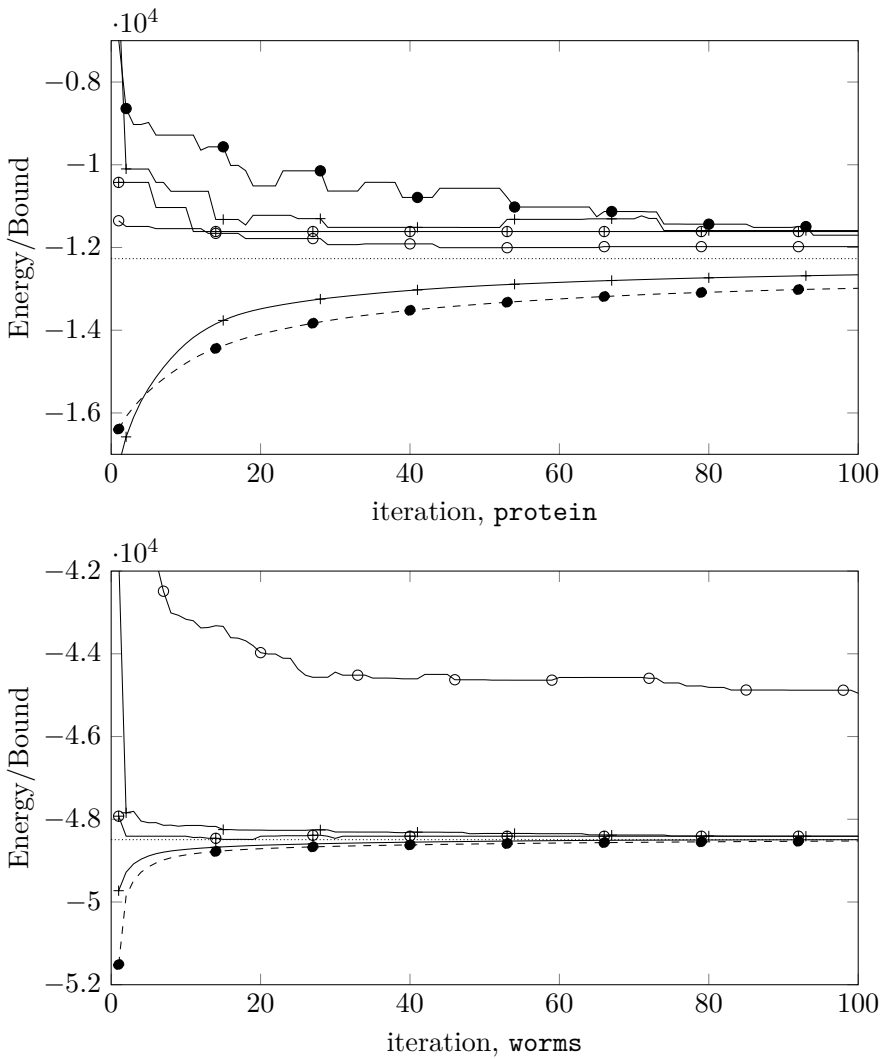

**Figure 10.7:** Primal-dual plots for the **(top)** protein and **(bottom)** worms problem instances. The legend for the plots is the same as in the **(top)** plot of Figure 10.6. Primal estimates of diffusion+naïve are not shown for worms, since they are very high because of infinite values presented in pairwise costs. This is significantly improved with the rounding used in diffusion+ICM. Note also that for protein the primal estimates of diffusion+ICM are better than those of TRWS+ICM. See the main text for the analysis.

The generalized algorithm obtained the name *Sequential Reweighted Message Passing (SRMP)*.

In this chapter we considered the block-coordinate ascent method, which optimizes over the dual variables associated with a given node at a time. As an alternative approach, one could optimize over the dual variables associated with a whole (acyclic) subproblem, one subproblem during one step of the algorithm. The algorithms suggested in early works in this direction [117, 137, 142] typically cannot compete with TRW-S and SRMP. However, a very efficient method based on a similar idea was proposed recently in [112]. It seems that this algorithm outperforms TRW-S, although a vast experimental comparison is lacking. Importantly, for grid graphs a parallelization of the algorithm was proposed in the same work. It shows a very high performance on the stereo problem instances used for testing.

A dual block-coordinate ascent framework that generalizes TRW-S, max-sum diffusion and MPLP methods to a fairly broad class of combinatorial optimization problems was recently proposed in [122]. Encouraging results using this framework have been obtained for the quadratic assignment [123] and the multicut problems [121].

# 11

## Min-Cut/Max-Flow Based Inference

So far, we considered a single class of graphical models with exactly solvable MAP-inference problem, that is, the acyclic problems introduced in Chapter 2. In this case the MAP-inference is exactly solvable due to the special property of the local polytope, which contains only integer vertexes.

However, there is another way to enforce integrality. Instead of restricting the graph structure, one can restrict the cost vectors such that they point to integer vertexes of the local polytope. A large class of such costs are called *submodular*, since the MAP-inference can be treated as *submodular minimization* in this case, a well-studied polynomially solvable type of problems in combinatorial optimization.

The chapter consists of three parts. We start with a general definition of submodular functions. Then we give a detailed treatment of a subclass of such functions, corresponding to the pairwise graphical model energies. In the second part we will show equivalence of the submodular MAP-inference and a polynomially solvable *min-st-cut* problem. Finally, the third part is devoted to one of the most efficient type of primal heuristics, commonly known as *graph cuts*. Although this type of algorithms is related to the submodularity property, it is also applicable to a large class of non-submodular problems.

## 11.1   Submodular functions

Let $X$ be a finite set and $2^X$ be its powerset. Mappings of the form $f \colon 2^X \to \mathbb{R}$ are called *set-functions*. Finite sets can be represented by binary indicator vectors with coordinates indexed by elements of $X$. A binary vector $\xi \in \{0,1\}^X$ corresponds to the set $A \subset X$ if for all $x \in X$ it holds that

$$\xi_x = \begin{cases} 0, & x \notin A, \\ 1, & x \in A. \end{cases} \tag{11.1}$$

**Example 11.1.** Let $(\mathcal{G}, \mathcal{Y}_\mathcal{V}, \theta)$ be a binary graphical model, i.e. $\mathcal{Y}_\mathcal{V} = \{0,1\}^\mathcal{V}$. Then the energy function

$$E(y) = \sum_{u \in \mathcal{V}} \theta_u(y_u) + \sum_{uv \in \mathcal{E}} \theta_{uv}(y_u, y_v) \tag{11.2}$$

can be seen as a set-function w.r.t. $y$, since each labeling $y \in \mathcal{Y}_\mathcal{V}$ is a binary vector.

**Definition 11.2.** A set-function $f \colon 2^X \to \mathbb{R}$ is called *submodular* if

$$f(A) + f(B) \geq f(A \cap B) + f(A \cup B) \quad \forall A, B \in 2^X. \tag{11.3}$$

Expressing this definition in terms of the corresponding graphical model energy, we obtain

$$E(y) + E(y') \geq E(y \wedge y') + E(y \vee y') \quad \forall y, y' \in \mathcal{Y}_\mathcal{V} \tag{11.4}$$

with logical *and* $\wedge$ and *or* $\vee$ operations defined as in Chapter 6, see Figure 11.1 for illustration.

Similarly to convex/concave functions, a function $f$ is called *supermodular*, if $(-f)$ is submodular.

The *submodular minimization problem* consists in computing

$$A^* = \arg \min_{A \in 2^X} f(A), \tag{11.5}$$

where the set-function $f$ is submodular. This problem is known to be polynomially solvable (assuming that the value of $f(A)$ can be computed in polynomial time for any $A \in X$), although the order of the polynomial is quite high in general. However, different subclasses of submodular functions have lower computational complexity. One of such sub-classes, namely, submodular pairwise energy minimization problems, we will consider below.

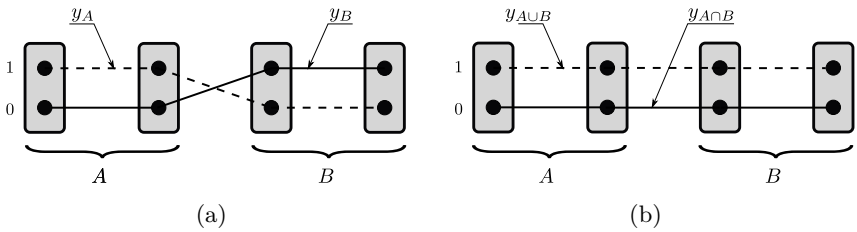

**Figure 11.1:** Representation of sets intersection and union with binary labelings. The sets $A$, $B$, $A \cup B$ and $A \cap B$ are represented with binary labelings $y_A$, $y_B$, $y_{A \cup B} = y_A \vee y_B$ and $y_{A \cap B} = y_A \wedge y_B$ respectively.

### 11.1.1 Submodularity on a lattice

Submodularity can also be defined for multilabel graphical model energies. To do so, we assume that the set of labels in each node is totally ordered, e.g. the relations $\leq$ and $\geq$ are naturally defined for any two labels $s, l \in \mathcal{Y}_u$, $u \in \mathcal{V}$. The set of labelings $\mathcal{Y}_\mathcal{V}$ in this case is partially ordered, where the operation $\geq$ ($\leq$) is defined point-wise, which is $y \leq y'$ for $y, y' \in \mathcal{Y}_\mathcal{V}$, if $y_u \leq y'_u$ for all $u \in \mathcal{V}$, see Figure 11.2. Moreover, the set of labelings $\mathcal{Y}_\mathcal{V}$ is a *lattice*, since it has its *supremum* and *infimum*, i.e. the *highest* labeling $\check{y}$ such that $\check{y} \geq y$ for all $y \in \mathcal{Y}_\mathcal{V}$ and the *lowest* one $\hat{y}$ such that $\hat{y} \leq y$ for all $y \in \mathcal{Y}_\mathcal{V}$.

In general, a non-empty partially ordered set $A$ equipped with operations $\vee$ and $\wedge$ is called *a lattice* if $x \wedge z$ and $x \vee z$ are defined for all $x, z \in A$.

For every two labels $s, l$ from the same node their maximum $\vee$ and minimum $\wedge$ are naturally defined, i.e. if $s \geq l$, then $s \vee l = s$ and $s \wedge l = l$. This generalizes to any two labelings $y, y'$, where operations $\wedge$ and $\vee$ are applied point-wise, see Figure 11.2.

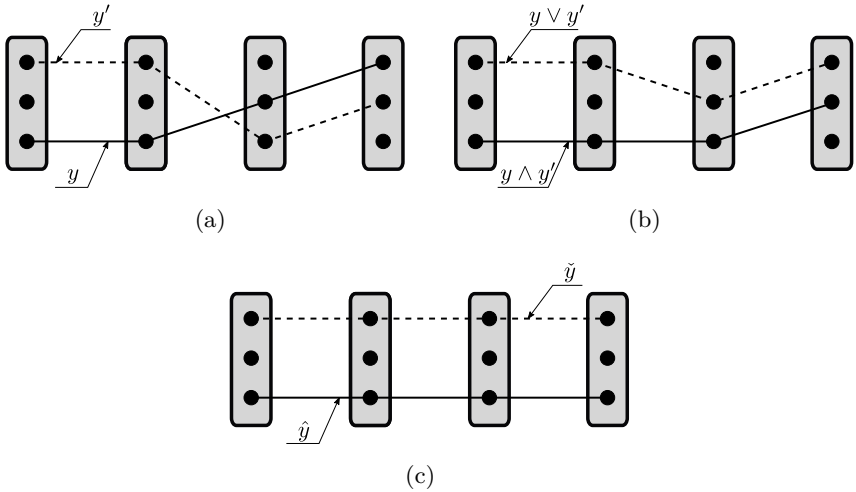

**Figure 11.2:** The multilabel energy as a special case of submodularity on a lattice. Labels are totally ordered in each node. **(a)** Labelings $y$ and $y'$ are incomparable, since neither $y'_u \le y_u$ nor $y_u \le y'_u$ holds for all nodes $u$. **(a)-(b)** For any two labelings $y$ and $y'$ their node-wise maximum $y \vee y'$ and minimum $y \wedge y'$ are always well-defined. It always holds $(y \vee y') \ge y \wedge y'$ **(c)** The supremum and infimum of the set of all labelings.

A classical example of a lattice is the *set lattice*, i.e. a family of subsets of a given set with the "subset" relation $\subseteq$ in place of $\le$ and closed w.r.t. the operations $\cap$ and $\cup$ standing for $\wedge$ and $\vee$. A lattice is called *finite*, if the underlying partially ordered set is finite.

Moreover, it is known (see [24]) that any finite lattice is isomorphic to some finite set lattice. This substantiates the generalized definition of submodularity:

**Definition 11.3.** Let $(X, \wedge, \vee)$ be a finite lattice. A function $f \colon X \to \mathbb{R}$ is called *submodular*, if

$$f(x) + f(z) \ge f(x \wedge z) + f(x \vee z) \tag{11.6}$$

holds for any $x, z \in X$.

For the supermodular function the inequality holds with the opposite sign $\le$.

Translated to the energy function of a graphical model with a lattice structure the submodularity requires that

$$E(y) + E(y') \geq E(y \wedge y') + E(y \vee y') \tag{11.7}$$

holds for any two labelings $y, y' \in \mathcal{Y}_\mathcal{V}$.

Although this is a valid definition, it can not be checked explicitly in practice, as it would require evaluating the relation (11.7) for all possible pairs of labelings. Fortunately, there is no need to do so, due to the following remarkable fact:

**Theorem 11.1** ([100]). *Let $(\mathcal{G}, \mathcal{Y}_\mathcal{V}, \theta)$ be a graphical model and $E$ be its energy. Let also label sets corresponding to the nodes of the graphical model be totally ordered and operations $\vee$ and $\wedge$ be defined as above. Then $E$ is submodular if and only if all its pairwise cost functions $\theta_{uv} \colon \mathcal{Y}_u \times \mathcal{Y}_v \to \mathbb{R}$, $uv \in \mathcal{E}$, are submodular.*

We will get back to this theorem and its proof in §11.2.2 after considering the class of submodular functions of two variables in §11.2.

**Remark 11.4.** Submodular functions on a lattices $(W, \min, \max)$ with $W = \{1, \ldots, n_1\} \times \cdots \times \{1, \ldots, n_d\}$ and min and max applied coordinate-wise are equivalent to *Monge arrays*. The latter have a number of applications in optimization, see e.g. the overview of [18].

### 11.1.2 General properties of submodular functions

**Definition 11.5.** Let $(X, \wedge, \vee)$ be a finite lattice. A function $f \colon X \to \mathbb{R}$ is called *modular*, if it is both sub- and supermodular, i.e.:

$$f(x) + f(x') = f(x \vee x') + f(x \wedge x'). \tag{11.8}$$

**Lemma 11.6.** *Let $X$ be a totally ordered set. Then any function $f \colon X \to \mathbb{R}$ defined on this set is modular.*

*Proof.* Consider any $x, x' \in X$. W.l.o.g. assume $x \leq x'$. Then $x \wedge x' = x$ and $x \vee x' = x'$. Therefore $f(x \wedge x') + f(x \vee x') = f(x) + f(x')$, which finalizes the proof. $\square$

**Lemma 11.7.** *Let $f$ and $g$ be (sub/super)modular and $\alpha, \beta \geq 0$. Then $\alpha f + \beta g$ is (sub/super)modular.*

The proof follows trivially from the definition of (sub/super)modularity.

In particular, Lemma 11.7 implies, that adding a modular function to a submodular one results in a submodular function.

## 11.2   Submodular pairwise energies

### 11.2.1   Submodular functions of two variables

Let us consider submodular functions of only two variables. In particular, we are interested when a pairwise cost function $\theta_{uv}$ is submodular.

We start with the definition (11.7), which can be reformulated as follows:

**Definition 11.8.** The pairwise cost function $\theta_{uv}$ is *submodular* if for all $s^1, s^2 \in \mathcal{Y}_u$ and all $l^1, l^2 \in \mathcal{Y}_v$ such that $s^1 \leq s^2$ and $l^1 \leq l^2$, it holds that

$$\theta_{uv}(s^1, l^1) + \theta_{uv}(s^2, l^2) \leq \theta_{uv}(s^2, l^1) + \theta_{uv}(s^1, l^2). \tag{11.9}$$

Submodular functions of two variables are also known in the literature as *Monge matrices*, see e.g. the overview of [18].

Figure 11.3(a) illustrates Definition 11.8. Below we consider several important examples.

**Example 11.9** (Binary problem)**.** In the binary case, when $\mathcal{Y}_u = \{0, 1\}$ for all $u \in \mathcal{V}$, the submodularity condition (11.9) simplifies to

$$\theta_{uv}(0, 1) + \theta_{uv}(1, 0) \geq \theta_{uv}(0, 0) + \theta_{uv}(1, 1). \tag{11.10}$$

It is also easy to see that in this case each factor is either submodular or supermodular, see Figure 11.3(b) for illustration.

**Example 11.10** (Ising model)**.** Let the energy function be binary and $\theta_{uv}(s, t) = \lambda_{uv} [\![ s \neq t ]\!]$ for some constants $\lambda_{uv}$. If $\lambda_{uv} > 0$, the corresponding cost function is submodular, for $\lambda_{uv} < 0$ it is supermodular. It can be shown, see Exercise 11.28, that using reparametrization any pairwise binary energy can be transformed into the form of Ising model.

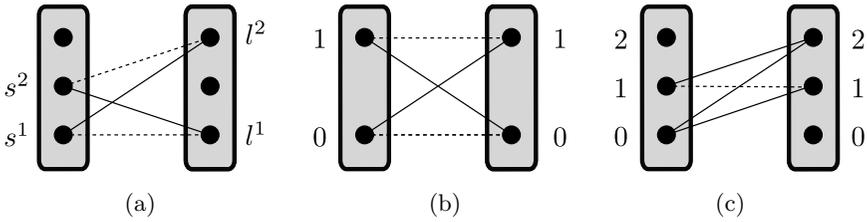

**Figure 11.3:** **(a)** Illustration of Definition 11.8 of a submodular pairwise cost function. The total cost of the label pairs corresponding to the dashed lines should not exceed the total cost of the label pairs corresponding to the solid lines; **(b)** Binary pairwise cost function. If the sum of the costs corresponding to the dashed lines does not exceed the sum of the costs associated with the solid lines, the cost function is submodular. Otherwise it is supermodular. **(c)** Potts cost function: solid lines are assigned the cost $\lambda > 0$, the dashed lines have a zero cost. The shown configuration is the one, where the submodularity definition (11.9) is not fulfilled, see Example 11.11.

**Example 11.11** (Potts model). The Potts model is the generalization of the Ising model into the multi-label case. That is, the pairwise factors have the same format as in the Ising model $\theta_{uv}(s,t) = \lambda_{uv}[\![s \neq t]\!]$, however the label sets $\mathcal{Y}_u$ and $\mathcal{Y}_v$ may have more than two elements. The Potts model is not submodular even for $\lambda > 0$, because already for 3 labels (i.e. $\mathcal{Y}_u = \{0, 1, 2\}$) it holds that

$$\theta(0,1) + \theta(1,2) = 2\lambda \geq \lambda = \theta(0,2) + \theta(1,1) \,, \tag{11.11}$$

see Figure 11.3(c) for illustration.

Since labels in each node are ordered, one may speak about the *next* label. More precisely, the label $l$ is called *next* for the label $s$, if $l \geq s$ and there is no other label $l'$ such that $l \geq l' \geq s$. The *previous* label is defined similarly. We will use the notations $s + 1$ and $s - 1$ for the label next and previous for $s$ respectively. This is motivated by the important special case $\mathcal{Y}_u = \{1, 2 \ldots, n\}$, when this notation becomes a natural choice.

Checking submodularity with the defining criterion (11.9) requires $O(|\mathcal{Y}_{uv}|^2)$ operations, which can be computationally quite expensive for large label sets. An equivalent definition provided by the following proposition allows us to check submodularity in linear time $O(|\mathcal{Y}_{uv}|)$.

**Proposition 11.12.** The cost function $\theta_{uv}$ is submodular if and only if

$$\theta_{uv}(s,l) + \theta_{uv}(s+1,l+1) - \theta_{uv}(s+1,l) - \theta_{uv}(s,l+1) \leq 0. \quad (11.12)$$

for all $s$ and $l$ such that $s+1$ and $t+1$ exist, see Figure 11.4(a). In particular, when $\mathcal{Y}_u = \{1, 2, \ldots, |\mathcal{Y}_u|\}$, then $s \in \{1, \ldots, |\mathcal{Y}_u| - 1\}$ and $l \in \{1, \ldots, |\mathcal{Y}_v| - 1\}$.

The technique used in the following proof plays an important role for this section and will also be used several times later on. First, we introduce the notation

$$\frac{\partial^2 \theta(s,l)}{\partial s \partial l} := \theta(s+1,l+1) + \theta(s,l) - \theta(s,l+1) - \theta(s+1,l). \quad (11.13)$$

for the mixed second derivative of the pairwise cost function. Note that the derivative is non-positive for submodular and non-negative for supermodular pairwise costs. The following lemma constitutes the core of the proof of Proposition 11.12:

**Lemma 11.13.** Let $\mathcal{Y}_u$ and $\mathcal{Y}_v$ be finite ordered sets and $\theta \colon \mathcal{Y}_u \times \mathcal{Y}_v \to \mathbb{R}$ be a function of two discrete variables. Then for $y^2 \geq y^1$, $z^2 \geq z^1$ the following equality holds

$$\theta(y^2,z^2) + \theta(y^1,z^1) - \theta(y^1,z^2) - \theta(y^2,z^1) = \sum_{s=y^1}^{y^2-1} \sum_{l=z^1}^{z^2-1} \frac{\partial^2 \theta(s,l)}{\partial s \partial l}. \quad (11.14)$$

Note that the sign of the left-hand-side of (11.14) determines the submodularity of $\theta$ according to Definition 11.8. Therefore, Proposition 11.12 directly follows from Lemma 11.13.

*Proof of Lemma 11.13.* In a nutshell, the proof is based on the following fact: For a function $f \colon \mathbb{Z} \to \mathbb{R}$ and $a > b$ it holds that

$$f(a) - f(b) = \sum_{i=b}^{a-1} f(i+1) - f(i), \quad (11.15)$$

which in the continuous case reads

$$f(a) - f(b) = \int_b^a df = \int_b^a f'(t)dt. \quad (11.16)$$

Here, the derivative $f'$ is represented by its discrete analog $\frac{f(t+dt)-f(t)}{dt}$ with $dt = 1$.

A similar property holds also for a function of two variables:

$$f(a,b) + f(a',b') - f(a',b) - f(a,b') = \int_a^{a'} \int_b^{b'} \frac{\partial^2 f(s,t)}{\partial s \partial t} ds dt. \quad (11.17)$$

Nevertheless, we provide a rigorous proof below.

Let us consider the differences $\theta(y^2, z^1) - \theta(y^1, z^1)$ and $\theta(y^2, z^2) - \theta_{uv}(y^1, z^2)$. Comparing them to (11.15) one obtains

$$\theta(y^2, z^2) - \theta(y^1, z^2) = \sum_{s=y^1}^{y^2-1} \theta(s+1, z^2) - \theta(s, z^2), \quad (11.18)$$

$$\theta(y^2, z^1) - \theta(y^1, z^1) = \sum_{s=y^1}^{y^2-1} \theta(s+1, z^1) - \theta(s, z^1). \quad (11.19)$$

Subtracting the second equality from the first one results in

$$\theta(y^2, z^2) + \theta(y^1, z^1) - \theta(y^1, z^2) - \theta(y^2, z^1)$$
$$= \sum_{s=y^1}^{y^2-1} \theta(s+1, z^2) + \theta(s, z^1) - \theta(s, z^2) - \theta(s+1, z^1). \quad (11.20)$$

The similar procedure w.r.t. each summand of the right-hand-side of (11.20) gives:

$$\theta(s+1, z^2) + \theta(s, z^1) - \theta(s, z^2) - \theta(s+1, z^1) =$$
$$\sum_{l=z^1}^{z^2-1} \theta(s+1, l+1) + \theta(s, l) - \theta(s, l+1) - \theta(s+1, l). \quad (11.21)$$

Plugging (11.21) into (11.20) one obtains the implication of the lemma. $\qquad \square$

Similarly to Proposition 11.12 which characterizes submodularity of a pairwise cost function through the sign of its second derivative, the following theorem gives a criterion for submodularity of the energy function through the sign of its second derivative:

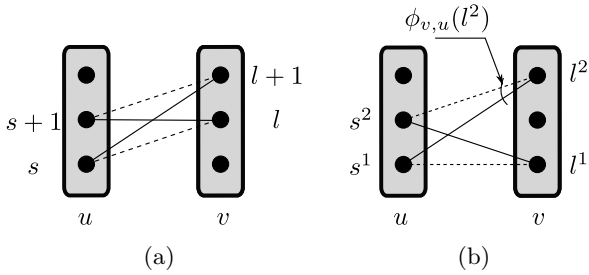

**Figure 11.4: (a)** Illustration of the claim of Proposition 11.12. To check submodularity it is sufficient to consider only pairs of the labels placed next to each w.r.t. the ordering within each node. For submodularity the sum of the costs corresponding to the dashed lines should not exceed those of the solid lines. **(b)** Illustration of the claim of Proposition 11.16: the coordinate $\phi_{v,u}(l^2)$ of the dual vector is added to the costs associated with both, dashed and solid lines, therefore the submodularity inequality (11.9) is not affected by the reparametrization.

**Theorem 11.2** ([1, 18, 62])**.** The graphical model energy $E$ is submodular if and only if

$$\begin{aligned}
\frac{\partial^2 E(y)}{\partial y_i \partial y_j} := & \ E(s_1, \ldots, s_i, \ldots, s_j, \ldots, s_n) \\
& + E(s_1, \ldots, s_i + 1, \ldots, s_j + 1, \ldots, s_n) \\
& - E(s_1, \ldots, s_i + 1, \ldots, s_j, \ldots, s_n) \\
& - E(s_1, \ldots, s_i, \ldots, s_j + 1, \ldots, s_n) \le 0
\end{aligned} \qquad (11.22)$$

for any labeling $y$ and any two neighboring nodes $i$ and $j$.

*Proof.* Necessity of the condition (11.22) follows directly from (11.7) as the inequality (11.22) is a special case thereof.

The proof of sufficiency is more involved, therefore, we omit it here and refer to the original manuscript. $\qquad \square$

Importantly, Theorem 11.2 remains also correct for higher order models, since it is a specialization of the corresponding fact for general submodular functions. Note that the non-positive second derivative is the property of concave functions in continuous case.

**Convex and truncated convex pairwise cost functions**

Proposition 11.12 can be efficiently used to check submodularity of certain functions of two variables.

**Example 11.14** (Convex functions). Let the labels sets $\mathcal{Y}_u$ and $\mathcal{Y}_v$ be equal to $\{1, 2 \ldots, n\}$, i.e. they represent a range of integer numbers. Let the subtraction operation "$-$" be defined in a natural way. Let $f \colon \mathbb{R} \to \mathbb{R}$ be a convex function and $\theta_{uv}(s, t) = f(s - t)$. Let us prove that $\theta_{uv}$ is submodular. According to (11.12) it is necessary to show that

$$f(s - t) + f((s + 1) - (t + 1)) \le f(s + 1 - t) + f(s - t - 1) \quad (11.23)$$

Indeed, the convexity of $f$ implies

$$f(s+1-t)+f(s-t-1) \ge 2f(\frac{1}{2}(s+1-t+s-t-1)) = 2f(s-t) \quad (11.24)$$

Note that the left-hand-side of (11.23) is equal to $2f(s - t)$ as well, which finalizes the proof.

In particular, the following functions are submodular:

$$\theta_{uv}(s, t) = |s - t| \qquad (11.25)$$
$$\theta_{uv}(s, t) = (s - t)^2. \qquad (11.26)$$

**Example 11.15** (Truncated convex). Consider the *truncated convex* pairwise cost functions

$$\theta_{uv}(s, t) = \min\{a, |s - t|\}, \ a > 0$$

and

$$\theta_{uv}(s, t) = \min\{a, (s - t)^2\}, \ a > 0.$$

In the binary case they coincide with the Ising model and therefore are submodular. However, as soon as three and more labels are considered, these functions are not submodular anymore. In particular, for $a = 1$ these functions reduce to the Potts cost function 11.11, which is non-submodular.

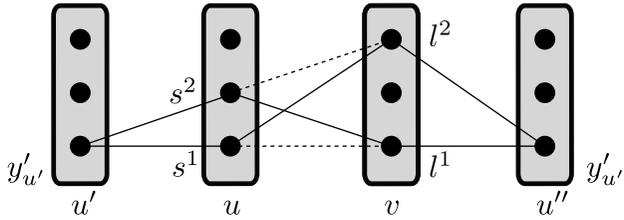

**Figure 11.5:** Illustration to the proof of Theorem 11.1. The solid lines correspond to the labelings $y^1$ and $y^2$, the dashed lines correspond to $y^1 \vee y^2 = (s^2, l^2, y')$ and $y^1 \wedge y^2 = (s^1, l^1, y')$.

### 11.2.2 Sum-of-submodular functions

Consider the energy function of a graphical model. Since modularity of unary costs (see Lemma 11.6) and submodularity of all pairwise costs implies submodularity of the energy according to Lemma 11.7, Theorem 11.1 states that the inverse claim holds as well. Below, we prove this inverse statement.

*Proof of Theorem 11.1.* Consider any edge $uv \in \mathcal{E}$, then the corresponding cost function is denoted by $\theta_{uv}$. It is sufficient to prove that for any four labels $s^1, s^2 \in \mathcal{Y}_u$ and $l^1, l^2 \in \mathcal{Y}_v$ such that $s^2 \geq s^1$ and $l^2 \geq l^1$ the submodularity inequality (11.9) holds.

Consider two labelings $y^1$ and $y^2$ such that

$$y_u^1 = s^1, \ y_u^2 = s^2, \ y_v^2 = l^1, \ y_v^1 = l^2\,,$$
and $y_{u'}^1 = y_{u'}^2 = y_{u'}'$ for $u' \in \mathcal{V}\backslash\{u,v\}$ and some $y' \in \mathcal{Y}_{\mathcal{V}\backslash\{u,v\}}\,. \quad (11.27)$

See Figure 11.5 for illustration.

W.l.o.g. and for brevity, we assume that the nodes $u$ and $v$ correspond to the first two coordinates in the vector representation of the labelings on $\mathcal{Y}_\mathcal{V}$ and the labelings $y^1$ and $y^2$ can be written as

$$y^1 = (s^1, l^2, y')\,, \quad y^2 = (s^2, l^1, y')\,. \quad (11.28)$$

Submodularity of the energy $E$ implies that the inequality (11.7) holds for $y^1$ and $y^2$. Note that $y^1 \vee y^2 = (s^2, l^2, y')$ and $y^1 \wedge y^2 = (s^1, l^1, y')$.

Therefore, the inequality (11.7) reduces to

$$\begin{aligned}
0 \leq E(y^1) &+ E(y^2) - E(y^1 \vee y^2) - E(y^1 \wedge y^2) \\
&= \theta_{uv}(s^2, l^1) + \theta_{uv}(s^1, l^2) - \theta_{uv}(s^2, l^2) - \theta_{uv}(s^1, l^1), \quad (11.29)
\end{aligned}$$

which finalizes the proof. □

Functions representable as a sum of submodular functions of a small number of variables are important for two reasons. First of all, due to their practical value, and second, because algorithms for their minimization are much more efficient than those for general submodular functions. We refer to §11.6 for the corresponding references.

### 11.2.3 Reparametrization and node-edge agreement for submodularity

Recall the definition of reparametrization given in Chapter 6. Since reparametrization does not change the energy of any labeling, it does not influence the submodularity property (11.7) of the energy function. Together with Theorem 11.1 it implies that reparametrization does not influence the submodularity of the pairwise costs as well. This fact can be also easily shown directly.

**Proposition 11.16.** Let the pairwise costs $\theta_{uv}$ be submodular for all $uv \in \mathcal{E}$. Then for any reparametrization $\phi$ the reparametrized costs $\theta_{uv}^{\phi}$ are submodular as well.

*Proof.* According to Lemma 11.6 the functions $\phi_{u,v}(s)$, $s \in \mathcal{Y}_u$, are modular. Therefore, $\theta_{uv}^{\phi}$ is submodular as a sum of a submodular and modular functions according to Lemma 11.7. □

Note that a statement similar to Proposition 11.16 does not hold for higher order costs in general.

An important fact, which directly implies the polynomial solvability of submodular MAP-inference problems, is stated in the following proposition:

**Proposition 11.17.** Let $(\mathcal{G}, \mathcal{Y}_{\mathcal{V}}, \theta)$ be a graphical model and the pairwise costs $\theta_{uv}$ be submodular for all $uv \in \mathcal{E}$. Then the node-edge agreement of $\theta$ implies that there exists an integer labeling $y$ consisting of locally optimal labels and edges. Moreover, this labeling can be found by selecting the *highest* optimal label belonging to the arc-consistency closure $\mathrm{cl}(\mathrm{mi}[\theta])$ in each node independently.

*Proof.* The proof is illustrated in Figure 11.4(b). Let $s^2 \in \mathcal{Y}_u$ and $l^2 \in \mathcal{Y}_v$ be the highest labels of the closure $\mathrm{cl}(\mathrm{mi}[\theta])$ in the nodes $u$ and $v$ respectively. Assume that there is no label pair $(s^2, l^2)$ in the closure. Due to the arc-consistency, however, there are $s^1 \in \mathcal{Y}_u$ and $l^1 \in \mathcal{Y}_v$ such that $(s^1, l^2)$ and $(s^2, l^1)$ belong to the closure. It implies that

$$\theta_{uv}(s^1, l^1) \geq \theta_{uv}(s^1, l^2) = \theta_{uv}(s^2, l^1) =: \alpha \leq \theta_{uv}(s^2, l^2)\,. \qquad (11.30)$$

At the same time, the submodularity implies

$$\theta_{uv}(s^1, l^1) + \theta_{uv}(s^2, l^2) \leq \theta_{uv}(s^1, l^2) + \theta_{uv}(s^2, l^1) = 2\alpha\,. \qquad (11.31)$$

Taken together these inequalities imply

$$\theta_{uv}(s^1, l^1) = \theta_{uv}(s^2, l^2) = \alpha\,, \qquad (11.32)$$

and, therefore, $\theta_{uv}(s^1, l^1)$ and $\theta_{uv}(s^2, l^2)$ are locally optimal and belong to the closure. This finalizes the proof. $\qquad \square$

Proposition 11.17 implies that

- The local polytope relaxation is tight for submodular problems. Indeed, node-edge agreement is necessary for dual optimality. Similar to the case of acyclic graphical models, considered in §6.3, Proposition 11.17 together with Proposition 6.2 implies that node-edge agreement is also sufficient for the dual optimality and, moreover, that the Lagrange dual is tight for the non-relaxed MAP-inference problem.

- Algorithms such as TRW-S, min-sum diffusion and sub-gradient method converge to the optimum of the non-relaxed binary MAP-inference problem. Indeed, the former two converge to node-edge agreement, and, therefore, to the dual optimum. The sub-gradient method converges to the dual optimum independently of the submodularity property.

Although tightness of the local polytope relaxation implies the polynomial solvability of the submodular MAP-inference, it does not imply the existence of a practically efficient polynomial algorithm, due to the universality of the local polytope (see Theorem 4.1).

In §11.3 and §11.4 we will show that submodular pairwise MAP-inference reduces to a special, much more narrow class of linear programs, known as min-cut/max-flow problems. The latter has a number of practically efficient dedicated finite step algorithms for their solution. These algorithms can be therefore used to solve the submodular MAP-inference problem as well.

### 11.2.4 Permuted submodularity

Submodularity is a property dependent on the order of labels. Note, however, that the order does not play any role for the energy itself. It is only needed to define the submodularity property.

**Definition 11.18.** A MAP-inference problem $(\mathcal{G}, \mathcal{Y}_y, \theta)$ is called *permuted submodular* if there exists an ordering on each label set $\mathcal{Y}_u$, $u \in \mathcal{V}$ such that the corresponding problem is submodular.

**Example 11.19.** Binary problems on a tree are always permuted submodular. One should traverse the tree and change the order of labels if needed to switch from supermodular cost to submodular.

**Example 11.20.** Binary supermodular problems on bipartite graphs (e.g. acyclic graphs, 2D grid) are permuted submodular. To obtain a submodular problem, the order of the labels on one "side" of the bipartite graph must be reversed.

Importantly, there is an efficient algorithm, which either finds the permutation turning the problem into a submodular one, or proves, that it is impossible. See §11.6 for the respective reference.

Note also that Proposition 11.17 holds also for permuted submodular problems. And in turn, the corollaries of the proposition hold as well: The local polytope relaxation is tight for such problems and the algorithms converging to the node-edge agreement converge to the optimum of the non-relaxed MAP-inference problem. However, to reconstruct an integer solution from the optimal reparametrization, in the same way as
it is done in the proof of Proposition 11.17, one has to know (compute)
the order of the labels. At the same time, knowing the order is not
mandatory for the reconstruction in general, see e.g. [138, Thm. 16].

## 11.3  Binary submodular problems as min-cut

We start this section with a definition of the min-$st$-cut problem, which
is one of the combinatorial problems having efficient (small) polynomial
time algorithms. Further in the section we will show how submodular
MAP-inference problems can be reduced to the min-$st$-cut. Moreover,
in Chapter 12 we will show that not only (multilabel) submodular
problems are reducible to the min-$st$-cut, but also the local polytope
relaxation of the binary problems allows for such a reduction.

### 11.3.1  min-$st$-cut Problem

**Definition 11.21.** Let $\mathcal{G}' = (\mathcal{V}', \mathcal{E}')$ be a directed graph, $c \colon \mathcal{E}' \to \mathbb{R}$ be
the *weight function of the edges* and $s, t \in \mathcal{V}'$ be two different nodes,
which will be called *source* and *target* respectively. The *st-cut* $C = (S, T)$
is a partition of $\mathcal{V}'$ into two parts $S$ and $T$ such that

$$S \cap T = \emptyset \tag{11.33}$$

$$S \cup T = \mathcal{V}' \tag{11.34}$$

$$s \in S, \ t \in T. \tag{11.35}$$

In the following we will use term st-cuts simply as *cuts*.

The weight of the cut $C$ is the total weight of edges leading from $S$
to $T$:

$$c(S, T) := \sum_{\substack{u \in S, v \in T \\ (u,v) \in \mathcal{E}'}} c_{u,v}. \tag{11.36}$$

The *min-st-cut* (or shortly *min-cut*) problem consists in finding a
cut with the minimal total weight.

The min-cut problem is illustrated in Figure 11.6.

As many combinatorial problems, min-cut can be formulated as an
integer linear program, see Example 11.22. In general it is NP-hard, just like other ILPs. However, if all weights $c_{u,v}$ are non-negative, its linear programming relaxation is tight and therefore the problem is polynomially solvable.[1] The corresponding linear programming relaxation has a special structure, which allows us to find its minimum with polynomial finite-step algorithms with the worst-case complexity $O(|\mathcal{V}'|^3)$. Typically, the full Lagrange dual of the problem is solved instead of the primal problem. It constitutes the *max-flow* problem and has a number of efficient polynomial finite-step solvers.

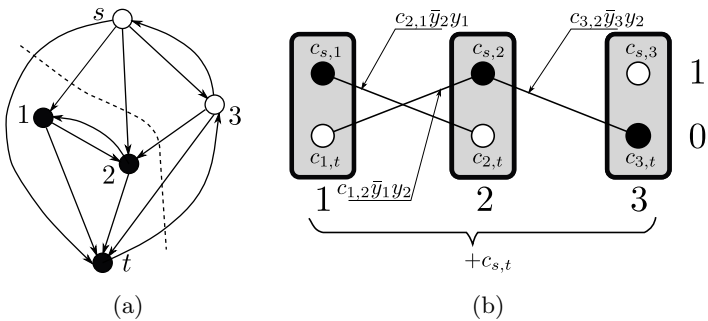

(a)          (b)

**Figure 11.6: (a)** An example of the min-$st$-cut problem. The shown problem can be equivalently written as minimization of the following pseudo-boolean function represented in the oriented form (see §11.3.2): $c_{s,t} + c_{s,1}x_1 + c_{s,2}x_2 + c_{s,3}x_3 + c_{1,t}\bar{x}_1 + c_{2,t}\bar{x}_2 + c_{3,t}\bar{x}_3 + c_{1,2}\bar{x}_1x_2 + c_{2,1}\bar{x}_2x_1 + c_{3,2}\bar{x}_3x_2$. The cut shown by a dashed line corresponds to the partitioning, where white nodes belong to the set $S$ and black ones to the set $T$. The cut corresponds to the partition $(x_1, x_2, x_3) = (1, 1, 0)$ and its weight is equal to $c_{s,t} + c_{s,1} + c_{s,2} + c_{3,2} + c_{3,t}$. Note that the edge $(t, 3)$ does not belong to the cut. Moreover, the edge $(t, 3)$ does not belong to *any* cut, as well as the edge $(3, s)$. Therefore, these edges can be ignored. **(b)** The MAP-inference problem is equivalent to the min-cut problem in (a). Black circles define the labeling corresponding to the cut shown in (a).

**Example 11.22** (ILP representation of the min-$st$-cut and its LP relaxation). The min-$st$-cut problem with non-negative edge weights has the following

---

[1]The combinatorial min-cut and max-cut (see Example 3.48) problems are closely related, therefore one typically speaks about the min-cut problem, if the edge weights are non-negative and the term *max-cut* is used, when this is not the case.

simple ILP representation:

$$\min_{x \in \{0,1\}^{\mathcal{V}' \cup \mathcal{E}'}} \sum_{(u,v) \in \mathcal{E}'} c_{u,v} x_{u,v} \tag{11.37}$$

$$\text{s.t. } x_u - x_v + x_{u,v} \geq 0 \quad \forall (u,v) \in \mathcal{E}' \tag{11.38}$$

$$x_s = 0, \ x_t = 1. \tag{11.39}$$

The coordinates of the binary vector $x$ are indexed by vertexes and edges of the graph. The constraint (11.38) essentially states that if $u \in S$ and $v \in T$, the edge connecting them is in the cut. The non-negativity of the weights assures that for $u \in T$ and $v \in S$ the edge variable $x_{u,v}$ is set to zero. Indeed, in this case the inequality (11.38) reduces to $1 + x_{u,v} \geq 0$, therefore $x_{u,v}$ may take both values, 0 and 1. However, if $c_{u,v} > 0$ the value 0 will be selected, as it corresponds to the smaller objective value.

The LP relaxation of (11.37)-(11.39) is obtained as usual, by convexifying the integrality constraints.

## 11.3.2 Oriented forms

Although the graphical form of combinatorial problems is of great importance for the development of the field, their algebraic forms have the advantage of a more compact problem presentation and help to give unambiguous proofs. Along with the ILP formulation for the min-cut problem, we will consider another algebraic form, the *quadratic pseudo-boolean minimization* problem.

For a binary variable $x \in \{0,1\}$ let $\bar{x}$ denote its negation, which is $\bar{x} = 1 - x$. Functions $f \colon \{0,1\}^n \to \mathbb{R}$ are called *pseudo-boolean*. Note that the class of pseudo-boolean functions coincides with the class of *set-functions* as defined in §11.1. In the rest of the book we will prefer the latter term as it is done in e.g. [14]. An important example of pseudo-boolean functions $f(x)$, $x \in \{0,1\}^n$, are those representable as a *quadratic* polynomial of $x_i$ and $\bar{x}_i$. A subclass of such functions, tightly related to the min-cut problem, is defined as follows:

**Definition 11.23.** Let $\nu$, $\alpha_i, \beta_i$ and $\gamma_{i,j}$, $i, j = 1, \ldots, n$ be arbitrary real numbers. For $x \in \{0,1\}^n$ a quadratic pseudo-boolean function of the

type

$$f(x) = \nu + \sum_{i=1}^{n} \alpha_i x_i + \sum_{i=1}^{n} \beta_i \bar{x}_i + \sum_{i=1}^{n} \sum_{j=1}^{n} \gamma_{ij} \bar{x}_i x_j \qquad (11.40)$$

will be called an *oriented form*.

Contrary to the other terms used in this book, the term *oriented form* has not been used in the literature before and is first introduced here.

**Example 11.24.** The following quadratic pseudo-boolean functions are oriented forms: $-5 + 7x_3\bar{x}_1 + 12x_3 + \bar{x}_1$; $2\bar{x}_1 x_2 + 3\bar{x}_2 x_1 - 6x_3\bar{x}_4$.

**Example 11.25.** The quadratic pseudo-boolean functions below are not oriented forms: $-5 + 7x_3 x_1 + 12x_3 + \bar{x}_1$; $2x_1 x_2 + 3\bar{x}_2 x_1 - 6x_3\bar{x}_4$.

### 11.3.3 Min-cut as oriented form

Oriented forms are well suited for representation of graph partitions as the one to be found in the min-cut problem. Let the coordinates of a binary vector $x \in \{0,1\}^{\mathcal{V}'}$ be the indicator variables, denoting the subset each vertex belongs to. In other words, if $u \in S$, $x_u = 0$, and, on the other side, $x_u = 1$ is equivalent to $u \in T$. In this way, the vector $x$ represents a partition of the graph $\mathcal{G}'$.

Let the monomial $c_{u,v}\bar{x}_u x_v$ correspond to the directed edge $(u,v)$ with the weight $c_{u,v}$. Then the total weight of the cut (11.36) corresponding to the partition associated with the vector $x$ is equal to

$$\sum_{(u,v)\in\mathcal{E}'} c_{u,v}\bar{x}_u x_v . \qquad (11.41)$$

Indeed, the items of the sum are non-zero only if $x_u = 0$ and $x_v = 1$, meaning that $u \in S$ and $v \in T$.

Note that since $s \in S$ and $t \in T$ it follows that $x_s = 0$ and $x_t = 1$. Therefore, $c_{s,v}\bar{x}_s x_v = c_{s,v}x_v$, $c_{u,t}\bar{x}_u x_t = c_{u,t}\bar{x}_u$ and $c_{s,t}\bar{x}_s x_t = c_{s,t}$. It also implies that (directed) edges $(u,s)$ and $(t,u)$, $u \in \mathcal{V}$ can be ignored, since they do not belong to any cut and the corresponding monomials $c_{u,s}\bar{x}_u x_s$ and $c_{t,u}\bar{x}_t x_u$ vanish.

**Table 11.1:** Mapping between the min-cut problem, minimization of the oriented forms and the MAP-inference in graphical models.

| min-cut | oriented form | $E(y; \theta)$ |
|---|---|---|
| $V'$ | $x \in \{0,1\}^{\mathcal{V}'}$ | $V = V' \backslash \{s, t\}$ |
| $u \in S$ | $x_u = 0$ | $y_u = 0$ |
| $u \in T$ | $x_u = 1$ | $y_u = 1$ |
| $c_{u,v}, (u, v) \in \mathcal{E}'$ | $c_{u,v} \bar{x}_u x_v$ | $\theta_{uv}(0, 1) = c_{u,v}$ |
| in particular: | | |
| $s \in S$ | $x_s = 0$ | — |
| $t \in T$ | $x_t = 1$ | — |
| $c_{s,v}, (s, v) \in \mathcal{E}'$ | $c_{s,v} x_v$ | $\theta_v(1) = c_{s,v}$ |
| $c_{u,t}, (u, t) \in \mathcal{E}'$ | $c_{u,t} \bar{x}_u$ | $\theta_u(0) = c_{u,t}$ |
| $c_{s,t}, (s, t) \in \mathcal{E}'$ | $c_{s,t}$ | constant |
| $c_{u,s}, (u, s) \in \mathcal{E}'$ | $0$ | — |
| $c_{t,v}, (t, v) \in \mathcal{E}'$ | $0$ | — |

Concluding, each min-cut problem can be equivalently represented as minimization of an oriented form

$$c_{s,t} + \sum_{v:\, (s,v) \in \mathcal{E}'} c_{s,v} x_v + \sum_{u:\, (u,t) \in \mathcal{E}'} c_{u,t} \bar{x}_u + \sum_{\substack{u,v \in \mathcal{E}' \\ u \neq s,\ v \neq t}} c_{u,v} \bar{x}_u x_v \qquad (11.42)$$

with the constant $c_{s,t}$ being zero if the edge $(s, t)$ does not exist. Conversely, minimization of an oriented form can be treated as a min-cut problem. Figure 11.6(a) illustrates this correspondence. In this representation weights of the graph edges become coefficients of the corresponding monomials. In particular it means that a min-cut problem with non-negative weights corresponds to an oriented form with non-negative coefficients and the other way around. Table 11.1 summarizes the mapping between the min-cut problem and minimization of the oriented forms.

### 11.3.4  Binary energy as oriented form

To show a close relation between binary MAP-inference and min-cut problems, we show that the former can be equivalently seen as minimization of an oriented form. We will start with an inverse transformation by showing that an oriented form can be seen as a binary graphical model energy up to a constant term. The latter does not influence the optimization, and, therefore, can be ignored.

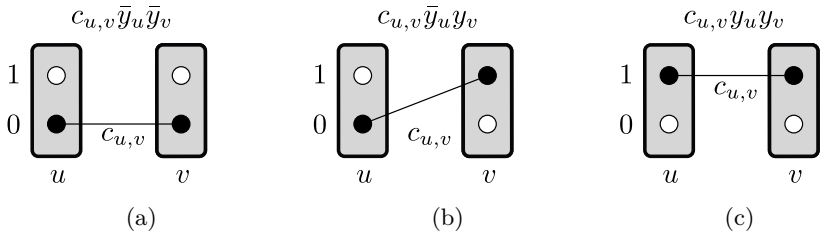

**Figure 11.7:** Correspondence between quadratic monomials and pairwise costs of binary graphical models. Solid lines corresponding to label pairs denote the only non-zero values of the cost vector of the pairwise factor.

For the beginning, consider Figure 11.7, which illustrates how different quadratic monomials can be represented as pairwise factors of a graphical model. So, for example, the monomial $c_{u,v}\bar{y}_u y_v$ for binary variables $y_u$ and $y_v$ can be represented as a pairwise factor with the cost vector

$$\theta_{uv}(0,1) = c_{u,v}, \qquad \theta_{uv}(0,0) = \theta_{uv}(1,1) = \theta_{uv}(1,0) = 0. \quad (11.43)$$

To specify a general transformation, consider the oriented form (11.40) with $n$ variables. Build the graphical model $(\mathcal{G}, \mathcal{Y}_\mathcal{V}, \theta)$ with the graph $\mathcal{G}$ consisting of $n$ nodes, $\mathcal{V} = \{1, 2, \ldots, n\}$, where $i$-th node corresponds to the $i$-th variable and the other way around. For each monomial $c_{i,j}\bar{x}_i x_j$ add an edge between nodes $i$ and $j$ to the graph $\mathcal{G}$, should it not exist yet. Set all costs $\theta$ initially to zero. For each monomial $\alpha_i x_i$ assign $\theta_i(1) := \alpha_i$, for $\beta_i \bar{x}_i$ assign $\theta_i(0) := \beta_i$ and for $c_{i,j}\bar{x}_i x_j$ assign $\theta_{ij}(0,1) := c_{i,j}$. Table 11.1 summarizes the construction, which is also illustrated in Figures 11.6(b) and 11.7.

To specify the inverse transformation of a graphical model to an oriented form note that only pairwise factors having the structure of (11.43), that is every entry is zero except one non-zero value on the diagonal, can be given by an oriented form according to the above construction. In the following we will adopt the name *diagonal form* for such factors. It turns out that any pairwise cost function can be transformed into the diagonal form by reparametrization. Figure 11.8 illustrates this transformation, given by the formulas below:

$$\phi_{u,v}(0) = \theta_{uv}(1,0) - \theta_{uv}(0,0) \tag{11.44}$$

$$\phi_{v,u}(0) = -\theta_{uv}(1,0) \tag{11.45}$$

$$\phi_{v,u}(1) = -\theta_{uv}(1,1) \tag{11.46}$$

Note that after the transformation (11.44)-(11.46) the only non-zero pairwise cost value is equal to

$$\theta_{uv}^{\phi}(0,1) = \theta_{uv}(0,1) + \theta_{uv}(1,0) - \theta_{uv}(1,1) - \theta_{uv}(0,0), \tag{11.47}$$

which is non-negative for submodular pairwise cost functions. To guarantee that all monomials in the oriented form are non-negative, one has to guarantee that the unary costs are non-negative as well. This can easily be achieved by subtracting the minimal value from the unary costs:

$$\forall s \in \mathcal{Y}_u \quad \theta_u(s) := \theta_u(s) - \min_{l \in \mathcal{Y}_u} \theta_u^{\phi}(l). \tag{11.48}$$

Since each labeling contains either the label 0 or 1 in each node, such operation only subtracts a constant value from energies of all labelings and therefore does not affect the optimum of the MAP-inference problem.

### 11.3.5 Equivalence: binary energy minimization = min-cut

Summarizing, transformations (11.44)-(11.46) performed for each edge $uv \in \mathcal{E}$ and followed by transformations (11.48) applied to each node $u \in \mathcal{V}$ translates the cost-vector into a form where it can be represented as an oriented form. Moreover, this oriented form has only non-negative coefficients, if the MAP-inference problem is submodular.

Combining this with the results of §11.3.3, one obtains the following statement:

**Proposition 11.26.** The binary MAP-inference problem reduces to the min-*st*-cut problem and the other way around. Should the costs of the MAP-inference problem be submodular, then it can be reduced to the min-*st*-cut problem with non-negative edge weights. Inversely, such min-cut problems are equivalently representable as binary submodular MAP-inference problems.

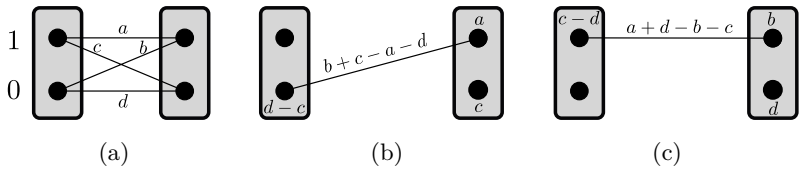

**Figure 11.8:** Transformation of an arbitrary pairwise cost vector **(a)** to the diagonal **(b)** and parallel **(c)** forms. For the use of the latter one see Chapter 12.

**Remark 11.27.** The transformation $(11.44)$-$(11.46)$ is asymmetric, since $\theta_{uv}(0,1) \neq 0$ and $\theta_{uv}(1,0) = 0$. A symmetric transformation can also be found, e.g. such that $\theta_{uv}(0,1) = \theta_{uv}(1,0) \neq 0$. Although it can be useful in some cases, the corresponding min-$st$-cut graph would contain two directed edges $(u,v)$ and $(v,u)$ for such a symmetric representation, whereas the diagonal representation $(11.43)$ requires only a single edge. The latter results in a smaller min-$st$-cut graph and faster optimization.

**Exercise 11.28.** Show that using reparametrization any pairwise binary energy can be transformed into the form of the Ising model (see Example 11.10). Moreover, submodular pairwise costs are transformed into submodular with $\lambda_{uv} \geq 0$ and supermodular into supermodular with $\lambda_{uv} \leq 0$.

## 11.4 Transforming multi-label problems into binary ones

As shown in §11.3 binary submodular MAP-inference problems can be reduced to the min-cut problem with non-negative weights. Interestingly, a similar reduction can be performed also for multi-label graphical models. This additionally supports the theoretical fact that generalization of the submodularity to general finite lattices from the set-lattice also preserves computational complexity of this class of functions.

Instead of providing a direct reduction of a multi-label problem to the min-cut problem, we will give a transformation of the general multi-label problem to the binary one. The latter can be reduced to the min-cut as described in §11.3. Importantly, our transformation will preserve the submodularity of the MAP-inference problem. As a result, submodular multi-label problems will be reduced to the min-cut problem with non-negative edge weights.

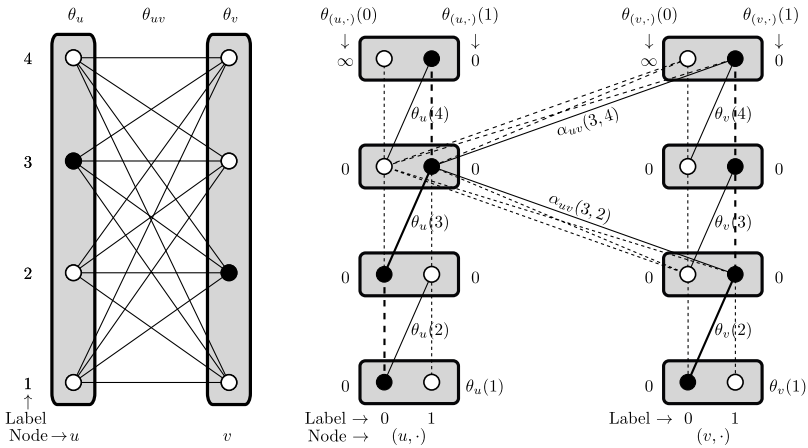

**Figure 11.9:** Transformation of a multi-label problem **(a)** into a binary one **(b)**. Dashed lines correspond to zero costs. Non-shown lines correspond to label pairs with infinitely large costs. Black circles denote a labeling for the multi-label problem and the corresponding one for the binary one. See also the main text for a detailed explanation.

Let the multi-label graphical model $(\mathcal{G} = (\mathcal{V}, \mathcal{E}), \mathcal{Y} = \prod_{u \in \mathcal{V}} \mathcal{Y}_u, \theta)$ be given. We will construct from it a binary graphical model $(\hat{\mathcal{G}} = (\hat{\mathcal{V}}, \hat{\mathcal{E}}), \hat{\mathcal{Y}}_{\hat{\mathcal{V}}} = \prod_{u \in \hat{\mathcal{V}}} \hat{\mathcal{Y}}_u, \hat{\theta})$ with $\hat{\mathcal{Y}}_u = \{0, 1\}$ for all $u \in \hat{\mathcal{V}}$, and build an injective mapping $\mathcal{Y}_{\mathcal{V}} \ni y \to \hat{y} \in \hat{\mathcal{Y}}_{\hat{\mathcal{V}}}$ such that the corresponding energies $E_{\mathcal{G}}(\theta, y)$ and $E_{\hat{\mathcal{G}}}(\hat{\theta}, \hat{y})$ coincide for the corresponding labeling pairs $y, \hat{y}$. The energy $E_{\hat{\mathcal{G}}}(\hat{\theta}, \hat{y})$ will be constructed such that it is equal to $\infty$ for those labelings $\hat{y} \in \hat{\mathcal{V}}$, which have no pre-image $y$ of the constructed mapping.

The binary graphical model is constructed according to the following rules (see illustration in Fig. 11.9):

- The node set contains one node for each label in each node of the multi-label problem: $\hat{\mathcal{V}} := \{(u, l) \colon u \in \mathcal{V}, l \in \mathcal{Y}_u\}$.

- The edges of the binary problem connect (i) nodes corresponding to consecutive labels and (ii) node pairs corresponding to the label

pairs of the multi-label problem:

$$\hat{\mathcal{E}} := \{((u,l),(u,l+1))\colon u \in \mathcal{V}, 1 \le l < |\mathcal{Y}_u|\}$$
$$\cup \{((u,s),(v,t))\colon uv \in \mathcal{E}, (s,t) \in \mathcal{Y}_{uv}\}. \quad (11.49)$$

- The set of labels contains two elements and is the same for each node: $\hat{\mathcal{Y}}_{\hat{u}} = \{0,1\}$, $\hat{u} \in \hat{\mathcal{V}}$.

- The mapping $y \to \hat{y}$ is constructed such that the label $y_u$ in the multi-label problem corresponds to the following labels $\hat{y}_{(u,l)}$ for all $l \in \mathcal{Y}_u$:
$$y_u \to \hat{y}_{(u,l)} = [\![l \ge y_u]\!]. \quad (11.50)$$

- Unary costs are set up to guarantee that the labelings $\hat{y} \in \hat{\mathcal{Y}}_{\hat{\mathcal{V}}}$, not conform to (11.50) get infinite energy:

$$\hat{\theta}_{(v,l)}(0) = \begin{cases} 0, & 1 \le l < |\mathcal{Y}_v| \\ \infty, & l = |\mathcal{Y}_v| \end{cases} \quad (11.51)$$

$$\hat{\theta}_{(v,l)}(1) = \begin{cases} 0, & 1 < l \le |\mathcal{Y}_v| \\ \theta_v(1), & l = 1 \end{cases} \quad (11.52)$$

- Pairwise costs for the pairs of nodes in the binary model, which correspond to the pairs of labels associated with the *same* node of the initial multi-label problem are constructed to take into account unary costs of the multi-label problem:

For $v \in \mathcal{V}$ and $1 \le s < |\mathcal{Y}_v|$: $\hat{\theta}_{(v,s)(v,s+1)} =$

|       | **0**    | **1**           |
|-------|----------|-----------------|
| **0** | 0        | $\theta_v(s+1)$ |
| **1** | $\infty$ | 0               |

- Pairwise costs for pairs of nodes corresponding to label pairs from different nodes of the initial multi-label problem are constructed as follows:

For $uv \in \mathcal{V}$ and $s \in \mathcal{Y}_u$, $t \in \mathcal{Y}_v$: $\hat{\theta}_{(u,s)(v,t)} =$

|       | **0** | **1**             |
|-------|-------|-------------------|
| **0** | 0     | 0                 |
| **1** | 0     | $\alpha_{uv}(s,t)$ |

Here

$$\alpha_{uv}(s,t) = \begin{cases} \theta_{uv}(s,t) + \theta_{uv}(s+1,t+1) \\ \qquad - \theta_{uv}(s,t+1) - \theta_{uv}(s+1,t), & s < |\mathcal{Y}_u|, t < |\mathcal{Y}_v| \\ \theta_{uv}(s,|\mathcal{Y}_v|) - \theta_{uv}(s+1,|\mathcal{Y}_v|), & s < |\mathcal{Y}_u|, t = |\mathcal{Y}_v| \\ \theta_{uv}(|\mathcal{Y}_u|,t) - \theta_{uv}(|\mathcal{Y}_u|,t+1), & s = |\mathcal{Y}_u|, t < |\mathcal{Y}_v| \\ \theta_{uv}(|\mathcal{Y}_u|,|\mathcal{Y}_v|), & s = |\mathcal{Y}_u|, t = |\mathcal{Y}_v| \end{cases}$$
(11.53)

Note that (i) for any $\hat{y} \in \hat{\mathcal{Y}}_{\hat{\mathcal{V}}}$ other than those representable as (11.50) for some labeling $y$ it holds that $E_{\hat{\mathcal{G}}}(\hat{\theta},\hat{y}) = \infty$ due to the pairwise costs $\hat{\theta}_{(v,s)(v,s+1)}$ corresponding to label pairs from the same node of the original multi-label problem; (ii) there is a bijective mapping between $\hat{y}$ which satisfy (11.50) and $\mathcal{Y}_{\mathcal{V}}$. Therefore, to prove that $E_{\mathcal{G}}(\theta,y) = E_{\hat{\mathcal{G}}}(\hat{\theta},\hat{y})$ it is sufficient to show that

$$\theta_u(y_u) = \sum_{l \in \mathcal{Y}_u} \hat{\theta}_{(u,l)}(\hat{y}_{(u,l)}) + \sum_{\substack{l \in \mathcal{Y}_u \\ l < |\mathcal{Y}_u|}} \hat{\theta}_{(u,l),(u,l+1)}(\hat{y}_{(u,l)},\hat{y}_{(u,l+1)}) \quad (11.54)$$

and

$$\theta_{uv}(s,t) = \sum_{l=s}^{|\mathcal{Y}_v|} \sum_{l'=t}^{|\mathcal{Y}_u|} \alpha_{uv}(l,l'). \tag{11.55}$$

hold for any $y \in \mathcal{Y}_{\mathcal{V}}$ and $\hat{y}$ defined as in (11.50).

The proof of (11.54) is straightforward. The proof of (11.55) follows from Lemma 11.13, which implies for $s < |\mathcal{Y}_u|, t < |\mathcal{Y}_v|$ (for the sake of notation we omit the indexes $uv$ in the next formula):

$$\sum_{l=s}^{|\mathcal{Y}_v|} \sum_{l'=t}^{|\mathcal{Y}_u|} \alpha(l,l')$$

$$\overset{(11.53)}{=} \sum_{l=s}^{|\mathcal{Y}_v|-1} \sum_{l'=t}^{|\mathcal{Y}_u|-1} \theta(s,t) + \theta(s+1,t+1) - \theta(s,t+1) - \theta(s+1,t)$$

$$+ \sum_{l=s}^{|\mathcal{Y}_v|-1} a(l,|\mathcal{Y}_v|) + \sum_{l'=t}^{|\mathcal{Y}_u|-1} a(|\mathcal{Y}_u|,l') + a(|\mathcal{Y}_u|,|\mathcal{Y}_v|)$$

$$\overset{(11.14)}{=} \theta(s,t) + \theta(|\mathcal{Y}_v|,|\mathcal{Y}_u|) - \theta(s,|\mathcal{Y}_u|) - \theta(|\mathcal{Y}_v|,t)$$

$$\overset{(11.53)}{+} \sum_{l=s}^{|\mathcal{Y}_v|-1} \theta(l,|\mathcal{Y}_v|) - \theta(l+1,|\mathcal{Y}_v|)$$

$$\overset{(11.53)}{+} \sum_{l'=t}^{|\mathcal{Y}_u|-1} \theta(|\mathcal{Y}_u|,l') - \theta(|\mathcal{Y}_u|,l'+1) + a(|\mathcal{Y}_u|,|\mathcal{Y}_v|)$$

$$\overset{(11.15)}{=} \theta(s,t) + \theta(|\mathcal{Y}_v|,|\mathcal{Y}_u|) - \theta(s,|\mathcal{Y}_u|) - \theta(|\mathcal{Y}_v|,t)$$

$$+ \theta(s,|\mathcal{Y}_v|) - \theta(|\mathcal{Y}_u|,|\mathcal{Y}_v|) + \theta(|\mathcal{Y}_u|,t) - \theta(|\mathcal{Y}_u|,|\mathcal{Y}_v|)$$

$$\overset{(11.53)}{+} \theta(|\mathcal{Y}_u|,|\mathcal{Y}_v|)$$

$$= \theta(s,t)\,. \tag{11.56}$$

**Elimination of the nodes corresponding to the highest labels**  The "top-most" nodes $(u,|\mathcal{Y}_u|)$ in the binary problem in Fig. 11.9 can be eliminated, as they contain a single allowed label (the one with a finite weight). The elimination changes only unary potentials as follows

$$\forall u \in \mathcal{V}, v \in \mathcal{N}_b(u)\colon \qquad \hat{\theta}_{(u,s)}(1) := \hat{\theta}_{(u,s)}(1) + \alpha_{uv}(s,|\mathcal{Y}_v|) \tag{11.57}$$

$$\forall u \in \mathcal{V}\colon \hat{\theta}_{(u,|\mathcal{Y}_u|-1)}(0) := \hat{\theta}_{(u,|\mathcal{Y}_u|-1)}(0) + \theta_u(|\mathcal{Y}_u|\,. \tag{11.58}$$

**Submodularity of the constructed binary problem**  The described transformation of multi-label problems into binary ones is valid in general, for any input pairwise graphical models. Note, however, that if the pairwise costs of the multi-label problem are submodular, so are the pairwise costs of the corresponding binary problem. Indeed,

- $\alpha_{uv}(s,t) \leq 0$ for $s < |\mathcal{Y}_u|$ and $t < |\mathcal{Y}_v|$ due to the submodularity of $\theta_{uv}$.

- Due to elimination of the "to-most" node $s < |\mathcal{Y}_u|$ and $t < |\mathcal{Y}_v|$ hold for any $uv \in \mathcal{E}$ above.

- $\hat{\theta}_{(v,s)(v,s+1)}$ are submodular since $\hat{\theta}_{(v,s)(v,s+1)}(1,0) = \infty$ according to the pairwise costs definition.

This shows equivalence of the min-cut and submodular pairwise energy minimization in general, including multi-label energies as well.

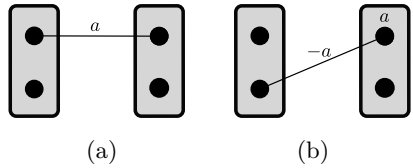

**Figure 11.10:** Transformation of the pairwise cost in a parallel form **(a)** to the diagonal form **(b)**. Edges corresponding to label pairs with zero cost are not shown.

**Diagonal form**   To transform a binary MAP-inference problem obtained from a multi-label one into min-cut one has to additionally transform all pairwise factors into diagonal form as described in §11.3.4. Note that the only factors in the non-diagonal form are $\hat{\theta}_{(u,s)(v,t)}$, $uv \in \mathcal{E}$, $s \in \mathcal{Y}_u$, $t \in \mathcal{Y}_v$, which correspond to labels from different nodes of the initial graph. However, these factors have a special parallel form, which has an especially simple reparametrization to the diagonal form, as shown in Figure 11.10.

**Computational complexity**   Let a binary graphical model be the reduction of an $L$-label one. Then the number of pairwise factors (and therefore, the size of the problem) grows as $O(|\mathcal{E}|L^2)$, where $\mathcal{E}$ is the number of edges in the multi-label model. Up to a constant, this coincides with the size of the multi-label model. Therefore, the considered reduction has a *linear* complexity if we assume the unary and pairwise costs to be defined explicitly as number arrays. In a special case, however, when the pairwise costs constitute an $\ell_1$-norm (11.25), an equivalent binary problem can be constructed, whose size grows as $O(|\mathcal{E}|L)$, only linearly with the number of labels $L$. We refer to [36] for details.

## 11.5   Move-making algorithms

Although the class of multi-label submodular energies is rich, it does not include a number of important pairwise costs, such as Potts (Example 11.11), truncated linear and quadratic costs (Example 11.15). Additionally, the size of the min-cut graph obtained as a reduction of a multi-label submodular problem grows quadratically with the number of labels, which results in high computational time for the min-cut algorithms.

Fortunately, there is an alternative class of min-cut-based inference methods, such that

- they are applicable to submodular problems and beyond, which includes Potts and truncated convex costs.

- on each iteration they compute an approximate integer labeling and energies of these labeling are monotonically decreasing over iterations;

- the size of the min-cut graph does not depend on the number of labels, although a number of such min-cut subproblems must be solved until the method terminates.

These methods are only approximative, which is the price for their favorable properties, although some of them provide approximation guarantees. In practice the obtained approximate solutions are often very good, much better than it is guaranteed by the theory.

### 11.5.1 $\alpha$-expansion algorithm

One of the most popular and efficient min-cut-based move-making algorithms is $\alpha$-*expansion*. Assume the label sets in all nodes are the same and equal to $\mathcal{Y}_0$, i.e. $\mathcal{Y}_u = \mathcal{Y}_0$, $u \in \mathcal{V}$. The algorithm starts with some initial labeling $y$ (e.g. consisting of locally optimal labels). On each iteration it selects a label $\alpha$ to be "expanded". That is, the binary auxiliary subproblem

$$\min_{y' \,:\, \forall u \ y'_u \in \{y_u, \alpha\}} E(y') \tag{11.59}$$

is considered, where the two labels in each node $u$ are $y_u$ and $\alpha$. (If $y_u = \alpha$ the number of labels in the node $u$ reduces to one.) The solution of this binary subproblem is considered as the current labeling for the next iteration, see Algorithm 17. Iterations proceed until the current labeling is not changing anymore.

---

**Algorithm 17** $\alpha$-expansion

---

1: **Init:** $y^0 \in \mathcal{Y}_\mathcal{V}$, $t = 0$
2: **repeat**
3:  **for** $\alpha \in \mathcal{Y}_0$ **do**
4:   $y^{t+1} = \arg\min_{y:\ \forall u\ y_u \in \{y_u^t, \alpha\}} E(y)$
5:   $t = t + 1$
6:  **end for**
7: **until** $E(y^t) = E(y^{t-|\mathcal{Y}_0|})$
8: **return** $y^t$

---

Let us assume that the minimization in line 4 of Algorithm 17 can be done efficiently. Then $E(y^{t+1}) \leq E(y^t)$, that is, Algorithm 17 monotonously decreases the energy $E$.

Let us now figure out, when the minimization (11.59) can be solved efficiently. This is obviously the case when the pairwise costs of the problem are submodular and the order of the labels in the auxiliary problem (11.59) does not change w.r.t. to the order of labels in the original label set $\mathcal{Y}_0$. That is, when $\alpha \geq y_u$ in some node $u$, then $\alpha$ corresponds to the label 1 and $y_u$ to 0 and the other way around:

**Proposition 11.29.** If the energy $E$ is submodular, the minimization problem (11.59) is permuted submodular. It is submodular when the original order of labels is kept.

However, the $\alpha$-expansion algorithm can be used also with non-submodular energies. It is sufficient that each pairwise cost function defines a metric:

**Definition 11.30.** For any set $L$ the function $f: L \times L \to \mathbb{R}$ is called *a semi-metric* if for all $x, x' \in L$ it holds that: (i) $f(x, x') \geq 0$; (ii) $f(x, x') = 0$ iff $x = x'$; (iii) $f(x, x') = f(x', x)$.

**Definition 11.31.** Function $f: L \times L \to \mathbb{R}$ is called *a metric* if it is a semi-metric and additionally the triangle inequality holds: $f(x, x') + f(x', x'') \geq f(x, x'')$, $\forall x, x', x'' \in L$.

**Proposition 11.32.** Let $\theta_{uv}(s, t)$ be a metric, then the binary optimization problem (11.59) is submodular when the label $\alpha$ is assumed to be lower than any other label. In other words, the label $\alpha$ is substituted with 0 and $y_u$ with 1.

*Proof.* The submodularity follows from the inequality based on the metric properties of $\theta_{uv}$:

$$\underbrace{\theta_{uv}(\alpha, \alpha)}_{=0} + \underbrace{\theta_{uv}(y_u, y_v) - \theta_{uv}(y_u, \alpha) - \theta_{uv}(\alpha, y_v)}_{\leq 0} \leq 0. \qquad (11.60)$$

$\square$

**Example 11.33.** It is easy to show that the following pairwise potentials are metrics:

- Potts potentials $\theta_{uv}(s, t) = [\![s \neq t]\!]$.

- Truncated $\ell_1$-norm: $\theta_{uv}(s, t) = \min\{|s - t|, M\}$.

- Truncated $\ell_2$-norm: $\theta_{uv}(s, t) = \min\{|s - t|^2, M\}$.

From Example 11.33 it follows that a metric need not be submodular. Conversely, submodularity does not imply metric properties. That is, if $\theta_{uv}$ is metric, then for any constant $a$ the cost function $\theta_{uv} + a$ is not a metric anymore. However, if $\theta_{uv}$ is submodular, $\theta_{uv} + a$ is submodular as well.

## 11.5.2 Optimality guarantees of the $\alpha$-expansion algorithm

The $\alpha$-expansion Algorithm 17 does not solve energy minimization optimally (even submodular) in general. It performs in a certain sense a "local search", similarly to the ICM algorithm considered in Chapter 8. The crucial difference between these two methods is in the size of the "vicinity" of the current labeling, where the search is performed. This vicinity for the $\alpha$-expansion algorithm is much larger than the one for ICM, which explains the typically much better quality of the approximate solutions obtained by $\alpha$-expansion.

The following theorem provides the worst-case optimality bound for Algorithm 17. In practice, however, the algorithm attains significantly better energy values than it is guaranteed by the bound.

**Theorem 11.3** (Multiplicative optimality bound Boykov *et al.* [16])**.** Let $\theta_{uv}$ be a metric for all $uv \in \mathcal{E}$, $y^*$ be an optimal solution and $y'$ be the labeling obtained by Algorithm 17. Then for

$$c := \max_{uv \in \mathcal{E}} \frac{\max_{y_v \neq y_u} \theta_{uv}(y_u, y_v)}{\min_{y_v \neq y_u} \theta_{uv}(y_u, y_v)} \tag{11.61}$$

it holds that:

$$2c \left( \sum_{v \in \mathcal{V}} \theta_v(y_v^*) + \sum_{uv \in \mathcal{V}} \theta_{uv}(y_u^*, y_v^*) \right) \geq \sum_{v \in \mathcal{V}} \theta_v(y_v') + \sum_{uv \in \mathcal{V}} \theta_{uv}(y_u', y_v') . \tag{11.62}$$

In particular, for a Potts model the constant $c$ is equal to 1 and for truncated linear and pairwise costs it is equal to the truncation value $a$, see Example 11.15.

### 11.5.3   $\alpha\beta$-swap algorithm

The triangle inequality of the metric properties (see Definition 11.31) can still be quite restrictive for certain applications. Therefore, we consider another popular min-cut-based move-making algorithm, applicable to semi-metric pairwise costs. This algorithm is called $\alpha\beta$-swap. On each iteration it considers only two labels of the current labeling: $\alpha$ and $\beta$ and finds an optimal configuration of those given that the rest are fixed.

Let the label set $\mathcal{Y}^{\alpha,\beta}(s)$ be defined as

$$\mathcal{Y}^{\alpha,\beta}(s) = \begin{cases} \{s\}, & s \in \mathcal{Y}_0 \backslash \{\alpha, \beta\} \\ \{\alpha, \beta\}, & s \in \{\alpha, \beta\} . \end{cases} \tag{11.63}$$

Algorithm 18 starts with an arbitrary labeling $y^0$ and tries to change to a new labeling with better energy on each iteration. Since the initial labeling can be an outcome of the binary optimization (5), the algorithm is monotonous, i.e. $E(y^{t+1}) \leq E(y^t)$.

Similar to $\alpha$-expansion, the $\alpha\beta$-swap algorithm can be applied to submodular functions, however it is also applicable to semi-metric, as we mentioned above. This follows from the fact, that on each iteration only labels $\alpha$ and $\beta$ are active. The following submodularity condition for this restricted problem follows directly from the semi-metric conditions:

---

**Algorithm 18** $\alpha\beta$-swap

---

1: **Init:** $y^0$, $t = 0$
2: **repeat**
3:   // The line below defines $\binom{\mathcal{Y}_0}{2}$ iterations
4:   **for** $(\alpha, \beta) \in \mathcal{Y}_0^2$ **do**
5:     $y^{t+1} = \arg \min\limits_{y \,:\, \forall u \; y_u \in \mathcal{Y}^{\alpha,\beta}(y_u^t)} E(y)$
6:     $t = t + 1$
7:   **end for**
8: **until** $E(y^t) = E(y^{t-\binom{\mathcal{Y}_0}{2}})$
9: **return** $y^t$

---

$$\underbrace{\theta_{uv}(\alpha, \beta)}_{>0} + \underbrace{\theta_{uv}(\beta, \alpha)}_{>0} > \underbrace{\theta_{uv}(\alpha, \alpha)}_{=0} + \underbrace{\theta_{uv}(\beta, \beta)}_{=0} \tag{11.64}$$

Each outer iteration of $\alpha\beta$-swap has a quadratic complexity with respect to the number of labels, which is more than the linear complexity of $\alpha$-expansion. However, since only those nodes participate in optimization, which have labels $\alpha$ or $\beta$, each min-cut computation is cheaper than those of $\alpha$-expansion.

Contrary to $\alpha$-expansion, $\alpha\beta$-swap does not provide any optimality bound. Moreover, the following example shows that the algorithm can return arbitrary bad results. However, in practice it often performs only somewhat worse than $\alpha$-expansion.

**Example 11.34** (Arbitrary bad optimality bound [16])**.** Consider an graphical model consisting of 3 nodes: $\mathcal{V} = \{1, 2, 3\}$. The edges connect the nodes 1,2 and 2,3: $\mathcal{E} = \{\{1, 2\}, \{2, 3\}\}$. Each node is associated with 3 labels: $\mathcal{Y}_0 = \{a, b, c\}$. The values of unary costs are defined by the following table:

|       | **1** | **2** | **3** |
|-------|-------|-------|-------|
| **a** | 0     | $K$   | $K$   |
| **b** | $K$   | 0     | $K$   |
| **c** | 2     | 2     | 0     |

Pairwise costs are homogeneous ($\theta_{12} = \theta_{23} = \theta$) and equal to $\theta(a, b) = \theta(b, c) = \frac{K}{2}$ and $\theta(a, c) = K$. It is easy to see that configuration $(a, b, c)$ is a fix-point for Algorithm 18. Its energy is $K$, while the optimal configuration is $(c, c, c)$ and has energy 4. Since $K$ is not bounded from above, the energy of a fix-point of the $\alpha\beta$-swap algorithm can be arbitrary worse than the energy of the optimal labeling.

**Exercise 11.35.** Show that $\alpha$-expansion Algorithm 17 obtains the optimal solution when started from the fix-point of $\alpha\beta$-swap in Example 11.34.

## 11.6    Bibliography and further reading

An overview of combinatorial algorithms for general submodular minimization can be found in [74].

The first work formulating the submodular binary energy minimization problem as min-cut was probably [29]. The classical paper, which describes reduction of the binary energy minimization problem to the min-cut problem is the seminal work [54].

Somewhat different transformations of multi-label problems into the min-cut problem were proposed independently in [2] and [36] for convex pairwise potentials (like in Exercise 11.14). Later this result was generalized in [99] to the submodularity-preserving transformation of general multi-label problems into binary ones. An algorithm determining whether a given problem is permuted submodular is due to [98]. Energies with metric potentials are analyzed in details in [45]. The $\alpha$-expansion and $\alpha\beta$-swap algorithms were proposed and analyzed in [16]. Recently some sufficient conditions guaranteeing global optimality of the $\alpha$-expansion algorithm were provided in [66]. Min-cut-based move-making methods attaining improved multiplicative bounds for truncated convex potentials were proposed and analyzed by [65]. Finally, the comparative studies [126] and [41] show how the considered $\alpha$-expansion and $\alpha\beta$-swap algorithms compare to some of the dual methods considered in the previous chapters.

# 12

## Relaxed Binary Energy as $st$-Min-Cut

In Chapter 11 we have shown an equivalence between min-$st$-cut and binary MAP-inference problems. That is, each binary MAP-inference problem can be reduced to min-$st$-cut of asymptotically the same size and the other way around: Each min-$st$-cut is reducible to a binary MAP-inference problem.

In particular, we have shown that submodular MAP-inference problems correspond to min-cut problems with non-negative edge weights. Hence, both problems have polynomial complexity and are moreover very efficiently solvable.

However, the non-submodular MAP-inference problems can be reduced only to min-cut with some edge weights being negative. This conforms to the fact that both problems are $\mathcal{NP}$-hard. Therefore, in order to deal with this kind of problems, one has to resort to approximate algorithms such as those addressing its LP relaxation.

Interestingly, the local polytope relaxation of the binary MAP-inference problems can be optimized very efficiently by reducing it to the min-$st$-cut with non-negative edge weights. This contrasts with the local polytope relaxation of the multi-label problems, which is known to be as difficult as general linear programs, see Theorem 4.1. In this chapter we will provide the above reduction. Additionally, we will analyze the local polytope corresponding to binary MAP-inference problems. In particular, we will show that it possesses the so called *persistency* or *partial* optimality property, which often allows us to reduce the size of the non-relaxed problem and therefore makes it possible to obtain its exact solution with combinatorial solvers. Another important use of this property are primal heuristics similar to the $\alpha$-expansion algorithm (see Chapter 11), however, in this case also applicable to problems with *arbitrary* pairwise costs. We address this type of algorithms in §12.3.2.

## 12.1 Half-integrality of local polytope

We start with the most important property of the local polytope of binary MAP-inference problems (we will refer to it as *binary local polytope*, although it is not a commonly used term in the literature). This is the characterization of its vertices, which all turn out to have the *half-integer* coordinates, i.e. for any vertex $\mu$ of the binary local polytope it holds that $\mu \in \{0, 1, \frac{1}{2}\}^{\mathcal{I}}$.

For the following statement recall the definition of the set $\mathcal{L}(\phi)$ from (6.10):

$$\mathcal{L}(\phi) = \{\mu \in \mathcal{L} \colon \mu_w(s) = 0 \text{ if } \theta_w^\phi(s) > \min_{s' \in \mathcal{Y}_w} \theta_w^\phi(s'), \ w \in \mathcal{V} \cup \mathcal{E}, \ s \in \mathcal{Y}_w\},$$

(12.1)

which is the set of all vectors in the local polytope such that their non-zero coordinates correspond only to the locally minimal entries of the cost vector.

**Theorem 12.1** ([100, 140, 53]). Let $(\mathcal{G}, \mathcal{Y}_{\mathcal{V}}, \theta)$ be a binary graphical model. Let $\phi \in \mathbb{R}^{\mathcal{J}}$ be a reparametrization and $\theta^\phi$ be the corresponding reparametrized costs such that the corresponding arc-consistency closure is not empty, i.e. $\mathrm{cl}(\mathrm{mi}[\theta^\phi]) \neq \bar{0}$. Then $\mathcal{L}(\phi) \cap \{0, \frac{1}{2}, 1\}^{\mathcal{I}} \neq \emptyset$.

In other words, there is a half-integral point in the local polytope, such that its non-zero coordinates correspond to locally minimal labels and label pairs of the cost vector $\theta^\phi$.

Note that due to Proposition 6.2, Theorem 12.1 implies that a non-empty arc-consistency closure is sufficient for dual optimality in case of binary problems.

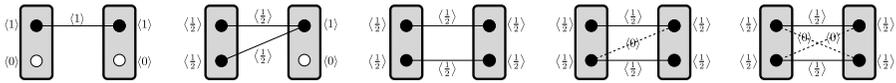

**Figure 12.1:** Possible configurations of the arc-consistency closure $\xi = \mathrm{cl}(\mathrm{mi}[\theta^\phi])$ up to the labels flip, see proof of Theorem 12.1. Black circles for a label $s$ in a node $u$ corresponds to $\xi_u(s) = 1$, white circle to $\xi_u(s) = 0$. Let $n_u$ be the number of black circles in the node $u$. Both solid and dashed lines connecting labels $s$ and $t$ in nodes $u$ and $v$ correspond to the values $\xi_{uv}(s, t) = 1$, otherwise $\xi_{uv}(s, t) = 0$. Angular brackets $\langle \rangle$ are used for coordinates of primal solutions $\mu$. To construct a feasible relaxed primal solution $\mu$ assign $\mu_u(s) = \frac{1}{n_u}$, if label $s$ corresponds to the black circle and $\mu_u(s) = 0$ if to a white one. Similarly, assign $\mu_{uv}(s, t) = \frac{1}{\max\{n_u, n_v\}}$, if $(s, t)$ corresponds to a solid line, otherwise assign $\mu_{uv}(s, t) = 0$. The assigned coordinates of $\mu$ are shown in the corner brackets $\langle \cdot \rangle$.

*Proof.* Consider $\xi = \mathrm{cl}(\mathrm{mi}[\theta^\phi])$. Denote $n_u$ the number of non-zero entries in $\xi_u$, i.e. $n_u = \xi_u(0) + \xi_u(1)$, $u \in \mathcal{V}$. Consider three incident factors $\xi_u$, $\xi_{uv}$ and $\xi_v$ for any $uv \in \mathcal{E}$. Up to the label flipping they can have only 5 possible configurations, depicted in Figure 12.1, where black circles, solid and dashed lines correspond to coordinates of $\xi$ equal to 1 and the rest to 0.

Initialize $\mu := 0$. Assign $\mu_u(s) = \frac{1}{n_u}$ for $u \in \mathcal{V}$. Assign $\mu_{uv}(s, t) = \frac{1}{\max\{n_u, n_v\}}$, if $(s, t)$ corresponds to a solid line. It is straightforward to check that such $\mu$ belong to the binary local polytope. It contains only half-integer coordinates and its non-zero coordinates correspond to the non-zero coordinates of $\xi$. Concluding, $\mu \in \mathcal{L}(\phi) \cap \{0, \frac{1}{2}, 1\}$, which finalizes the proof. □

**Corollary 12.1.** All vertices of the binary local polytope have half-integral coordinates.

*Proof.* Let $\theta$ be a cost vector such that the relaxed problem $\min_{\mu \in \mathcal{L}} \langle \theta, \mu \rangle$ has a unique solution. Let also $\phi$ be an optimal reparametrization, which implies that $\mathrm{cl}(\mathrm{mi}[\theta^\phi]) \neq \bar{0}$. According to Theorem 12.1 the relaxed solution is half-integral, which proves the statement according to Definition 3.19. □

## 12.2  LP relaxation as min-cut

A goal of this section is to reduce the LP relaxation of the binary MAP-inference problem to the min-$st$-cut problem. Such a transformation would imply existence of a polynomial practically efficient algorithm for the LP relaxation.

**Separability of the local polytope**    As noted in §6.1, the local polytope (for arbitrary, not only binary problems) can be written in the following compact form (see (6.2)):

$$
\mathcal{L} = \begin{cases} \mu_u \in \Delta^{\mathcal{Y}_u}, & \forall u \in \mathcal{V} \\ \mu_{uv} \in \Delta^{\mathcal{Y}_{uv}}, & \forall uv \in \mathcal{E} \\ \sum_{t \in \mathcal{Y}_v} \mu_{uv}(s,t) = \mu_u(s), & \forall u \in \mathcal{V}, \ v \in \mathcal{N}_b(u), \ s \in \mathcal{Y}_u . \end{cases} \tag{12.2}
$$

Let us fix the "unary" vectors $\mu_u$ for all nodes $u \in \mathcal{V}$. Let $\mathcal{L}_{uv}(\mu_u, \mu_v)$ be the set of all "pairwise" vectors $\mu_{uv}$, such that the whole vector $\mu$ belongs to $\mathcal{L}$, i.e.:

$$
\mathcal{L}_{uv}(\mu_u, \mu_v) = \{\mu_{uv} \in \Delta^{\mathcal{Y}_{uv}} : \sum_{t \in \mathcal{Y}_v} \mu_{uv}(s,t) = \mu_u(s) \quad \forall s \in \mathcal{Y}_u\} \tag{12.3}
$$

This separable structure implies that given the "unary" vectors $\mu_u$, $u \in \mathcal{V}$, optimization with respect to the "pairwise" ones decomposes into $|\mathcal{E}|$ small independent optimization problems, one for each pairwise factor:

$$
\min_{\mu \in \mathcal{L}} \langle \theta, \mu \rangle = \min_{\substack{\mu_u \in \Delta^{\mathcal{Y}_u} \\ u \in \mathcal{V}}} \sum_{uv \in \mathcal{E}} \min_{\mu_{uv} \in \mathcal{L}_{uv}(\mu_u, \mu_v)} \langle \theta_{uv}, \mu_{uv} \rangle . \tag{12.4}
$$

We will use this fact below to reduce the binary local polytope relaxation to the min-$st$-cut problem.

**Idea of the transformation**    As follows from §11.3.4, binary energy can be represented as an oriented form. The submodular pairwise costs $\theta_{uv}$ can be turned into diagonal form (see Figure 11.8(b)) equivalent to a quadratic monomial $c_{u,v}\bar{y}_u y_v$ with non-negative coefficients $c_{u,v}$, whereas for supermodular costs these coefficients are non-positive.

Alternatively, supermodular pairwise costs can be turned into the parallel form, as shown in Figure 11.8(c), corresponding to the quadratic monomial $c_{u,v}y_u y_v$, see Figure 11.7(c). Although the coefficient of this monomial is non-negative for supermodular costs, they can not be directly represented as edges of the min-$st$-cut graph. To do so, we will introduce an additional binary vector $z$ such that $y = \bar{z}$ and rewrite $c_{u,v}y_u y_v = \frac{1}{2}c_{u,v}y_u \bar{z}_v + \frac{1}{2}c_{u,v}\bar{z}_u y_v$. The right-hand side is an oriented form, and, therefore, can be represented as edges of a min-cut graph.

**Binary energy as oriented form with a double number of variables**
In general, we will assume that the binary energy is transformed into a quadratic pseudo-boolean function with *non-negative* coefficients. By adding an additional binary vector $z$ such that $y = \bar{z}$ we rewrite each term of the pseudo-boolean function as a symmetric oriented form containing both $y$ and $z$ vectors.

- Unary costs are transformed as follows

$$c_u y_u \to \frac{1}{2}c_u y_u + \frac{1}{2}c_u \bar{z}_u \, , \tag{12.5}$$

$$c_u \bar{y}_u \to \frac{1}{2}c_u \bar{y}_u + \frac{1}{2}c_u z_u \, . \tag{12.6}$$

- Submodular pairwise terms are transformed as follows

$$c_{u,v}\bar{y}_u y_v \to \frac{1}{2}c_{u,v}\bar{y}_u y_v + \frac{1}{2}c_{u,v}z_u \bar{z}_v \, . \tag{12.7}$$

- Supermodular pairwise terms are transformed as follows

$$c_{u,v}y_u y_v \to \frac{1}{2}c_{u,v}y_u \bar{z}_v + \frac{1}{2}c_{u,v}\bar{z}_u y_v \, . \tag{12.8}$$

As noted above, we assume all coefficients $c_u$, $c_v$ and $c_{u,v}$ above to be non-negative.

Let $E$ be an energy and $g(y, z)$ be the oriented form, constructed according to (12.5)-(12.8). By construction $E(x) = g(y, \bar{y})$, therefore it holds that:

$$\min_{y,z \in \{0,1\}^{\mathcal{V}}} g(y, z) \leq \min_{\substack{y,z \in \{0,1\}^{\mathcal{V}} \\ z=\bar{y}}} g(y, z) = \min_{y \in \{0,1\}^{\mathcal{V}}} E(y) \, . \tag{12.9}$$

The following theorem relates the first term in (12.9) to the LP relaxation of the energy minimization:

**Theorem 12.2.** Let $(\mathcal{G}, \mathcal{Y}_\mathcal{V}, \theta)$ be a binary MAP-inference problem and its submodular and supermodular pairwise terms be represented in the "diagonal" and "parallel" forms respectively. Let $g$ be constructed according to (12.5)-(12.8). Then it holds that

$$\min_{y,z \in \{0,1\}^\mathcal{V}} g(y,z) = \min_{\mu \in \mathcal{L}} \langle \theta, \mu \rangle . \qquad (12.10)$$

Moreover, a minimizer $\mu'$ of the right-hand-side can be constructed from a minimizer $(y', z')$ of the left-hand-side as follows:

$$\mu'_u(1) = \begin{cases} y_u, & y_u = \bar{z}_u \\ \frac{1}{2}, & y_u \neq \bar{z}_u , \end{cases} \qquad (12.11)$$

$$\mu'_u(0) = 1 - \mu'_u(1) \qquad (12.12)$$

$$\mu'_{uv} = \underset{\mu_{uv} \in \mathcal{L}_{uv}(\mu'_u, \mu'_v) \cap \{0, 1, \frac{1}{2}\}^{\mathcal{Y}_{uv}}}{\arg\min} \langle \theta_{uv}, \mu_{uv} \rangle , \qquad (12.13)$$

where $\mathcal{L}_{uv}$ is defined by (12.3).

*Proof.* The proof of the theorem is quite straightforward: except for a special case noted separately, for each factor $w \in \mathcal{V} \cup \mathcal{E}$ we will show that if $\mu'_w$ is constructed from *arbitrary* $(y'_w, z'_w)$ the relaxed energy $\langle \theta_w, \mu'_w \rangle$ is equal to the corresponding term of $g(y', z')$.

1. If $y_w = \bar{z}_w$ the condition above holds by construction.

2. Unary terms (12.5)-(12.6). If $y_u \neq \bar{z}_u$ then the corresponding term of $g$ is equal to $\frac{1}{2} c_u$, which is equal to $\langle \theta_u, \mu'_u \rangle$.

3. Submodular terms (12.7). Let $y_v = \bar{z}_v$ and $y_u \neq \bar{z}_u$.

   (a) if $y_v = \bar{z}_v = 1$ the corresponding term (12.7) of $g$ turns into $\frac{1}{2} c_{u,v}(y_u + \bar{z}_u)$, which is $\frac{1}{2} c_{u,v}$ (independently on the values of $y_u \neq \bar{z}_u$) and is equal to $\langle \theta_{uv}, \mu'_{uv} \rangle$ according to Figure 12.2(a).

   (b) if $x_v = \bar{y}_v = 0$ the corresponding term (12.7) of $g$ is equal to 0 (independently on the values of $x_u \neq \bar{y}_u$), which is equal to $\langle \theta_{uv}, \mu'_{uv} \rangle$, see Figure 12.2(a) for illustration.

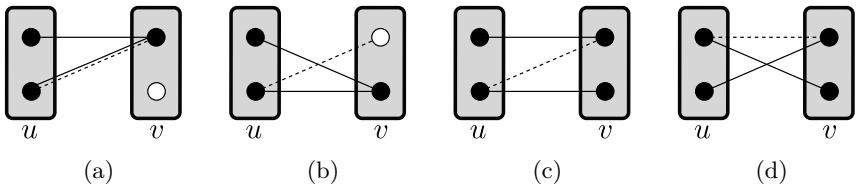

**Figure 12.2:** Optimal primal solutions for pairwise factors, see the proof of Theorem 12.2. Solid lines correspond to the value $\frac{1}{2}$ of the primal variable $\mu$, black labels - either to 1 (in case of a single black label per node) or to $\frac{1}{2}$ (for two black labels per node). Dashed lines correspond to the diagonal and parallel forms of the submodular and supermodular pairwise factors respectively.

4. Supermodular terms (12.8) with $y_v = \bar{z}_v$ and $y_u \neq \bar{z}_u$ turn into 3 by changing $y_v$ to $\bar{z}_v$.

5. Submodular (12.7) and supermodular (12.8) terms. Let $y_v \neq \bar{z}_v$ and $y_u \neq \bar{z}_u$. This is the mentioned above exceptional case, it holds that $\langle \theta_{uv}, \mu'_{uv} \rangle = 0$ (see Figure 12.2(c) and 12.2(d) for illustration). The same 0 value the terms (12.7) and (12.8) have when $y_v = z_v = y_u = z_u = 0$. There are combinations of values of $y_v$ and $z_u$ (e.g. $y_v = z_v = 1$ and $y_u = z_u = 0$ for (12.7)), corresponding to bigger values of the considered terms, however they are never optimal. This is due to the fact that as soon as $y_u \neq \bar{z}_u$ the change from $y_u = z_u = 1$ to $y_u = z_u = 0$ can either decrease the value of the corresponding components of $g(x, y)$ (as in the current case) or leave it unchanged (as in all other situations described above). This precisely means that considering $y_u = z_u = 0$ a soon as $y_u \neq \bar{z}_u$ gives smaller or equal value of $g$.

$\square$

## 12.3 Persistency of the binary local polytope

In this section we get back to the question, how integer coordinates of a relaxed solution are related to an exact non-relaxed solution. As shown in Example 4.6, in general, one can not assume that the integer coordinates are kept unchanged in some exact non-relaxed solution.

However, for the binary local polytope this property indeed holds, as we show below.

First, we will give the definition of the considered property and a sufficient condition for its existence for general graphical models. Afterwards we concentrate on the binary MAP-inference and show that the sufficient condition always holds in this case.

**Persistency definition and criterion for general graphical models**

**Definition 12.2** (Persistency of an integer coordinate). Let

$$\mu' \in \arg\min_{\mu \in \mathcal{L}} \langle \theta, \mu \rangle \tag{12.14}$$

be a relaxed solution of the MAP-inference problem for a graphical model $(\mathcal{G}, \mathcal{Y}_\mathcal{V}, \theta)$. For some $u \in \mathcal{V}$ its coordinate $\mu'_u(s)$ is called *(weakly) persistent* or *(weakly) partially optimal*, if (i) $\mu'_u(s) \in \{0,1\}$ and (ii) there exists an exact integer solution

$$\mu^* \in \arg\min_{\mu \in \mathcal{L} \cap \{0,1\}^\mathcal{I}} \langle \theta, \mu \rangle \tag{12.15}$$

such that $\mu^*_u(s) = \mu'_u(s)$. It is *strongly persistent* if this property holds for *any* exact integer solution $\mu^*$.

A somewhat more general definition of persistency can also be given without relation to any relaxation.

**Definition 12.3** (Persistency of a partial labeling). A partial labeling $y' \in \mathcal{Y}_\mathcal{A}$ on a subset $\mathcal{A} \subseteq \mathcal{V}$ is *(weakly) partially optimal* or *(weakly) persistent* if $y' = y^*|_\mathcal{A}$ for some $y^* \in \arg\min_{y \in \mathcal{Y}_\mathcal{V}} \langle \theta, \delta(y) \rangle$.

The following proposition formulates a sufficient persistency condition:

**Proposition 12.4** (Persistency criterion [124]). A partial labeling $y' \in \mathcal{Y}_\mathcal{A}$ is persistent, if for all $\tilde{y} \in \mathcal{Y}_{\mathcal{V} \setminus \mathcal{A}}$ the following holds

$$y' \in \arg\min_{y \in \mathcal{Y}_\mathcal{A}} E((y, \tilde{y})), \tag{12.16}$$

where $(\cdot, \cdot)$ stands for the labeling concatenation, i.e.

$$(y, \tilde{y}) = \begin{cases} y_u, & u \in \mathcal{A}, \\ \tilde{y}_u, & \mathcal{V} \setminus \mathcal{A}. \end{cases} \tag{12.17}$$

*Proof.* Consider the equation

$$\min_{y \in \mathcal{Y}_\mathcal{V}} E(y) = \min_{\tilde{y} \in \mathcal{Y}_{\mathcal{V} \setminus \mathcal{A}}} \min_{y \in \mathcal{Y}_\mathcal{A}} E((y, \tilde{y})) . \tag{12.18}$$

Let $\tilde{y} \in \mathcal{Y}_{\mathcal{V} \setminus \mathcal{A}}$ be such that it leads to a minimal value on the right hand side of (12.18). Then $\tilde{y}$ is part of an optimal solution. By the assumption (12.16), $y'$ is an optimal solution to the inner minimization problem of (12.18), hence $(y', \tilde{y})$ minimizes $E$. $\square$

### 12.3.1 Persistency of the binary local polytope

Let $\mathcal{I}$ be the set of indexes associated with primal solutions of a binary graphical model. For an arbitrary vector $\mu \in \mathbb{R}^\mathcal{I}$ let $\mathrm{Int}(\mu) = \{u \in \mathcal{V} \colon \mu_u \in \{0, 1\}^2\}$ be the set of nodes corresponding to integer coordinates of $\mu$.

**Theorem 12.3.** Let $\mu' \in \arg\min_{\mu \in \mathcal{L}} \langle \theta, \mu \rangle$ be a relaxed solution of a binary MAP-inference problem. Then $\mu'|_{\mathrm{Int}(\mu')}$ is persistent.

*Proof.* The proof is illustrated in Figure 12.3. Let us denote $\mathcal{A} := \mathrm{Int}(\mu')$. Let also $y' \in \mathcal{X}_\mathcal{A}$ be selected such that $\delta(y') = \mu'|_\mathcal{A}$. The criterion of Proposition 12.4 does not depend on reparametrization of $\theta$. W.l.o.g. we will assume the latter to be (dual) optimal and such that for all $w \in \mathcal{V} \cup \mathcal{E}$ it holds that

$$\min_{s \in \mathcal{Y}_w} \theta_w(s) = 0 . \tag{12.19}$$

Let us consider the persistency criterion (12.16). For a fixed $\tilde{y} \in \mathcal{Y}_{\mathcal{V} \setminus \mathcal{A}}$ it holds that

$$\arg\min_{y \in \mathcal{Y}_\mathcal{A}} E((y, \tilde{y})) = \arg\min_{y \in \mathcal{Y}_\mathcal{A}} \left( E_\mathcal{A}(y) + \sum_{\substack{uv \in \mathcal{E} \colon \\ u \in \mathcal{A}, v \in \mathcal{V} \setminus \mathcal{A}}} \theta_{uv}(y_u, \tilde{y}_v) \right) \geq 0 , \tag{12.20}$$

where the energy $E_\mathcal{A}$ is restricted to the induced subgraph of $\mathcal{A}$ and the inequality holds due to (12.19).

Since $\mu'$ is primal (relaxed) optimal and $\theta$ is a dual optimal, it holds that $E_\mathcal{A}(y') = 0$. Moreover, since the solution $\mu'$ outside $\mathcal{A}$ is non-integer, it holds that $\mu'_u(s) > 0$ for all $s \in \mathcal{X}_u$, $u \in \mathcal{V} \setminus \mathcal{A}$. Therefore, due to the local polytope constraints, $\mu_{uv}(y'_u, y_v) > 0$ for all $u \in \mathcal{A}$, $x_u \in \mathcal{Y}_u$. Statement 5 of Proposition 6.2 implies $\theta_{uv}(y'_u, y_v) = 0$, in particular $\theta_{uv}(y'_u, \tilde{y}_v) = 0$. Hence, we obtain

$$\left( E_{\mathcal{A}}(y') + \sum_{\substack{uv \in \mathcal{E}_{\partial \mathcal{A}}: \\ u \in \mathcal{A}}} \theta_{uv}(y'_u, \tilde{y}_v) \right) = 0 \,. \qquad (12.21)$$

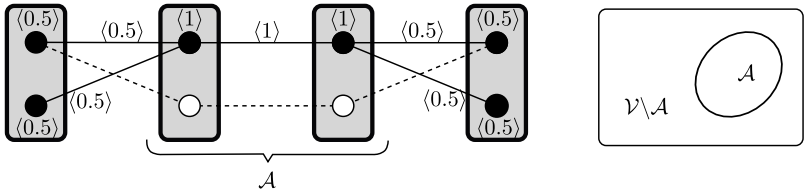

**Figure 12.3:** Illustration to the proof of Theorem 12.3. Black labels and solid lines correspond to the non-zero values of the primal variables of a given relaxed solution $\mu'$. An example of such values is written in the angular brackets $\langle \rangle$. These non-zero values correspond to locally optimal labels and label pairs in an optimal reparametrization $\theta$. We assume these locally optimal labels and label pairs to have weight 0. Dashed lines correspond to an arbitrary selected alternative labeling. Since it may contain locally non-minimal labels and label pairs, it is not minimal in general.

Comparing to (12.20) we obtain $y' \in \arg\min_{y \in \mathcal{Y}_{\mathcal{A}}} E((y', \tilde{y}))$, which finalizes the proof. $\qquad \square$

**Remark 12.5.** An important property of $y'$ follows from the fact that it minimizes (12.20). It reads

$$E((y, \tilde{y})) \geq E((y', \tilde{y})), \quad \forall y \in \mathcal{Y}_{\mathrm{Int}(\mu)} \text{ and } \tilde{y} \in \mathcal{Y}_{\mathcal{V} \setminus \mathrm{Int}(\mu)} \,, \qquad (12.22)$$

where $\mu$ is the relaxed solution. This property constitutes a basis for a powerful class of min-cut-based primal heuristics applicable to general graphical models, which we consider in §12.3.2.

### 12.3.2 Fusion moves

The property (12.22) allows algorithms like $\alpha$-expansion (or $\alpha\beta$-swap) to work with arbitrary pairwise potentials. Consider the $\alpha$-expansion Algorithm 17. Let $\mathcal{L}(y^t, \alpha)$ be the local polytope of the auxiliary binary problem in line 4 of the algorithm, such that each node $u$ contains only two labels, $\alpha$ and $y_u^t$. Let

$$\mu^* \in \underset{\mu \in \mathcal{L}(y^t, \alpha)}{\arg\min} \langle \theta, \mu \rangle \qquad (12.23)$$

be a solution of such relaxed binary problem. Construct the new labeling as

$$\forall u \in \mathcal{V} \quad y_u^{t+1} := \begin{cases} \alpha, & \mu_u^*(\alpha) = 1, \\ y_u^t, & \text{otherwise} \,. \end{cases} \qquad (12.24)$$

The property (12.22) implies that $E(y^{t+1}) \leq E(y^t)$, which guarantees monotonicity of the algorithm. Since the new current solution $y^{t+1}$ is "fused" from the old one $y^t$ and the *proposal* consisting of a constant labeling $\alpha$, the modified algorithm is usually referred to as *fusion moves*.

Note several important differences between Algorithm 17 and the fusion moves:

- Fusion moves do not require any metric properties of the initial multi-label problem to be applicable, since the key computational step – the relaxed problem (12.23) is efficiently solvable for any costs.

- Due to the above property there is no need to restrict the proposal labeling to a constant one. The best results are often obtained, when an application-specific knowledge is used for proposal generation. E.g. for stereo reconstruction of smooth surfaces the labelings corresponding to such surfaces can be used as proposals. The constant proposals like in the $\alpha$-expansion algorithm correspond to the piece-wise constant compactness prior, more typical for the image segmentation problem. Note that if the set of proposals becomes infinite, or at least exponentially large, there is no clear stopping criterion for the fusion moves algorithm anymore, except by prescribing the number of algorithm's iterations in advance.

- Similar to $\alpha\beta$-swap, the fusion moves algorithm does not have any optimality guarantees in general, which conforms to the fact that such guarantees can not be obtained by polynomial algorithms for general graphical models. In particular, there is no guarantee that the set of nodes with integer coordinates $\text{Int}(\mu^*)$ is non-empty for the solution of any of the auxiliary problems (12.23).

## 12.4   Bibliography and further reading

Although there is a fair amount of literature on quadratic pseudo-Boolean optimization and its relation to binary labeling problems and min-cut, the presentation of the chapter is quite new. Typically the proof of Theorem 12.2 is done for dual objectives, which is more practical for applications. The primal exposition, like the one in this monograph, seems to be more illustrative and intuitive. Moreover, it is more straightforwardly related to persistency. The typical proof for the persistency property is done by reducing the problem to the vertex cover problem and operating by known results in that domain. It turned out that this result can easily be proven directly, see Theorem 12.3, although similar proofs exist in the literature, e.g. [53].

Related works on quadratic pseudo-boolean optimization are [13, 14, 34]. Extensions and applications of the persistency property of the binary local polytope are discussed in [93, 50, 46].

In a recent line of research [125, 108, 124, 110, 120, 113, 109, 111] persistency of multi-label problems have been studied and efficient polynomial algorithms for its computation were proposed.

The fusion moves were proposed and discussed in [68, 52, 69].

# Acknowledgements

The initial idea for such a monograph belongs to my PhD supervisor, Prof. Michail Schlesinger in whose Image Processing and Recognition Department I spent ten years from 1999 until 2009. He was one of the first, who have gave rise to the powerful modeling technique known now as probabilistic graphical models. Already in the 70s he published a number of papers about it, formulated the underlying combinatorial problem and analyzed its linear programming relaxation. As the power of this tool was recognized by the wider pattern recognition and computer vision community in the early 2000s, and a number of other researchers got into the field, he suggested to write a book about the topic. Such a work was obviously beyond the capabilities and skills of a PhD student as I was at that time.

After moving to Germany in 2009 I worked mainly in the field this monograph is related to, and in 2015 I started to give lectures on graphical models. At this point another person important for this work appeared in my life. Prof. Carsten Rother, in whose Computer Vision Lab Dresden I was working, suggested to write a book on the topic of my lectures. It was not just a void suggestion, as my boss he fully supported me in this work, from the very beginning until the end. Without exaggeration, this monograph would have never been written without him.

There were many other people, who motivated me and helped with writing. First of all, these are my PhD students, Stefan Haller, Siddharth Tourani and Lisa Kruse. They have done a huge amount of work by proofreading the manuscript, preparing the lectures and exercises for my course, implementing the algorithms and performing their experimental evaluation. Yet, Lisa's contribution to the monograph is outstanding, even when compared to others. She has found a number of mathematical flaws in my writing and suggested corrections to them. Moreover, a number of proofs in the monograph belong to her, since my own variant was notably less rigorous and understandable. Without her work this monograph would have been much more erroneous and way less readable.

I also appreciate the feedback given by my colleagues from other labs and institutions: Florian Jug and Mangal Prakash from the Max Planck Institute of Molecular Cell Biology and Genetics, as well as Jan Kuske, Fabrizio Savarino and Ruben Hühnerbein from Heidelberg University.

A selection of the content of the monograph as well as its presentation appeared as a result of joint work and fruitful discussions with my colleagues actively working in the field. In particular, I would like to thank Paul Swoboda, Oleksandr Shekhovtsov, Tomas Werner, Boris Flach, Dmitrij Schlesinger, Christoph Schnörr, Jörg Kappes, Stefan Schmidt, Alexander Kirillov and Björn Andres. Moreover, corrections of Oleksandr and Tomas notably improved readability of the monograph and enriched its reference list.

Finally, students of the Technical University of Dresden and the Heidelberg University, which attended my course since 2015 helped to improve the book by their questions and comments as well as by pointing out errors and typos in the draft. The help of Daniel Gonsalez, Konstantin Klumpp, Ino Schrot and Jakob Schnell was the most valuable in this respect.

I am also thankful to my daughter and my wife for their patience during numerous weekends, which have been spent on working on the monograph.

During the work on the monograph my research on graphical models was supported by the DFG grant "Exact Relaxation-Based Inference in Graphical Models" (SA 2640/1-1).